\documentclass[12pt]{article}
\usepackage{amsfonts, amsmath, latexsym, vmargin, epic, eepic, pifont}
\usepackage{epsf}
\usepackage{color}
\usepackage{graphics}
\usepackage{epsfig}
\title{Elementary elliptic $(R,q)$-polycycles}
\setpapersize{custom}{21cm}{29.7cm}
\setmarginsrb{1.7cm}{2cm}{1.7cm}{3.2cm}{0pt}{0pt}{0pt}{5mm}
\author{Michel DEZA\\
CNRS/ENS, Paris, and Institute of Statistical Mathematics, Tokyo,\\
\ Mathieu DUTOUR \\
ENS, Paris, and Institut Rudjer Boskovi\'c, Zagreb\\
\ Mikhail SHTOGRIN
\thanks{Third author acknowledges financial support of the Russian 
Foundation of Fundamental Research (grant 05-01-00170),
the Russian Foundation for Scientific Schools (grant NSh.2185.2003.1)
and the Program OMN (Division of Mathematical Sciences) 
"Modern problems of theoretical mathematics"
of Russian Academy of Sciences.}\\
Steklov Mathematical Institute, Moscow.
}

\begin{document}
\newcommand{\RR}{\ensuremath{\mathbb{R}}}
\newcommand{\NN}{\ensuremath{\mathbb{N}}}
\newcommand{\QQ}{\ensuremath{\mathbb{Q}}}
\newcommand{\CC}{\ensuremath{\mathbb{C}}}
\newcommand{\ZZ}{\ensuremath{\mathbb{Z}}}
\newcommand{\TT}{\ensuremath{\mathbb{T}}}
\newtheorem{proposition}{Proposition}
\newtheorem{theorem}{Theorem}
\newtheorem{corollary}{Corollary}
\newtheorem{lemma}{Lemma}
\newtheorem{conjecture}{Conjecture}
\newtheorem{claim}{Claim}
\newtheorem{remark}{Remark}
\newtheorem{definition}{Definition}
\newcommand{\qed}{\hfill $\Box$ }
\newcommand{\proof}{\noindent{\bf Proof.}\ \ }

\maketitle

\begin{abstract}
We consider the
following generalization of the decomposition theorem for
polycycles.
A {\em $(R,q)$-polycycle} is, roughly, a plane graph,
whose faces, besides some disjoint {\em holes}, are $i$-gons, $i \in R$,
and whose vertices, outside of holes, are $q$-valent. 
Such polycycle is called {\em elliptic}, {\em parabolic} or
{\em hyperbolic} if $\frac{1}{q} + \frac{1}{r} - \frac{1}{2}$
(where $r={max_{i \in R}i}$)
is positive, zero or negative, respectively.

An edge on the boundary of a hole in such polycycle is called 
{\em open} if both its end-vertices have degree less than $q$. We enumerate
all elliptic {\em elementary} polycycles, i.e. those that any elliptic 
$(R,q)$-polycycle can be obtained from them by agglomeration along some 
open edges.

\end{abstract}

\section{Introduction}
Given $q\in \NN$ and $R\subset \NN$, a $(R,q)$-{\em polycycle} is
a non-empty $2$-connected plane, locally finite (i.e. any circle
contain only finite number of its vertices) graph $G$
with faces partitioned in two non-empty sets $F_1$ and $F_2$, so that:

(i) all elements of $F_1$ (called {\em proper faces}) are
combinatorial $i$-gons with $i\in R$;

(ii) all elements of $F_2$ (called {\em holes}, the exterior face(s)
\footnote{Any finite plane graph has an unique exterior face; any 
infinite plane graph can have any number of exterior faces, including $0$ 
and infinity ($2$-gonal faces are permitted).}
are amongst them) are pair-wisely disjoint, i.e. have no common vertices;

(iii) all vertices have degree within $\{2,\dots,q\}$ and all 
{\em interior} (i.e. not on the boundary of a hole)
vertices are $q$-valent.

The plane graph $G$ can be finite or infinite and some of the faces
of the set $F_2$ can be $i$-gons with $i\in R$.
Two $(R,q)$-polycycles, which are isomorphic as plane graphs,
but have different pairs $(F_1, F_2)$, will be considered non-isomorphic
in our context. The symmetry group $Aut(P)$ of a polycycle $P$,
considered below, consists of all automorphisms of plane
graph $G$ preserving the pair $(F_1,F_2)$.
Note that it is different from $Aut(G)$, the full automorphism group of 
the plane graph $G$, i.e. the group of transformations,
which preserves the edge-set
and the face-set of the plane graph $G$. In fact, $Aut(P)$ is the
stabilizer of the pair $(F_1,F_2)$ in $Aut(G)$.

%TWO POLICIES, NEED TO CHOICE (I BEND TOWARDS SECOND CHOICE BUT I AM NOT COMPLETELY SURE, THE QUESTION WAS VOID FOR usual $(q,r)$-polycycles but it become CRUCIAL HERE): of its vertices
%\begin{enumerate}
%\item Also, if a $(R,q)$-polycycle has one exterior face of
%gonality $h$ with $h\in P$ and if all its incident vertices are of
%degree $r$, then we can consider this face as interior. There are
%two choices for every such face and we will consider the corresponding
%$(R,q)$-polycycles to be different.
%\item Note that we forbide having an exterior face of gonality $h\in P$
%incident with only vertices of degree $r$. Those faces are fake exterior
%faces and should be consider to be interior.
%\end{enumerate}

The notion of $(R,q)$-polycycle is a large generalization
of the case $|R|=|F_2|=1$,
i.e. {\em $(r,q)$-polycycle} introduced by Deza and Shtogrin in \cite{DSp1}
and studied in their papers \cite{DSp1}, \cite{DSp2}, \cite{DSp3}, \cite{DS5},
\cite{7}, \cite{DS8}, \cite{DS11}, \cite{DS10}, \cite{DS12},
\cite{DS13a}, \cite{DS13b}, \cite{S1}, \cite{S2}. The case $|R|=1$, i.e.
{\em $(r,q)$-polycycles with holes}, was considered in \cite{DDS2}.

%has only one exterior face and if $P=\{p\}$ then
%it is called in ???? a $(p,3)$-polycycle.

A {\em boundary} of a $(R,q)$-polycycle $P$ is the boundary of any of
its holes.

A {\em bridge} of a $(R,q)$-polycycle is an edge, which is not on a boundary
and goes from a hole to a hole (possibly, the same).

An $(R,q)$-polycycle is called {\em elementary} if it has no bridges.

An {\em open edge} of a $(R,q)$-polycycle is an edge on a boundary,
such that each of its end-vertices have degree less than $q$.

%Given this notion, one can state the following theorem:
\begin{theorem}\label{Fundamental_Decomposition_theorem}
Every $(R,q)$-polycycle is uniquely formed by the agglomeration of 
elementary ones along open edges or, in other words, it can be 
uniquely cut,
along the bridges, into the elementary ones.
\end{theorem}

Hence, the interesting question is to enumerate those elementary
$(R,q)$-polycycles. Call a $(R,q)$-polycycle {\em elliptic},
{\em parabolic} or {\em hyperbolic} if the number
$\frac{1}{q} + \frac{1}{r} - \frac{1}{2}$ (where $r=max_{i \in R}i$)
is positive, zero or
negative, respectively. The number of elementary
$(r,q)$-polycycles is uncountable for any parabolic or hyperbolic
pairs $(r,q)$. But in \cite{DS8} and \cite{DS10} all elliptic
elementary $(r,q)$-polycycles were determined. Namely, the countable
set of all elementary $(5,3)$- and $(3,5)$-polycycles was 
described; the cases of $(3,3)$-, $(4,3)$- and
$(3,4)$-polycycles are easy. In \cite{DDS2} all elliptic
elementary $(\{r\},q)$-polycycles (i.e. $(r,q)$-polycycles with holes)
were determined as a part of a more 
elaborate classification.

The purpose of this paper is to generalize it for elliptic
$(R,q)$-polycycles. In fact, we will consider the main case 
$R=\{i\mbox{~:~} 2\leq i\leq r\}$ and so, all such elliptic
possibilities are $(\{2,3,4,5\}, 3)$-, $(\{2,3\},4)$- and
$(\{2,3\},5)$-polycycles.

Given a $(R,q)$-polycycle $P$, one can define another
$(R,q)$-polycycle $P'$ by removing a face $f$ 
from $F_1$, i.e. by considering it as a hole.
If $f$ has no common vertices with other faces from $F_1$, then 
removing of it leaves unchanged the plane
graph $G$ and only changes the pair $(F_1, F_2)$.
If $f$ has some edges in common with a hole, then we remove them.
If $f$ has a common vertex with a hole, then we split it in two.
See below this operation for a $2$-gon, which is incident to a 
vertex on the boundary.
\begin{center}
\epsfig{height=30mm, file=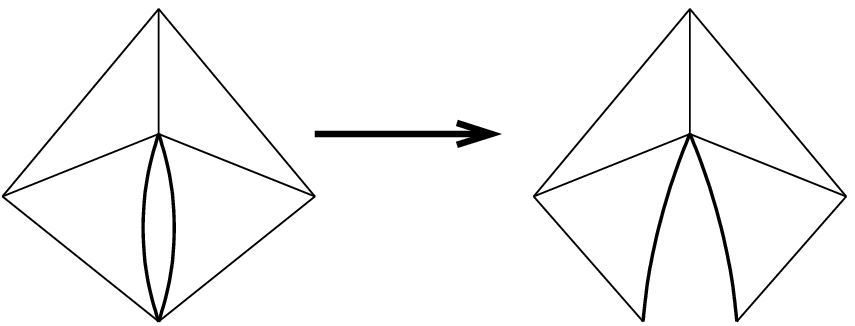}\par
\end{center}
A $(R,q)$-polycycle $P$ is called {\em extensible} if there exists
another $(R,q)$-polycycle $P'$, such that the elimination of a face
of $P'$ yields $P$.

The same plane graph $G$ can admit several polycycle realizations.
For example, 
$Prism_3$ admits two realizations as a $(\{3,4\},3)$-polycycle (see 
$(4,20-22, D_{3h})$ below). 
For elementary $(\{2,3,4,5\},3)$-polycycles,
we give the list in the form 
$(n_f,nb_1-\dots-nb_p,Aut(G))$ with $n_f$ being the number of proper
faces (from $F_1$),
and $nb_1, \dots, nb_p$ being the position in the corresponding lists
of elementary polycycles with $n_f$ faces and $Aut(G)$ being 
the automorphism group of plane graph $G$:
$(4,12-13,C_s)$, $(4,20-22,D_{3h})$, 
$(4,23-5^{th}\, of\,\, Lemma\, 1, C_{2\nu})$, $(5,26-28-30, C_{2\nu})$, 
$(5,27-29, C_s)$, $(6,21-22-31, C_{3\nu})$, $(6,28-34,D_{5h})$, $(7,24-36, D_{3d})$, 
$(7,29-31, C_{2\nu})$, $(8, 25-29, D_{3h})$, $(9, 11-16, D_{4d})$.
For elementary $(\{2,3\},5)$-polycycles, this concerns 
the 10th and 20th of Figure \ref{Sporadic5_valentFirst} ($Aut(G)=C_s$),
the 21st of Figure \ref{Sporadic5_valentFirst} and 1st
of Figure \ref{Sporadic5_valentSecond} ($Aut(G)=D_3$), 
the 14th and 22nd of Figure \ref{Sporadic5_valentSecond} ($Aut(G)=D_{2d}$).
If $Aut(P) \not= Aut(G)$ and no other polycycle with the same 
graph $G$ exists, we indicate it by putting the second group in parenthesis.

Theorem \ref{Fundamental_Decomposition_theorem}, together with the
determination of the elementary $(R,q)$-polycycles, is
especially useful in extremal problems, related to the number of
interior/exterior vertices, see \cite{DS8} and \cite{DS10}, and
in classification of {\em face-regular two-faced} maps,
see \cite{FaceRegular}.

%Any $2$-connected $s$-valent plane graph, whose gonalities of faces
%belong to the set $P$ is a $(R,q)$-polycycle according to our
%definition with $F_2=\emptyset$. Furthermore, those graphs are elementary,
%since there are no bridges. Since those graphs have already been
%enumerated, we exclude them in the enumeration process.

%An elementary $(P,3)$-polycycle is called {\em extensible} if there exists
%another elementary $(P,3)$-polycycle, which can be obtained by the
%addition of a face on the boundary of a hole; so, at least one boundary
%vertex will become an interior vertex.
%An elementary $(P,s)$-polycycle with $s=4$ or $5$ is called {\em extensible}
%if there exist another $(P,s)$-polycycle, which can be obtained by the
%addition of faces around a boundary vertex, so that it becomes an
%interior vertex.

We thank Gil Kalai for putting a question, which
lead us to this study.

\section{Classification of elementary $(\{2,3,4,5\},3)$-polycycles}

A $(R,3)$-polycycle is called {\em totally elementary} if
it is elementary and if, after removing any face adjacent to a hole,
one obtains a non-elementary $(R,3)$-polycycle. So, an elementary
$(R,3)$-polycycle is totally elementary if and only if it is not
the result of an extension of some elementary $(R,3)$-polycycle.

We will classify those polycycles in a number of steps. At first, we find
all totally elementary $(\{2,3,4,5\},3)$-polycycles.
Then, all other elementary ones are obtained by adding faces to the
existing elementary or totally elementary polycycles.

\begin{theorem}\label{TheoremTotalElementary}
The list of totally elementary $(\{3,4,5\}, 3)$-polycycles consists of:
\begin{enumerate}
\item[(i)] three isolated $i$-gons, $i\in \{3,4,5\}$\footnote{In fact, 
they belong to the infinite series of $(r,2)$-polycycles ({\em 
monocycles}) consisting of isolated $r$-gons, $r\geq 2$.}:
\begin{center}
%\begin{minipage}{3cm}
%\centering
%\resizebox{2.0cm}{!}{\includegraphics{ElementaryDrawing/2gon.eps}}\par
%\end{minipage}
\begin{minipage}{3cm}
\centering
\resizebox{2.0cm}{!}{\includegraphics{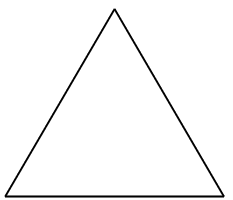}}\par
$C_{3\nu}$ $(D_{3h})$
\end{minipage}
\begin{minipage}{3cm}
\centering
\resizebox{2.0cm}{!}{\includegraphics{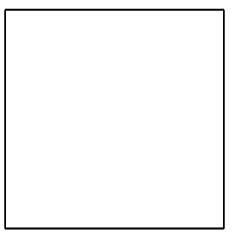}}\par
$C_{4\nu}$ $(D_{4h})$
\end{minipage}
\begin{minipage}{3cm}
\centering
\resizebox{2.0cm}{!}{\includegraphics{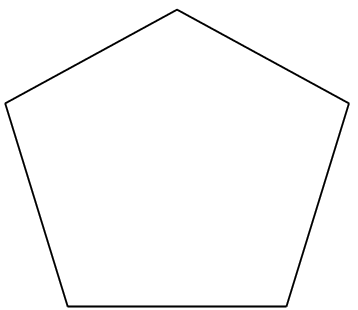}}\par
$C_{5\nu}$ $(D_{5h})$
\end{minipage}
\end{center}
\item[(ii)] all ten triples of $i$-gons, $i\in \{3,4,5\}$:
\begin{center}
\begin{minipage}{3cm}
\centering
\resizebox{2.0cm}{!}{\includegraphics{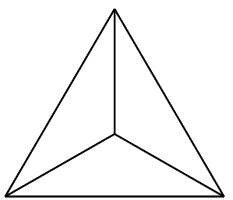}}\par
$C_{3\nu}$,~nonext.~$(T_d)$
\end{minipage}
\begin{minipage}{3cm}
\centering
\resizebox{2.0cm}{!}{\includegraphics{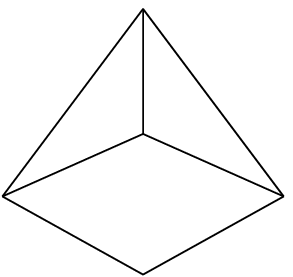}}\par
$C_s$,~nonext.~$(C_{2\nu})$
\end{minipage}
\begin{minipage}{3cm}
\centering
\resizebox{2.0cm}{!}{\includegraphics{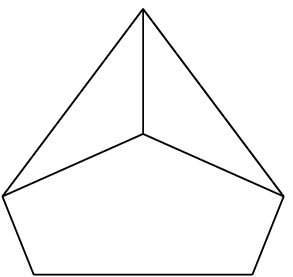}}\par
$C_s$ $(C_{2\nu})$
\end{minipage}
\begin{minipage}{3cm}
\centering
\resizebox{2.0cm}{!}{\includegraphics{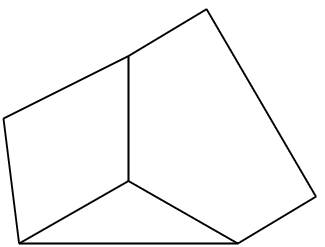}}\par
$C_{1}$
\end{minipage}
\begin{minipage}{3cm}
\centering
\resizebox{2.0cm}{!}{\includegraphics{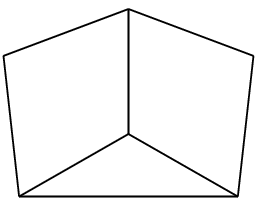}}\par
$C_s$
\end{minipage}
\begin{minipage}{3cm}
\centering
\resizebox{2.0cm}{!}{\includegraphics{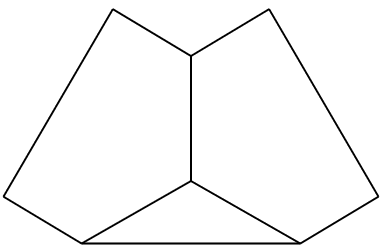}}\par
$C_s$
\end{minipage}
\begin{minipage}{3cm}
\centering
\resizebox{2.0cm}{!}{\includegraphics{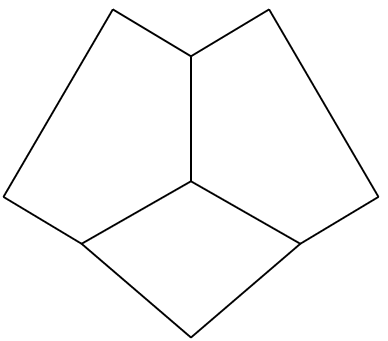}}\par
$C_s$
\end{minipage}
\begin{minipage}{3cm}
\centering
\resizebox{2.0cm}{!}{\includegraphics{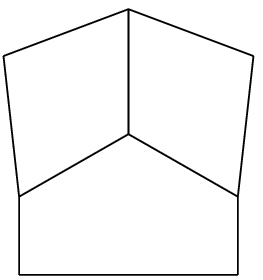}}\par
$C_s$
\end{minipage}
\begin{minipage}{3cm}
\centering
\resizebox{2.0cm}{!}{\includegraphics{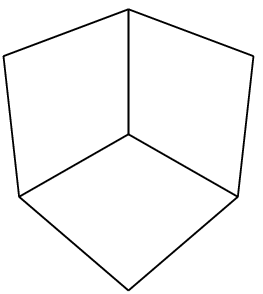}}\par
$C_{3\nu}$
\end{minipage}
\begin{minipage}{3cm}
\centering
\resizebox{2.0cm}{!}{\includegraphics{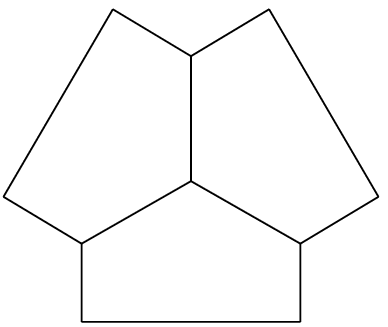}}\par
$C_{3\nu}$
\end{minipage}

\end{center}
\item[(iii)] the following doubly infinite $(\{5\},3)$-polycycle, denoted by $Barrel_{\infty}$:
\begin{center}
\begin{minipage}{14cm}
\centering
\epsfig{height=16mm, file=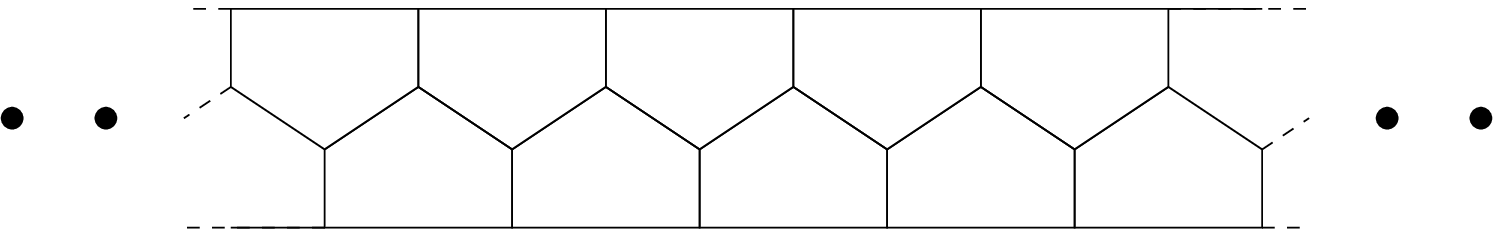}\par
pma2, nonext.
\end{minipage}
\end{center}
\item[(iv)] the infinite series\footnote{The $Barrel_m$ is a $3$-valent 
plane graph, consisting of two $m$-gons separated by two $m$-rings of 
$5$-gons.} of $Barrel_m$, $m\geq 2$ (two $m$-gonal holes, non-extensible for 
$m\not=3,4,5$, with symmetry $D_{md}$), represented below for 
$m=2,3,4,5$:
\begin{center}
\begin{minipage}{3cm}
\centering
\epsfig{height=20mm, file=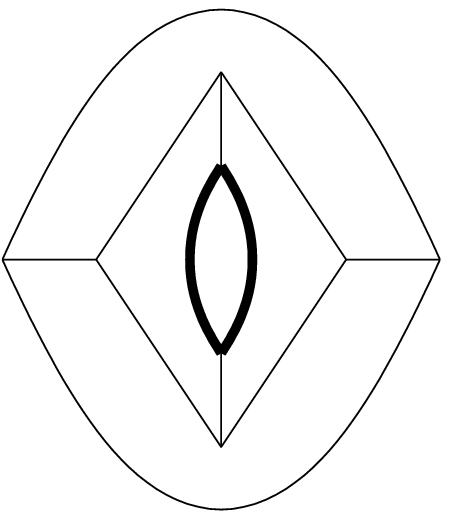}\par
$D_{2d}$
\end{minipage}
\begin{minipage}{3cm}
\centering
\epsfig{height=20mm, file=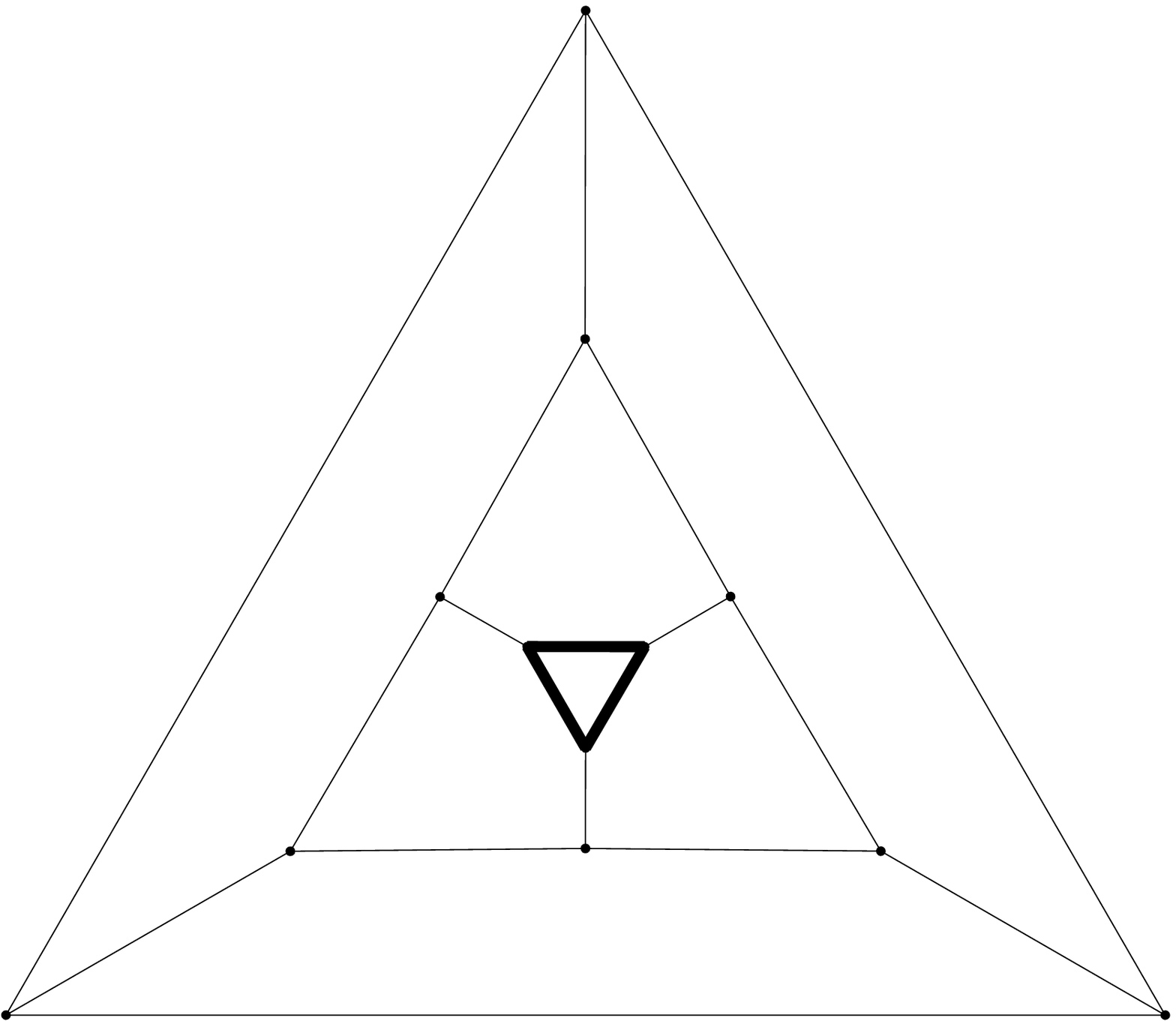}\par
$D_{3d}$
\end{minipage}
\begin{minipage}{3cm}
\centering
\epsfig{height=20mm, file=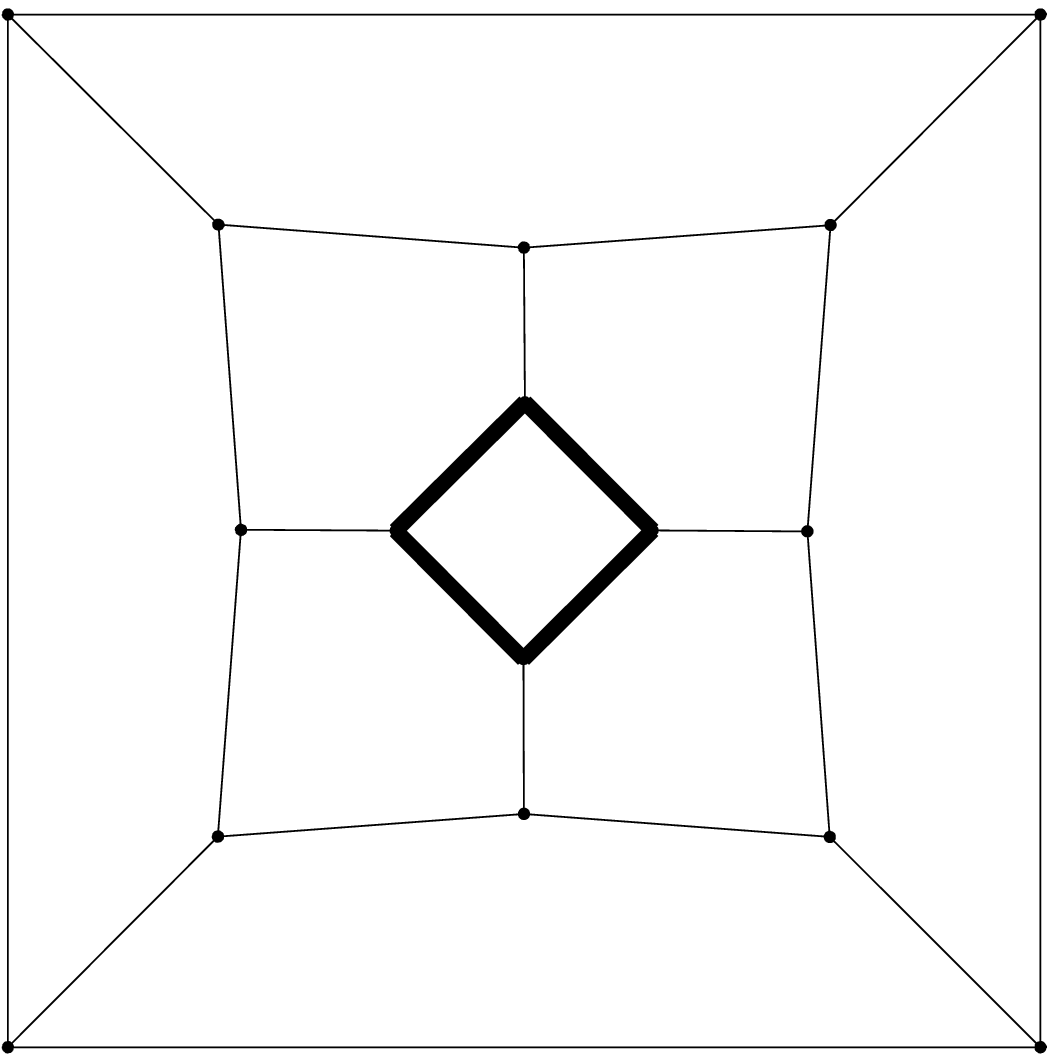}\par
$D_{4d}$
\end{minipage}
\begin{minipage}{3cm}
\centering
\epsfig{height=20mm, file=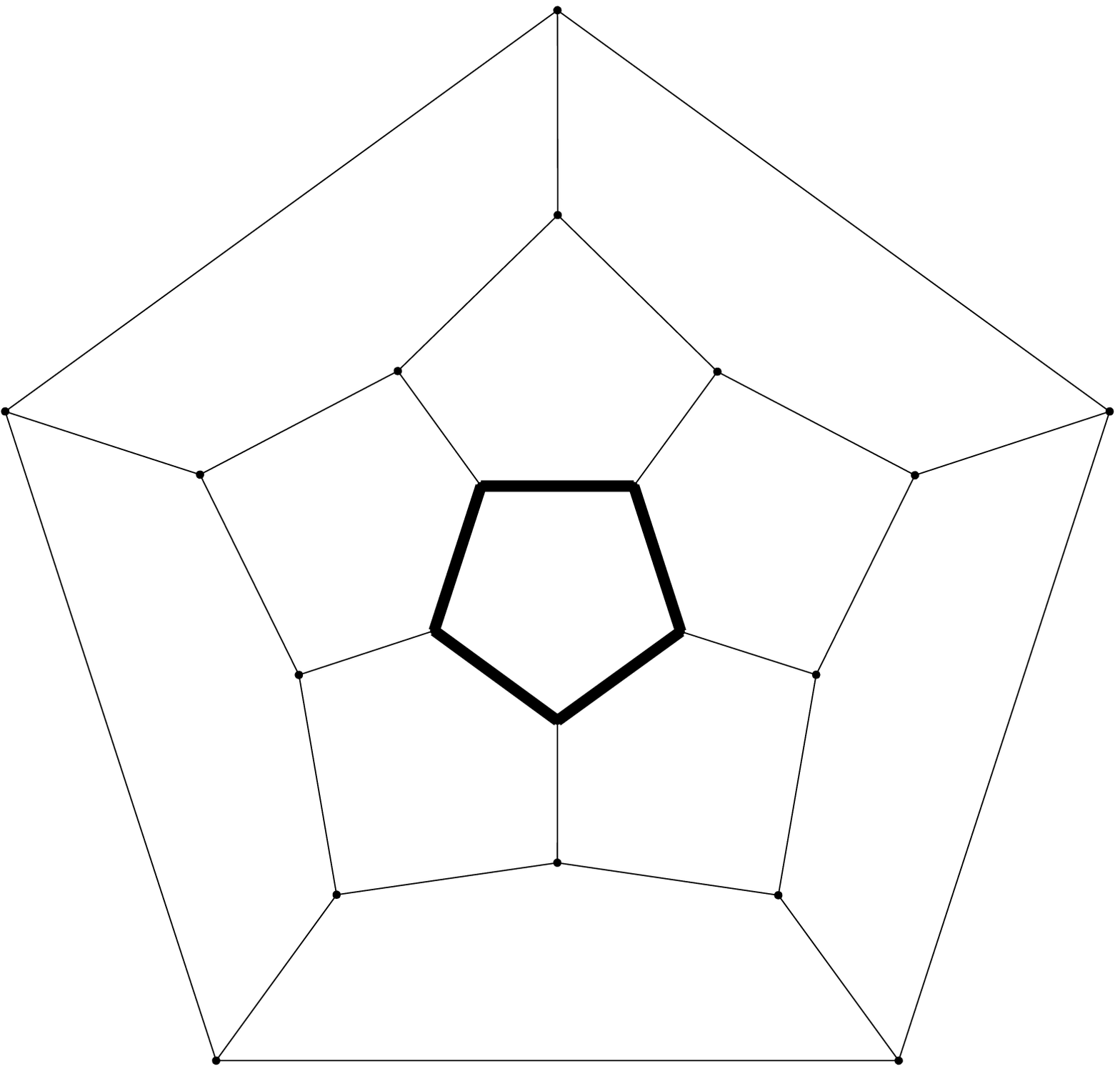}\par
$D_{5d}$ $(I_h)$
\end{minipage}

\end{center}

\end{enumerate}

\end{theorem}
\proof Take a totally elementary $(\{2,3,4,5\},3)$-polycycle $P$.
If $|F_1|=1$, then $P$ is, clearly, totally elementary; so,
let us assume that $|F_1|\geq 2$. If $|F_1|=2$,  
then it is, clearly, not elementary; so, assume $|F_1|\geq 3$.
Of course, $P$ has at least one interior vertex; let $v$ be such a vertex.
Furthermore, one can assume that $v$ is adjacent to a vertex $v'$, which
is incident to a hole.

The vertex $v$ is incident to three faces $f_1$, $f_2$, $f_3$. Let us
denote by $v_{ij}$ unique vertex incident to $f_i$, $f_j$ and
adjacent to $v$. Without loss of generality, one can suppose that $v'$
is incident to the faces $f_1$ and $f_2$, i.e. that $v'=v_{12}$.
%\begin{center}
%\epsfig{height=30mm, file=ElementaryDrawing/ProofDraw3reg_1.eps}\par
%\end{center}

\begin{center}
\centering
\begin{picture}(0,0)%
\includegraphics{ProofDraw3reg_1.pstex}%
\end{picture}%
\setlength{\unitlength}{3947sp}%
\begingroup\makeatletter\ifx\SetFigFont\undefined%
\gdef\SetFigFont#1#2#3#4#5{%
  \reset@font\fontsize{#1}{#2pt}%
  \fontfamily{#3}\fontseries{#4}\fontshape{#5}%
  \selectfont}%
\fi\endgroup%
\begin{picture}(1974,1797)(64,-1873)
\put(1126,-1111){\makebox(0,0)[lb]{\smash{\SetFigFont{12}{14.4}{\rmdefault}{\mddefault}{\updefault}{\color[rgb]{0,0,0}$v$}%
}}}
\put(1126,-661){\makebox(0,0)[lb]{\smash{\SetFigFont{12}{14.4}{\rmdefault}{\mddefault}{\updefault}{\color[rgb]{0,0,0}$v'=v_{12}$}%
}}}
\put(1651,-1561){\makebox(0,0)[lb]{\smash{\SetFigFont{12}{14.4}{\rmdefault}{\mddefault}{\updefault}{\color[rgb]{0,0,0}$v_{23}$}%
}}}
\put(226,-1561){\makebox(0,0)[lb]{\smash{\SetFigFont{12}{14.4}{\rmdefault}{\mddefault}{\updefault}{\color[rgb]{0,0,0}$v_{13}$}%
}}}
\put(1501,-961){\makebox(0,0)[lb]{\smash{\SetFigFont{12}{14.4}{\rmdefault}{\mddefault}{\updefault}{\color[rgb]{0,0,0}$f_2$}%
}}}
\put(976,-1561){\makebox(0,0)[lb]{\smash{\SetFigFont{12}{14.4}{\rmdefault}{\mddefault}{\updefault}{\color[rgb]{0,0,0}$f_3$}%
}}}
\put(451,-886){\makebox(0,0)[lb]{\smash{\SetFigFont{12}{14.4}{\rmdefault}{\mddefault}{\updefault}{\color[rgb]{0,0,0}$f_1$}%
}}}
\end{picture}

\end{center}

The removal of the face $f_1$ yields a non-elementary polycycle; so,
there is at least one bridge separating $P-f_1$ in two parts.
Such bridge should have an end-vertex incident to $f_1$. 
The same holds for $f_2$. The proof consists of a number of cases.

{\bf $1^{st}$ case:} If $e_1=\{v,v_{23}\}$ and $e_2=\{v,v_{13}\}$ are
bridges for $P-f_1$, $P-f_2$, respectively, then from the constraint
that faces $f_i$ are $p$-gons with $p\leq 5$, one sees that each 
face $f_i$ is adjacent in $P$ to at most one other face.
Furthermore, if $f_i$ is adjacent to another face, then this adjacence is
along a bridge, which is forbidden.
Hence, $F_1=\{f_1, f_2, f_3\}$.

{\bf $2^{nd}$ case:} Let us assume now that $e_1=\{v, v_{23}\}$ is a
bridge for $P-f_1$, but $e_2=\{v,v_{13}\}$ is not a bridge for $P-f_2$.
Then, since $f_2$ is a $p$-gon with $p\leq 5$, it is adjacent to at most
one other face and, if so, then along a bridge, which is impossible.
So, $f_2$ is adjacent to only $f_1$ and $f_3$ and, since $e_2$ is
not a bridge for $P-f_2$, one obtains that $P-f_2$ is elementary,
which contradicts the hypothesis.

{\bf $3^{rd}$ case:} Let us assume that neither $e_1$, nor $e_2$ 
are bridges for $P-f_1$ and $P-f_2$. From the consideration of
previous two cases, one has that every vertex $v$, adjacent to 
a vertex on the boundary, is in this $3^{rd}$ case.

The first subcase, which can happen only if $f_1$ is $5$-gon,
happens, when the over-lined edge $e'$, in the drawing below, is a bridge.
\begin{center}
\centering
\begin{picture}(0,0)%
\includegraphics{ProofDraw3reg_2.pstex}%
\end{picture}%
\setlength{\unitlength}{3947sp}%
\begingroup\makeatletter\ifx\SetFigFont\undefined%
\gdef\SetFigFont#1#2#3#4#5{%
  \reset@font\fontsize{#1}{#2pt}%
  \fontfamily{#3}\fontseries{#4}\fontshape{#5}%
  \selectfont}%
\fi\endgroup%
\begin{picture}(2874,1674)(-836,-1873)
\put(1126,-1111){\makebox(0,0)[lb]{\smash{\SetFigFont{12}{14.4}{\rmdefault}{\mddefault}{\updefault}{\color[rgb]{0,0,0}$v$}%
}}}
\put(1501,-961){\makebox(0,0)[lb]{\smash{\SetFigFont{12}{14.4}{\rmdefault}{\mddefault}{\updefault}{\color[rgb]{0,0,0}$f_2$}%
}}}
\put(1126,-661){\makebox(0,0)[lb]{\smash{\SetFigFont{12}{14.4}{\rmdefault}{\mddefault}{\updefault}{\color[rgb]{0,0,0}$v'=v_{12}$}%
}}}
\put(-74,-1786){\makebox(0,0)[lb]{\smash{\SetFigFont{12}{14.4}{\rmdefault}{\mddefault}{\updefault}{\color[rgb]{0,0,0}$h$}%
}}}
\put(976,-1786){\makebox(0,0)[lb]{\smash{\SetFigFont{12}{14.4}{\rmdefault}{\mddefault}{\updefault}{\color[rgb]{0,0,0}$f_3$}%
}}}
\put(526,-886){\makebox(0,0)[lb]{\smash{\SetFigFont{12}{14.4}{\rmdefault}{\mddefault}{\updefault}{\color[rgb]{0,0,0}$f_1$}%
}}}
\put(1651,-1561){\makebox(0,0)[lb]{\smash{\SetFigFont{12}{14.4}{\rmdefault}{\mddefault}{\updefault}{\color[rgb]{0,0,0}$v_{23}$}%
}}}
\put(226,-1486){\makebox(0,0)[lb]{\smash{\SetFigFont{12}{14.4}{\rmdefault}{\mddefault}{\updefault}{\color[rgb]{0,0,0}$v_{13}$}%
}}}
\put(-524,-736){\makebox(0,0)[lb]{\smash{\SetFigFont{12}{14.4}{\rmdefault}{\mddefault}{\updefault}{\color[rgb]{0,0,0}$g$}%
}}}
\put(-299,-1186){\makebox(0,0)[lb]{\smash{\SetFigFont{12}{14.4}{\rmdefault}{\mddefault}{\updefault}{\color[rgb]{0,0,0}$e'$}%
}}}
\end{picture}

\end{center}
The face $g$ is adjacent to the faces $h$ and $f_1$ and, possibly,
to another face $g'$. But if $g$ is adjacent to such a face
$g'$, it is along a bridge of $P$; hence, $g$ is adjacent only to
$h$ and $f_1$. So, $P-g$ is elementary, which is impossible.

So, the edge $e'$ is not a bridge and this forces the face $h$ to
be $5$-gonal. Hence, the vertex $v_{13}$ is in the same situation 
as the vertex $v$, described in the diagram below:
\begin{center}
\centering
\begin{picture}(0,0)%
\includegraphics{ProofDraw3reg_3.pstex}%
\end{picture}%
\setlength{\unitlength}{3947sp}%
\begingroup\makeatletter\ifx\SetFigFont\undefined%
\gdef\SetFigFont#1#2#3#4#5{%
  \reset@font\fontsize{#1}{#2pt}%
  \fontfamily{#3}\fontseries{#4}\fontshape{#5}%
  \selectfont}%
\fi\endgroup%
\begin{picture}(2874,1962)(-836,-2161)
\put(-299,-1711){\makebox(0,0)[lb]{\smash{\SetFigFont{12}{14.4}{\rmdefault}{\mddefault}{\updefault}{\color[rgb]{0,0,0}$h$}%
}}}
\put(1651,-1561){\makebox(0,0)[lb]{\smash{\SetFigFont{12}{14.4}{\rmdefault}{\mddefault}{\updefault}{\color[rgb]{0,0,0}$v_{23}$}%
}}}
\put(1126,-1111){\makebox(0,0)[lb]{\smash{\SetFigFont{12}{14.4}{\rmdefault}{\mddefault}{\updefault}{\color[rgb]{0,0,0}$v$}%
}}}
\put(151,-1111){\makebox(0,0)[lb]{\smash{\SetFigFont{12}{14.4}{\rmdefault}{\mddefault}{\updefault}{\color[rgb]{0,0,0}$w$}%
}}}
\put(-299,-1261){\makebox(0,0)[lb]{\smash{\SetFigFont{12}{14.4}{\rmdefault}{\mddefault}{\updefault}{\color[rgb]{0,0,0}$e'$}%
}}}
\put(1126,-661){\makebox(0,0)[lb]{\smash{\SetFigFont{12}{14.4}{\rmdefault}{\mddefault}{\updefault}{\color[rgb]{0,0,0}$v_{12}$}%
}}}
\put(-449,-811){\makebox(0,0)[lb]{\smash{\SetFigFont{12}{14.4}{\rmdefault}{\mddefault}{\updefault}{\color[rgb]{0,0,0}$g$}%
}}}
\put(226,-1486){\makebox(0,0)[lb]{\smash{\SetFigFont{12}{14.4}{\rmdefault}{\mddefault}{\updefault}{\color[rgb]{0,0,0}$v_{13}$}%
}}}
\put(901,-1711){\makebox(0,0)[lb]{\smash{\SetFigFont{12}{14.4}{\rmdefault}{\mddefault}{\updefault}{\color[rgb]{0,0,0}$f_3$}%
}}}
\put(451,-886){\makebox(0,0)[lb]{\smash{\SetFigFont{12}{14.4}{\rmdefault}{\mddefault}{\updefault}{\color[rgb]{0,0,0}$f_1$}%
}}}
\put(1501,-961){\makebox(0,0)[lb]{\smash{\SetFigFont{12}{14.4}{\rmdefault}{\mddefault}{\updefault}{\color[rgb]{0,0,0}$f_2$}%
}}}
\end{picture}

\end{center}

So, one can repeat the construction. If, at some point, $e_1$ is a bridge,
then the construction stops; otherwise, one can continue indefinitely.
If one obtains only different vertices, then it
means that we have the  $(\{5\},3)$-polycycle
$Barrel_{\infty}$; otherwise, we obtain a loop of vertices, i.e. a circuit
of vertices,
which appear over and over, i.e. $Barrel_m$ for some $m\geq 2$. \qed

\begin{lemma}\label{WeHaveA2gon}
All elementary $(\{2,3,4,5\}, 3)$-polycycles, containing a $2$-gon, are 
following eight ones:
\begin{center}
\begin{minipage}{14cm}
\begin{center}
\begin{minipage}{3cm}
\centering
\epsfig{height=20mm, file=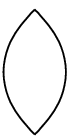}\par
$C_{2\nu}$ $(D_{2h})$
\end{minipage}
\begin{minipage}{3cm}
\centering
\epsfig{height=20mm, file=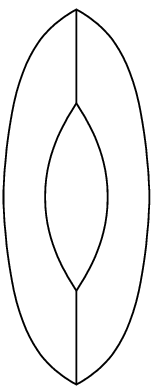}\par
$C_{2\nu}$,~nonext.~$(D_{2h})$
\end{minipage}
\begin{minipage}{3cm}
\centering
\epsfig{height=20mm, file=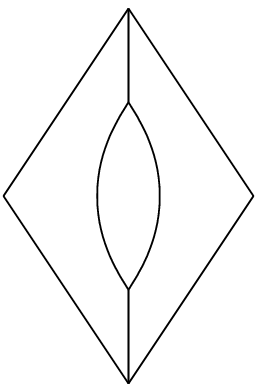}\par
$C_{2\nu}$
\end{minipage}
\begin{minipage}{3cm}
\centering
\epsfig{height=20mm, file=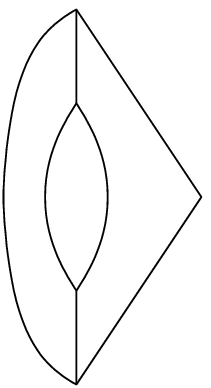}\par
$C_s$, nonext.
\end{minipage}
\begin{minipage}{3cm}
\centering
\epsfig{height=20mm, file=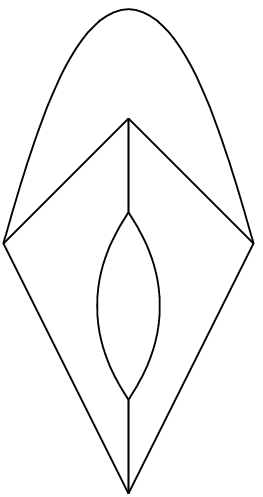}\par
$C_s$,~nonext.~$(C_{2\nu})$
\end{minipage}
\begin{minipage}{3cm}
\centering
\epsfig{height=20mm, file=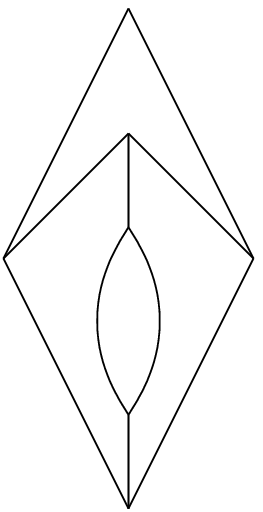}\par
$C_s$,~nonext.~$(C_{2\nu})$
\end{minipage}
\begin{minipage}{3cm}
\centering
\epsfig{height=20mm, file=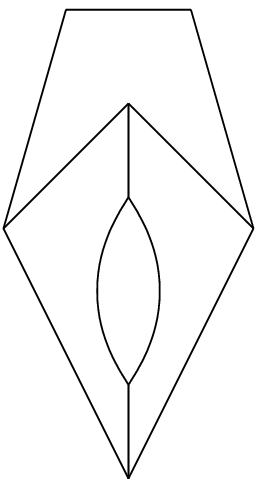}\par
$C_s$ $(C_{2\nu})$
\end{minipage}
\begin{minipage}{3cm}
\centering
\epsfig{height=20mm, file=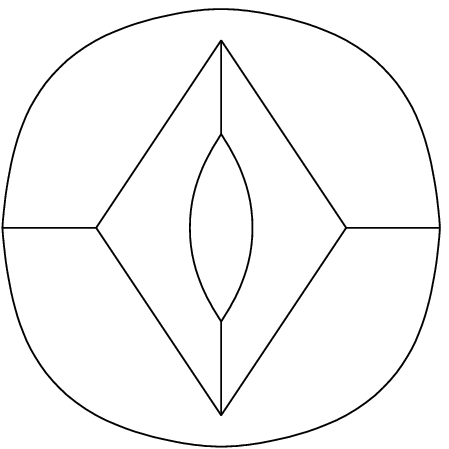}\par
$C_{2\nu}$,~nonext.~$(D_{2d})$
\end{minipage}
\end{center}
\end{minipage}
\end{center}

\end{lemma}
\proof Let $P$ be such polycycle. Clearly, the $2$-gon is
the only possibility if $|F_1|=1$. If $|F_1|=2$, then it is
not elementary. If $|F_1|\geq 3$, then the $2$-gon should be inside
of the structure. So, $P$ contain, as a subgraph, one of three
following graphs:
\begin{center}
\begin{minipage}{3cm}
\centering
\epsfig{height=20mm, file=ElementaryDrawing/El244.eps}\par
\end{minipage}
\begin{minipage}{3cm}
\centering
\epsfig{height=20mm, file=ElementaryDrawing/El245.eps}\par
\end{minipage}
\begin{minipage}{3cm}
\centering
\epsfig{height=20mm, file=ElementaryDrawing/El255.eps}\par
\end{minipage}
\end{center}
So, the only possibilities for $P$ are those given in above Lemma. \qed

We will enumerate now all elementary $(\{3,4,5\}, 3)$-polycycles.
If such a polycycle is not totally elementary, then it is obtained from
another elementary $(\{3,4,5\}, 3)$-polycycle (totally elementary or not) by
addition of another face.

\begin{figure}
\begin{center}
\begin{center}
Infinite series $\alpha\alpha$ of elementary $(\{2,3,4,5\},3)$-polycycles:
\begin{center}
\begin{minipage}{55mm}
\centering
\epsfig{height=14mm, file=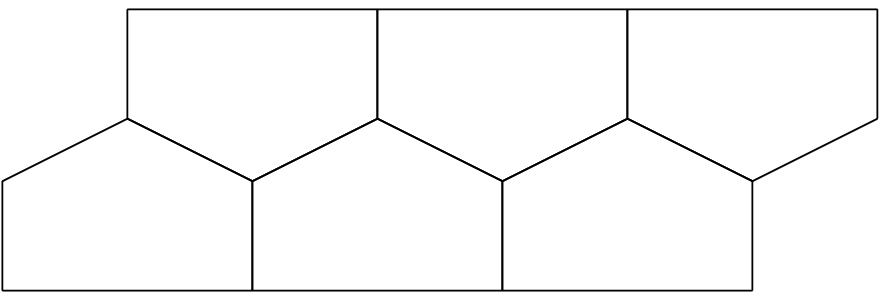}\par
\end{minipage}
\begin{minipage}{55mm}
\centering
\epsfig{height=14mm, file=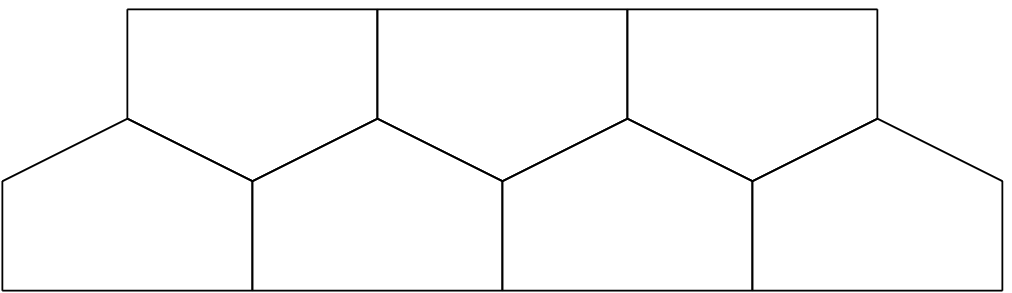}\par
\end{minipage}
\begin{minipage}{55mm}
\centering
\epsfig{height=14mm, file=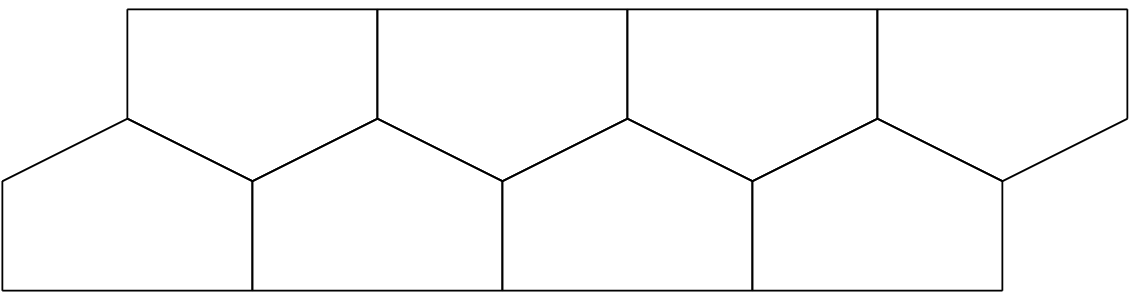}\par
\end{minipage}

\end{center}
\end{center}
\begin{center}
Infinite series $\beta\varepsilon$ of elementary $(\{2,3,4,5\},3)$-polycycles:
\begin{center}
\begin{minipage}{55mm}
\centering
\epsfig{height=17mm, file=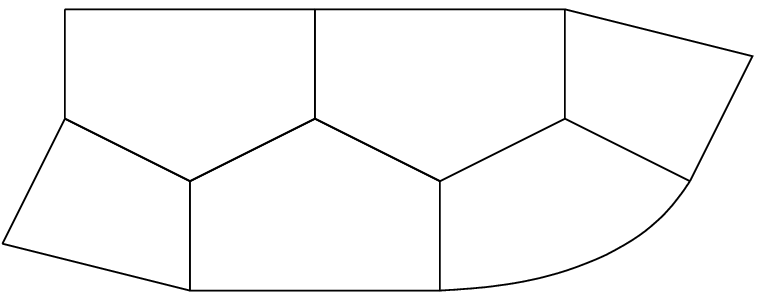}\par
\end{minipage}
\begin{minipage}{55mm}
\centering
\epsfig{height=17mm, file=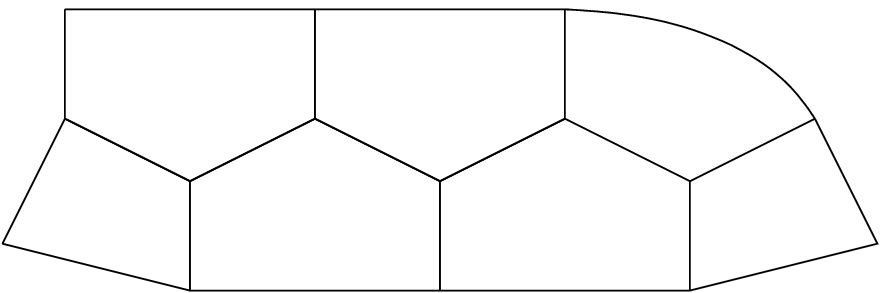}\par
\end{minipage}
\begin{minipage}{55mm}
\centering
\epsfig{height=17mm, file=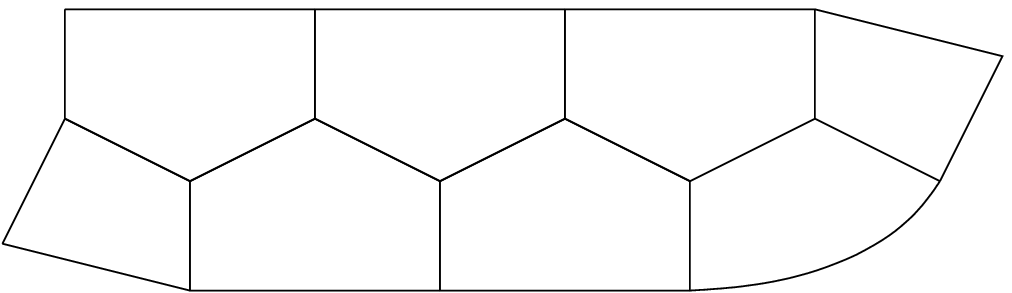}\par
\end{minipage}

\end{center}
\end{center}

\end{center}
\caption{The first $3$ members (starting with $6$ faces) of two infinite 
series, amongst $21$ series of $(\{2,3,4,5\},3)$-polycycles in 
Theorem \ref{Theorem345_3valent} (v)}
\label{ExampleOf21InfiniteSeries}
\end{figure}

\begin{theorem}\label{Theorem345_3valent}
The list of elementary $(\{2,3,4,5\},3)$-polycycles consists of:

\begin{enumerate}
\item[(i)] eight $(\{2,3,4,5\},3)$-polycycles, given in the list of
Lemma \ref{WeHaveA2gon}.

\item[(ii)] the sporadic $(\{3,4,5\},3)$-polycycles, given for $1$ and $3$
faces in Theorem \ref{TheoremTotalElementary} (i),(ii)  and the remainder,
after proof of this Theorem.

\item[(iii)] the infinite series of $(\{5\},3)$-polycycles $Barrel_m$ with 
$2\le m\le \infty$, given in Theorem 
\ref{TheoremTotalElementary} (iii), (iv).

\item[(iv)] six $(\{3,4,5\},3)$-polycycles, infinite in one direction:
\begin{center}
\begin{minipage}{8cm}
\centering
\epsfig{height=10mm, file=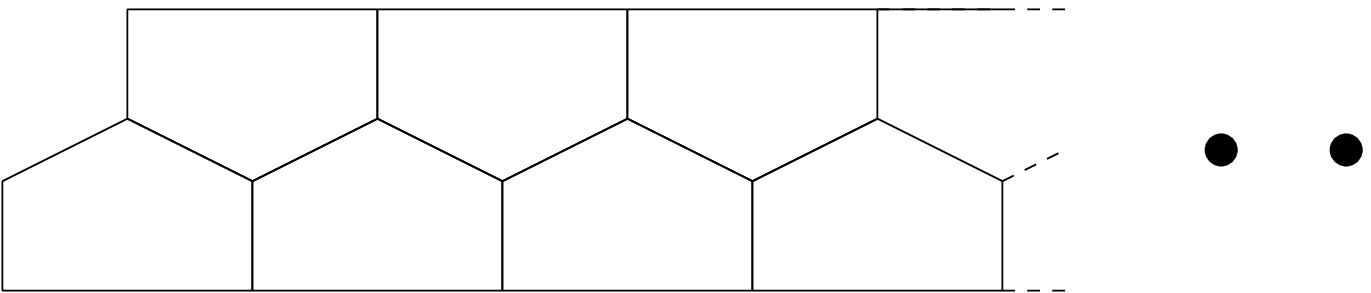}\par
$\alpha$: $C_1$
\end{minipage}
\begin{minipage}{8cm}
\centering
\epsfig{height=10mm, file=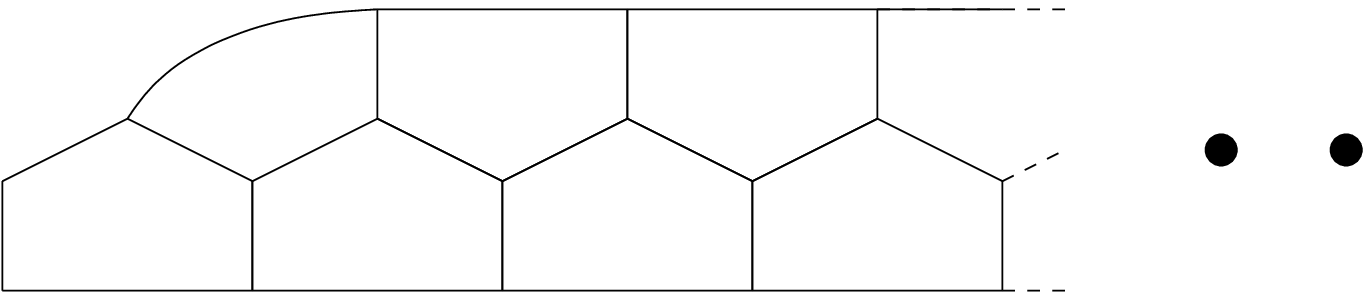}\par
$\delta$: $C_1$
\end{minipage}
\begin{minipage}{8cm}
\centering
\epsfig{height=10mm, file=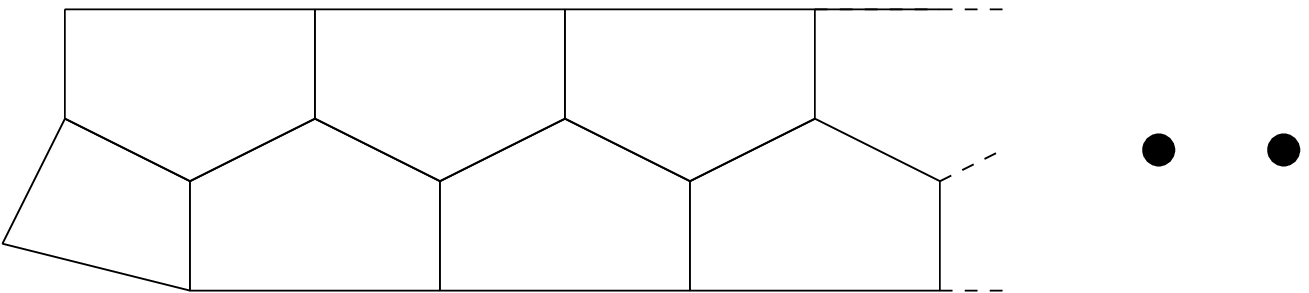}\par
$\beta$: $C_1$
\end{minipage}
\begin{minipage}{8cm}
\centering
\epsfig{height=10mm, file=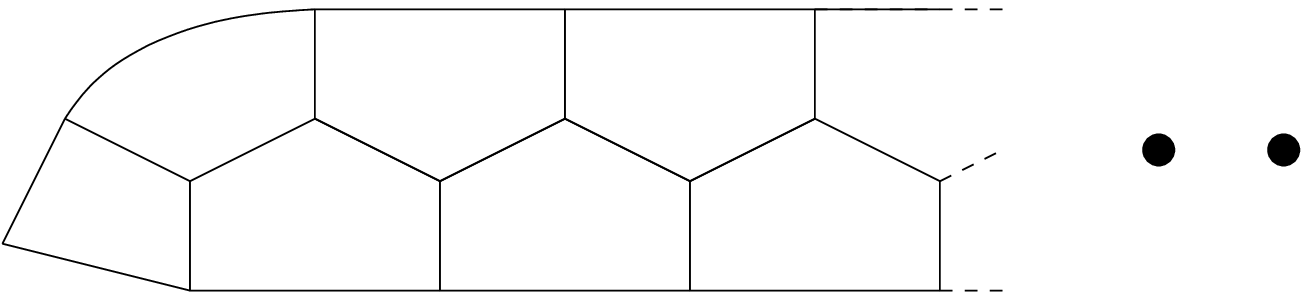}\par
$\varepsilon$: $C_1$, nonext.
\end{minipage}
\begin{minipage}{8cm}
\centering
\epsfig{height=10mm, file=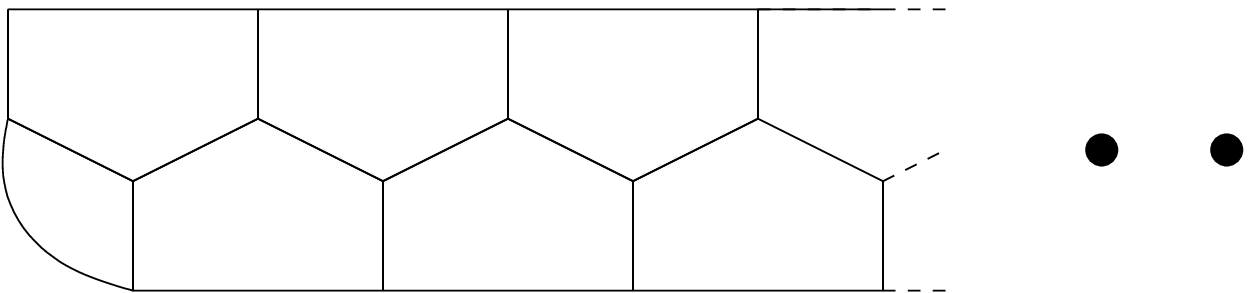}\par
$\gamma$: $C_1$, nonext.
\end{minipage}
\begin{minipage}{8cm}
\centering
\epsfig{height=10mm, file=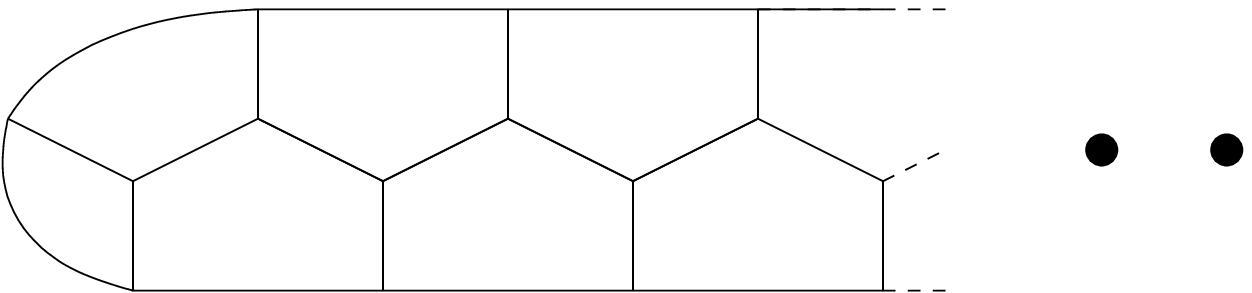}\par
$\mu$: $C_1$, nonext.
\end{minipage}
\end{center}

\item[(v)] $21={{6+1}\choose 2}$ infinite series obtained by taking two 
endings of the infinite polycycles of (iv) above and concatenating them.

For example, merging of $\alpha$ with itself produces the infinite 
series of elementary $(\{5\},3)$-polycycles, denoted by $E_n$ in
\cite{DS8}, $0 \le n \le \infty$.
For $n=0,1,2,3, \infty$, it is the $5$-gon
(Theorem \ref{TheoremTotalElementary} (i)), the triple of 
$5$-gons (Theorem \ref{TheoremTotalElementary} (ii)), 16th
in the list for $4$ faces below, $22th$ in the list for $5$ faces below,
$Barrel_{\infty}$. 
See Figure \ref{ExampleOf21InfiniteSeries} for the first $3$ 
members (starting with 6 faces) of 
two such series: $\alpha \alpha$ and $\beta \varepsilon$.

\end{enumerate}
\end{theorem}
\proof The proof consists of taking the totally elementary 
polycycles, given in Theorem \ref{TheoremTotalElementary}, adding
a face with right number of sides, which preserves the elementarity 
in all possible ways. Then we reduce, by isomorphism of $(R,q)$-polycycles,
and obtain the list of finite elementary $(\{2,3,4,5\},3)$-polycycles 
with one hole. If a $(\{3,4,5\},3)$-polycycle has two holes and is
not a $Barrel_m$, then it is not elementary. So, it can be
obtained from another elementary polycycle with one hole less,
by the addition of one face.
It is easy to see that this cannot happen. So, we have the 
complete list of finite $(\{2,3,4,5\},3)$-polycycles.

Take now an elementary $(\{3,4,5\},3)$-polycycle $P$, which is infinite.
Remove all $3$- or $4$-gonal faces of it. The result is a
$(\{5\},3)$-polycycle $P'$, which is not necessarily elementary. 
We will now use the classification of elementary $(\{5\},3)$-polycycles
(possibly, infinite) done in \cite{DDS2}. If the infinite $(\{5\},3)$-polycycle
$Barrel_{\infty}$ appears in the decomposition, then, clearly, $P$ 
is reduced to it. If the infinite polycycle $\alpha$ appear in the
decomposition, then there are two possibilities for extending it,
indicated below.
\begin{center}
\begin{minipage}{8cm}
\centering
\epsfig{height=15mm, file=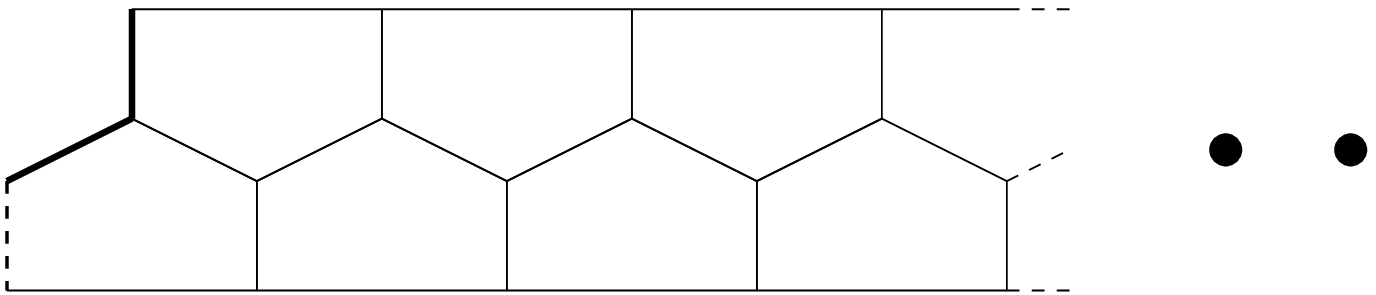}\par
\end{minipage}
\end{center}
If a $3$- or $4$-gonal face is adjacent on the dotted line, then there 
should
be another face on the over-line edges. So, in any case, there is a 
face,
adjacent on the over-line edges, and we can assume that it is a $3$- or
$4$-gonal face. Then, consideration of all possibilities to extend it,
yields $\beta$, \dots, $\mu$. Suppose now, that $P$ does not contain
any infinite $(\{5\},3)$-polycycles. Then we can find an infinite path
$f_0,\dots, f_i,\dots$ of distinct faces of $P$ in $F_1$, such that
$f_i$ is adjacent to $f_{i+1}$ and $f_{i-1}$ is not adjacent to $f_{i+1}$.
The condition on $P$ implies that an infinite number of faces are $3$-
or $4$-gons, but the condition of non-adjacency of $f_{i-1}$ with 
$f_{i+1}$
forbids $3$-gons. Take now a $4$-gon $f_i$ and assume that $f_{i-1}$
and $f_{i+1}$ are $5$-gons. The consideration of all possibilities of
extension around that face, lead us to an impossibility. If some of 
$f_{i-1}$
or $f_{i+1}$ are $4$-gons, then we have a path of $4$-gons and the case is
even simpler.  \qed

List of elementary $(\{3,4,5\},3)$-polycycles with $4$ faces:
\begin{center}
\begin{minipage}{3cm}
\centering
\epsfig{height=20mm, file=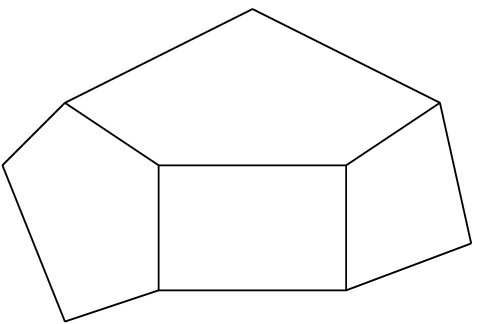}\par
$C_1$
\end{minipage}
\begin{minipage}{3cm}
\centering
\epsfig{height=20mm, file=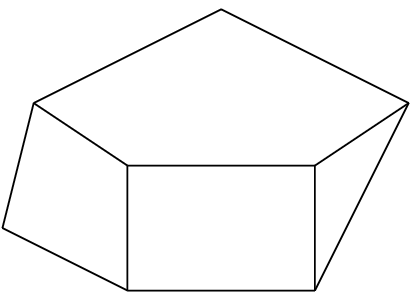}\par
$C_1$
\end{minipage}
\begin{minipage}{3cm}
\centering
\epsfig{height=20mm, file=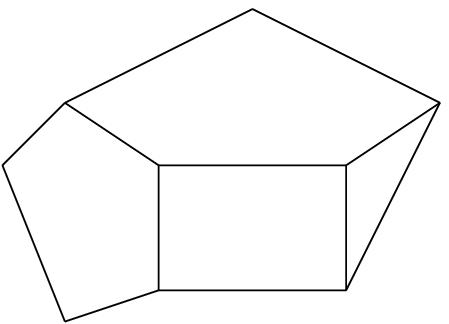}\par
$C_1$
\end{minipage}
\begin{minipage}{3cm}
\centering
\epsfig{height=20mm, file=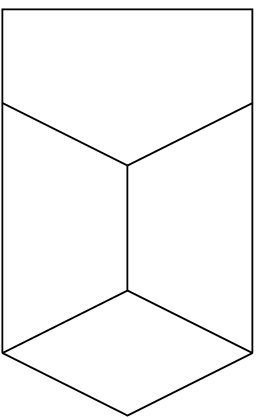}\par
$C_s$
\end{minipage}
\begin{minipage}{3cm}
\centering
\epsfig{height=20mm, file=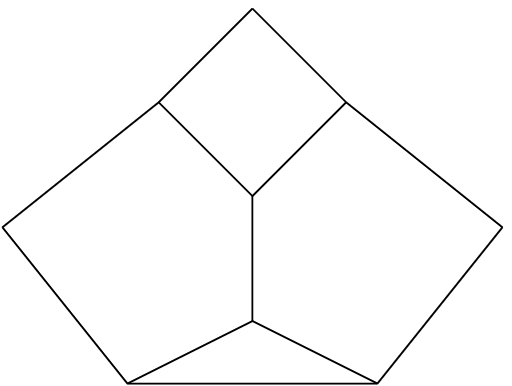}\par
$C_s$
\end{minipage}
\begin{minipage}{3cm}
\centering
\epsfig{height=20mm, file=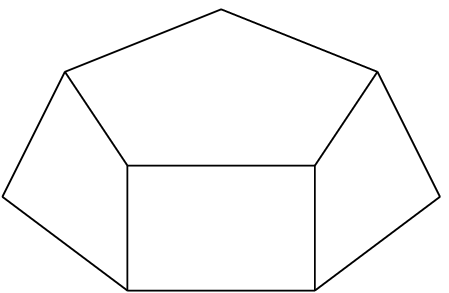}\par
$C_s$
\end{minipage}
\begin{minipage}{3cm}
\centering
\epsfig{height=20mm, file=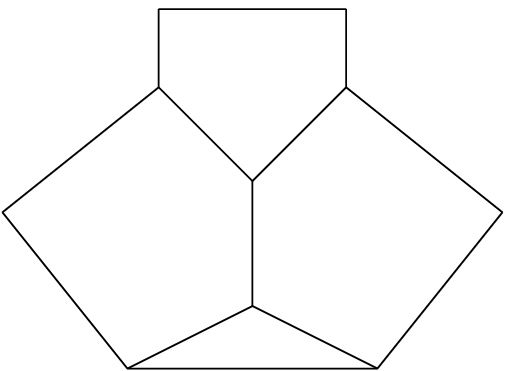}\par
$C_s$
\end{minipage}
\begin{minipage}{3cm}
\centering
\epsfig{height=20mm, file=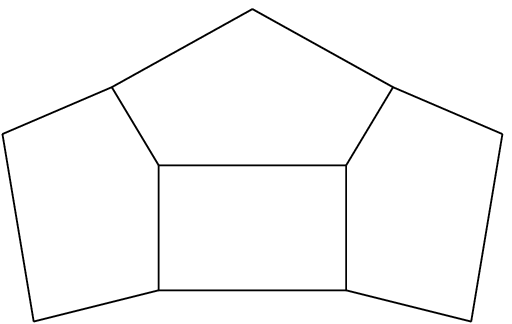}\par
$C_s$
\end{minipage}
\begin{minipage}{3cm}
\centering
\epsfig{height=20mm, file=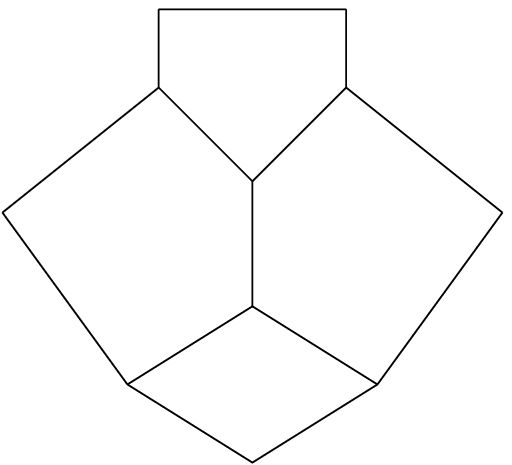}\par
$C_s$
\end{minipage}
\begin{minipage}{3cm}
\centering
\epsfig{height=20mm, file=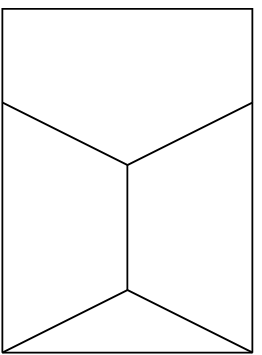}\par
$C_s$
\end{minipage}
\begin{minipage}{3cm}
\centering
\epsfig{height=20mm, file=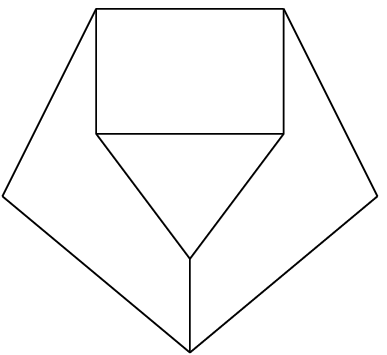}\par
$C_s$
\end{minipage}
\begin{minipage}{3cm}
\centering
\epsfig{height=20mm, file=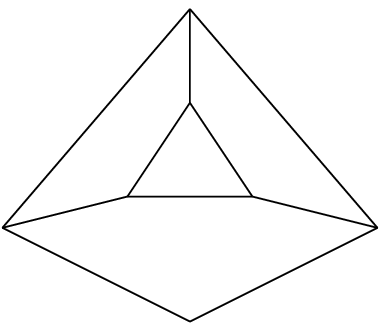}\par
$C_s$, nonext.
%PAIR1
\end{minipage}
\begin{minipage}{3cm}
\centering
\epsfig{height=20mm, file=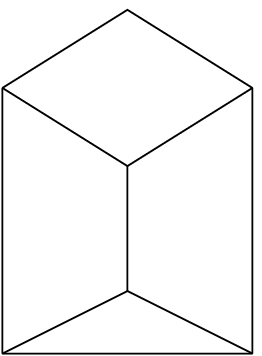}\par
$C_s$, nonext.
%PAIR1
\end{minipage}
\begin{minipage}{3cm}
\centering
\epsfig{height=20mm, file=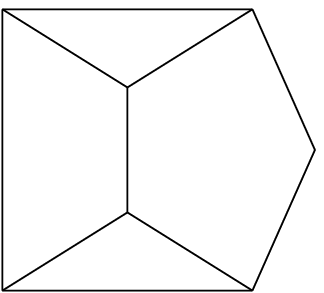}\par
$C_s$,~nonext.~$(C_{2\nu})$
\end{minipage}
\begin{minipage}{3cm}
\centering
\epsfig{height=20mm, file=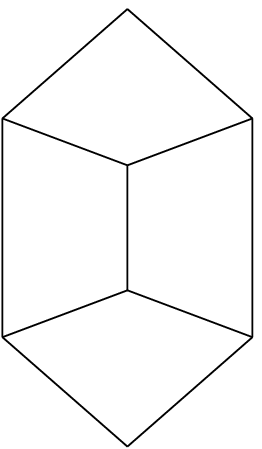}\par
$C_{2\nu}$
\end{minipage}
\begin{minipage}{3cm}
\centering
\epsfig{height=20mm, file=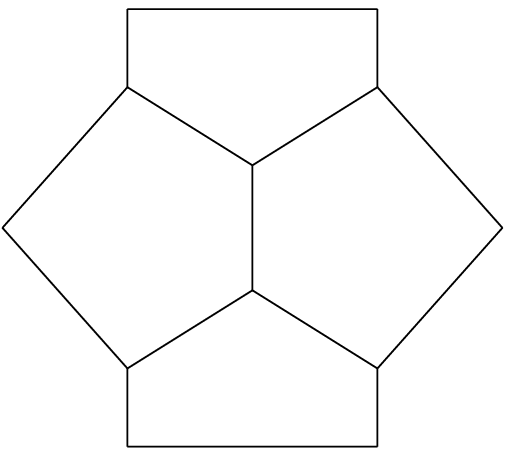}\par
$C_{2\nu}$
\end{minipage}
\begin{minipage}{3cm}
\centering
\epsfig{height=20mm, file=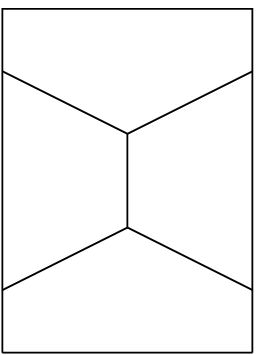}\par
$C_{2\nu}$
\end{minipage}
\begin{minipage}{3cm}
\centering
\epsfig{height=20mm, file=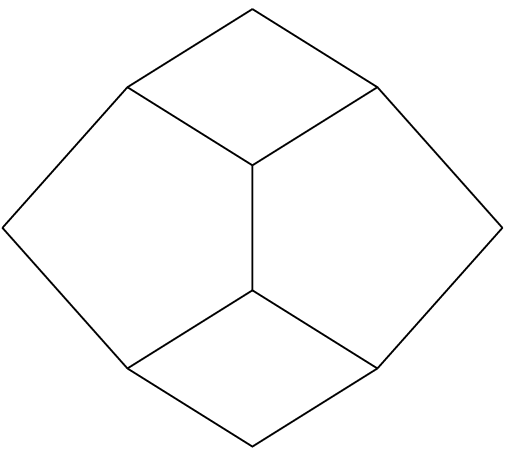}\par
$C_{2\nu}$
\end{minipage}
\begin{minipage}{3cm}
\centering
\epsfig{height=20mm, file=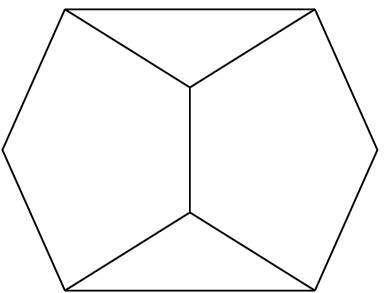}\par
$C_{2\nu}$
\end{minipage}
\begin{minipage}{3cm}
\centering
\epsfig{height=20mm, file=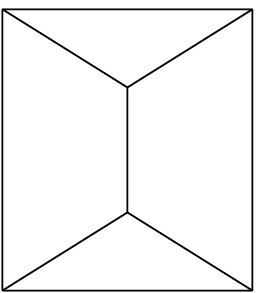}\par
$C_{2\nu}$, nonext.
%PAIR2
\end{minipage}
\begin{minipage}{3cm}
\centering
\epsfig{height=20mm, file=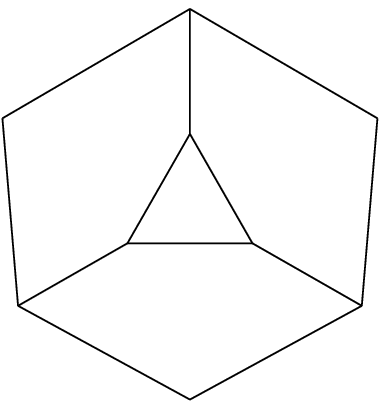}\par
$C_{3\nu}$
\end{minipage}
\begin{minipage}{3cm}
\centering
\epsfig{height=20mm, file=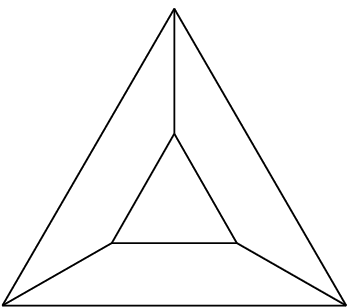}\par
$C_{3\nu}$, nonext.
%PAIR2
\end{minipage}
\begin{minipage}{3cm}
\centering
\epsfig{height=20mm, file=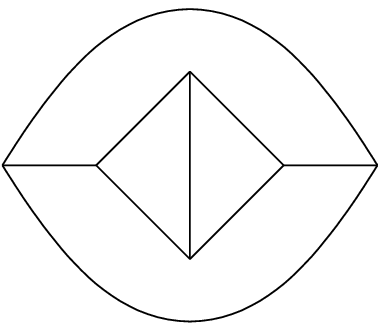}\par
$C_{2\nu}$, nonext.
%PAIR3
\end{minipage}

\end{center}
List of elementary $(\{3,4,5\},3)$-polycycles with $5$ faces:
\begin{center}
\begin{minipage}{3cm}
\centering
\epsfig{height=20mm, file=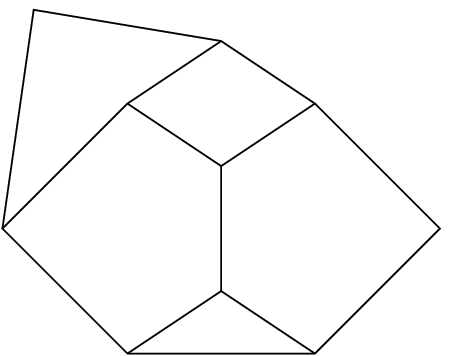}\par
$C_1$
\end{minipage}
\begin{minipage}{3cm}
\centering
\epsfig{height=20mm, file=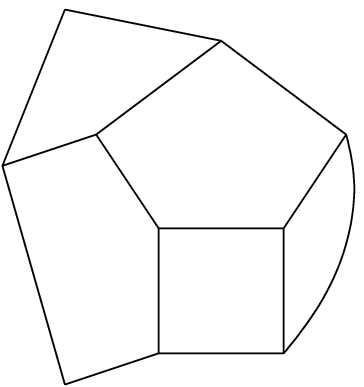}\par
$C_1$
\end{minipage}
\begin{minipage}{3cm}
\centering
\epsfig{height=20mm, file=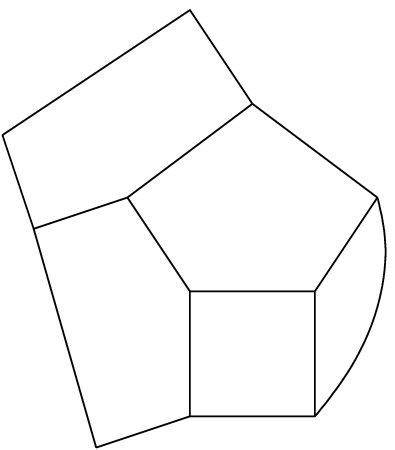}\par
$C_1$
\end{minipage}
\begin{minipage}{3cm}
\centering
\epsfig{height=20mm, file=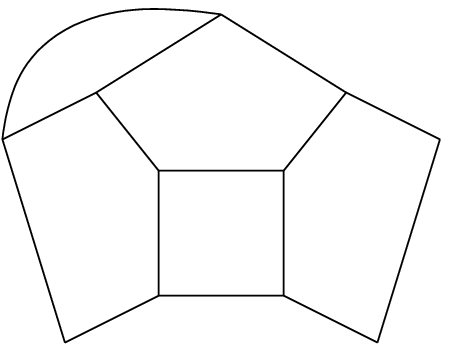}\par
$C_1$
\end{minipage}
\begin{minipage}{3cm}
\centering
\epsfig{height=20mm, file=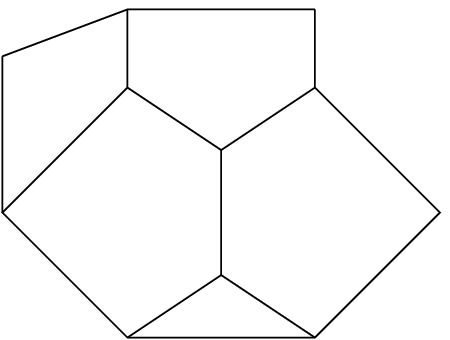}\par
$C_1$
\end{minipage}
\begin{minipage}{3cm}
\centering
\epsfig{height=20mm, file=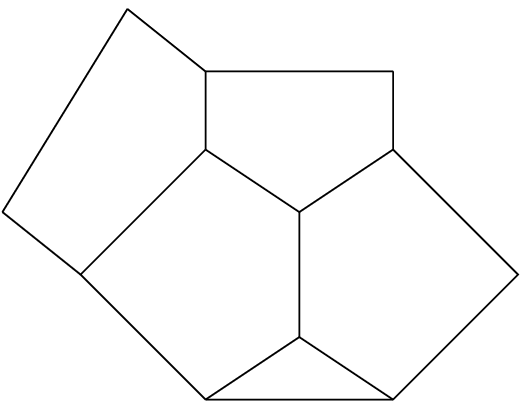}\par
$C_1$
\end{minipage}
\begin{minipage}{3cm}
\centering
\epsfig{height=20mm, file=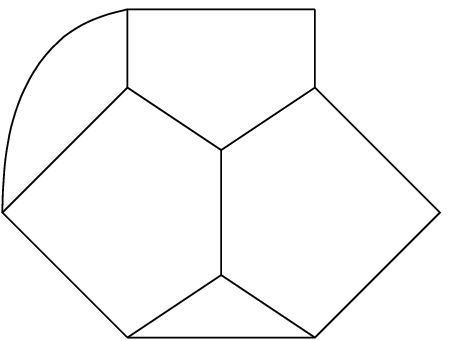}\par
$C_1$
\end{minipage}
\begin{minipage}{3cm}
\centering
\epsfig{height=20mm, file=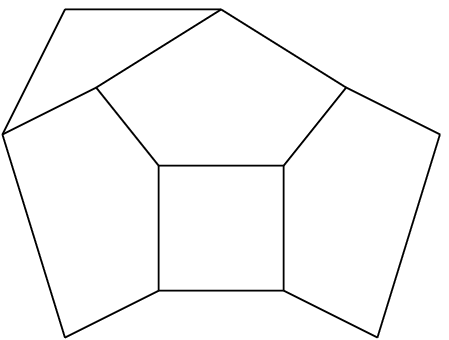}\par
$C_1$
\end{minipage}
\begin{minipage}{3cm}
\centering
\epsfig{height=20mm, file=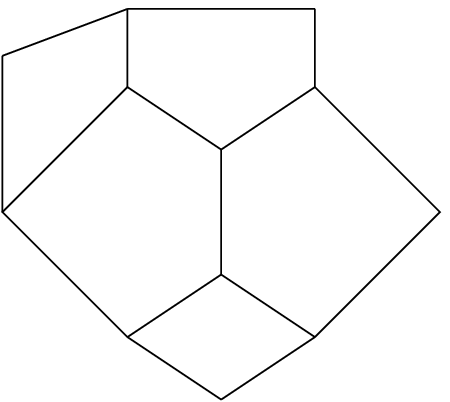}\par
$C_1$
\end{minipage}
\begin{minipage}{3cm}
\centering
\epsfig{height=20mm, file=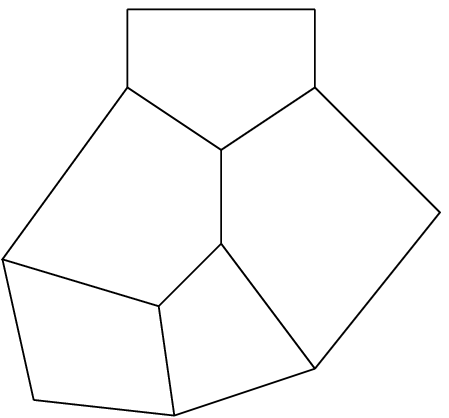}\par
$C_1$
\end{minipage}
\begin{minipage}{3cm}
\centering
\epsfig{height=20mm, file=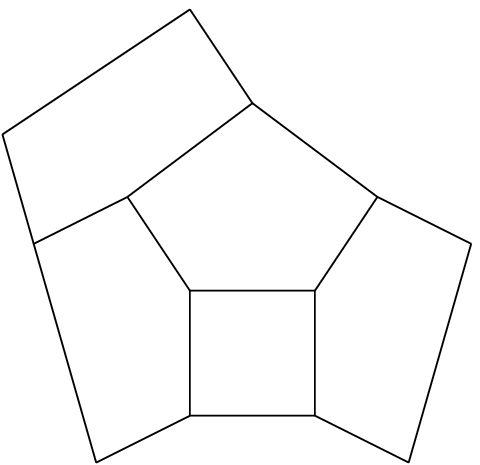}\par
$C_1$
\end{minipage}
\begin{minipage}{3cm}
\centering
\epsfig{height=20mm, file=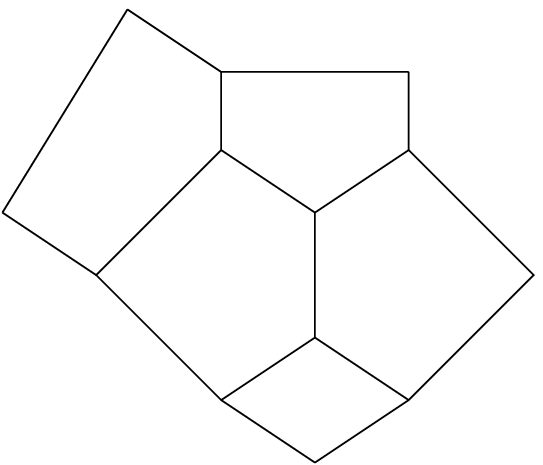}\par
$C_1$
\end{minipage}
\begin{minipage}{3cm}
\centering
\epsfig{height=20mm, file=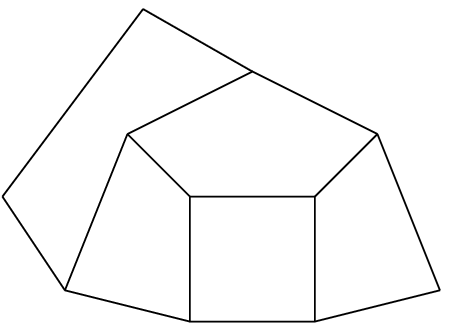}\par
$C_1$
\end{minipage}
\begin{minipage}{3cm}
\centering
\epsfig{height=20mm, file=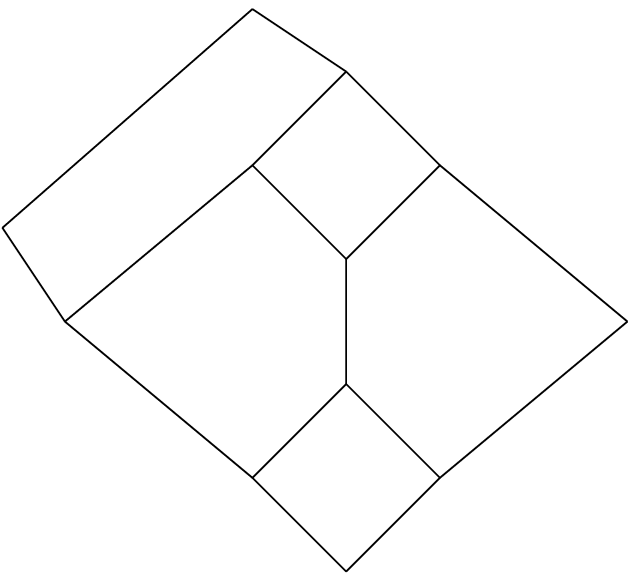}\par
$C_1$
\end{minipage}
\begin{minipage}{3cm}
\centering
\epsfig{height=20mm, file=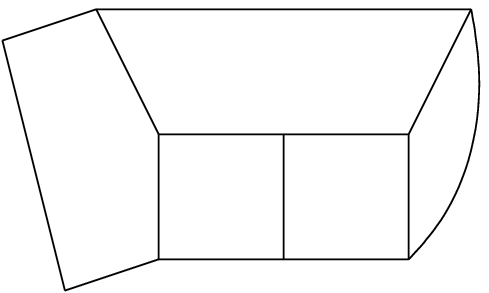}\par
$C_1$
\end{minipage}
\begin{minipage}{3cm}
\centering
\epsfig{height=18mm, file=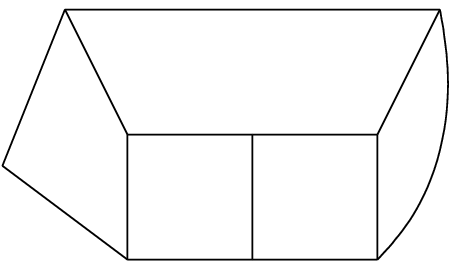}\par
$C_1$, nonext.
\end{minipage}
\begin{minipage}{3cm}
\centering
\epsfig{height=18mm, file=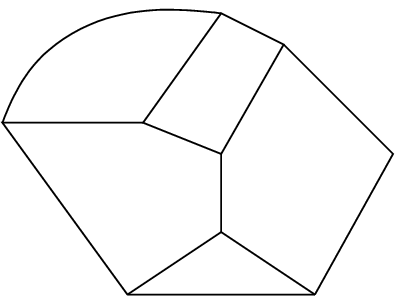}\par
$C_1$, nonext.
\end{minipage}
\begin{minipage}{3cm}
\centering
\epsfig{height=20mm, file=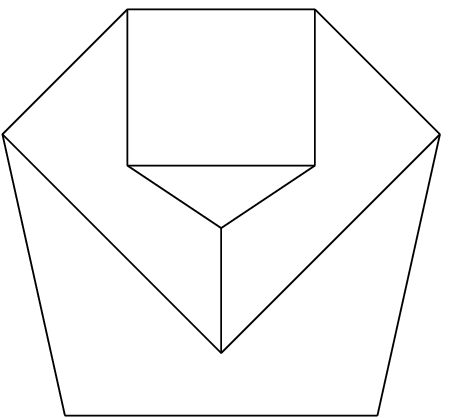}\par
$C_s$
\end{minipage}
\begin{minipage}{3cm}
\centering
\epsfig{height=20mm, file=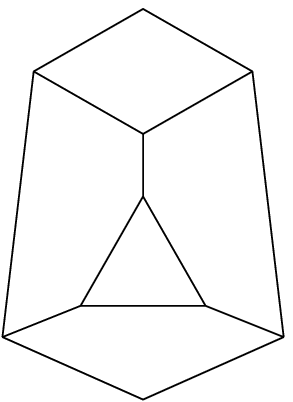}\par
$C_s$
\end{minipage}
\begin{minipage}{3cm}
\centering
\epsfig{height=20mm, file=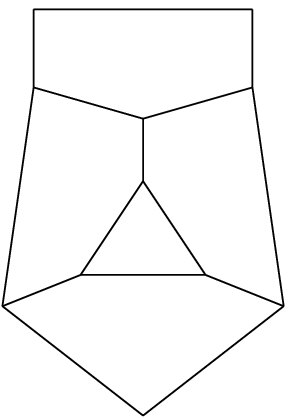}\par
$C_s$
\end{minipage}
\begin{minipage}{3cm}
\centering
\epsfig{height=20mm, file=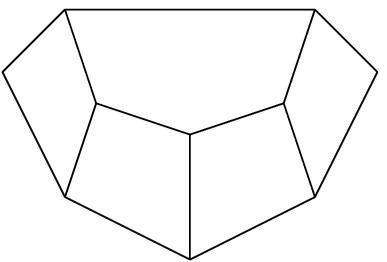}\par
$C_s$
\end{minipage}
\begin{minipage}{3cm}
\centering
\epsfig{height=14mm, file=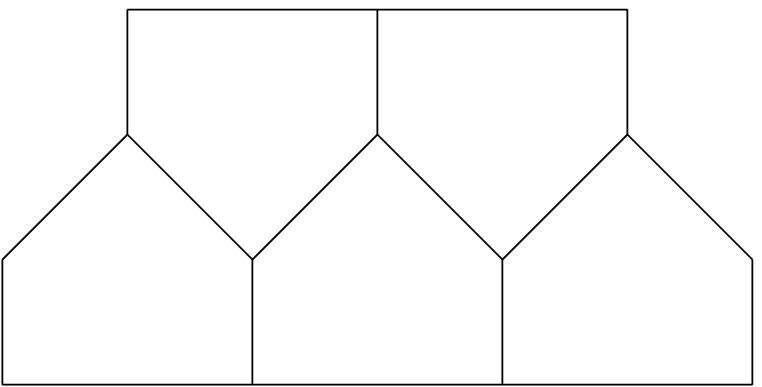}\par
$C_s$
\end{minipage}
\begin{minipage}{3cm}
\centering
\epsfig{height=20mm, file=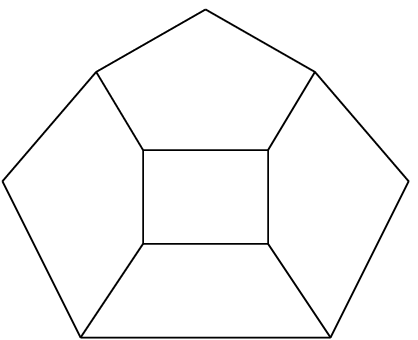}\par
$C_s$
\end{minipage}
\begin{minipage}{3cm}
\centering
\epsfig{height=17mm, file=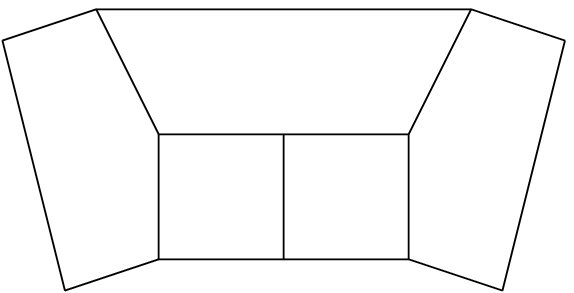}\par
$C_s$
\end{minipage}
\begin{minipage}{3cm}
\centering
\epsfig{height=20mm, file=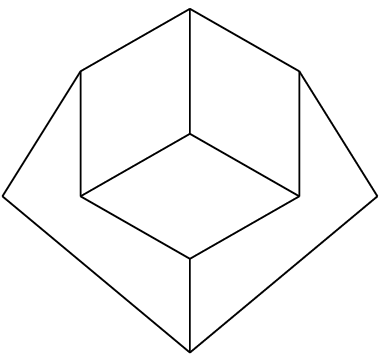}\par
$C_s$
\end{minipage}
\begin{minipage}{3cm}
\centering
\epsfig{height=20mm, file=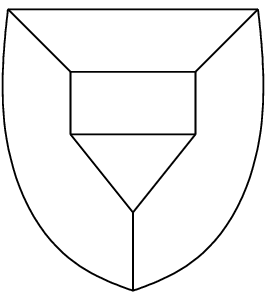}\par
$C_s$, nonext.
%PAIR1
\end{minipage}
\begin{minipage}{3cm}
\centering
\epsfig{height=20mm, file=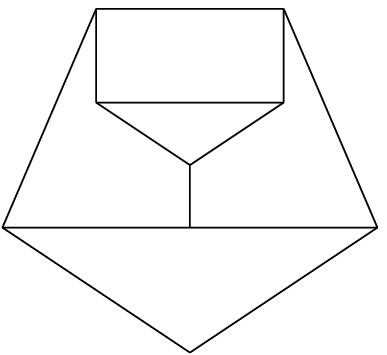}\par
$C_s$, nonext.
\end{minipage}
\begin{minipage}{3cm}
\centering
\epsfig{height=20mm, file=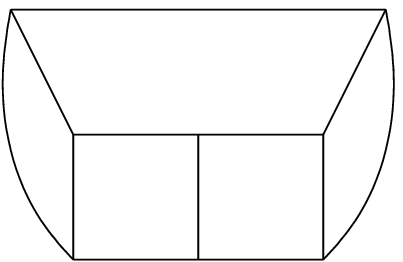}\par
$C_s$, nonext.
%PAIR1
\end{minipage}
\begin{minipage}{3cm}
\centering
\epsfig{height=20mm, file=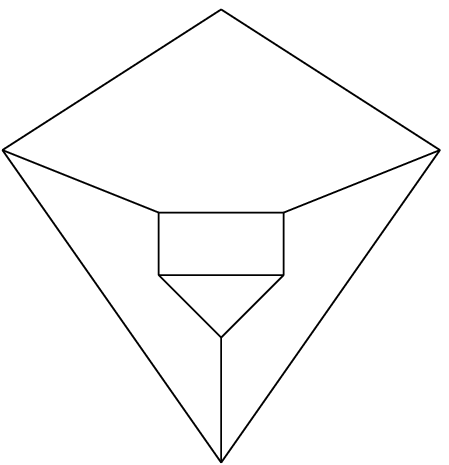}\par
$C_s$, nonext.
%PAIR2
\end{minipage}
\begin{minipage}{3cm}
\centering
\epsfig{height=20mm, file=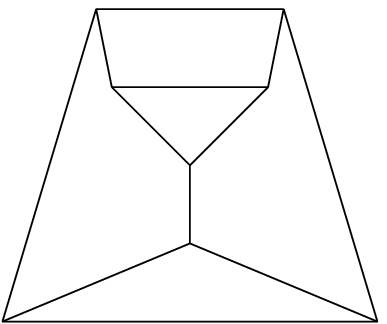}\par
$C_s$, nonext.
%PAIR1
\end{minipage}
\begin{minipage}{3cm}
\centering
\epsfig{height=20mm, file=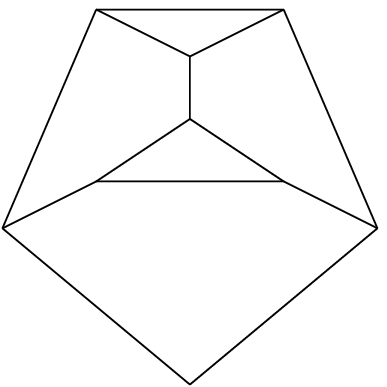}\par
$C_s$,~nonext.~$(C_{2\nu})$
%PAIR2
\end{minipage}
\begin{minipage}{3cm}
\centering
\epsfig{height=20mm, file=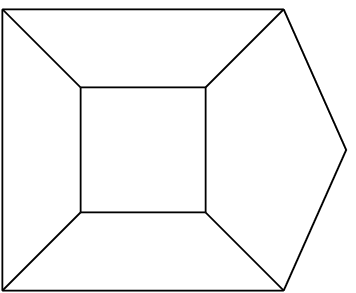}\par
$C_s$,~nonext.~$(C_{2\nu})$
\end{minipage}
\begin{minipage}{3cm}
\centering
\epsfig{height=20mm, file=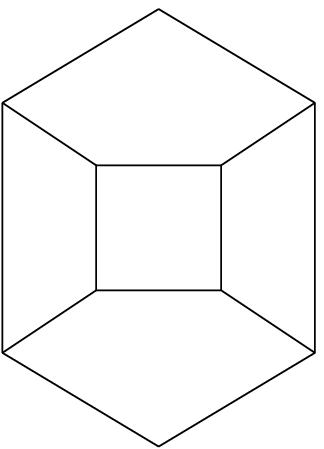}\par
$C_{2\nu}$
\end{minipage}
\begin{minipage}{3cm}
\centering
\epsfig{height=20mm, file=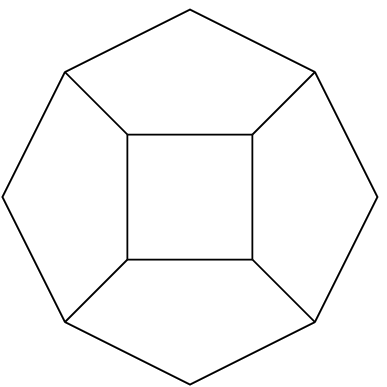}\par
$C_{4\nu}$
\end{minipage}
\begin{minipage}{3cm}
\centering
\epsfig{height=20mm, file=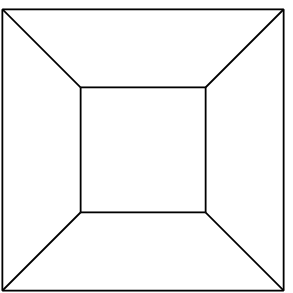}\par
$C_{4\nu}$,~nonext.~$(O_h)$
\end{minipage}

\end{center}

List of sporadic elementary $(\{3,4,5\},3)$-polycycles with $6$ faces:
\begin{center}
\begin{minipage}{3cm}
\centering
\epsfig{height=20mm, file=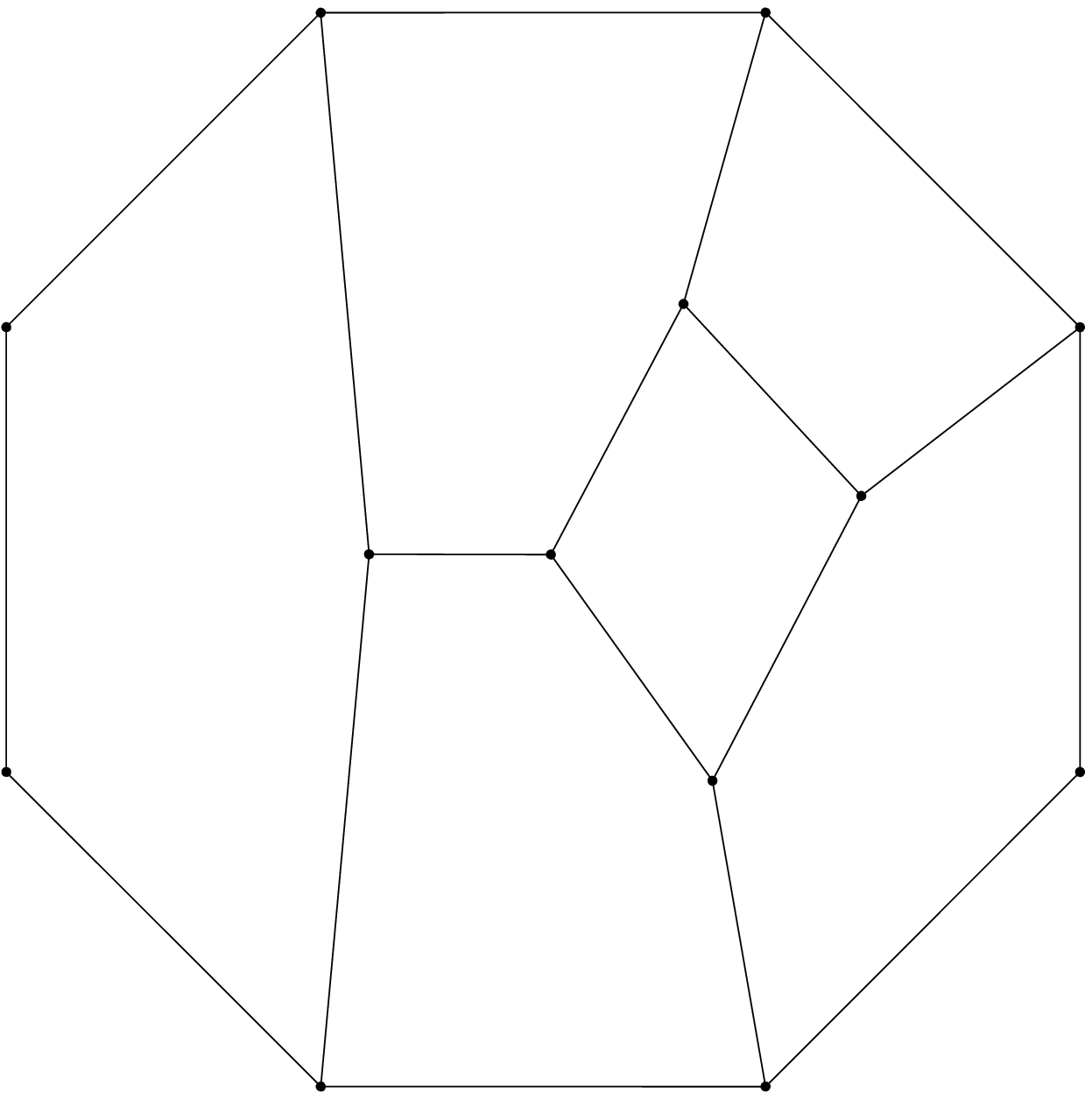}\par
$C_1$
\end{minipage}
\begin{minipage}{3cm}
\centering
\epsfig{height=20mm, file=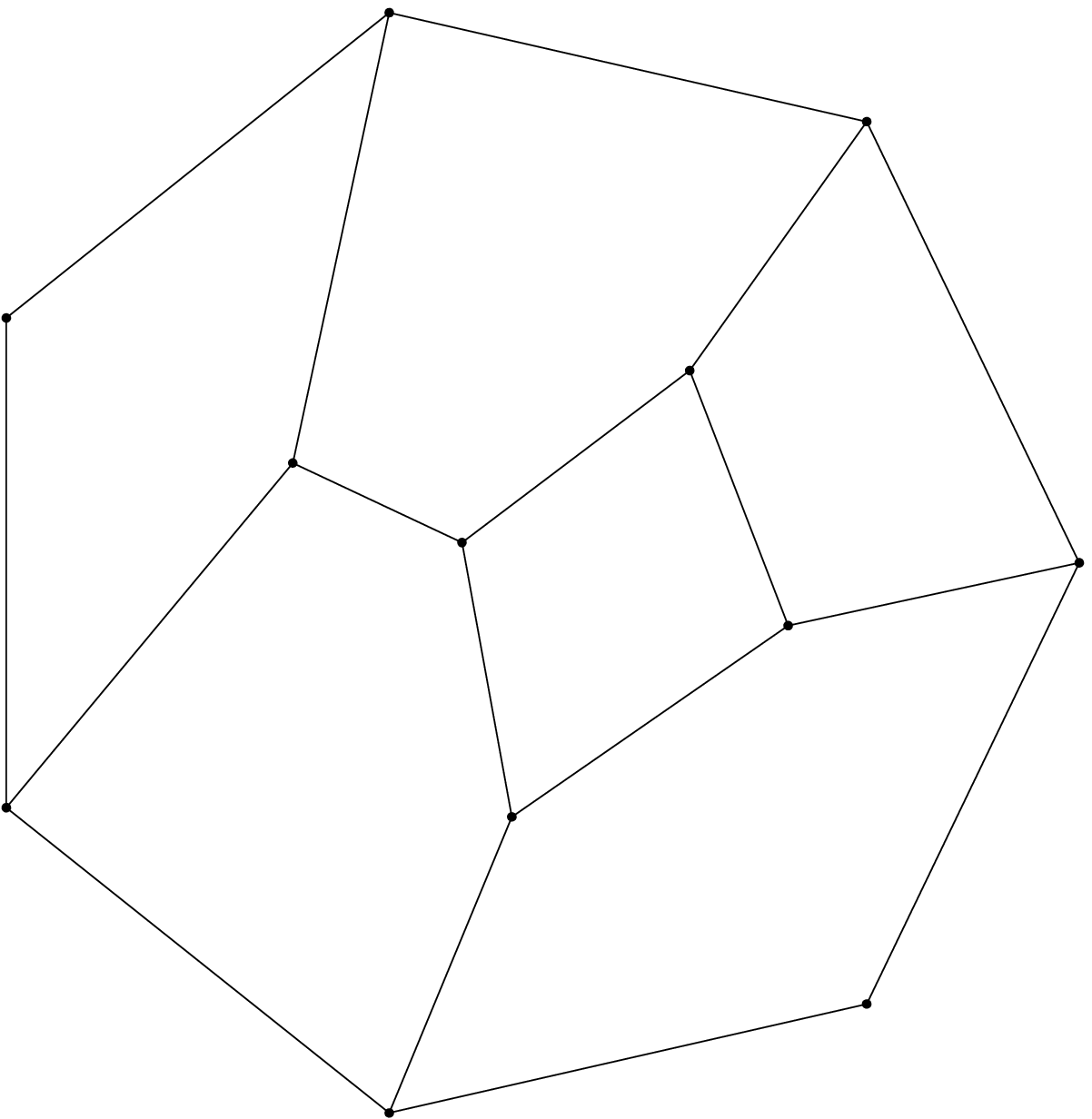}\par
$C_1$
\end{minipage}
\begin{minipage}{3cm}
\centering
\epsfig{height=20mm, file=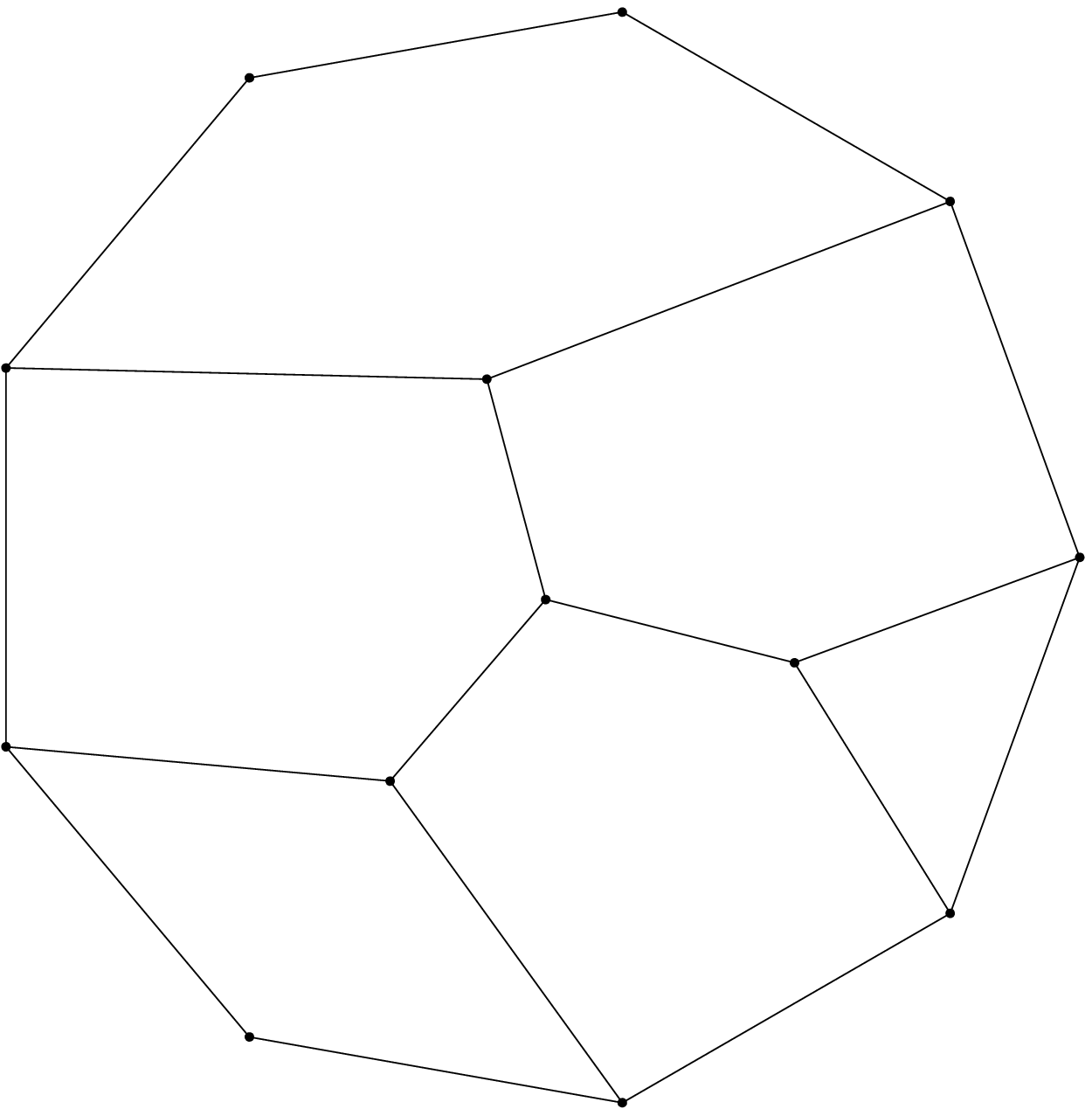}\par
$C_1$
\end{minipage}
\begin{minipage}{3cm}
\centering
\epsfig{height=20mm, file=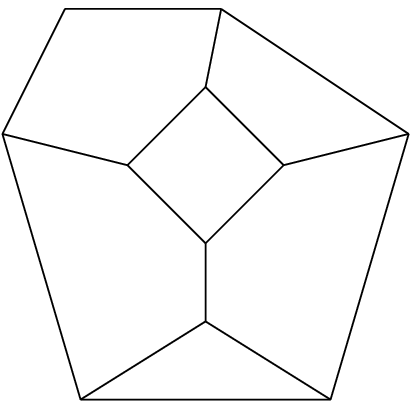}\par
$C_1$, nonext.
\end{minipage}
\begin{minipage}{3cm}
\centering
\epsfig{height=20mm, file=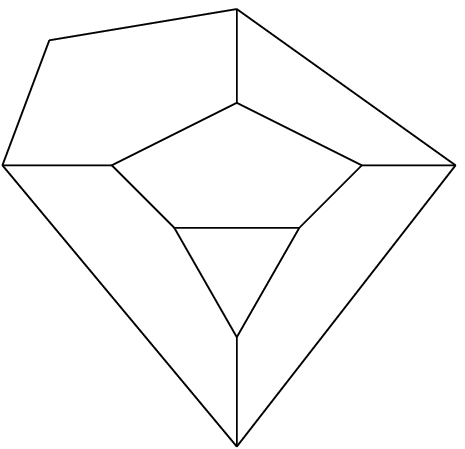}\par
$C_1$, nonext.
\end{minipage}
\begin{minipage}{3cm}
\centering
\epsfig{height=20mm, file=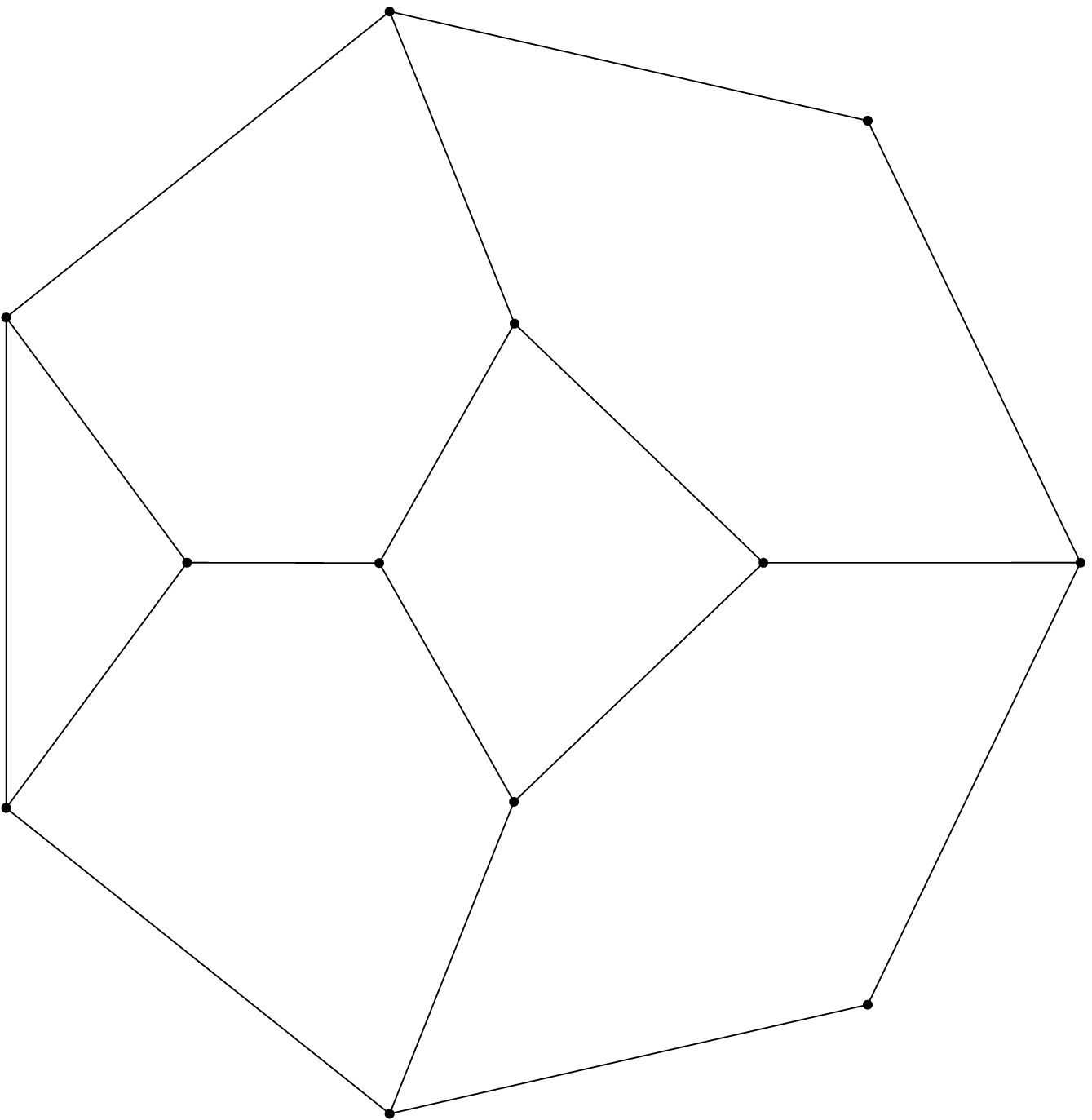}\par
$C_s$
\end{minipage}
\begin{minipage}{3cm}
\centering
\epsfig{height=20mm, file=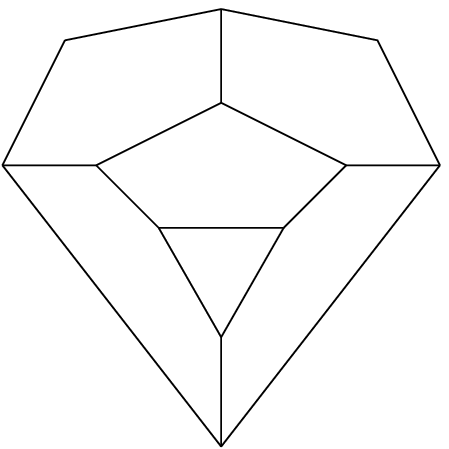}\par
$C_s$
\end{minipage}
\begin{minipage}{3cm}
\centering
\epsfig{height=20mm, file=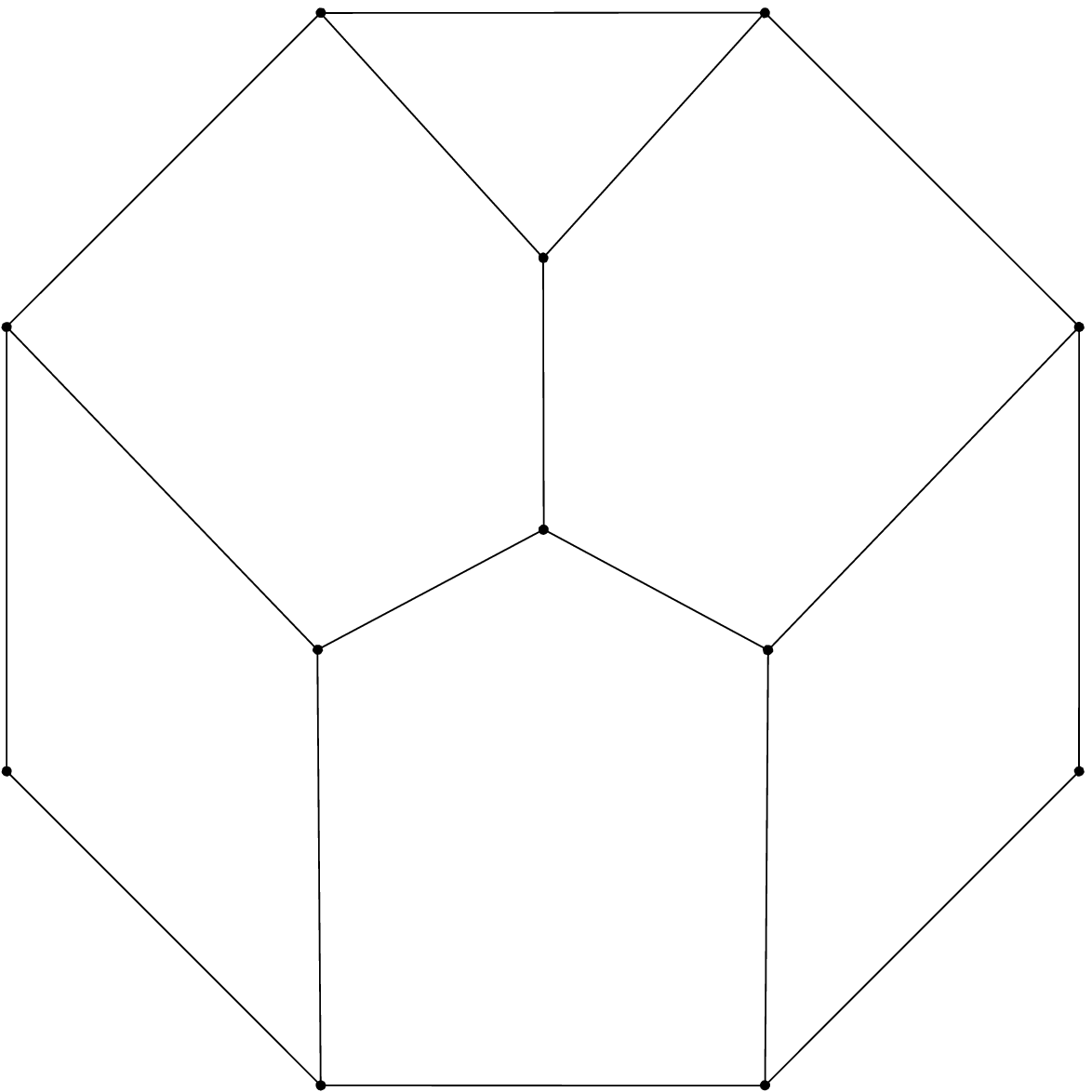}\par
$C_s$
\end{minipage}
\begin{minipage}{3cm}
\centering
\epsfig{height=20mm, file=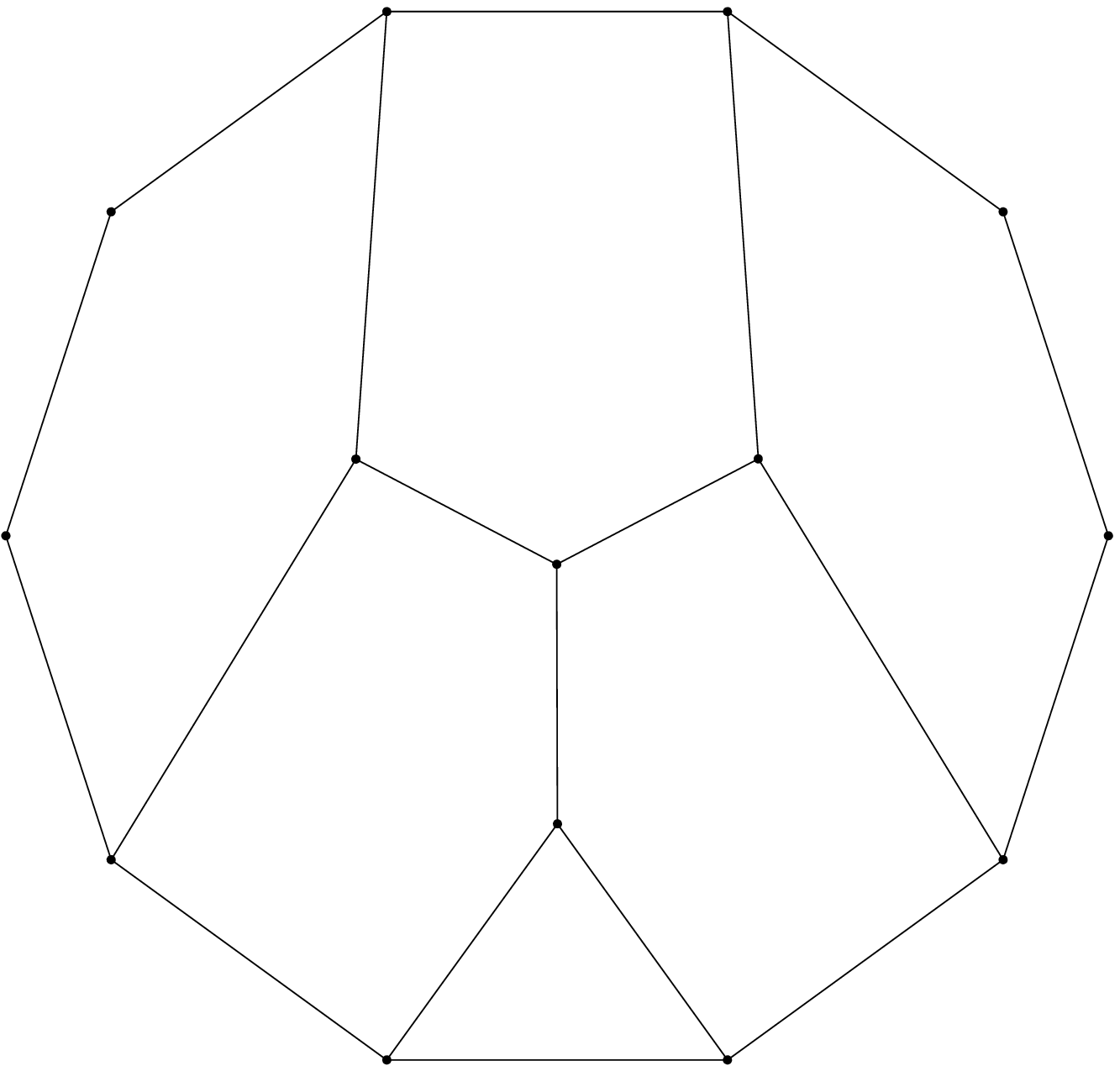}\par
$C_s$
\end{minipage}
\begin{minipage}{3cm}
\centering
\epsfig{height=20mm, file=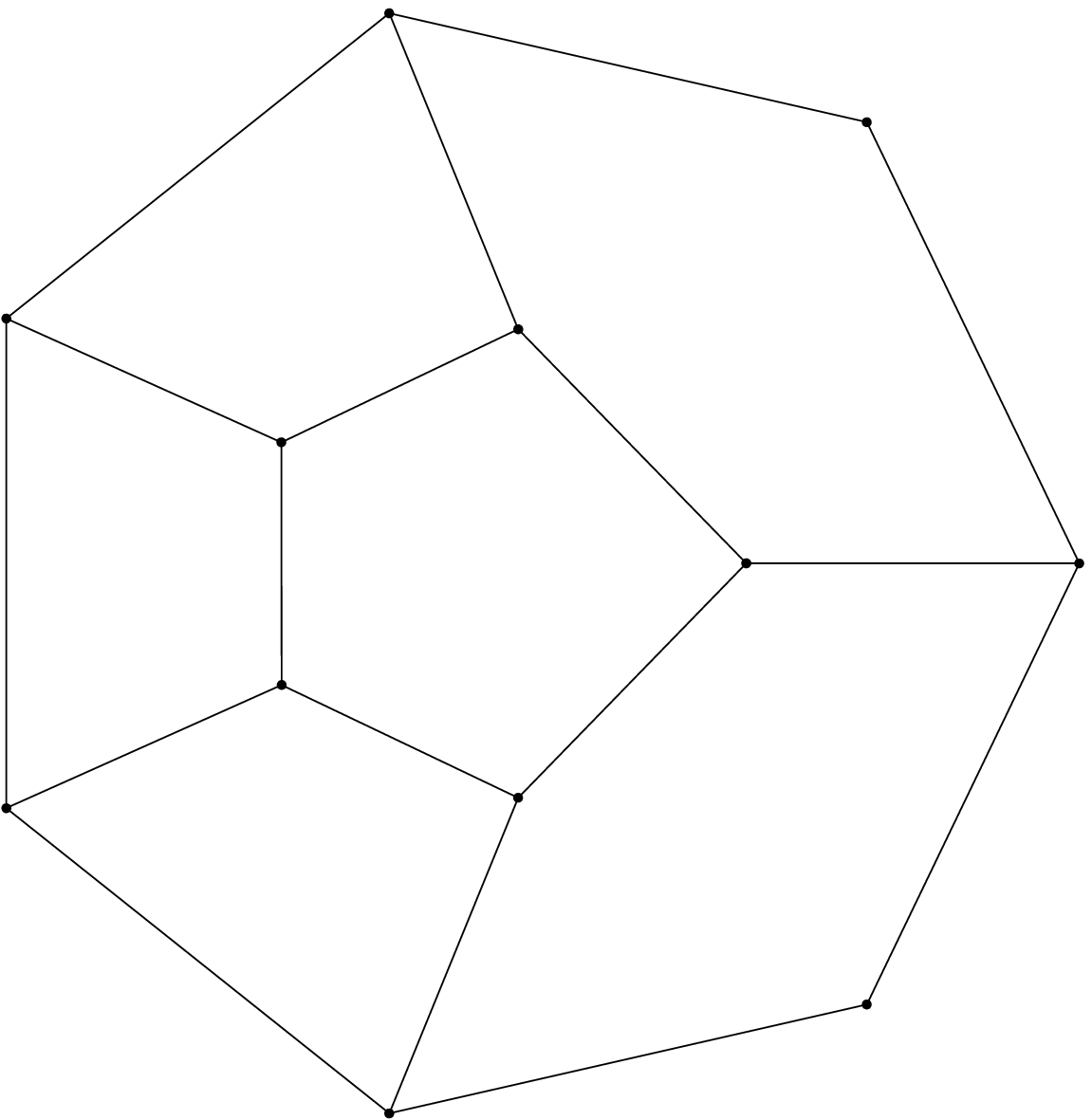}\par
$C_s$
\end{minipage}
\begin{minipage}{3cm}
\centering
\epsfig{height=20mm, file=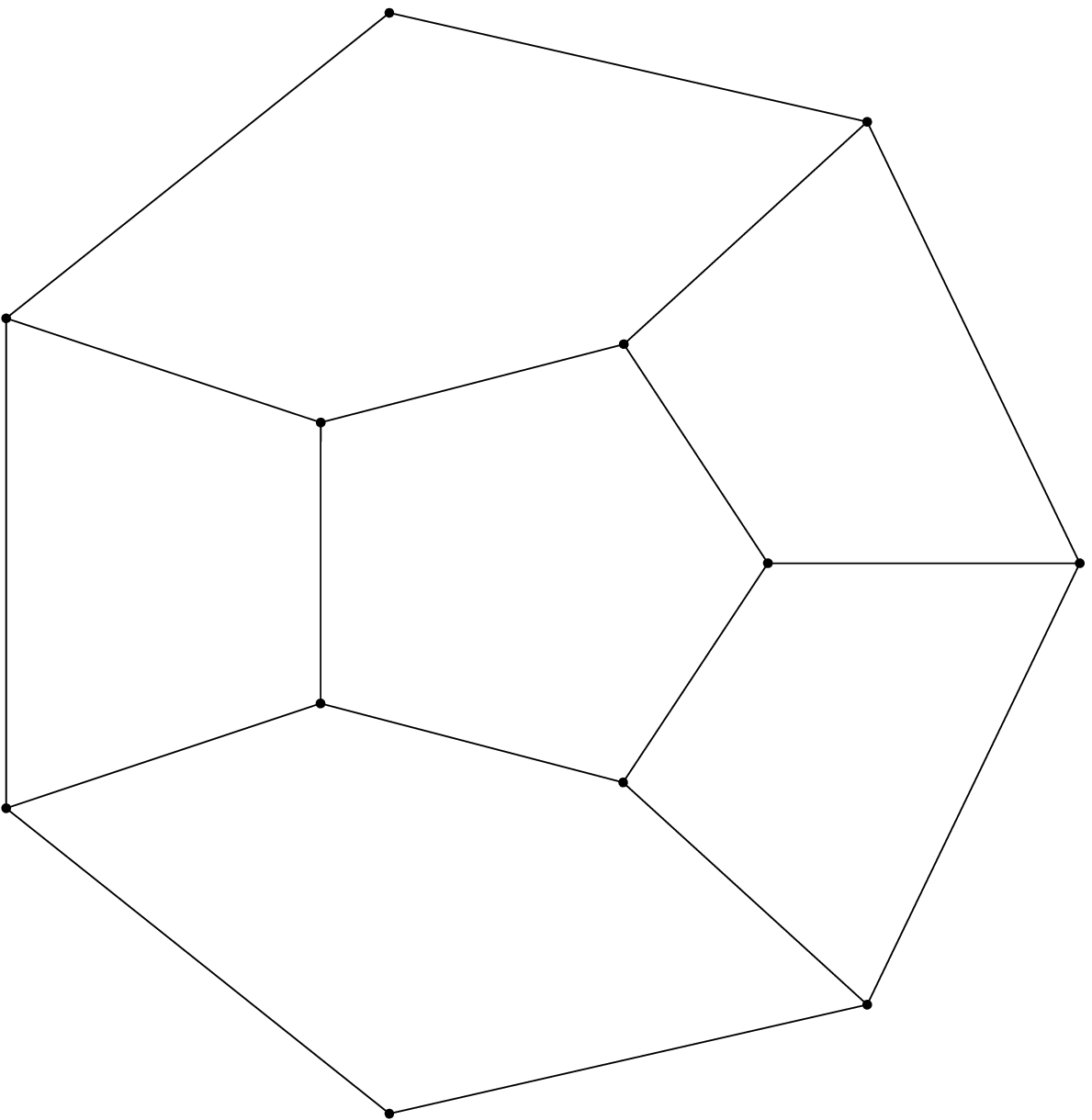}\par
$C_s$
\end{minipage}
\begin{minipage}{3cm}
\centering
\epsfig{height=20mm, file=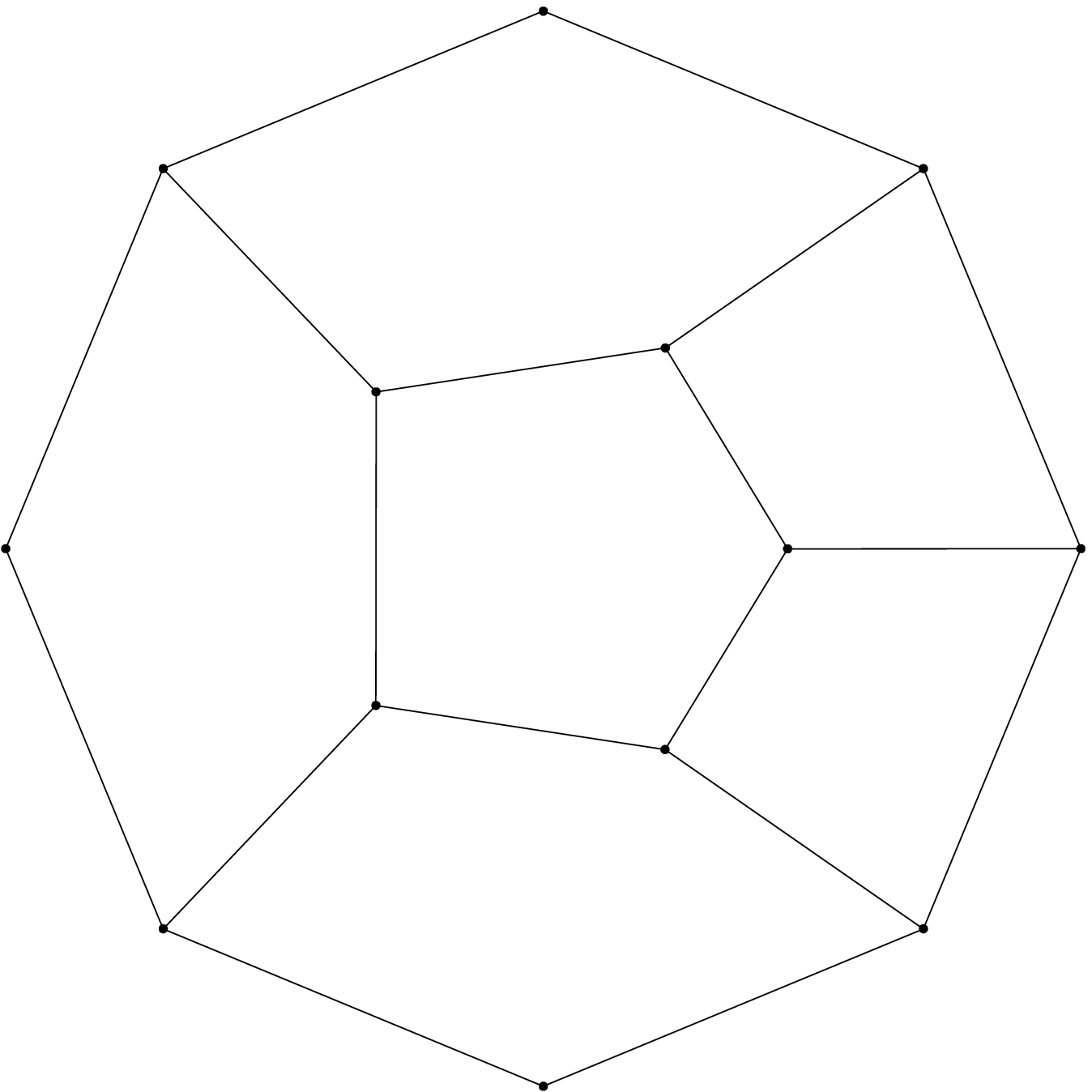}\par
$C_s$
\end{minipage}
\begin{minipage}{3cm}
\centering
\epsfig{height=20mm, file=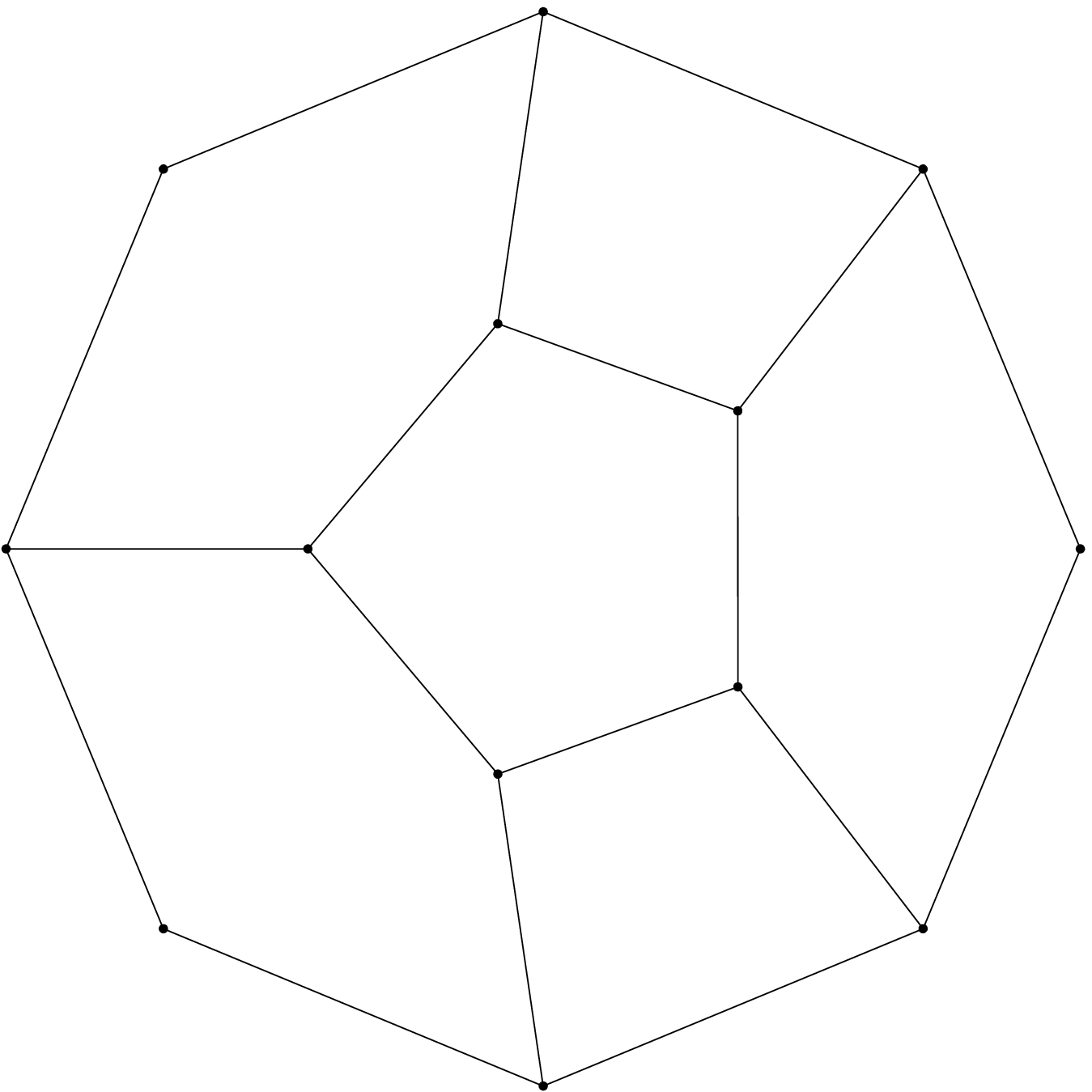}\par
$C_s$
\end{minipage}
\begin{minipage}{3cm}
\centering
\epsfig{height=20mm, file=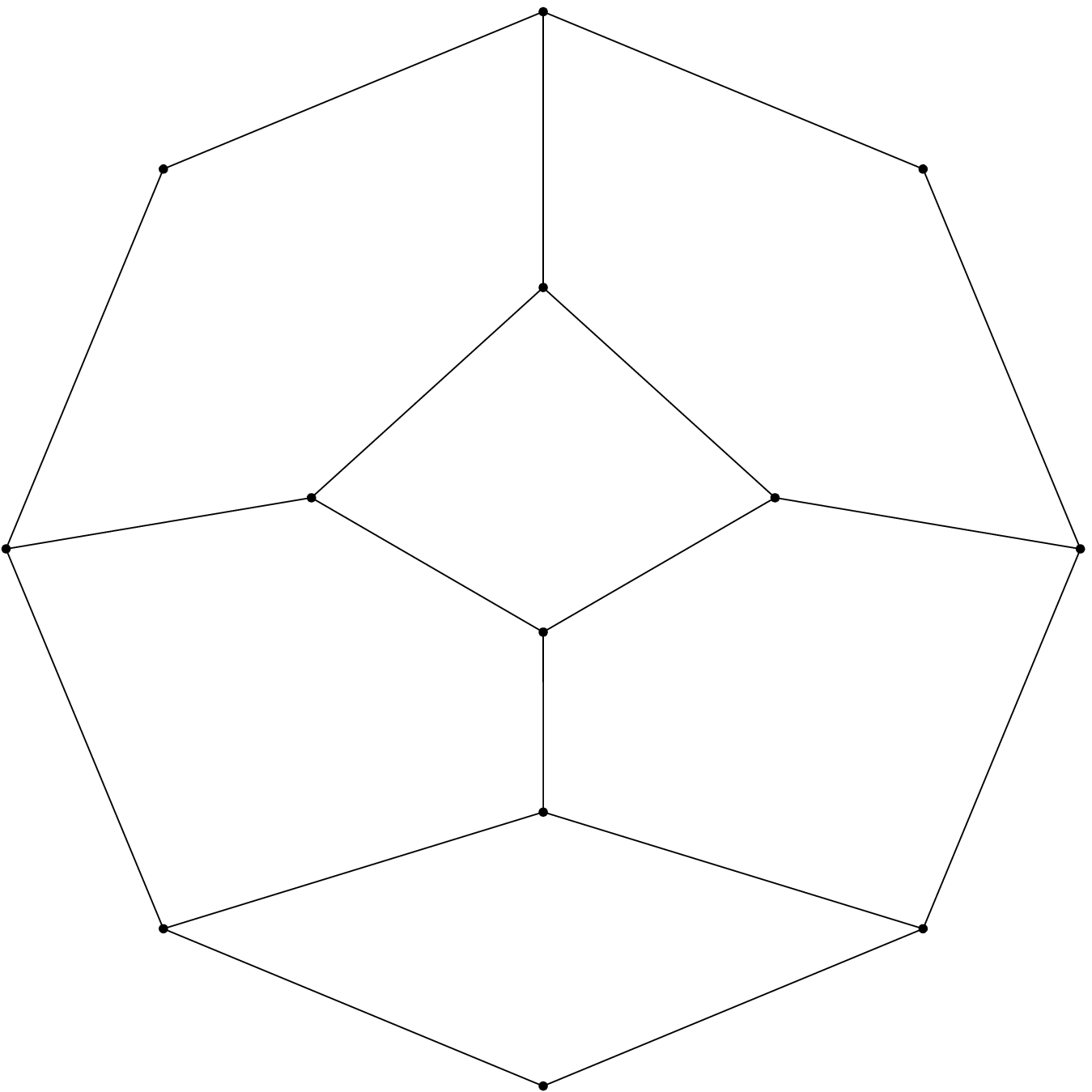}\par
$C_s$
\end{minipage}
\begin{minipage}{3cm}
\centering
\epsfig{height=20mm, file=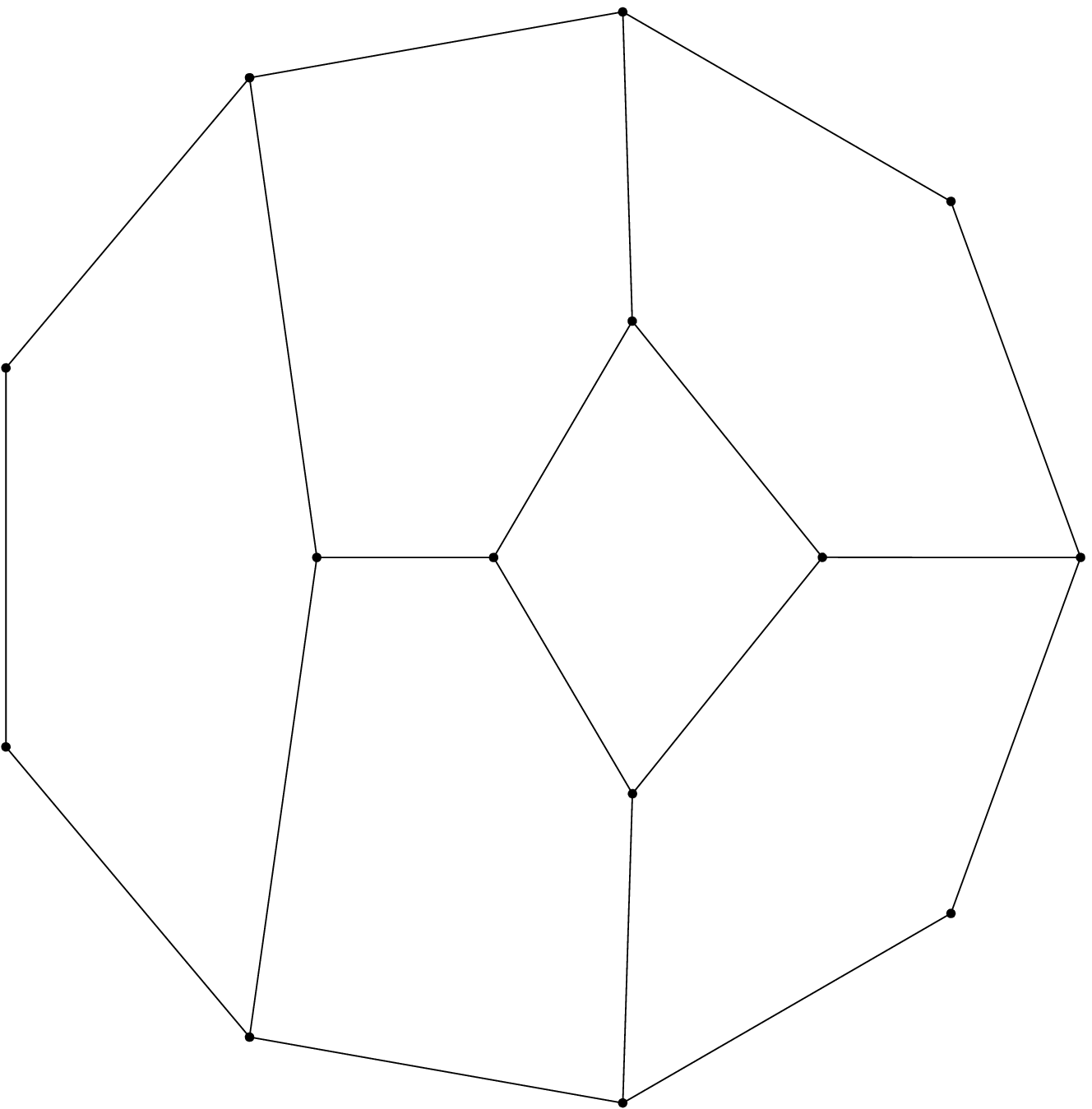}\par
$C_s$
\end{minipage}
\begin{minipage}{3cm}
\centering
\epsfig{height=20mm, file=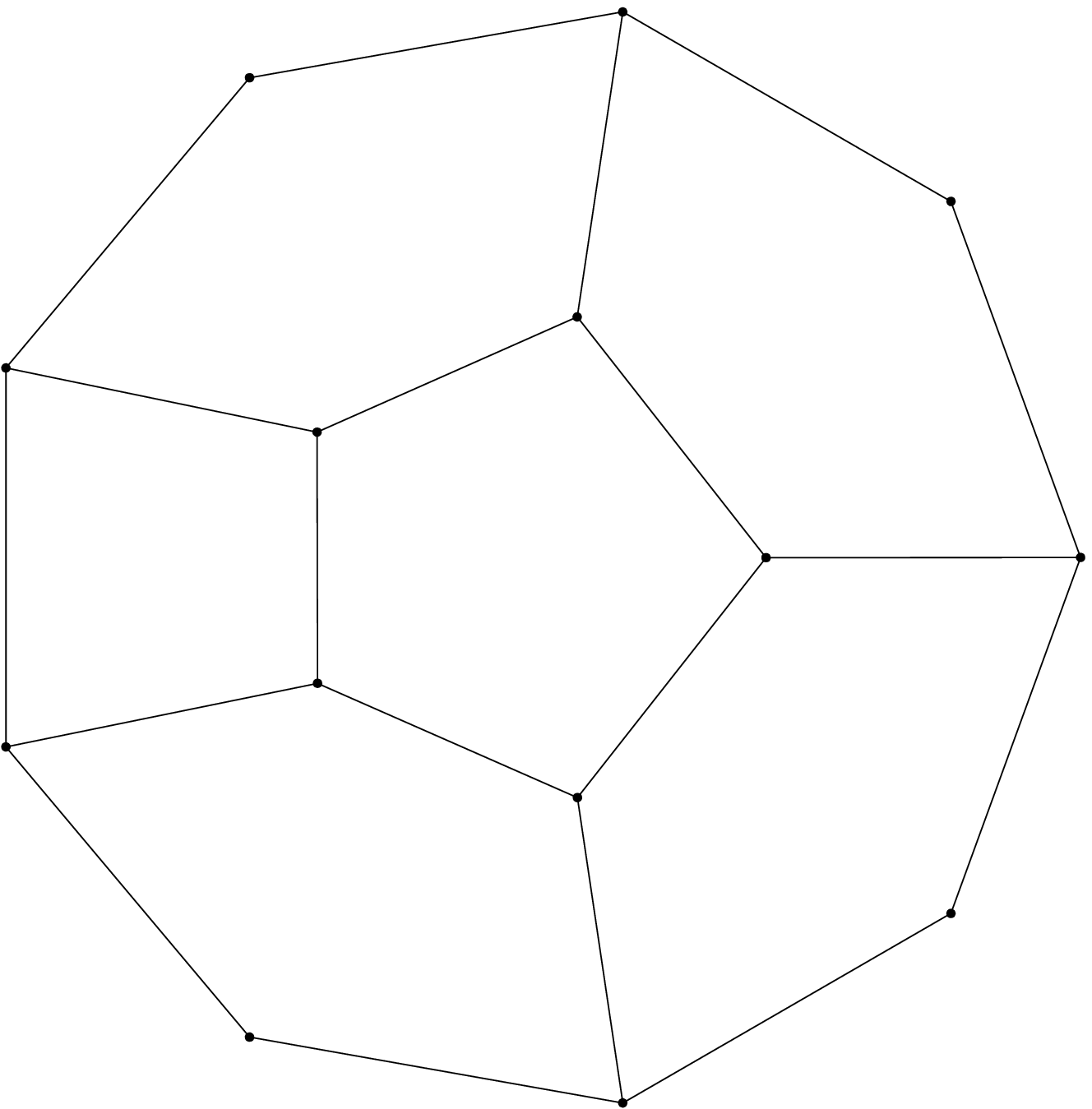}\par
$C_s$
\end{minipage}
\begin{minipage}{3cm}
\centering
\epsfig{height=20mm, file=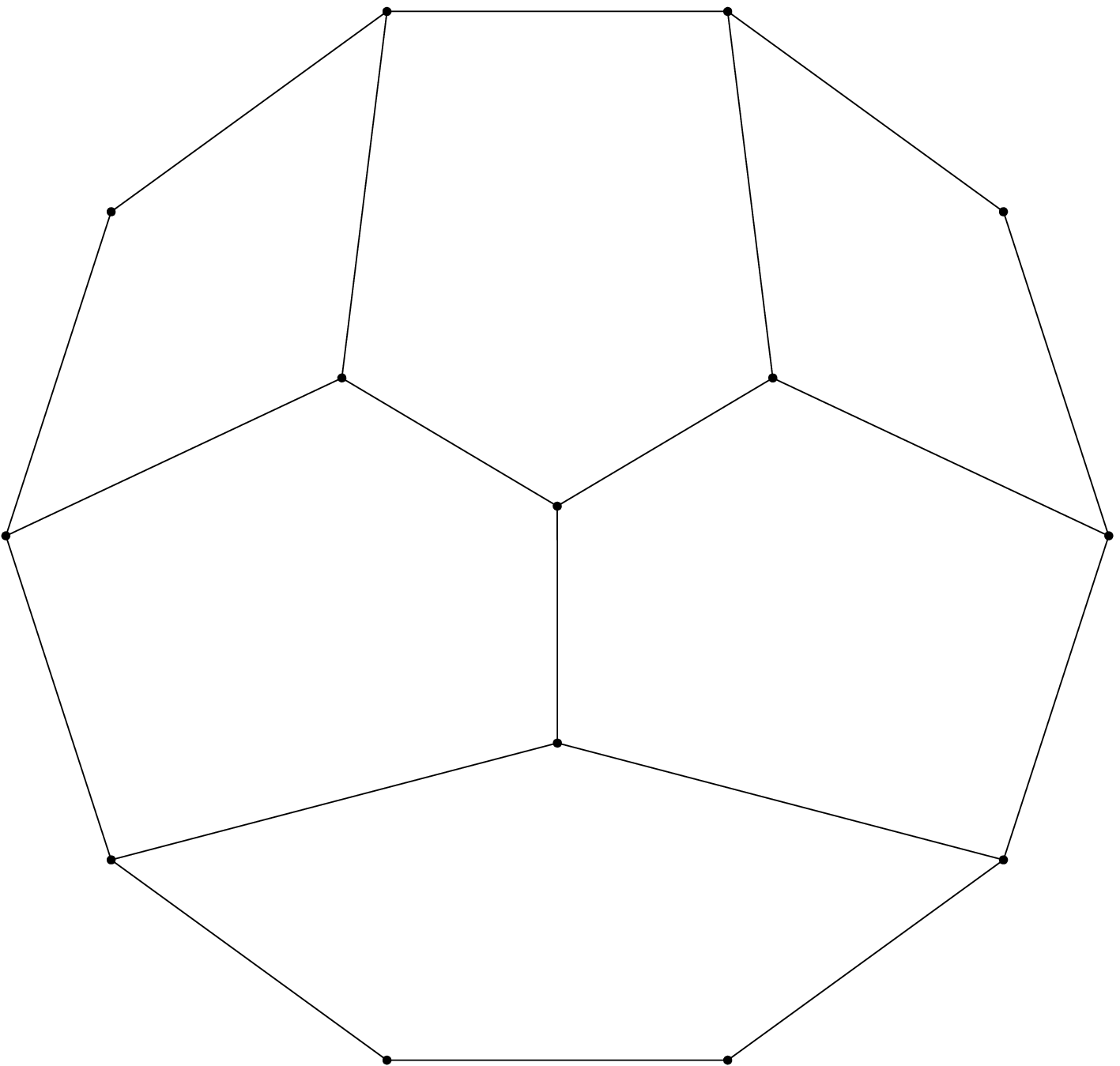}\par
$C_s$
\end{minipage}
\begin{minipage}{3cm}
\centering
\epsfig{height=20mm, file=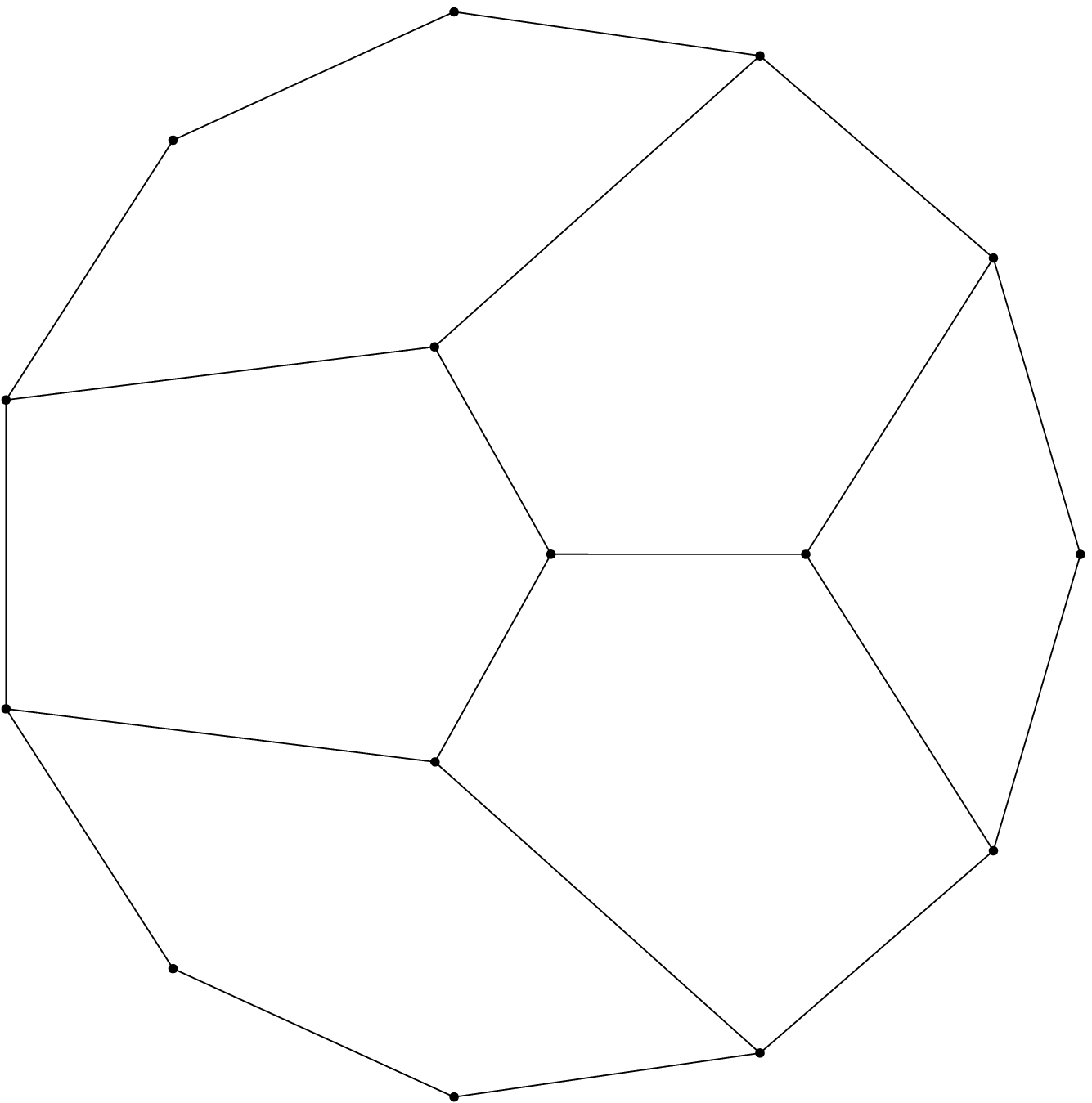}\par
$C_s$
\end{minipage}
\begin{minipage}{3cm}
\centering
\epsfig{height=20mm, file=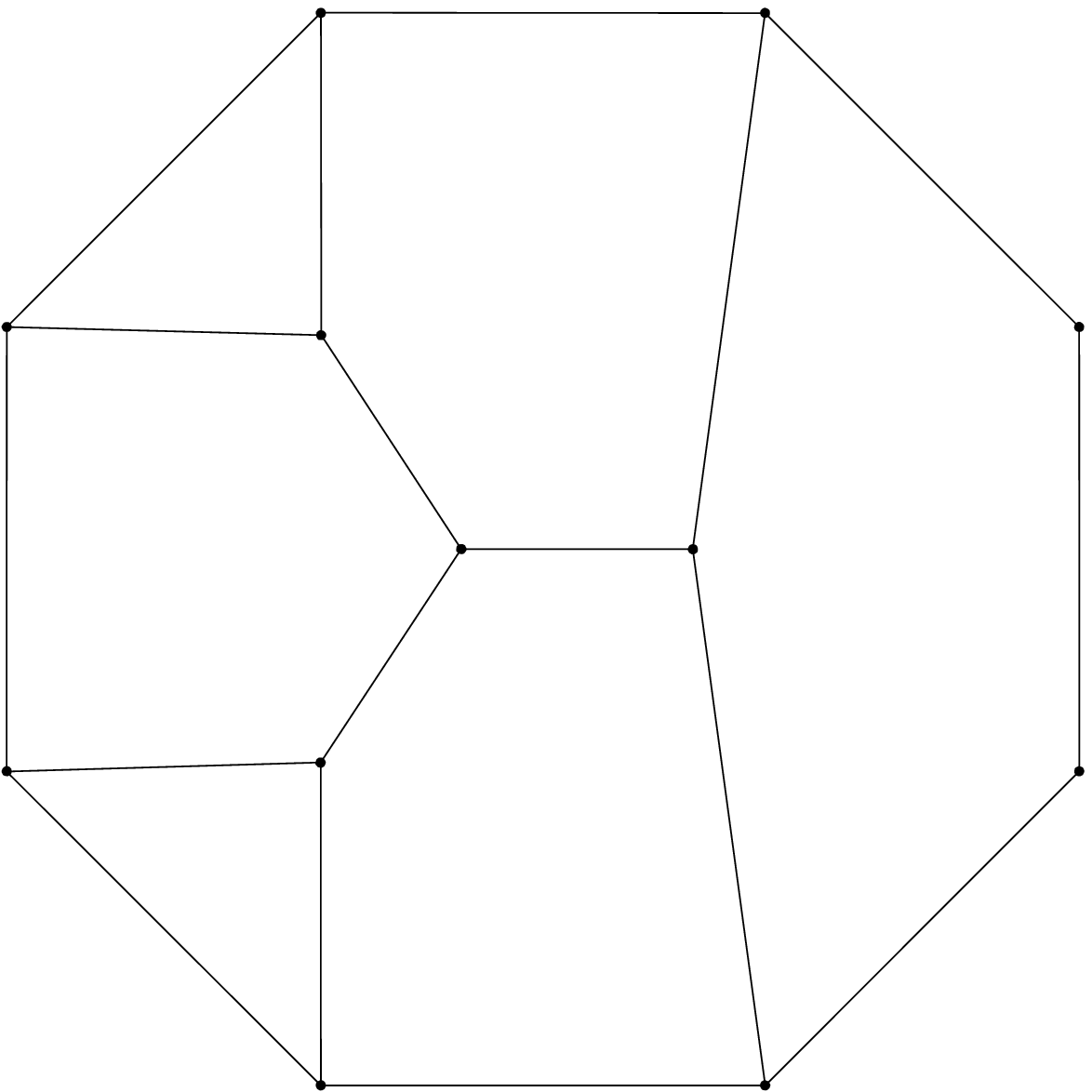}\par
$C_s$
\end{minipage}
\begin{minipage}{3cm}
\centering
\epsfig{height=20mm, file=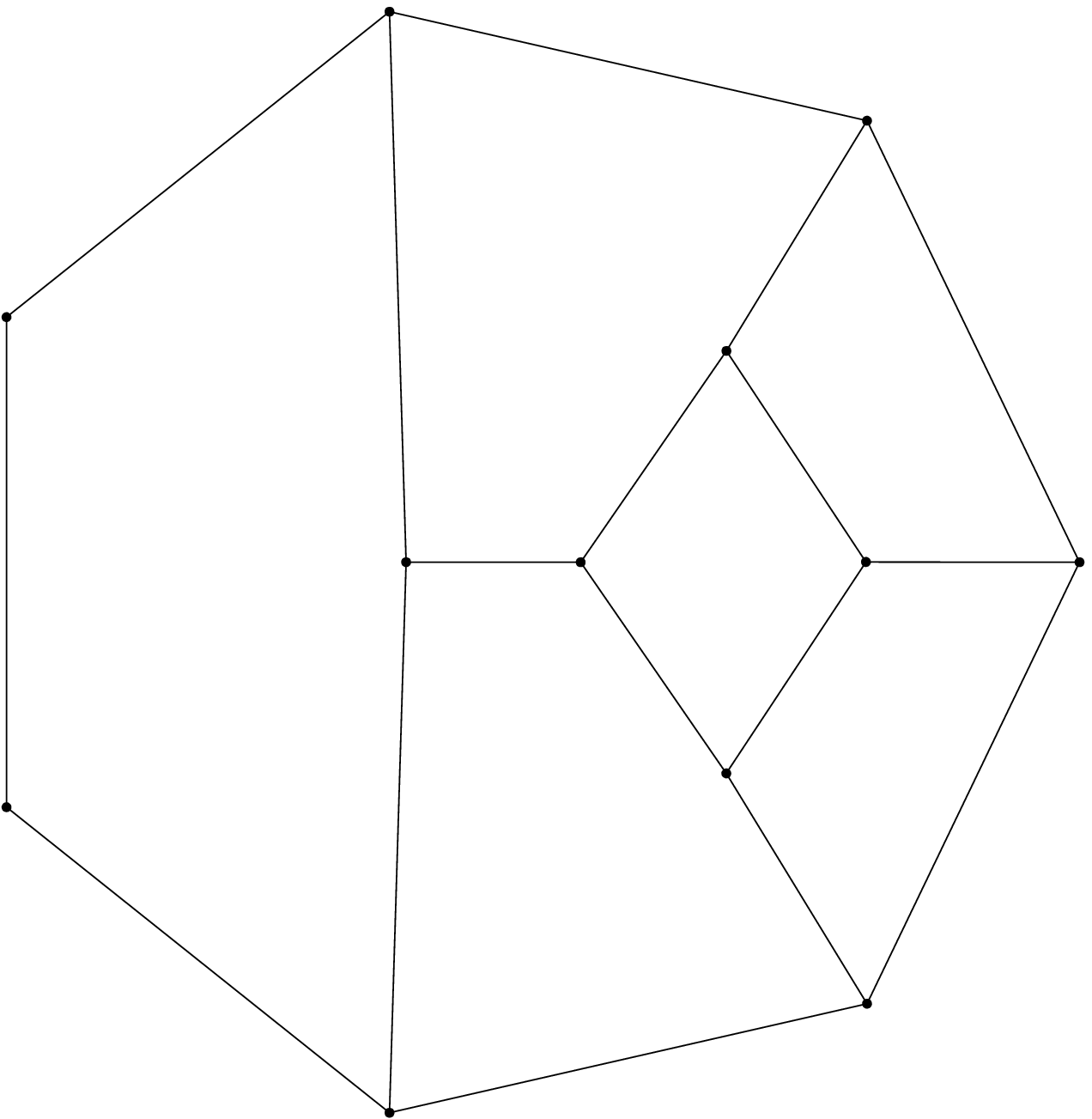}\par
$C_s$
\end{minipage}
\begin{minipage}{3cm}
\centering
\epsfig{height=20mm, file=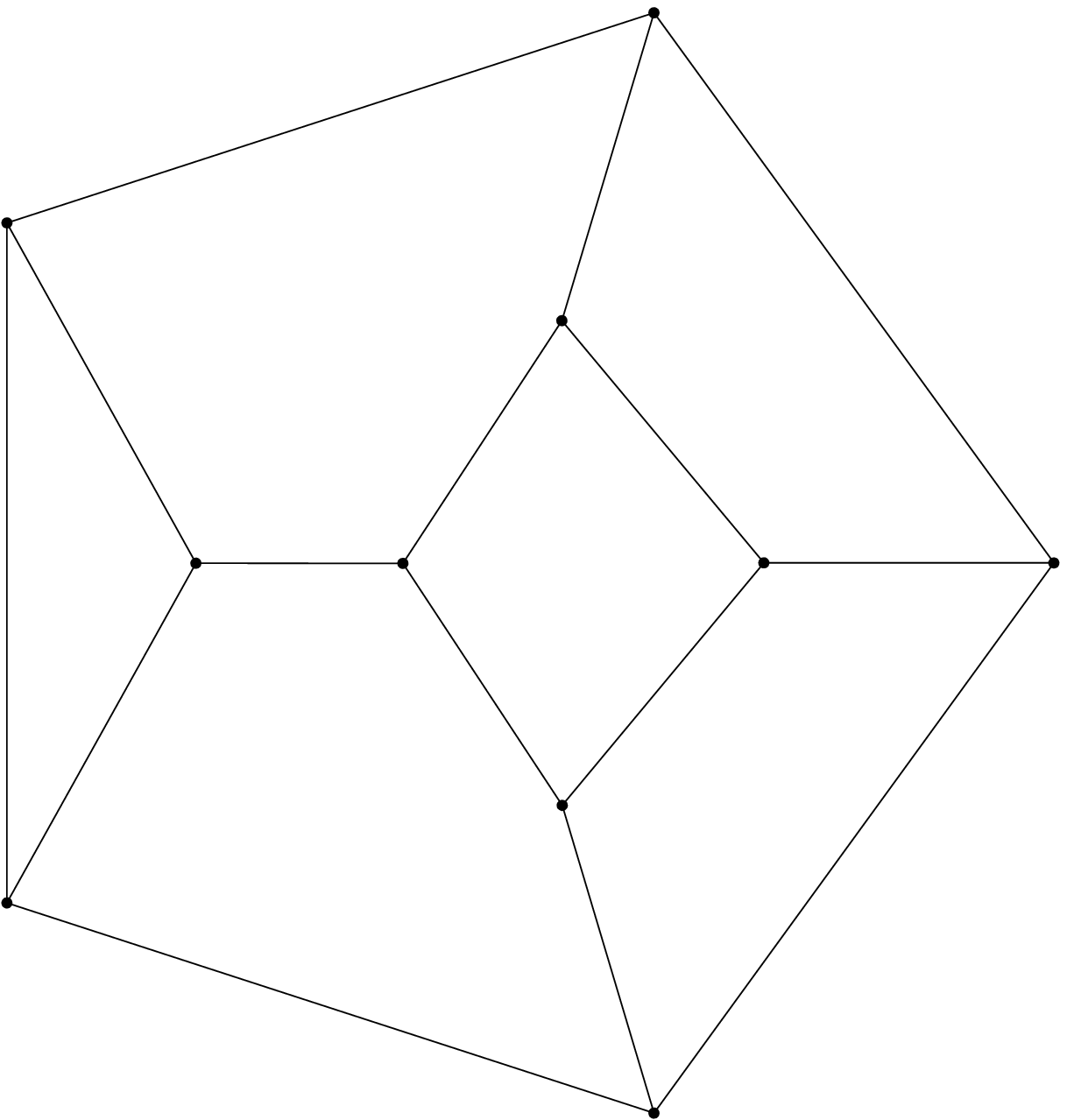}\par
$C_s$, nonext.
%PAIR1
\end{minipage}
\begin{minipage}{3cm}
\centering
\epsfig{height=20mm, file=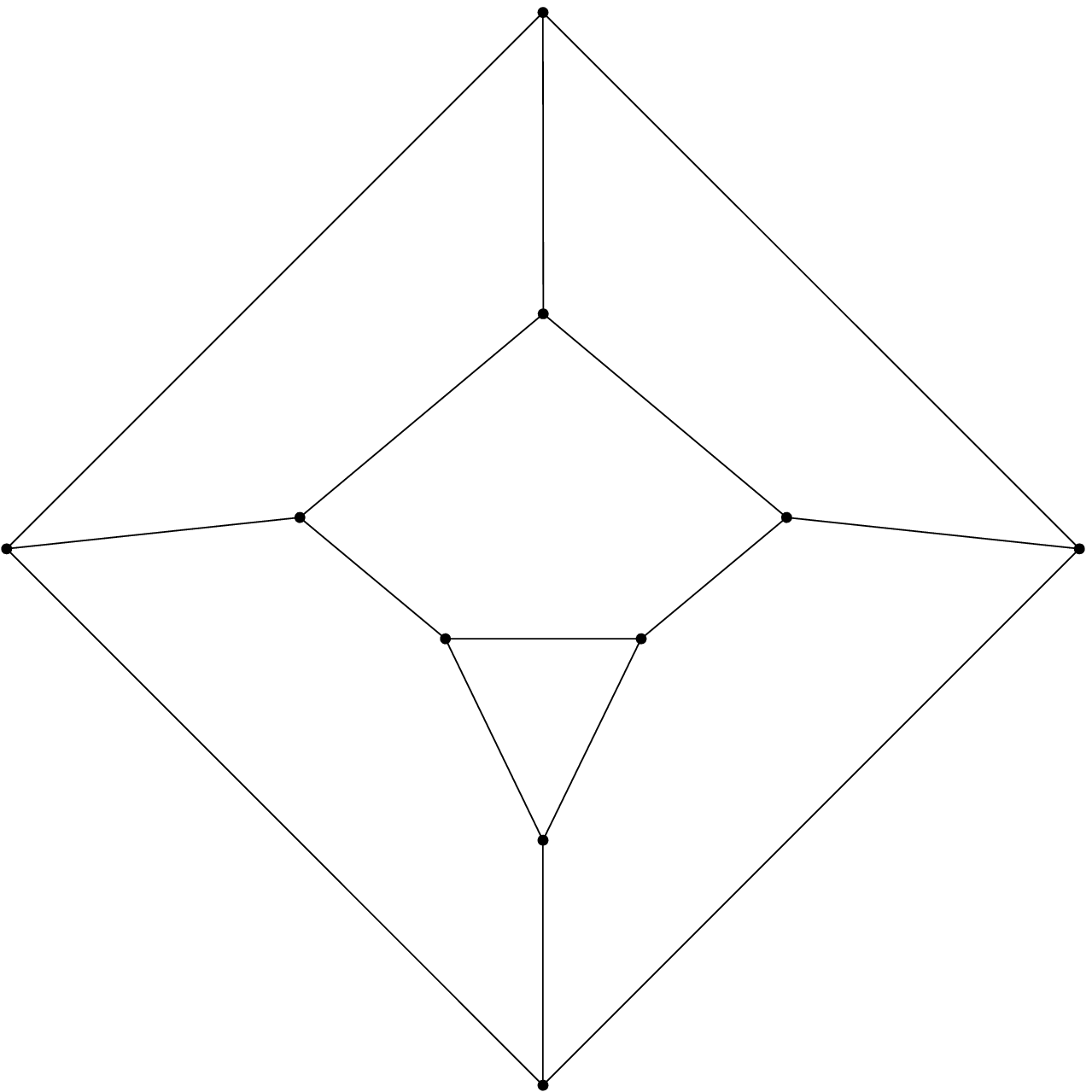}\par
$C_s$, nonext.
%PAIR1
\end{minipage}
\begin{minipage}{3cm}
\centering
\epsfig{height=20mm, file=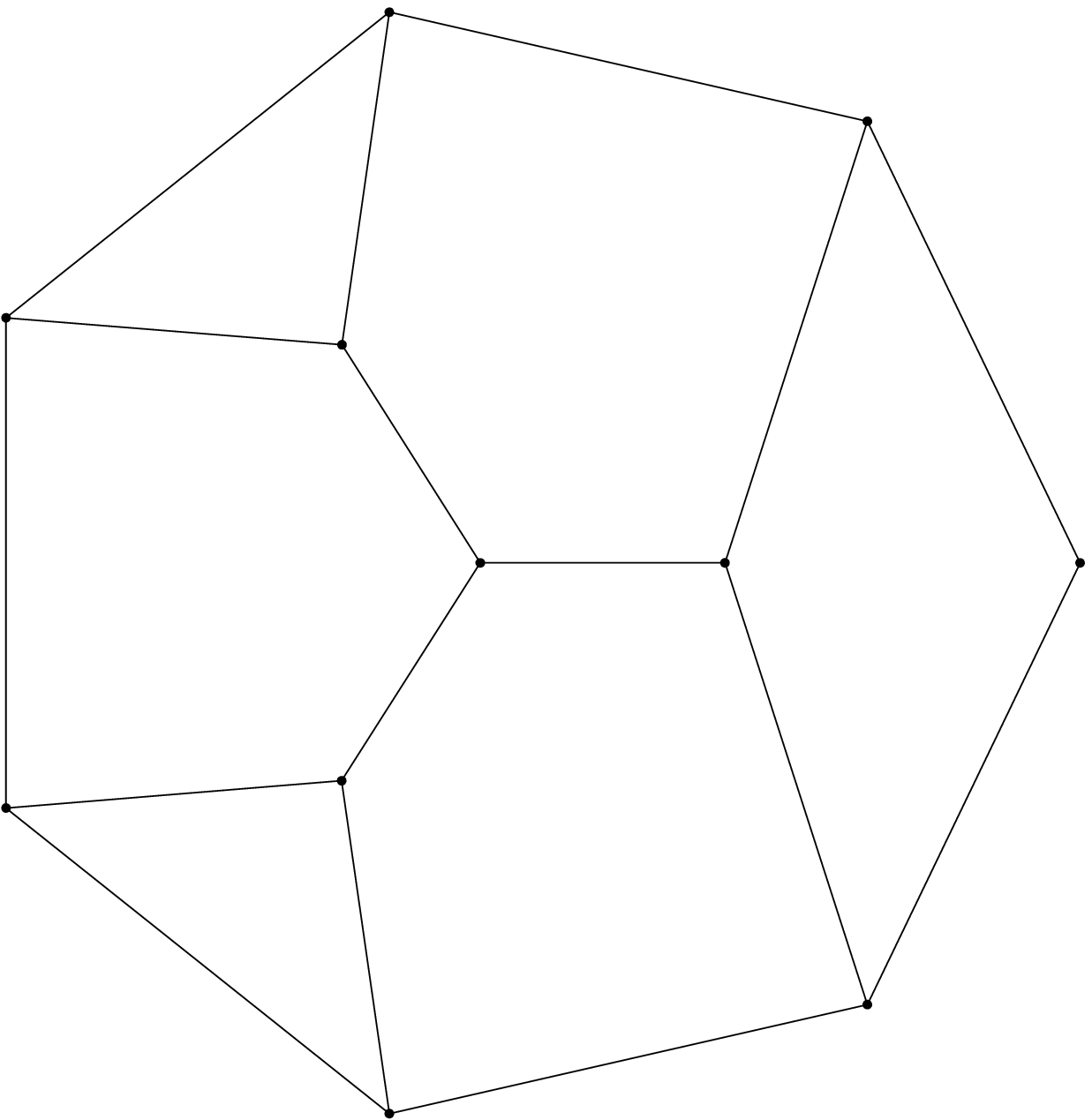}\par
$C_s$, nonext.
\end{minipage}
\begin{minipage}{3cm}
\centering
\epsfig{height=20mm, file=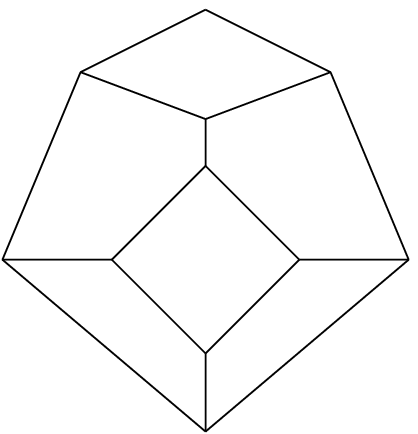}\par
$C_s$, nonext.
\end{minipage}
\begin{minipage}{3cm}
\centering
\epsfig{height=20mm, file=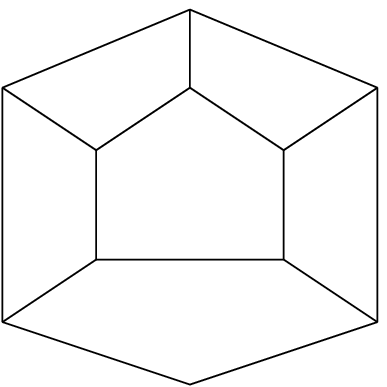}\par
$C_s$, nonext.
\end{minipage}
\begin{minipage}{3cm}
\centering
\epsfig{height=20mm, file=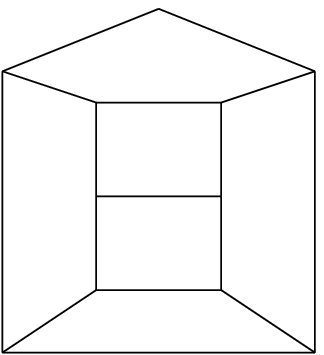}\par
$C_s$,~nonext.~$(C_{2\nu})$
\end{minipage}
\begin{minipage}{3cm}
\centering
\epsfig{height=20mm, file=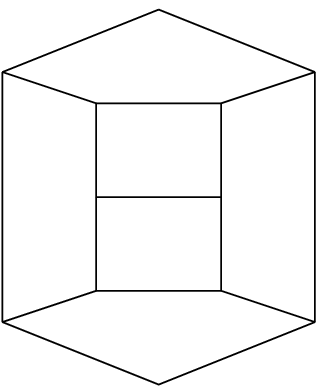}\par
$C_{2\nu}$
\end{minipage}
\begin{minipage}{3cm}
\centering
\epsfig{height=20mm, file=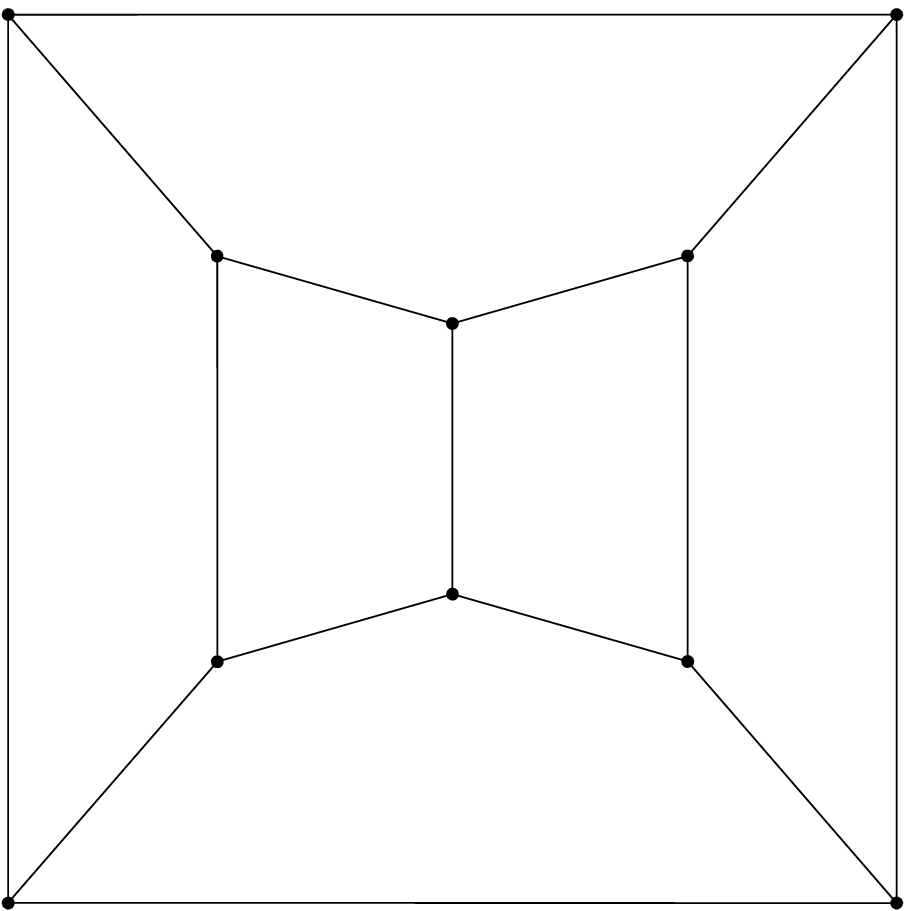}\par
$C_{2\nu}$, nonext.
%PAIR2
\end{minipage}
\begin{minipage}{3cm}
\centering
\epsfig{height=20mm, file=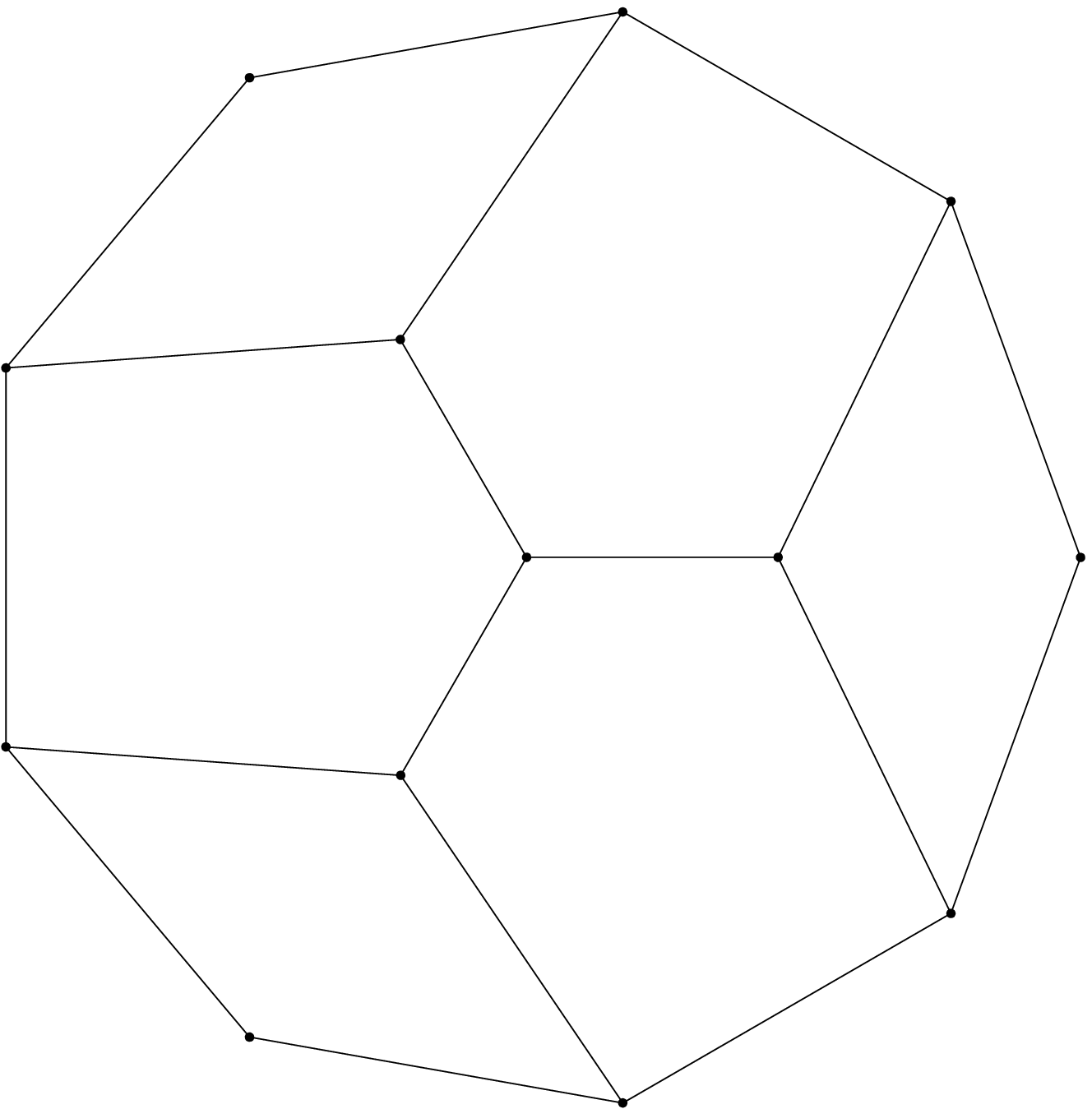}\par
$C_{3\nu}$
\end{minipage}
\begin{minipage}{3cm}
\centering
\epsfig{height=20mm, file=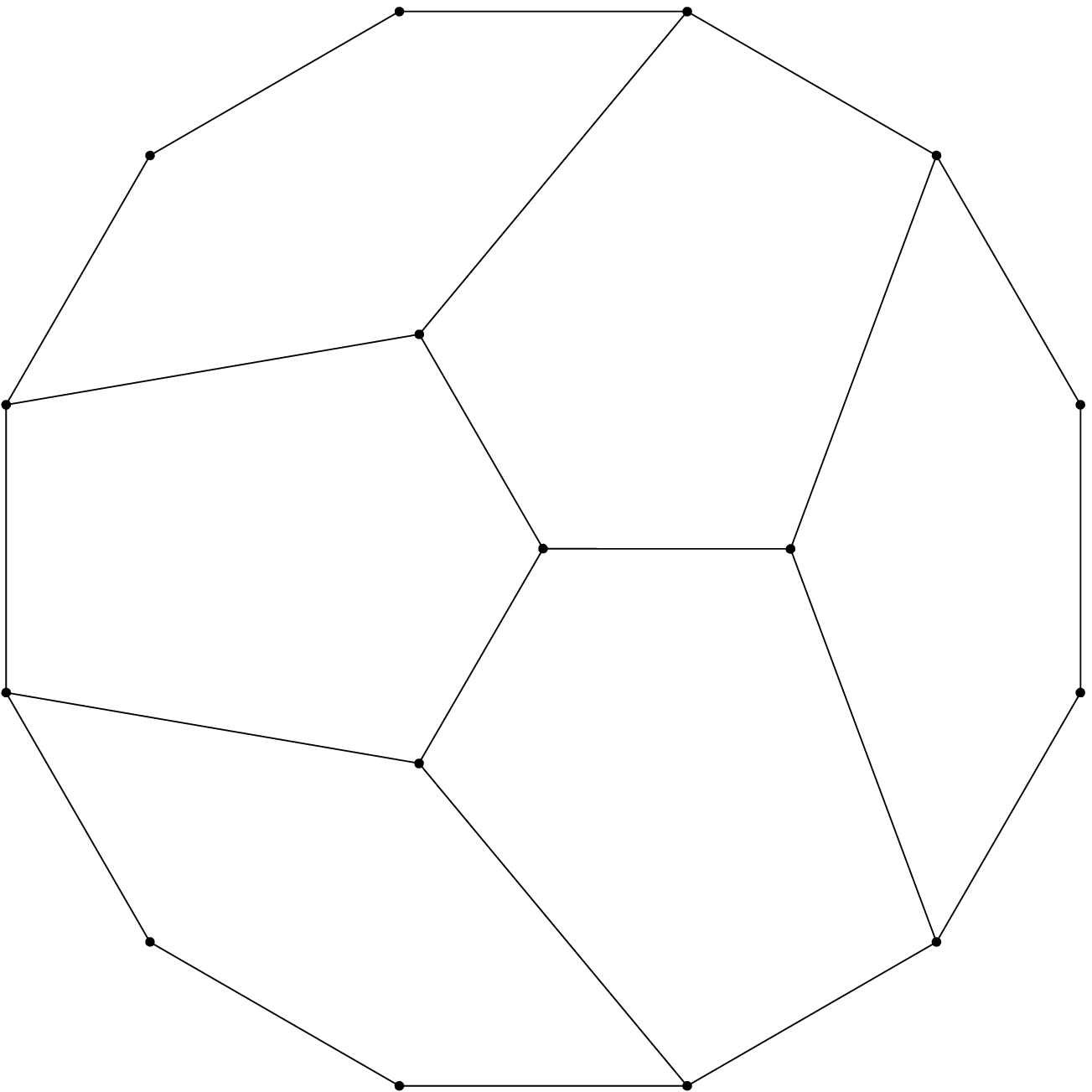}\par
$C_{3\nu}$
\end{minipage}
\begin{minipage}{3cm}
\centering
\epsfig{height=20mm, file=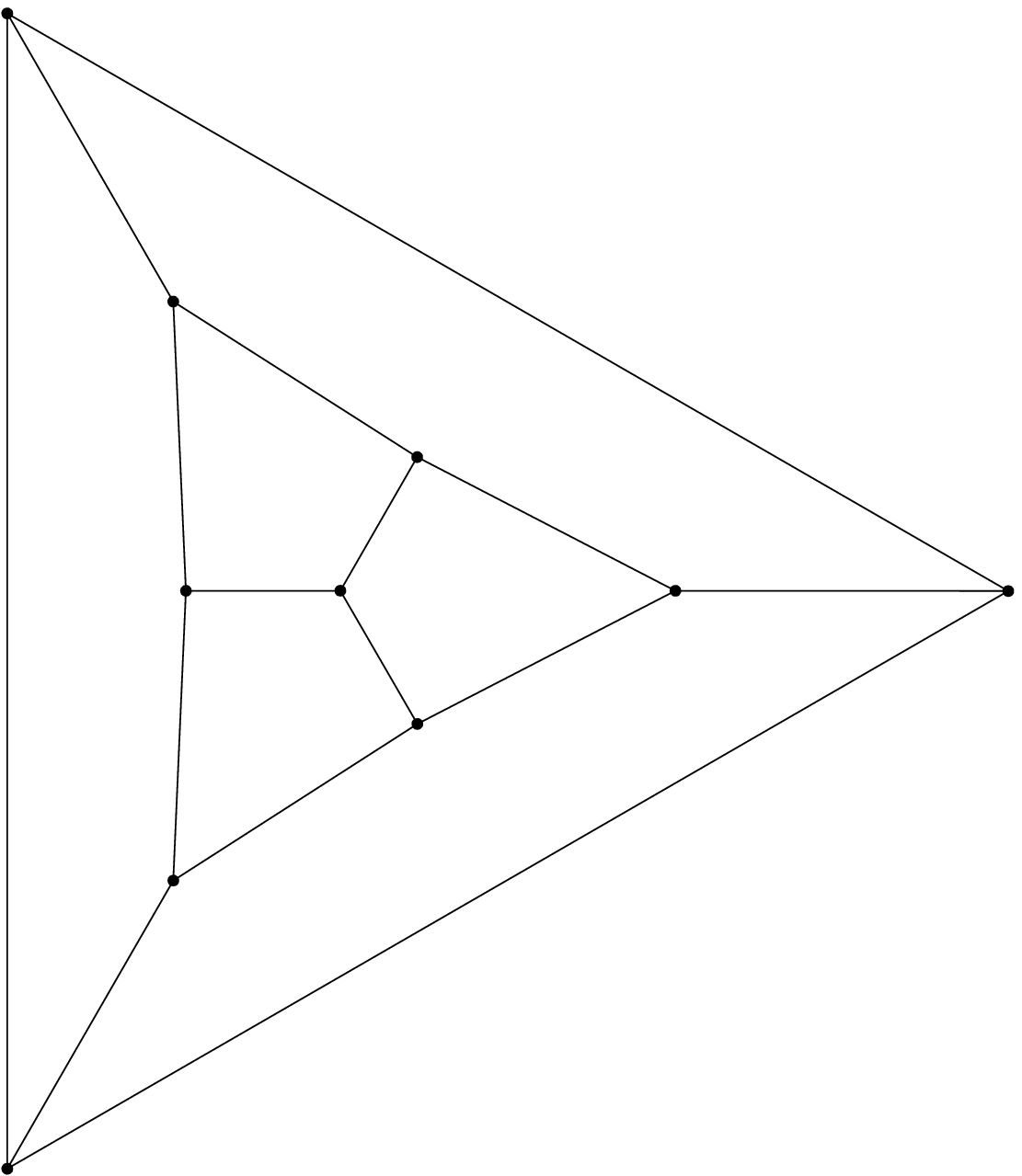}\par
$C_{3\nu}$, nonext.
%PAIR1
\end{minipage}
\begin{minipage}{3cm}
\centering
\epsfig{height=20mm, file=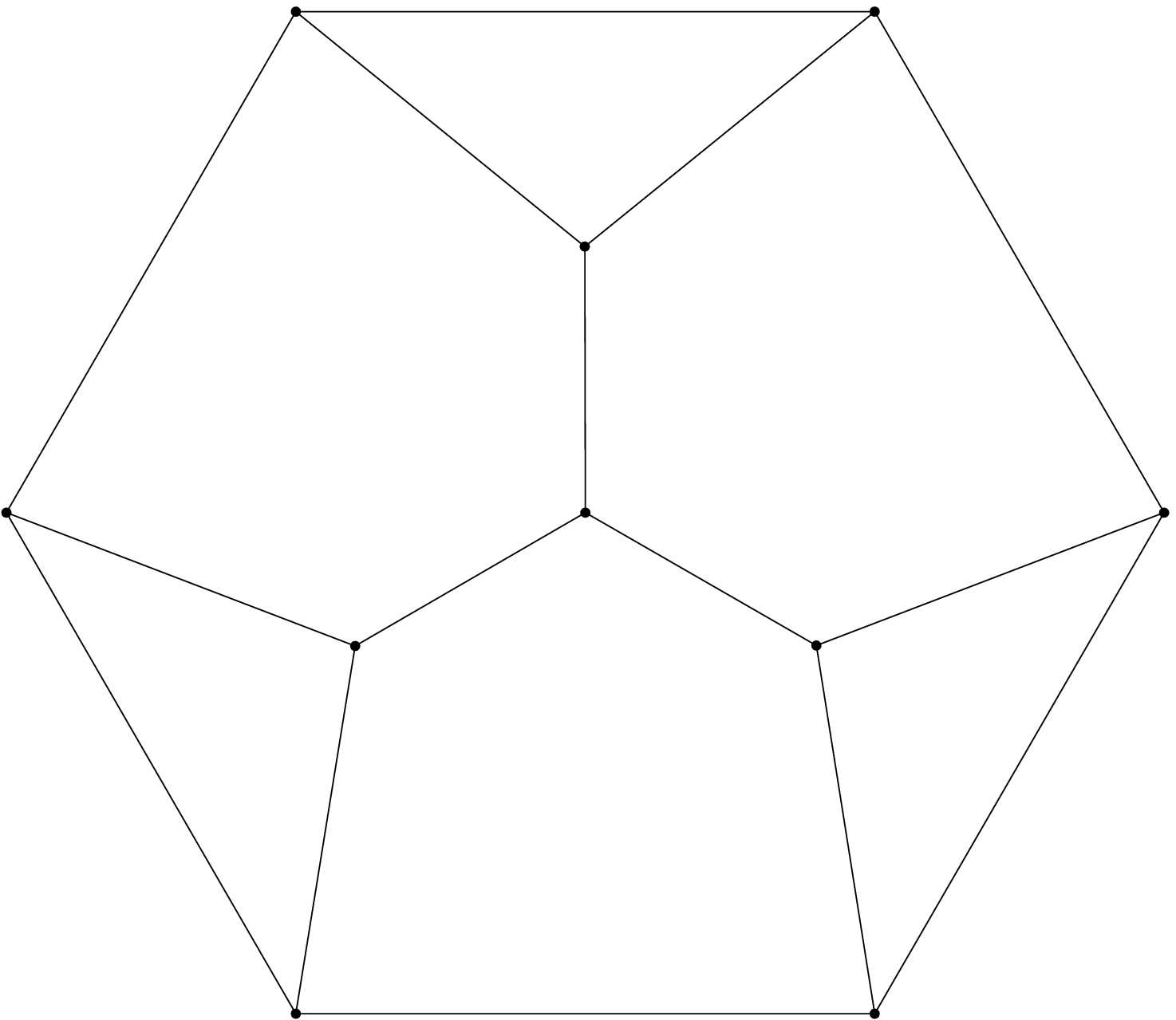}\par
$C_{3\nu}$, nonext.
\end{minipage}
\begin{minipage}{3cm}
\centering
\epsfig{height=20mm, file=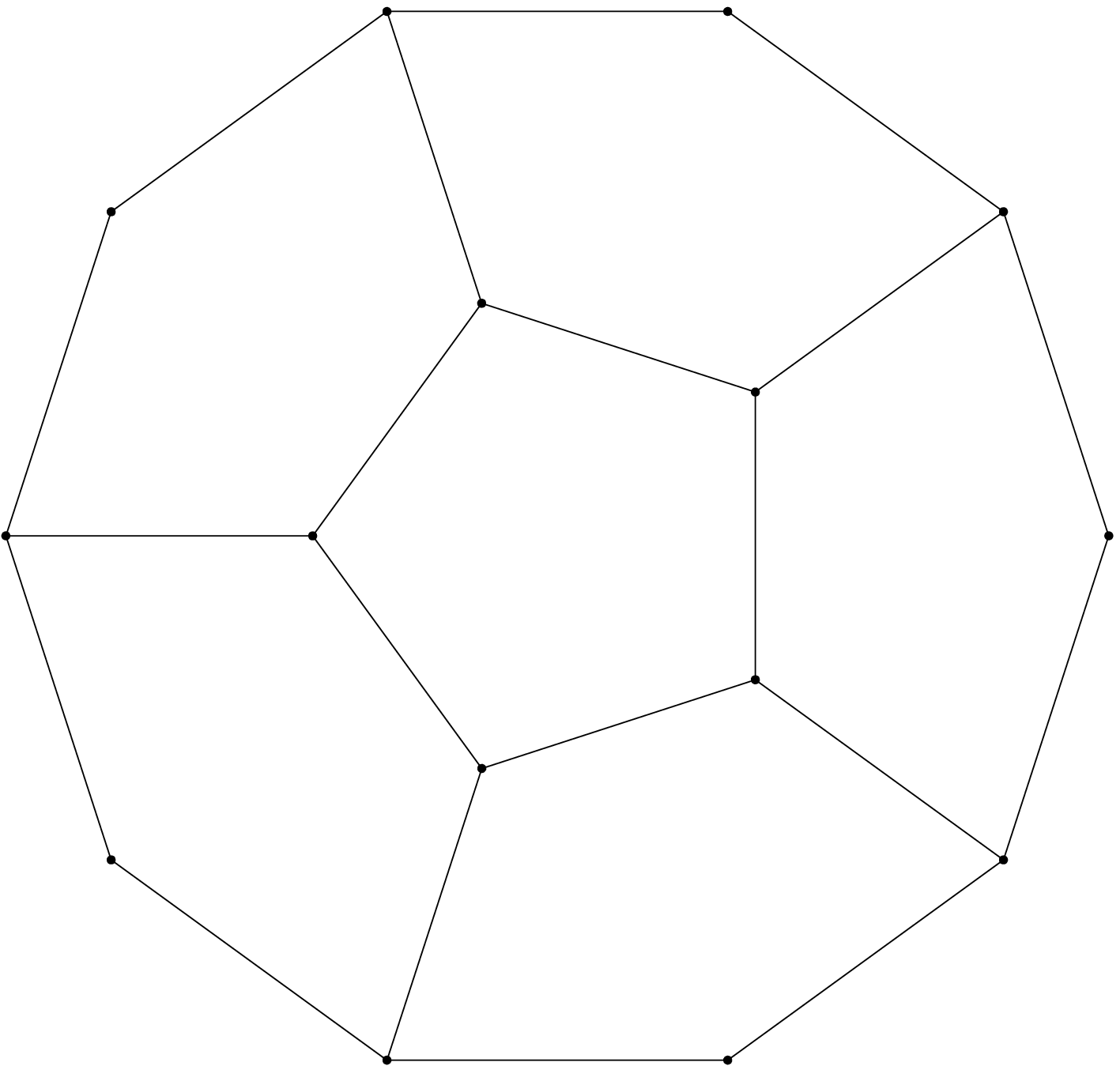}\par
$C_{5\nu}$
\end{minipage}
\begin{minipage}{3cm}
\centering
\epsfig{height=20mm, file=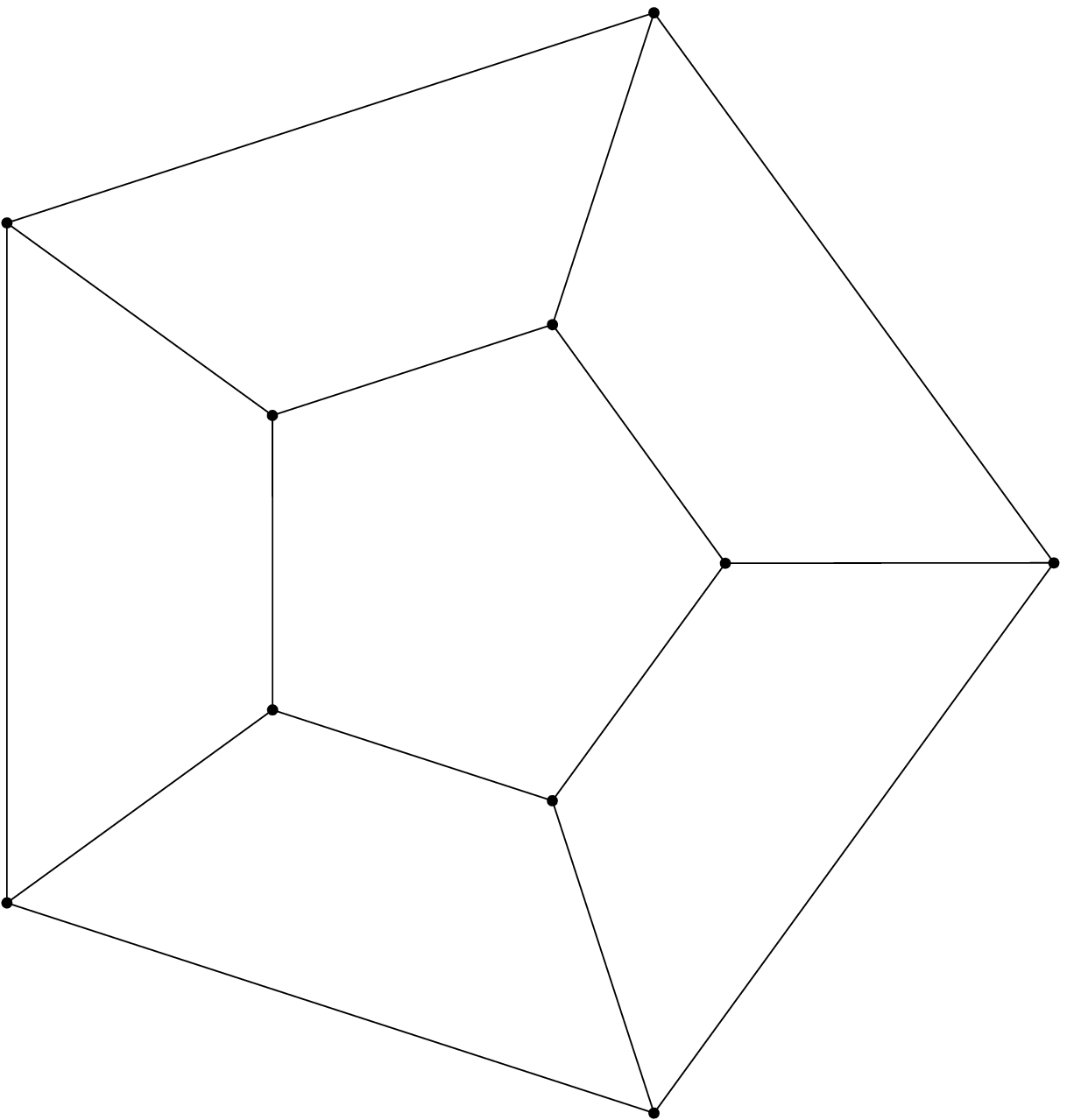}\par
$C_{5\nu}$, nonext.
%PAIR2
\end{minipage}

\end{center}
List of sporadic elementary $(\{3,4,5\},3)$-polycycles with $7$ faces:
\begin{center}
\begin{minipage}{3cm}
\centering
\epsfig{height=20mm, file=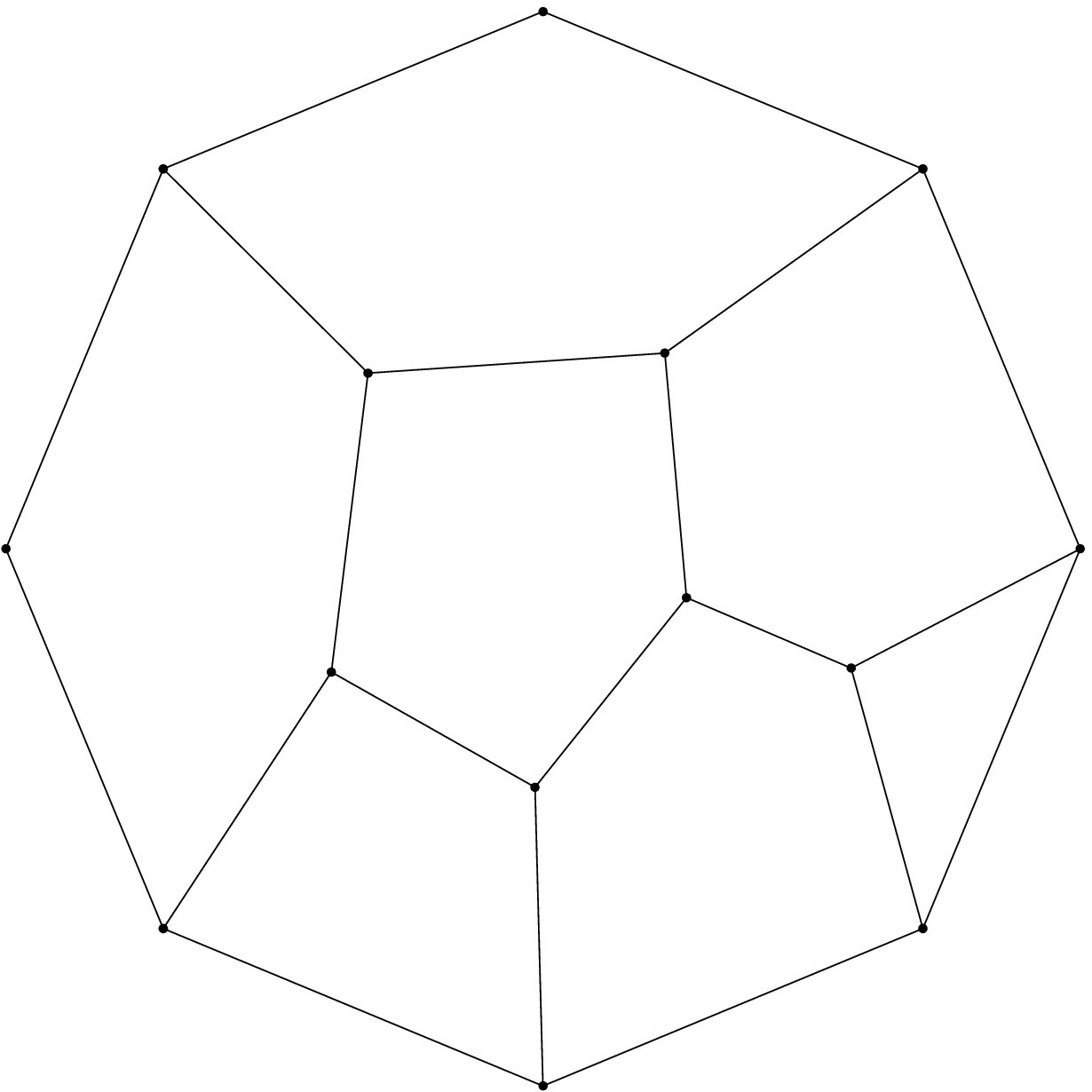}\par
$C_1$
\end{minipage}
\begin{minipage}{3cm}
\centering
\epsfig{height=20mm, file=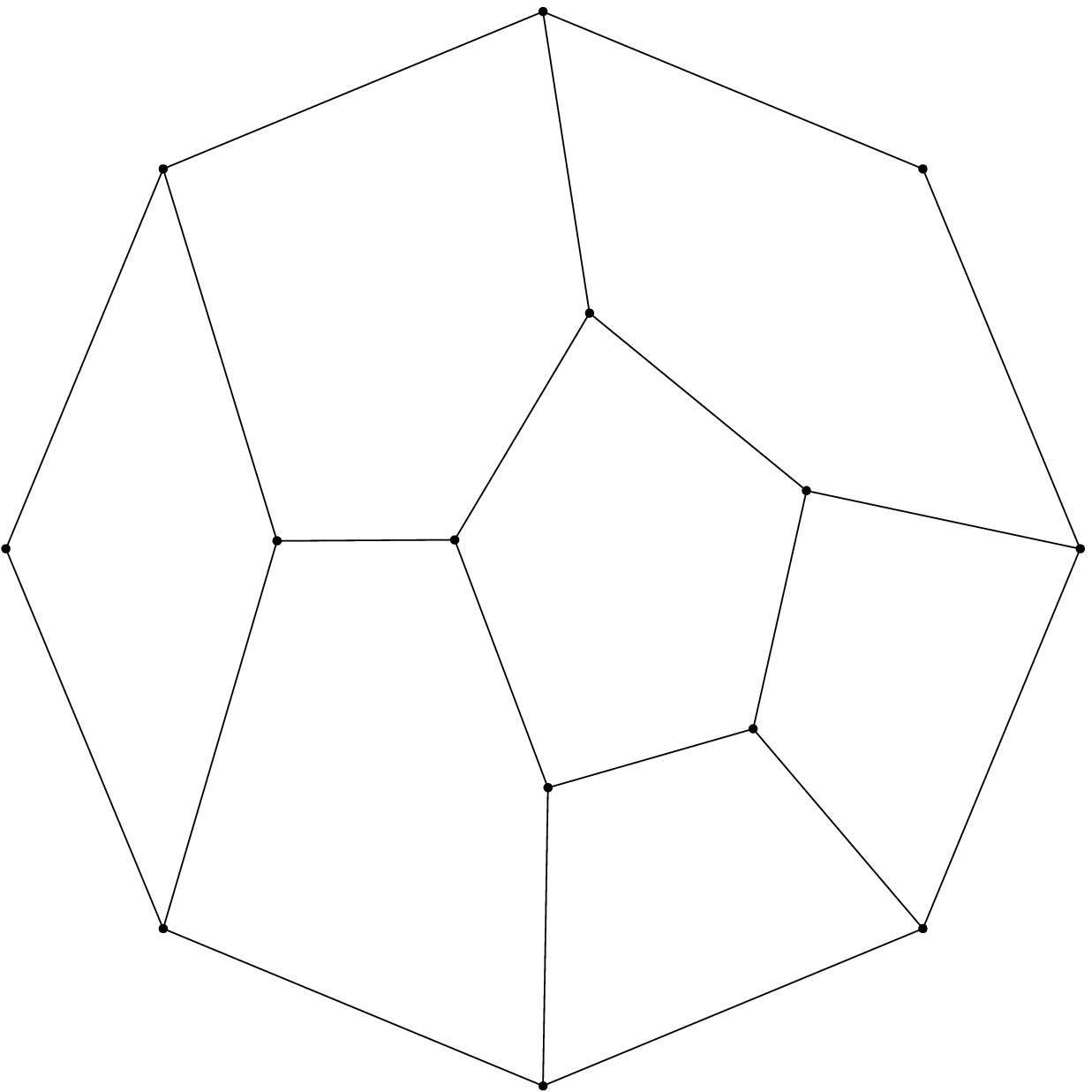}\par
$C_1$
\end{minipage}
\begin{minipage}{3cm}
\centering
\epsfig{height=20mm, file=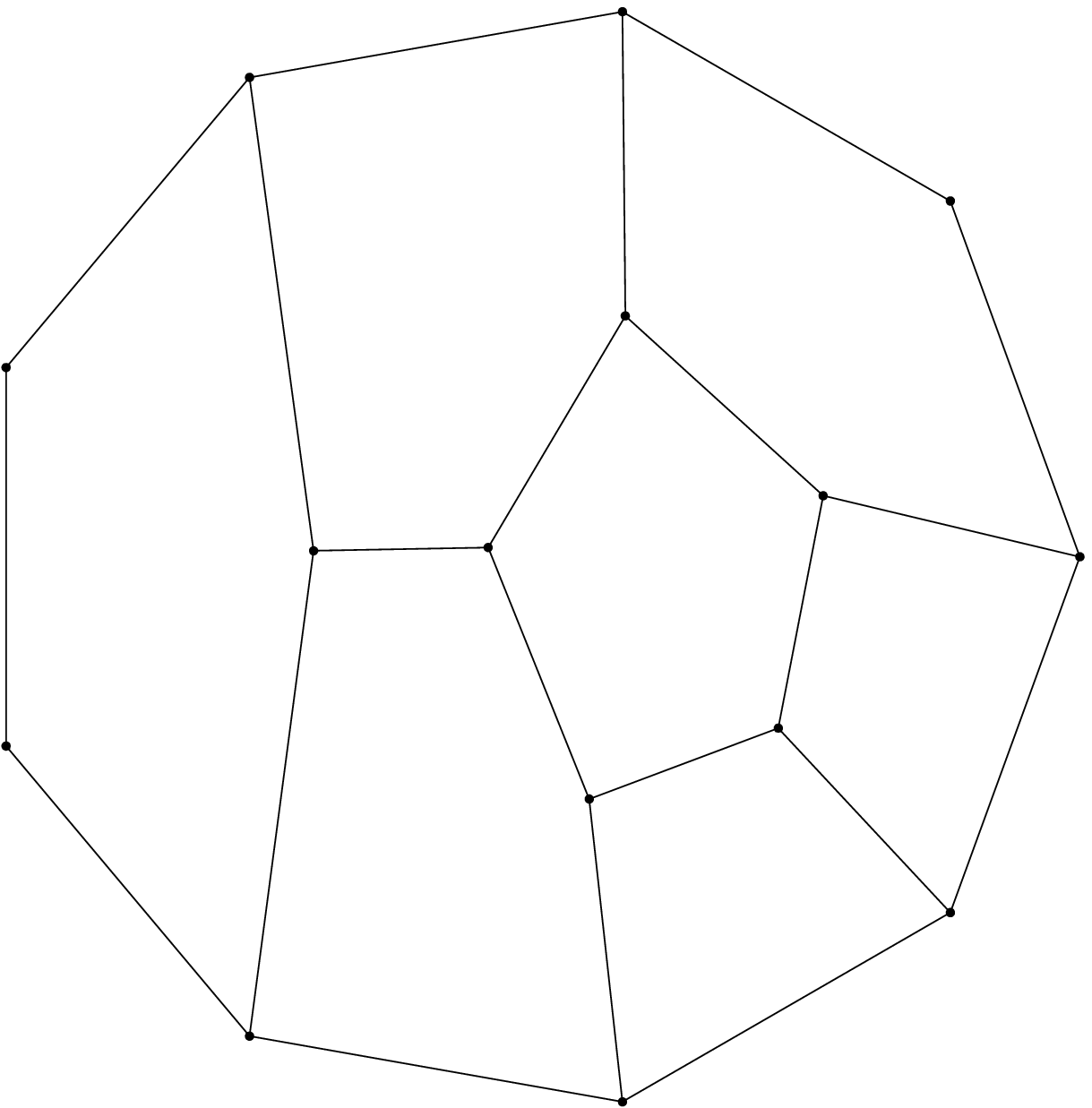}\par
$C_1$
\end{minipage}
\begin{minipage}{3cm}
\centering
\epsfig{height=20mm, file=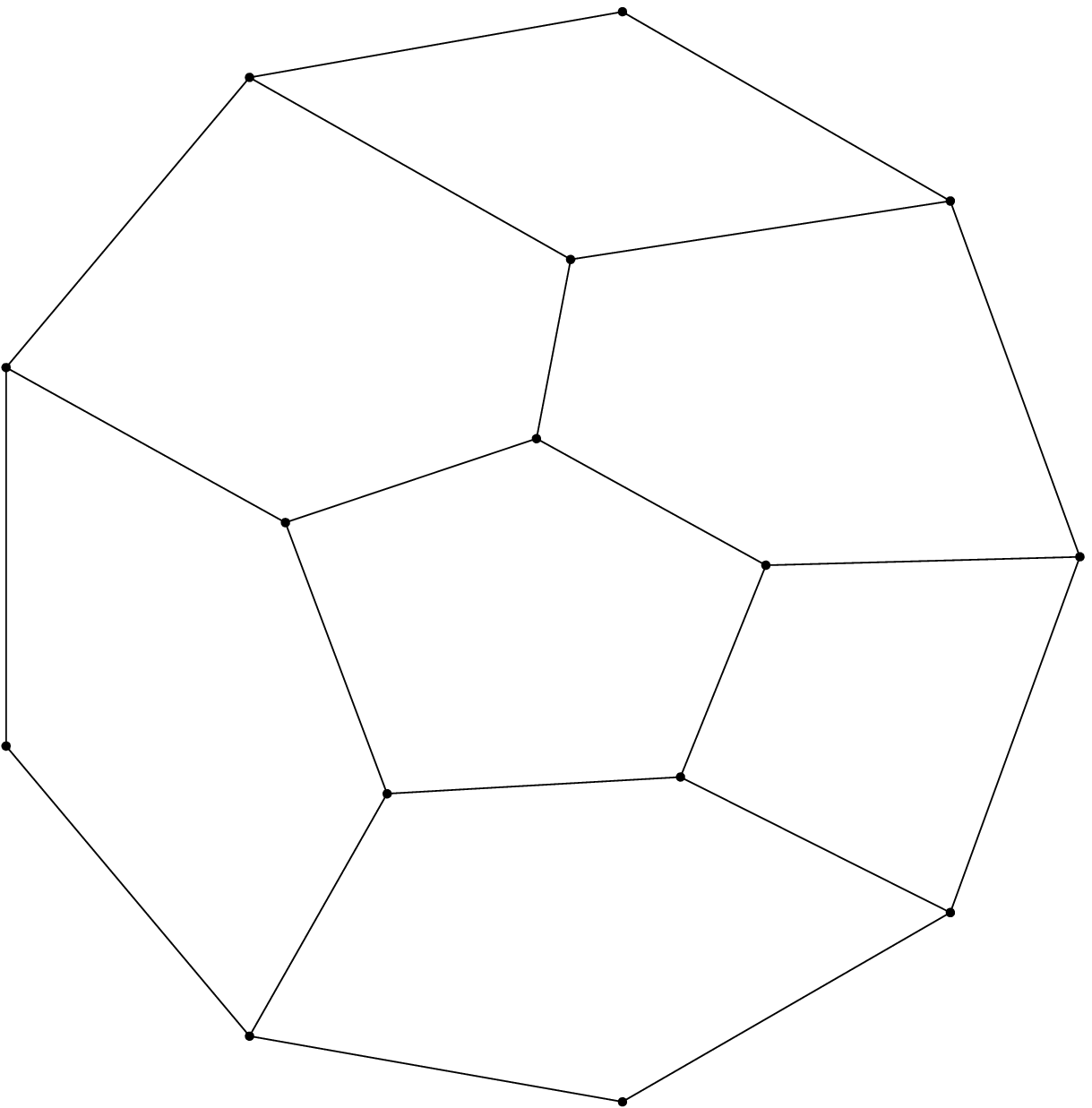}\par
$C_1$
\end{minipage}
\begin{minipage}{3cm}
\centering
\epsfig{height=20mm, file=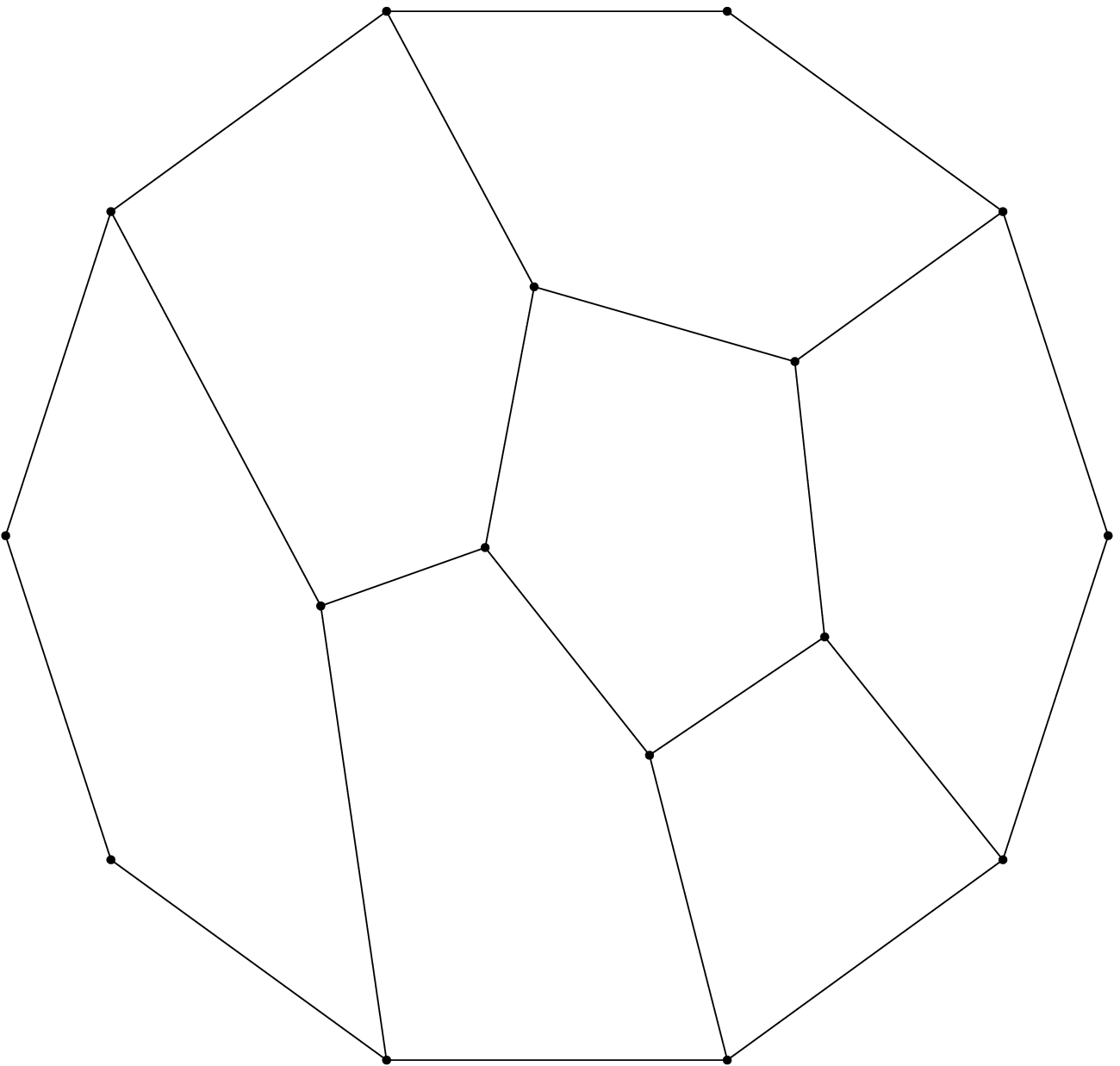}\par
$C_1$
\end{minipage}
\begin{minipage}{3cm}
\centering
\epsfig{height=20mm, file=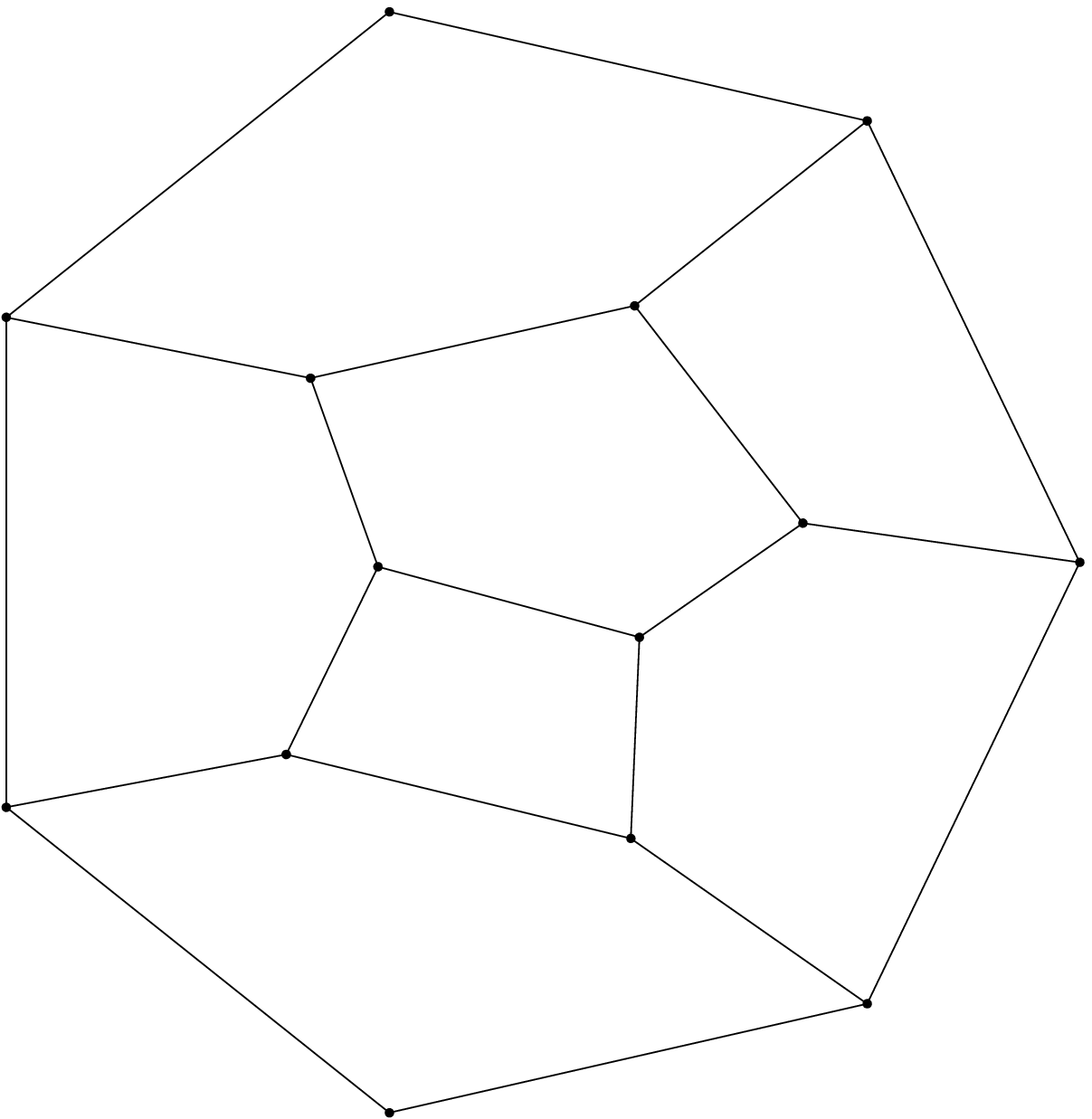}\par
$C_1$
\end{minipage}
\begin{minipage}{3cm}
\centering
\epsfig{height=20mm, file=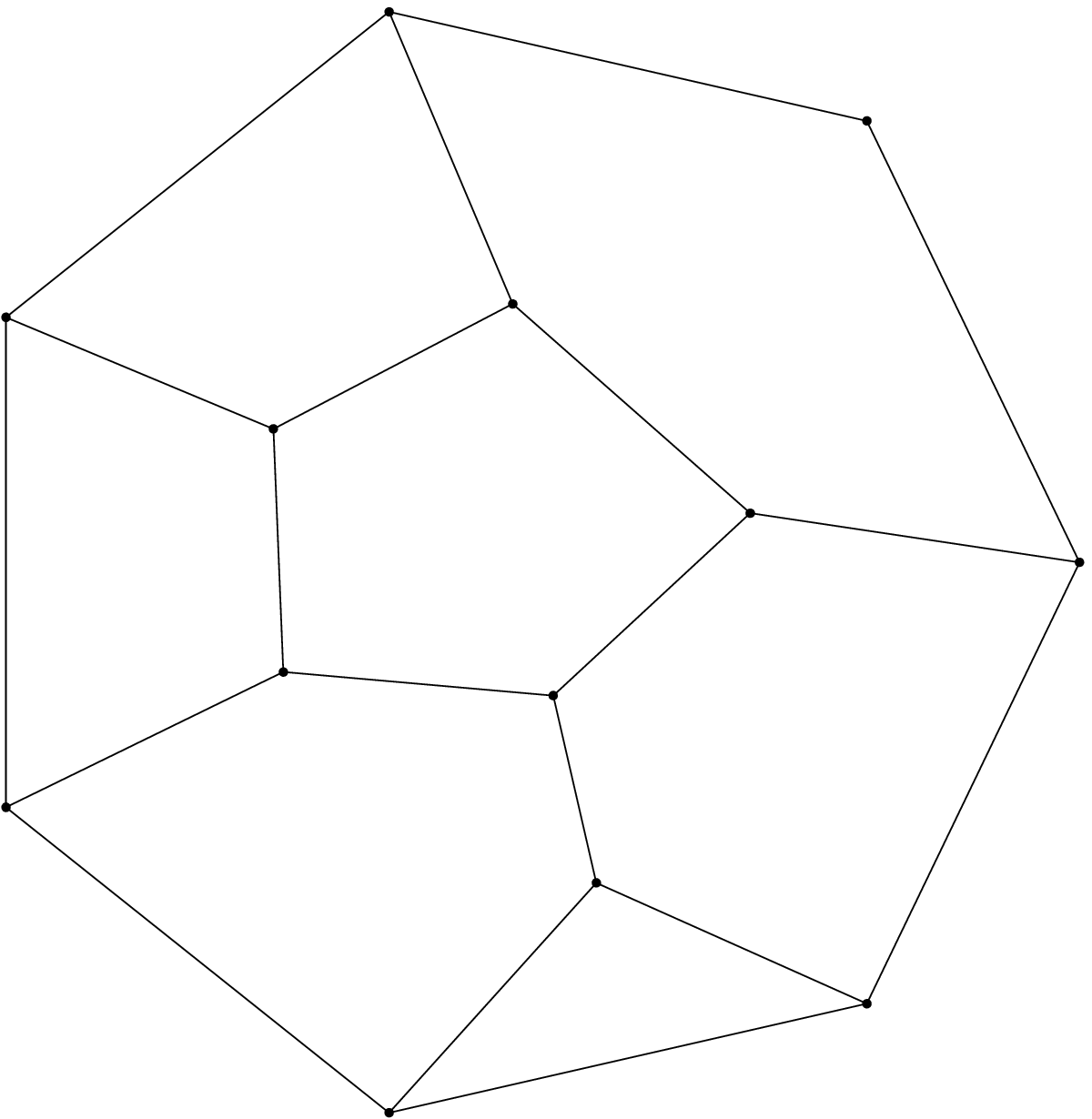}\par
$C_1$, nonext.
\end{minipage}
\begin{minipage}{3cm}
\centering
\epsfig{height=20mm, file=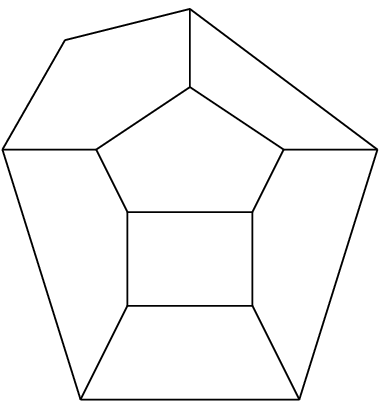}\par
$C_1$, nonext.
\end{minipage}
\begin{minipage}{3cm}
\centering
\epsfig{height=20mm, file=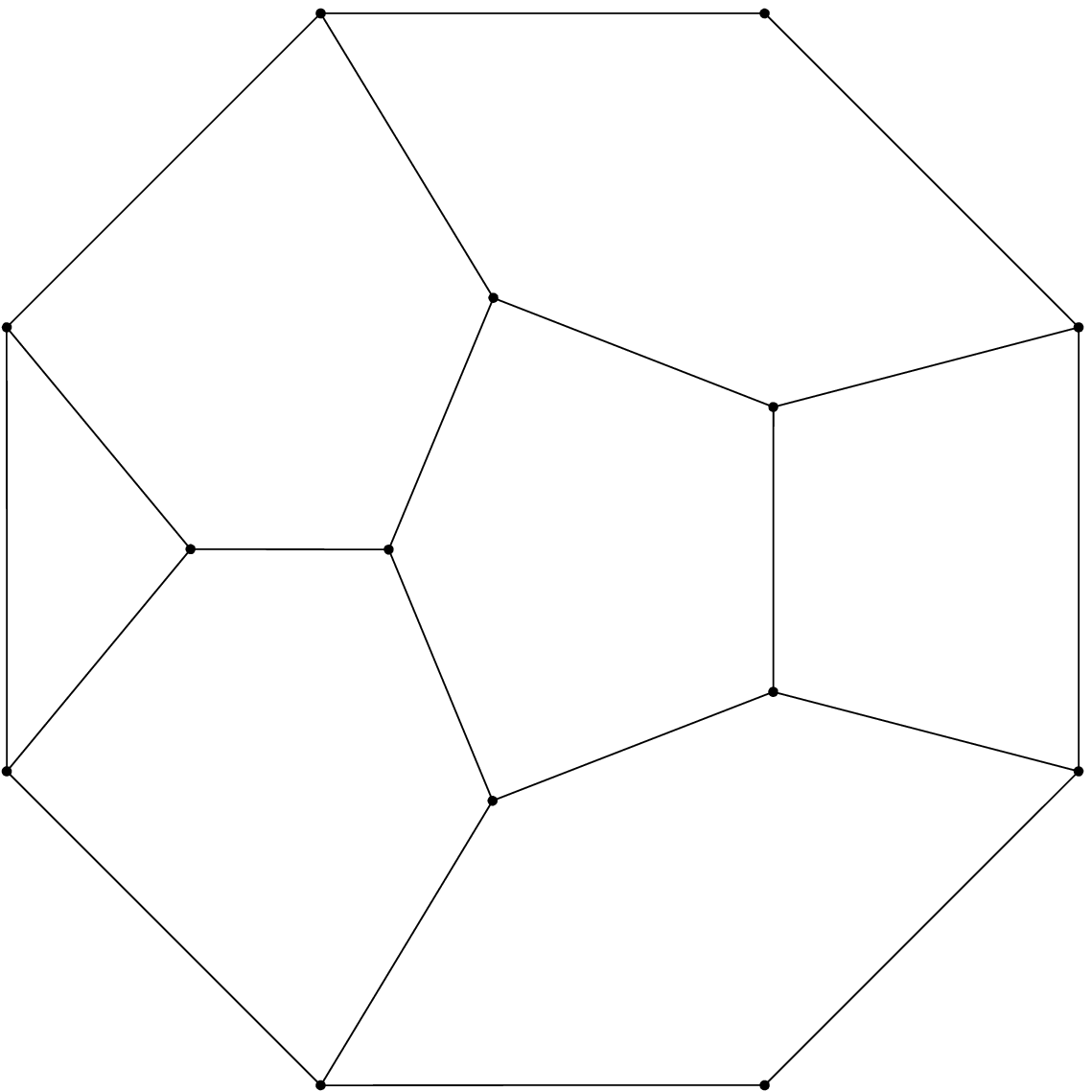}\par
$C_s$
\end{minipage}
\begin{minipage}{3cm}
\centering
\epsfig{height=20mm, file=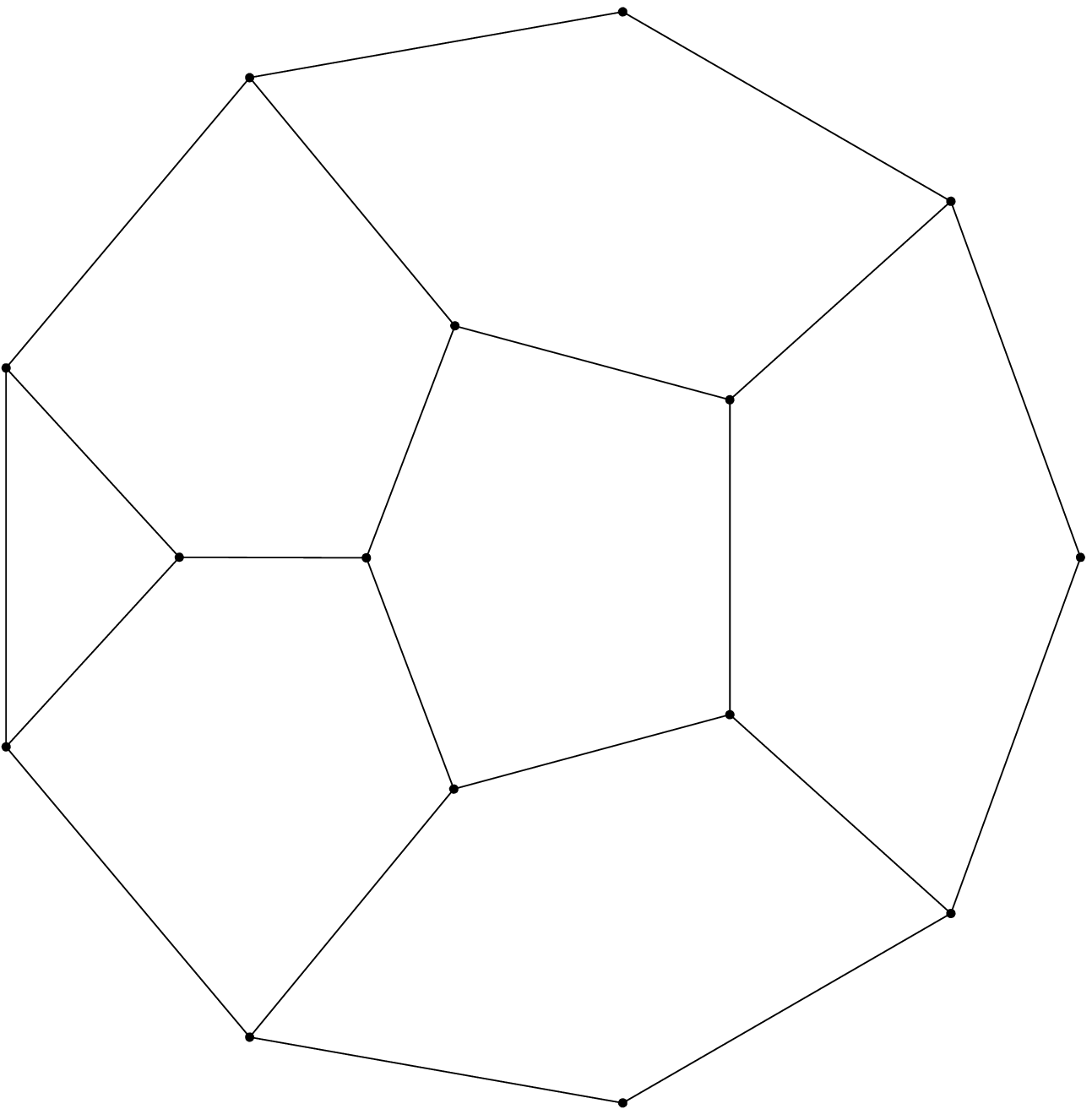}\par
$C_s$
\end{minipage}
\begin{minipage}{3cm}
\centering
\epsfig{height=20mm, file=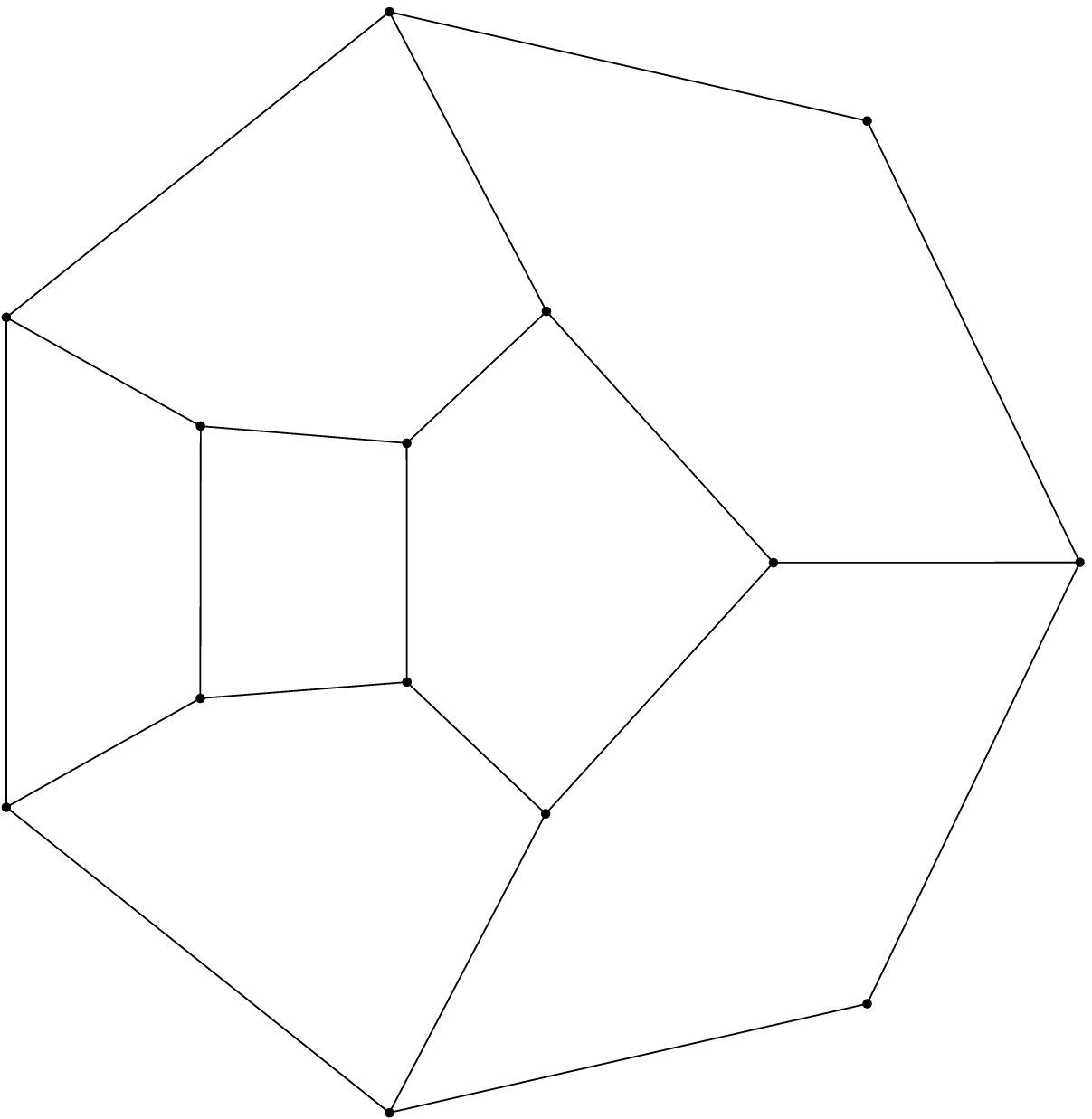}\par
$C_s$
\end{minipage}
\begin{minipage}{3cm}
\centering
\epsfig{height=20mm, file=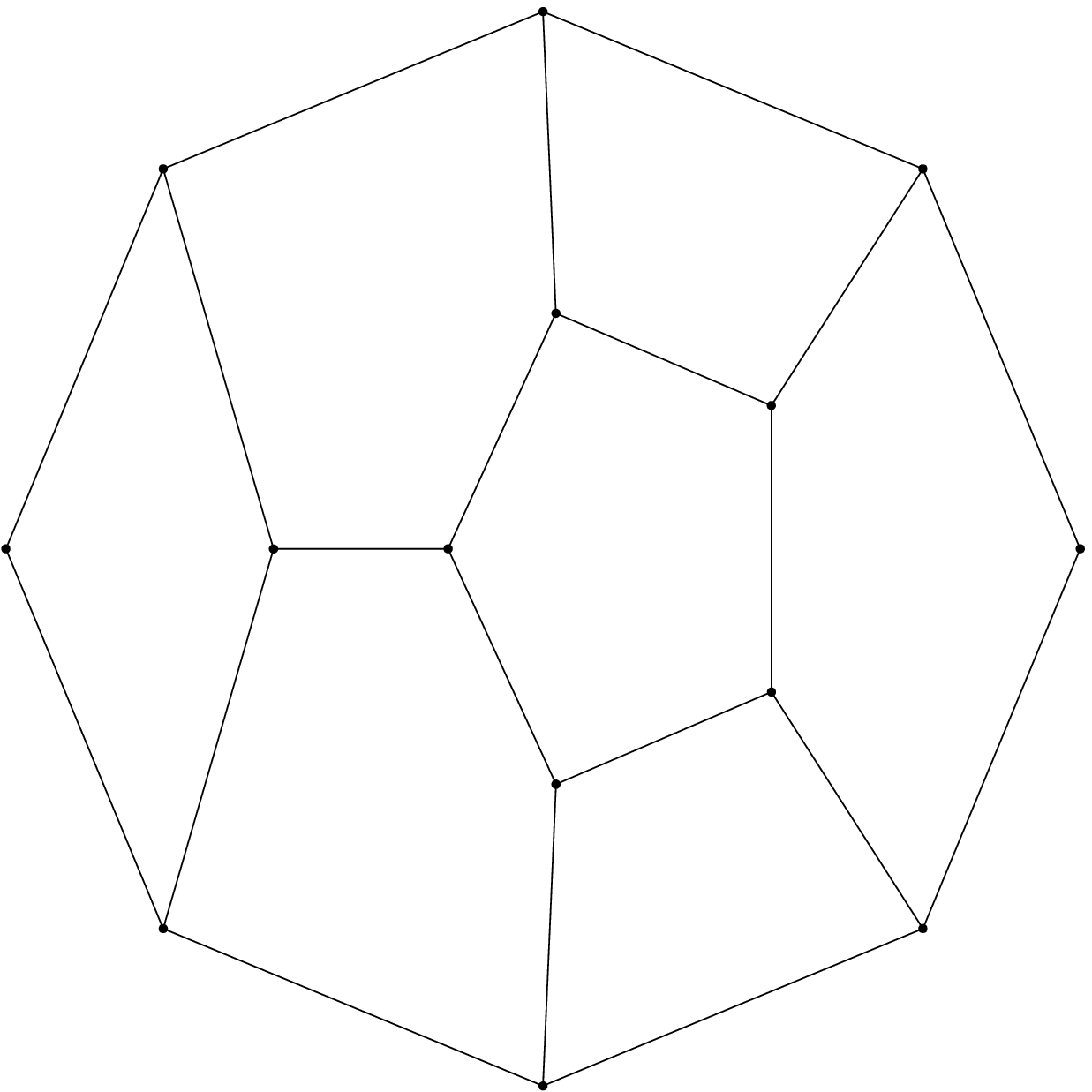}\par
$C_s$
\end{minipage}
\begin{minipage}{3cm}
\centering
\epsfig{height=20mm, file=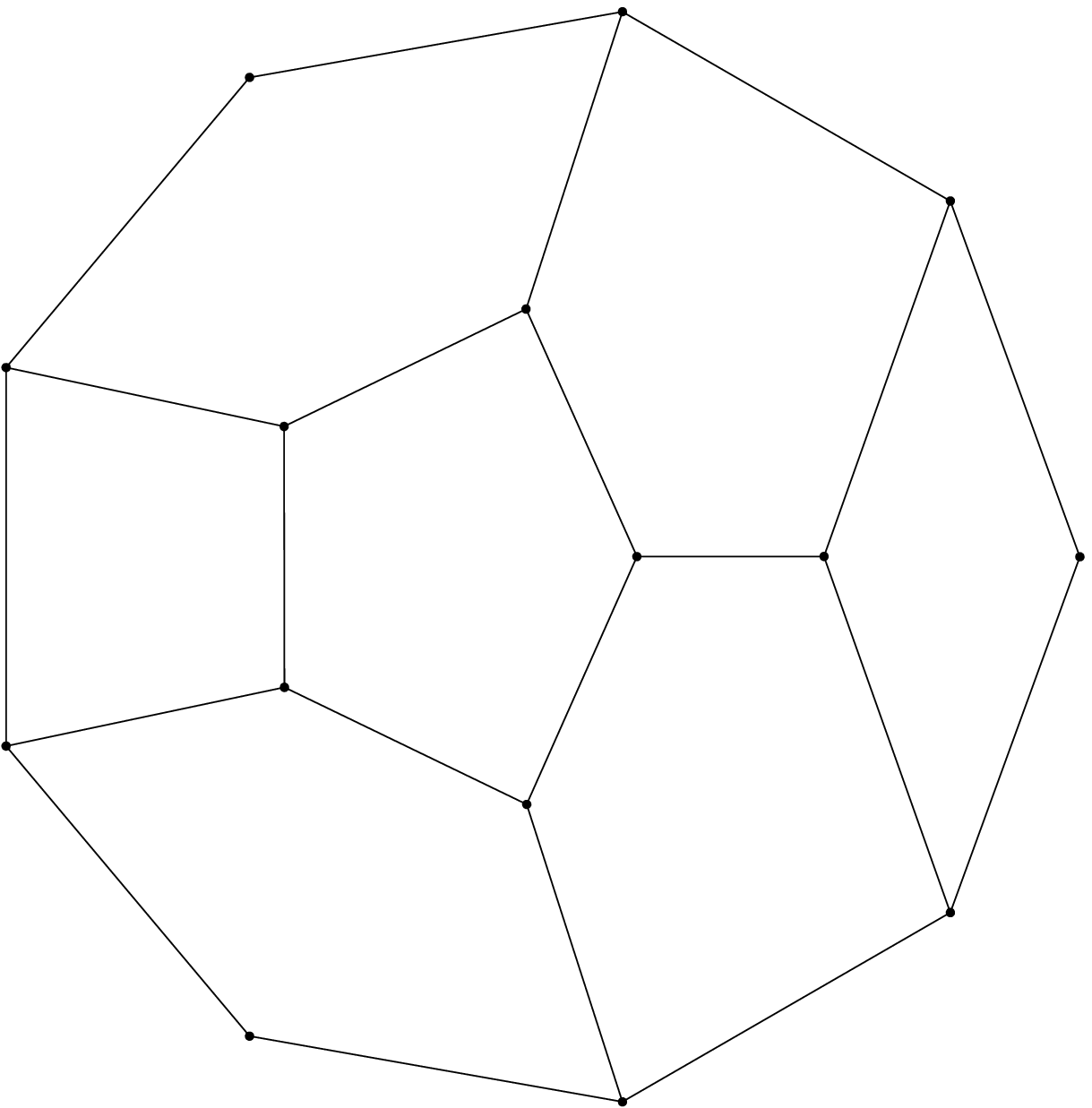}\par
$C_s$
\end{minipage}
\begin{minipage}{3cm}
\centering
\epsfig{height=20mm, file=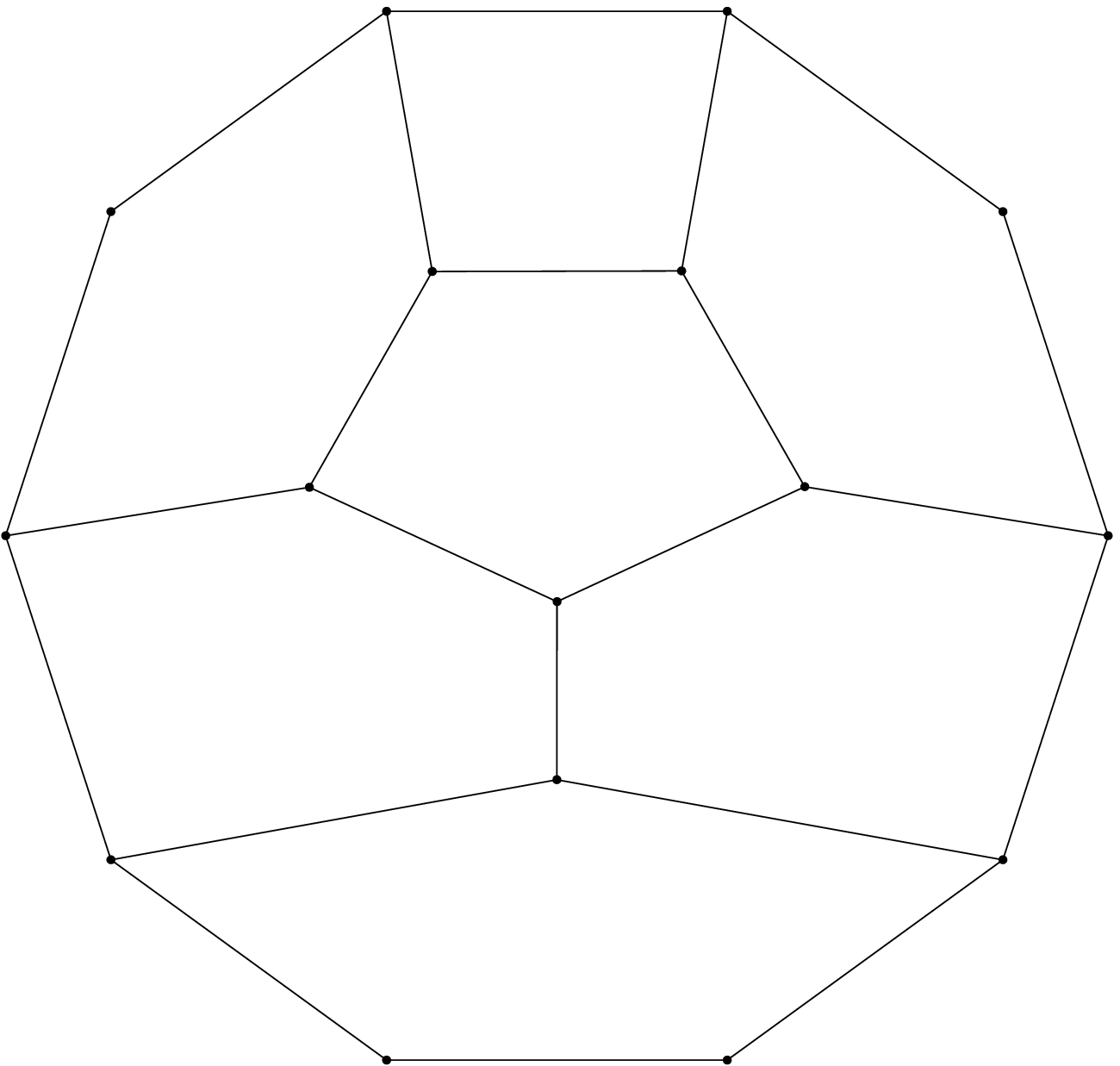}\par
$C_s$
\end{minipage}
\begin{minipage}{3cm}
\centering
\epsfig{height=20mm, file=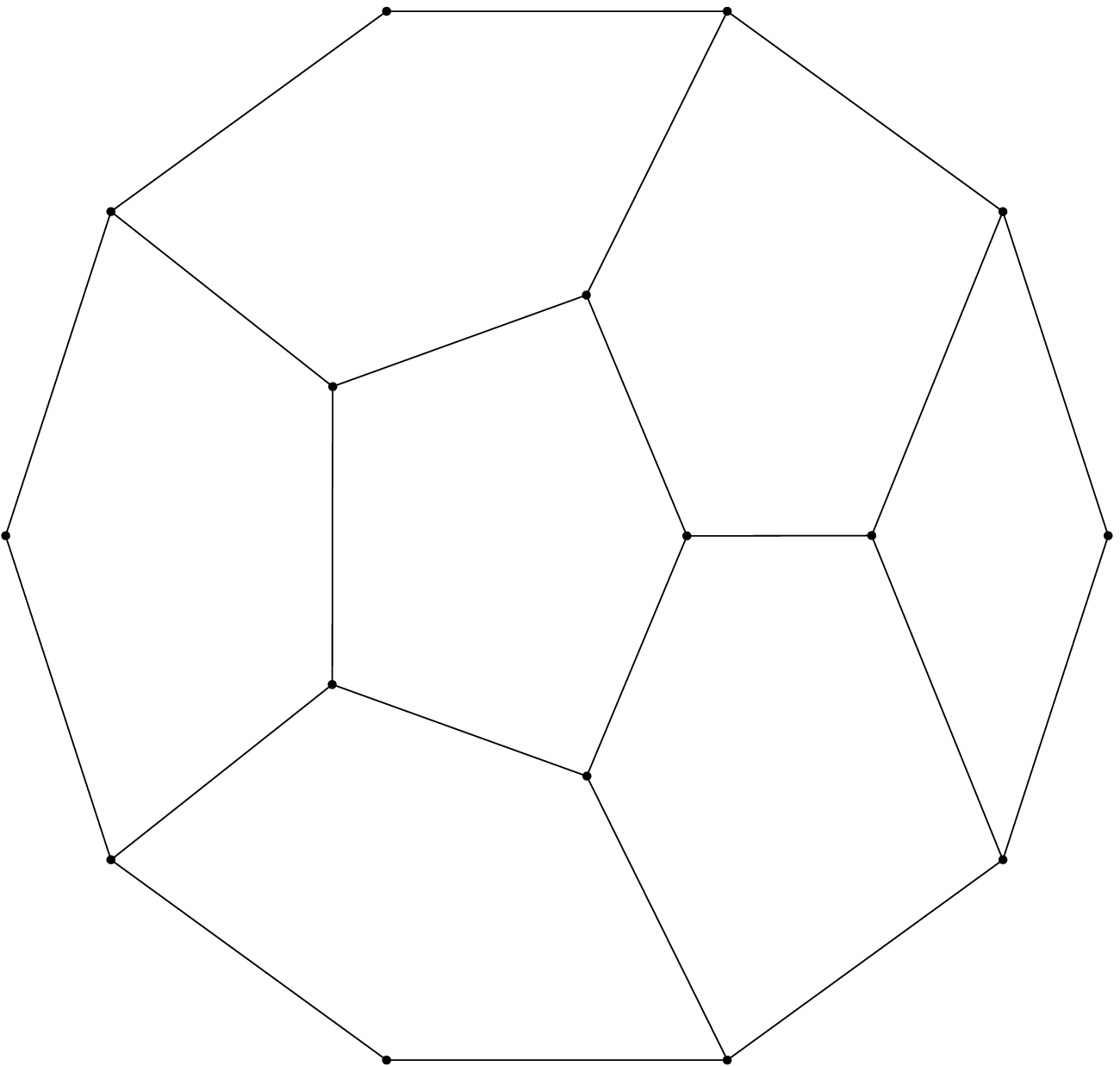}\par
$C_s$
\end{minipage}
\begin{minipage}{3cm}
\centering
\epsfig{height=20mm, file=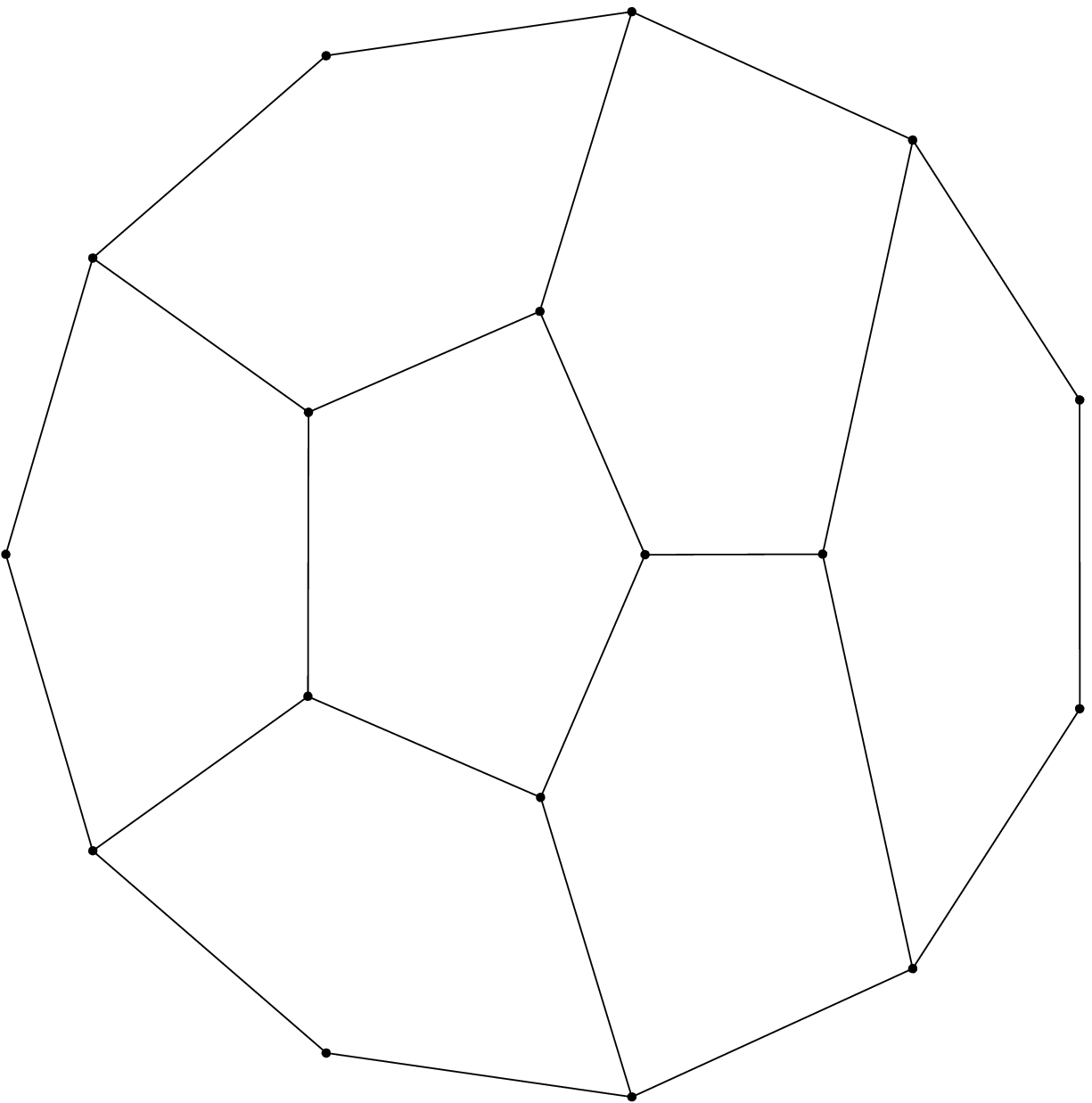}\par
$C_s$
\end{minipage}
\begin{minipage}{3cm}
\centering
\epsfig{height=20mm, file=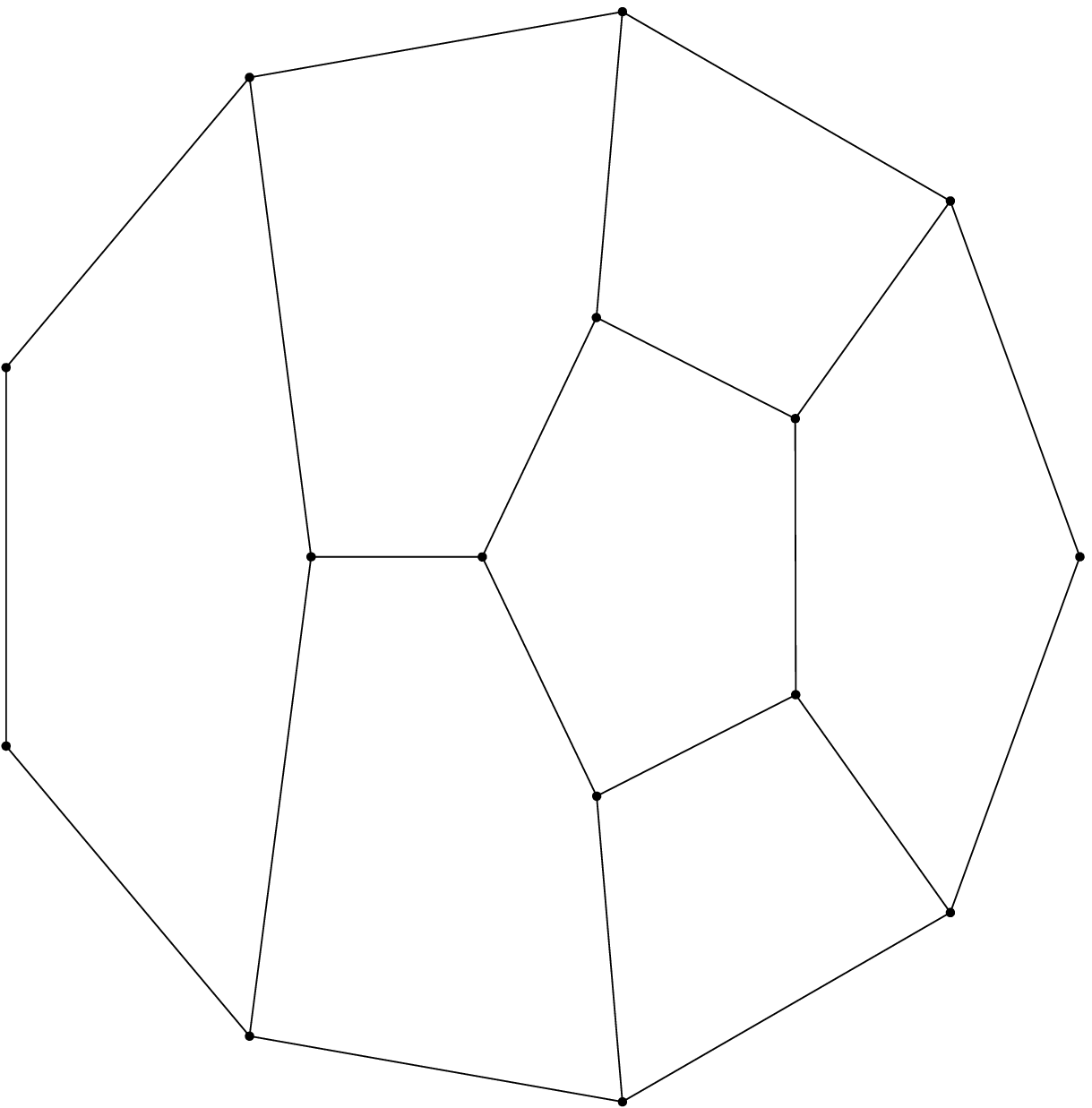}\par
$C_s$
\end{minipage}
\begin{minipage}{3cm}
\centering
\epsfig{height=20mm, file=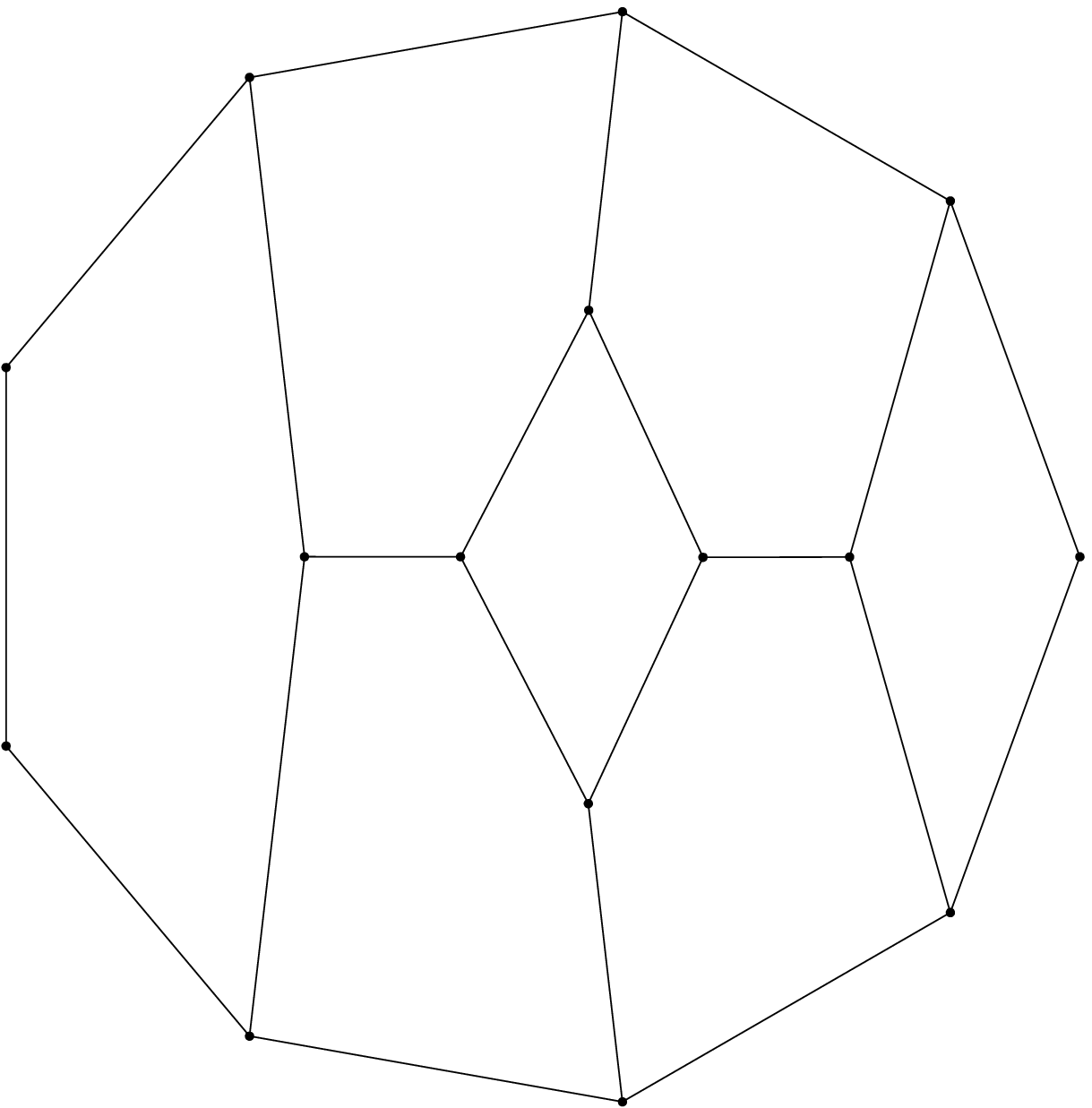}\par
$C_s$
\end{minipage}
\begin{minipage}{3cm}
\centering
\epsfig{height=20mm, file=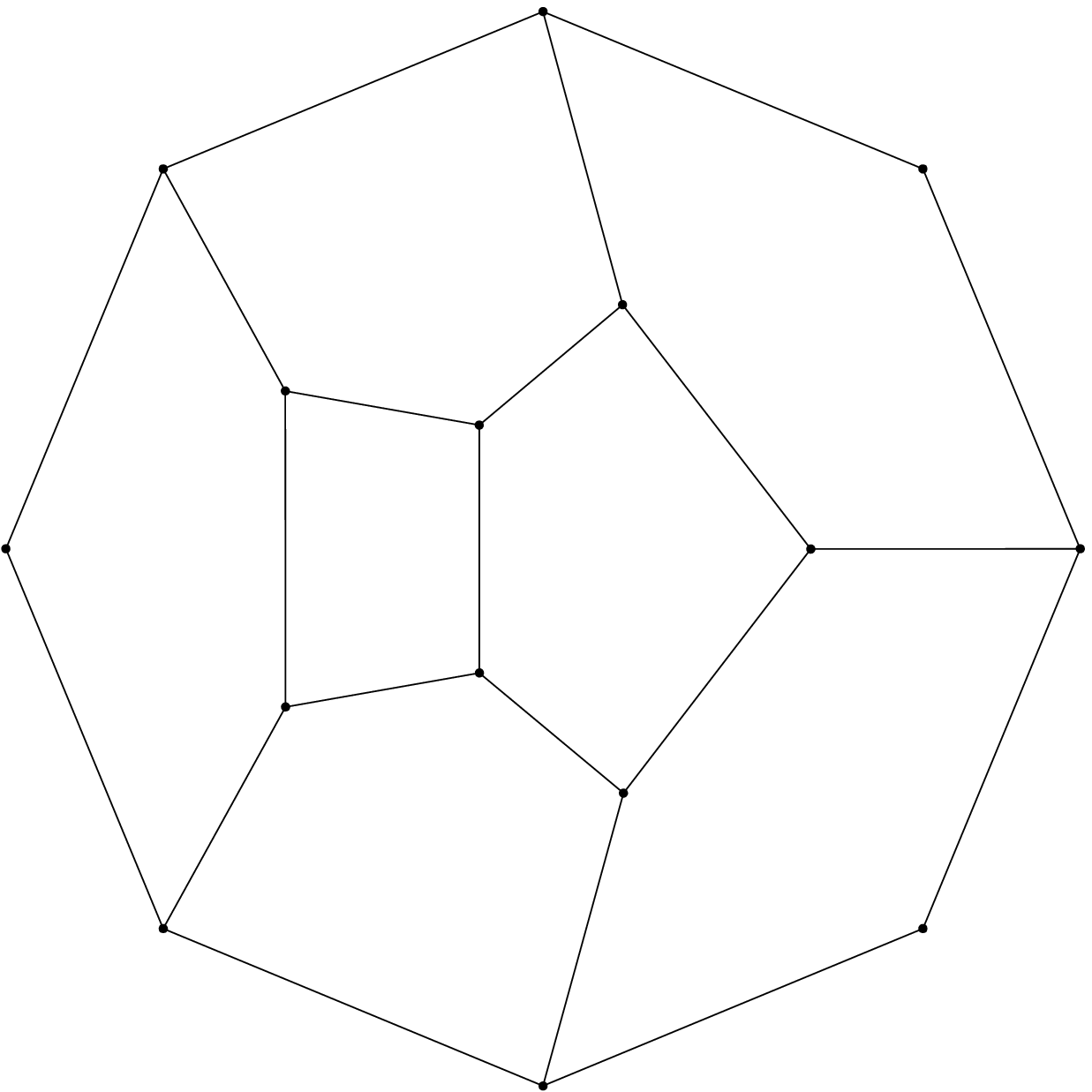}\par
$C_s$
\end{minipage}
\begin{minipage}{3cm}
\centering
\epsfig{height=20mm, file=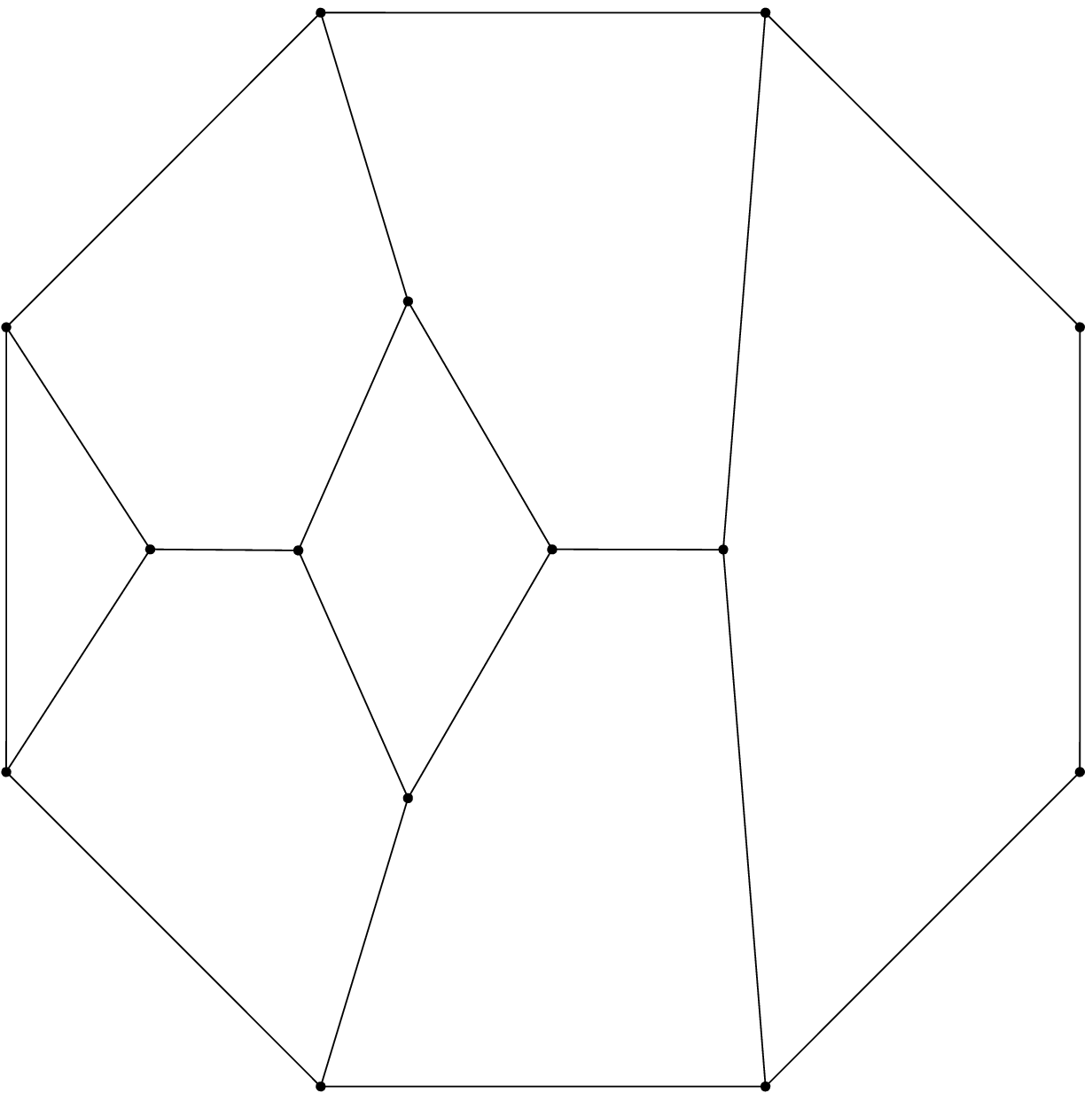}\par
$C_s$
\end{minipage}
\begin{minipage}{3cm}
\centering
\epsfig{height=20mm, file=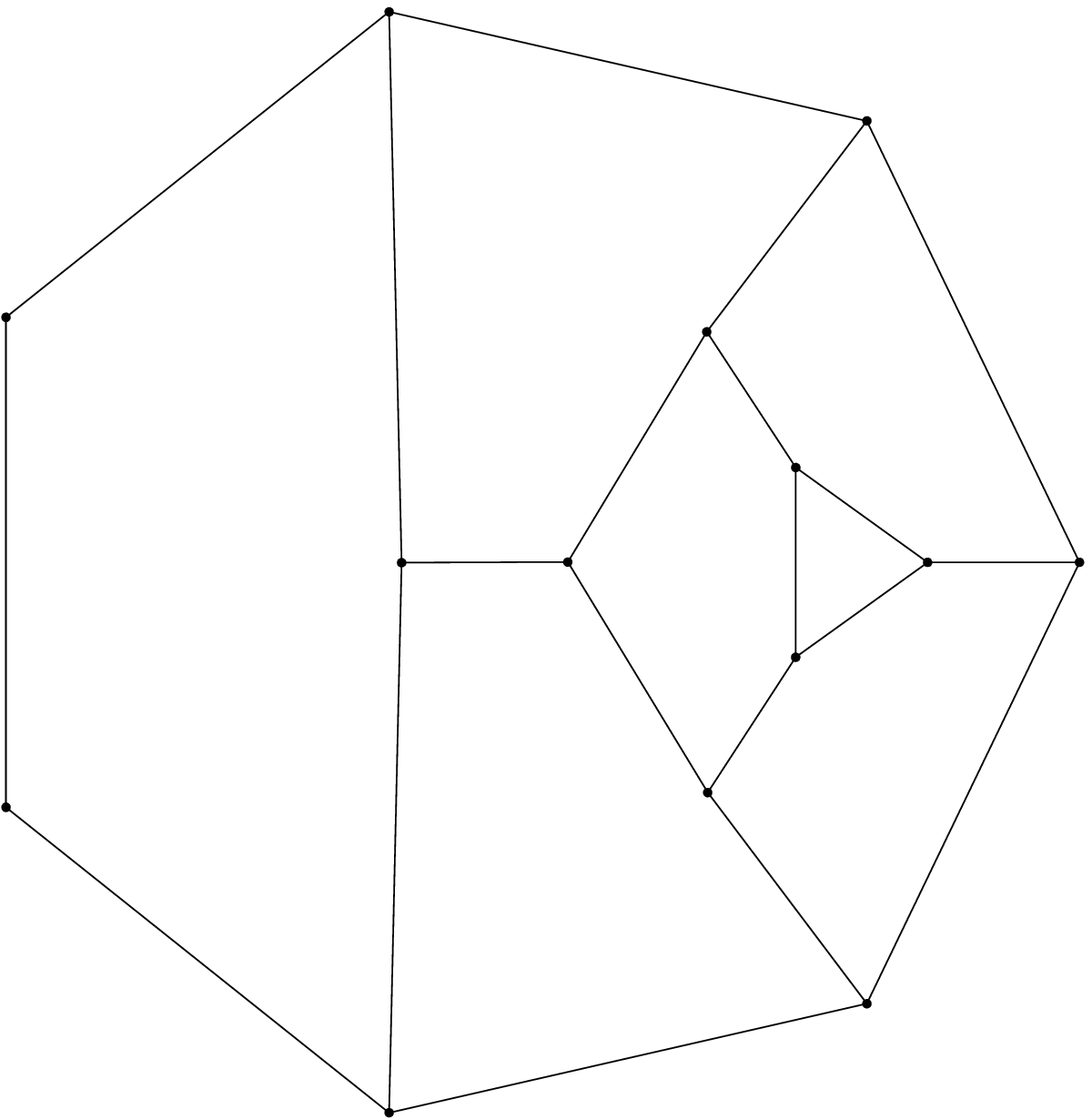}\par
$C_s$
\end{minipage}
\begin{minipage}{3cm}
\centering
\epsfig{height=20mm, file=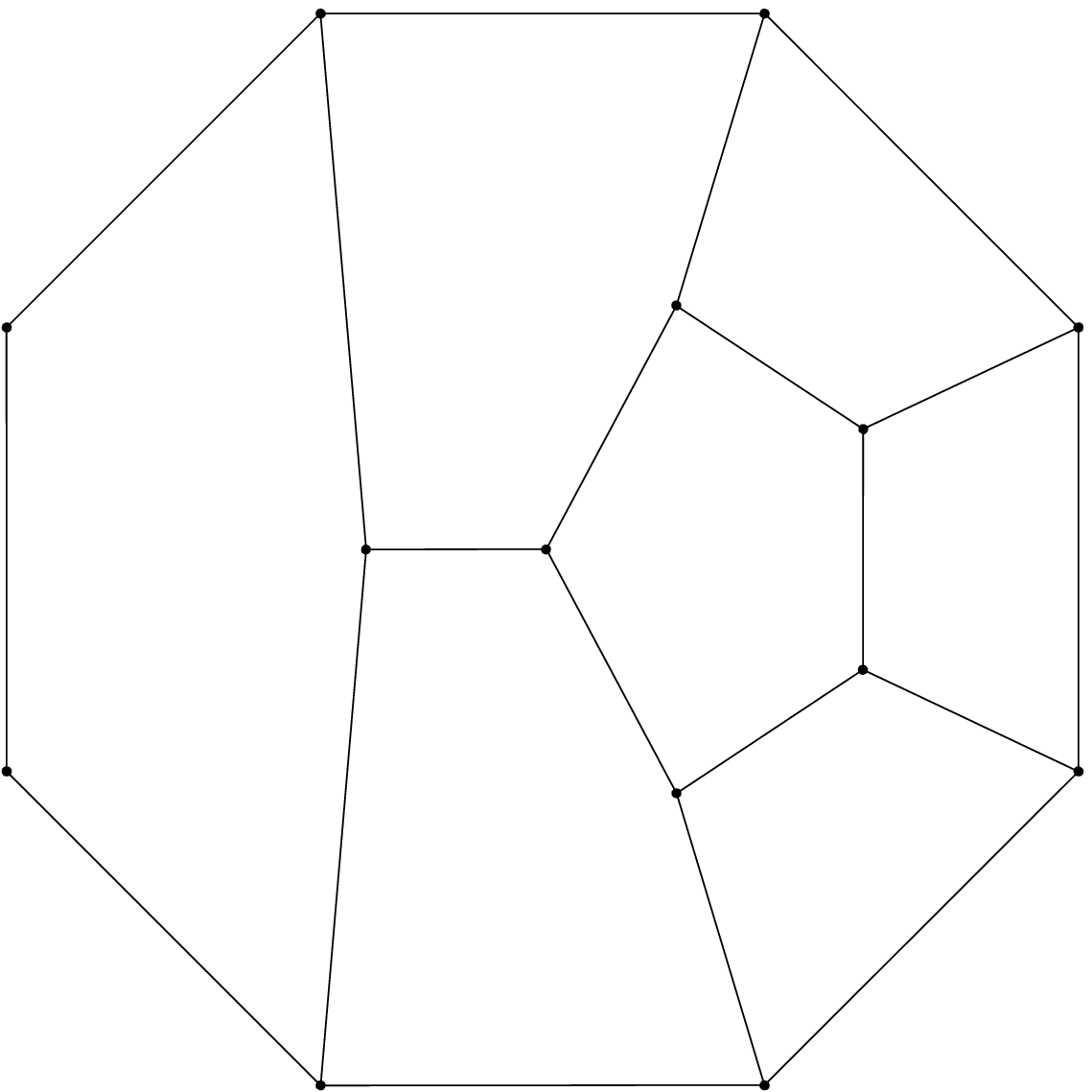}\par
$C_s$
\end{minipage}
\begin{minipage}{3cm}
\centering
\epsfig{height=20mm, file=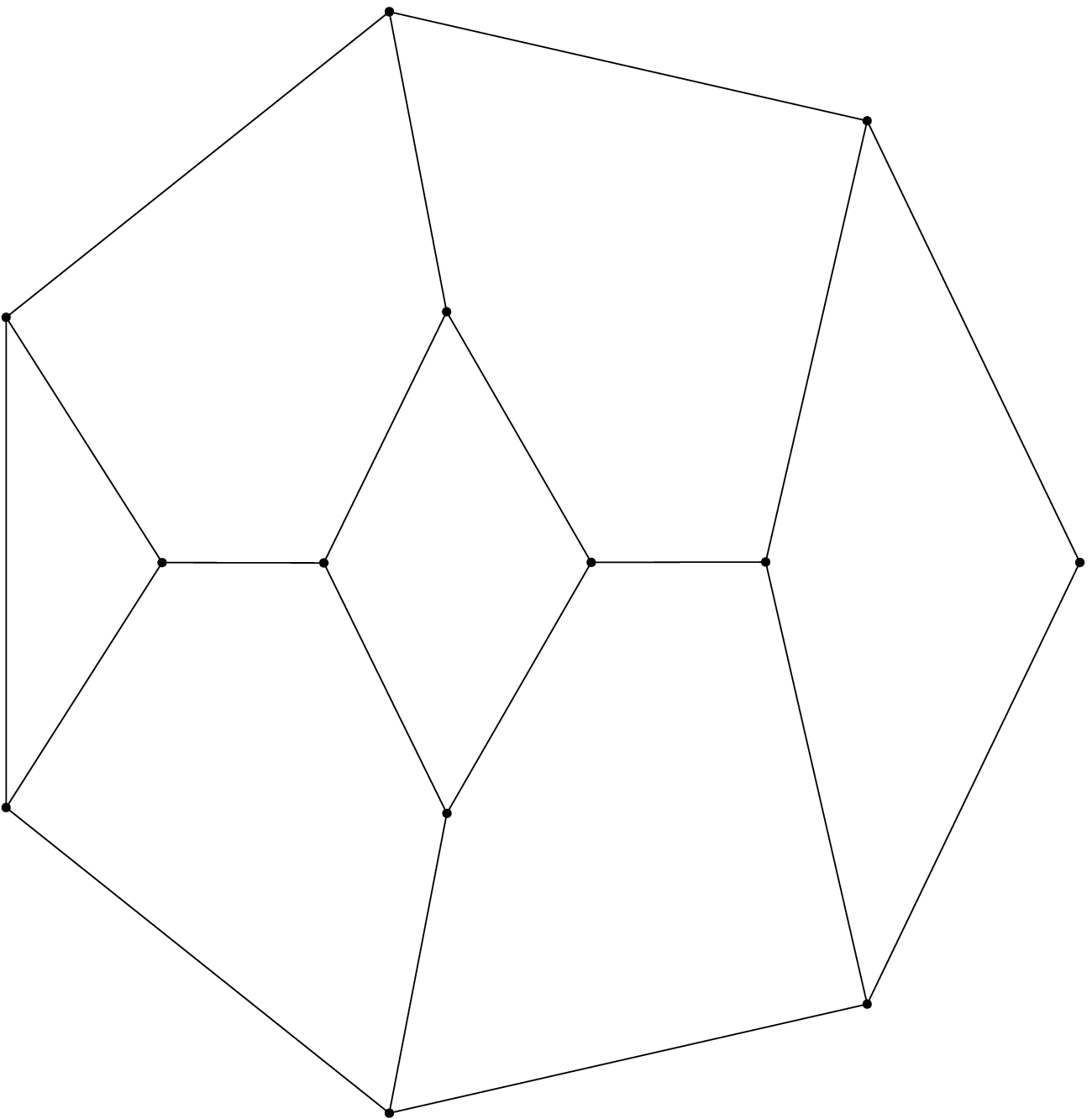}\par
$C_s$, nonext.
\end{minipage}
\begin{minipage}{3cm}
\centering
\epsfig{height=20mm, file=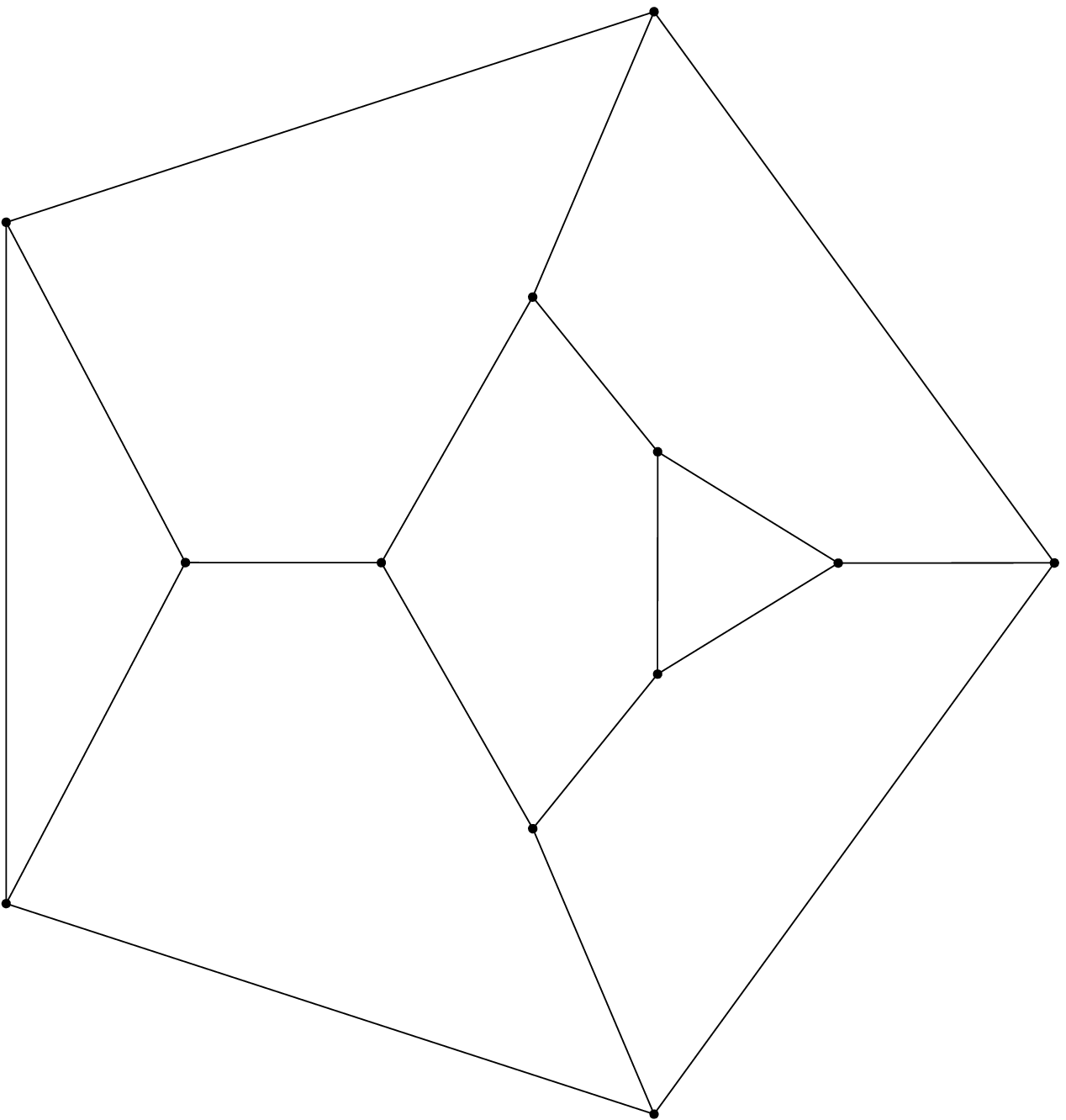}\par
$C_s$, nonext.
%PAIR1
\end{minipage}
\begin{minipage}{3cm}
\centering
\epsfig{height=20mm, file=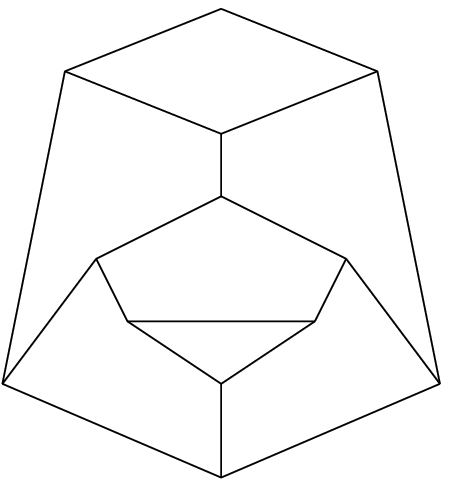}\par
$C_s$, nonext.
\end{minipage}
\begin{minipage}{3cm}
\centering
\epsfig{height=20mm, file=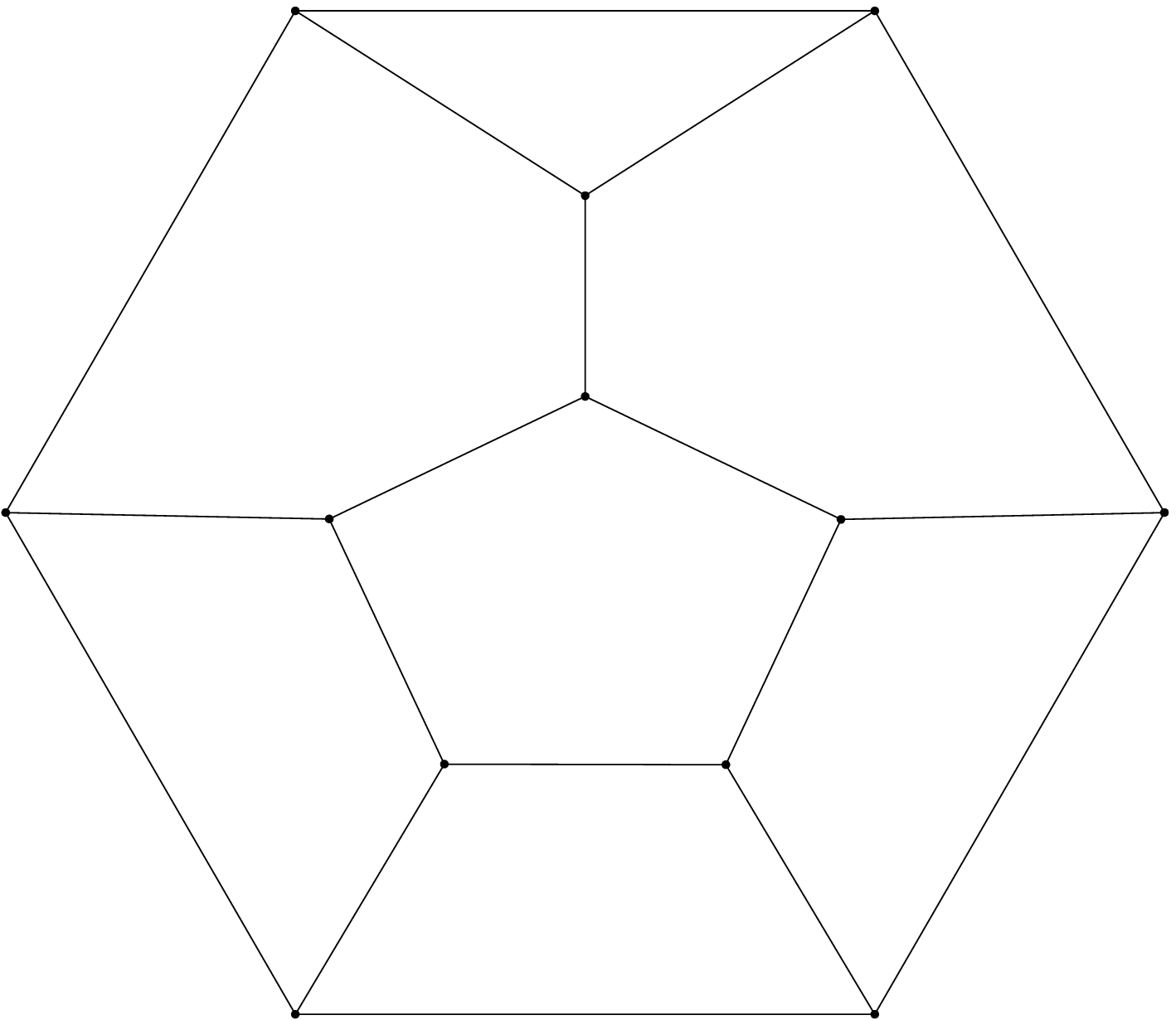}\par
$C_s$, nonext.
\end{minipage}
\begin{minipage}{3cm}
\centering
\epsfig{height=20mm, file=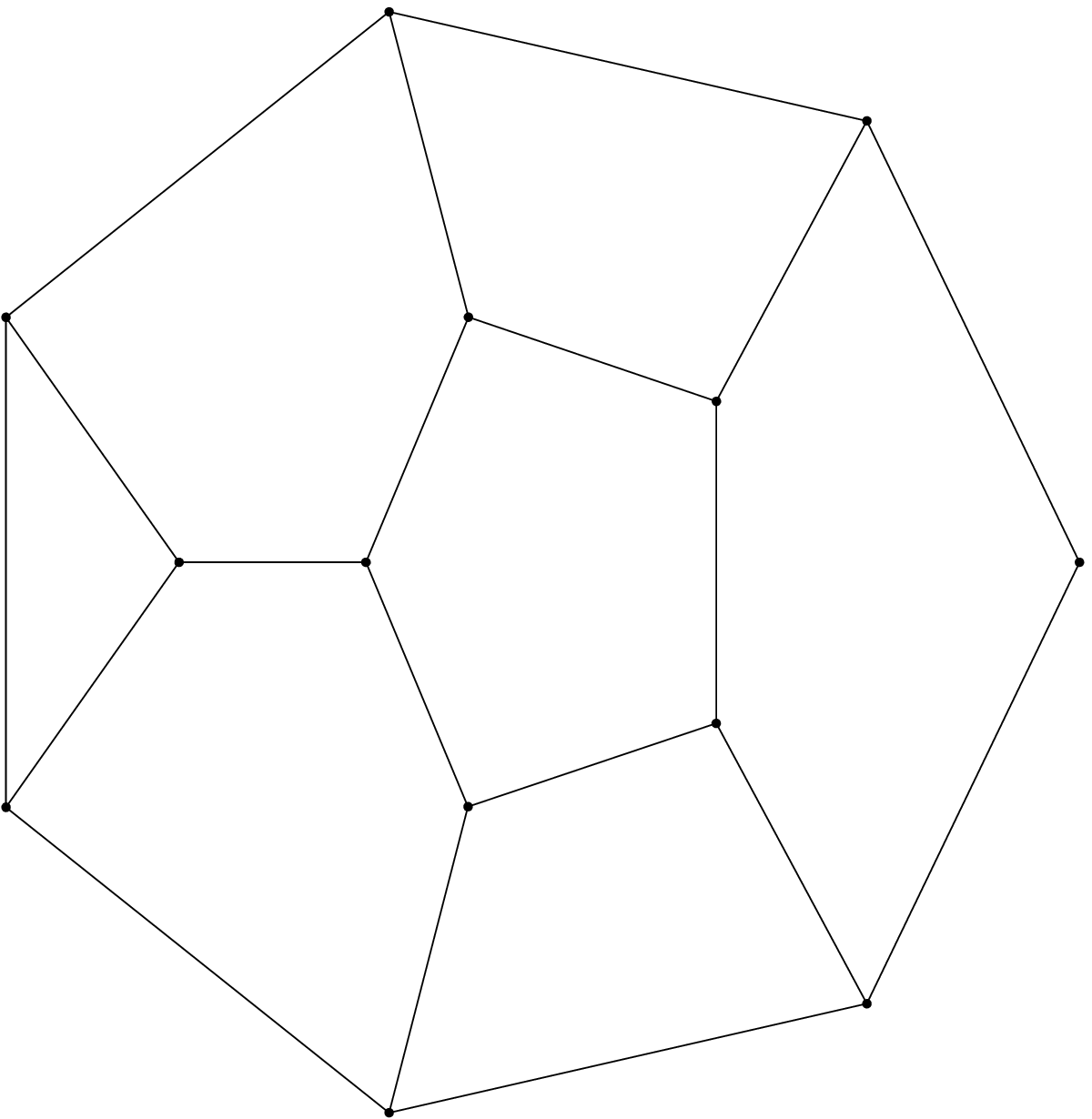}\par
$C_s$, nonext.
\end{minipage}
\begin{minipage}{3cm}
\centering
\epsfig{height=20mm, file=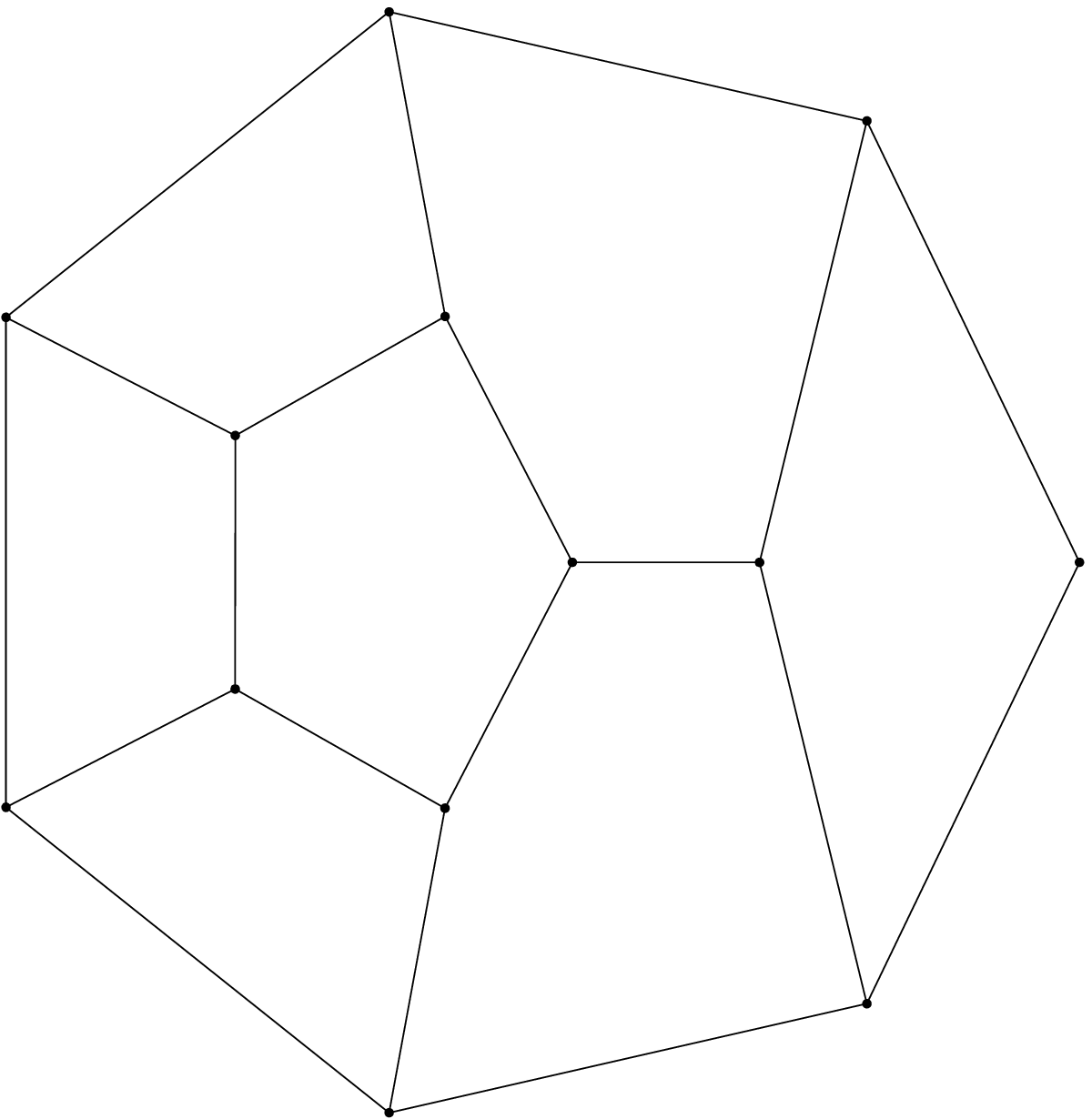}\par
$C_s$, nonext.
\end{minipage}
\begin{minipage}{3cm}
\centering
\epsfig{height=20mm, file=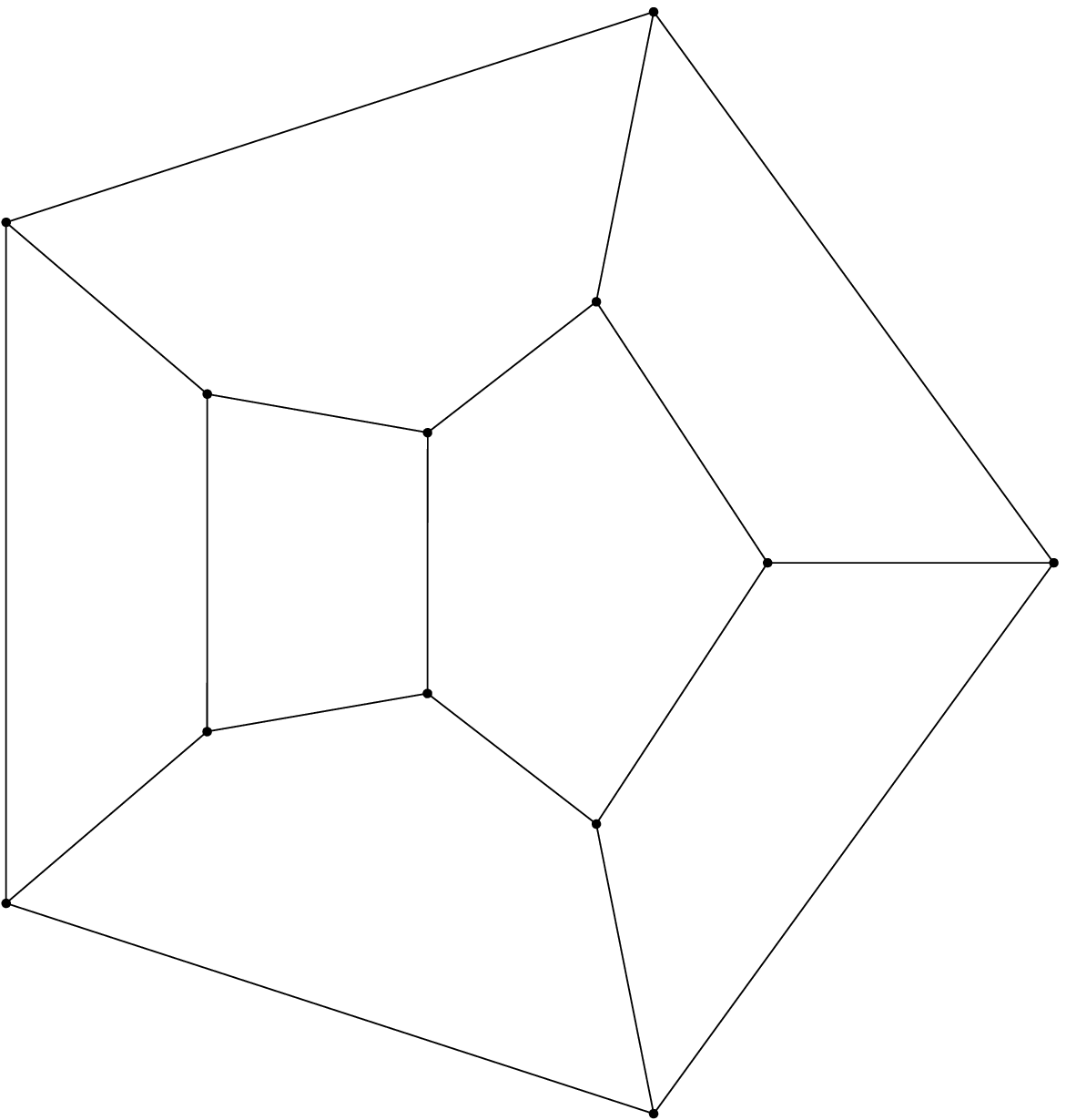}\par
$C_s$, nonext.
%PAIR2
\end{minipage}
\begin{minipage}{3cm}
\centering
\epsfig{height=20mm, file=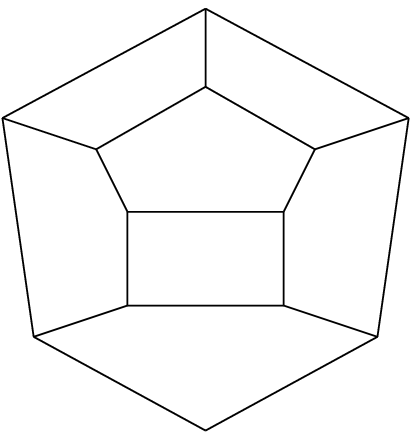}\par
$C_s$, nonext.
\end{minipage}
\begin{minipage}{3cm}
\centering
\epsfig{height=20mm, file=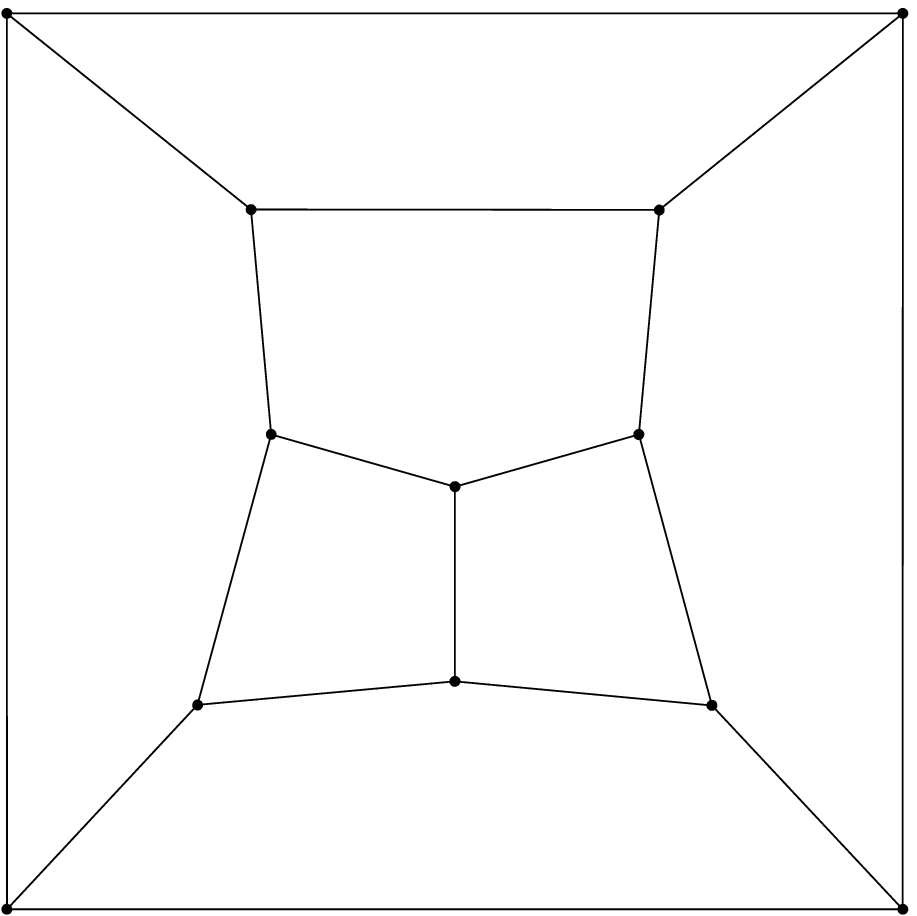}\par
$C_s$, nonext.
%PAIR2
\end{minipage}
\begin{minipage}{3cm}
\centering
\epsfig{height=20mm, file=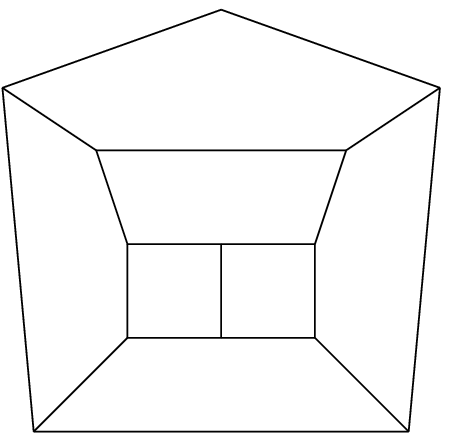}\par
$C_s$,~nonext.~$(C_{2\nu})$
\end{minipage}
\begin{minipage}{3cm}
\centering
\epsfig{height=20mm, file=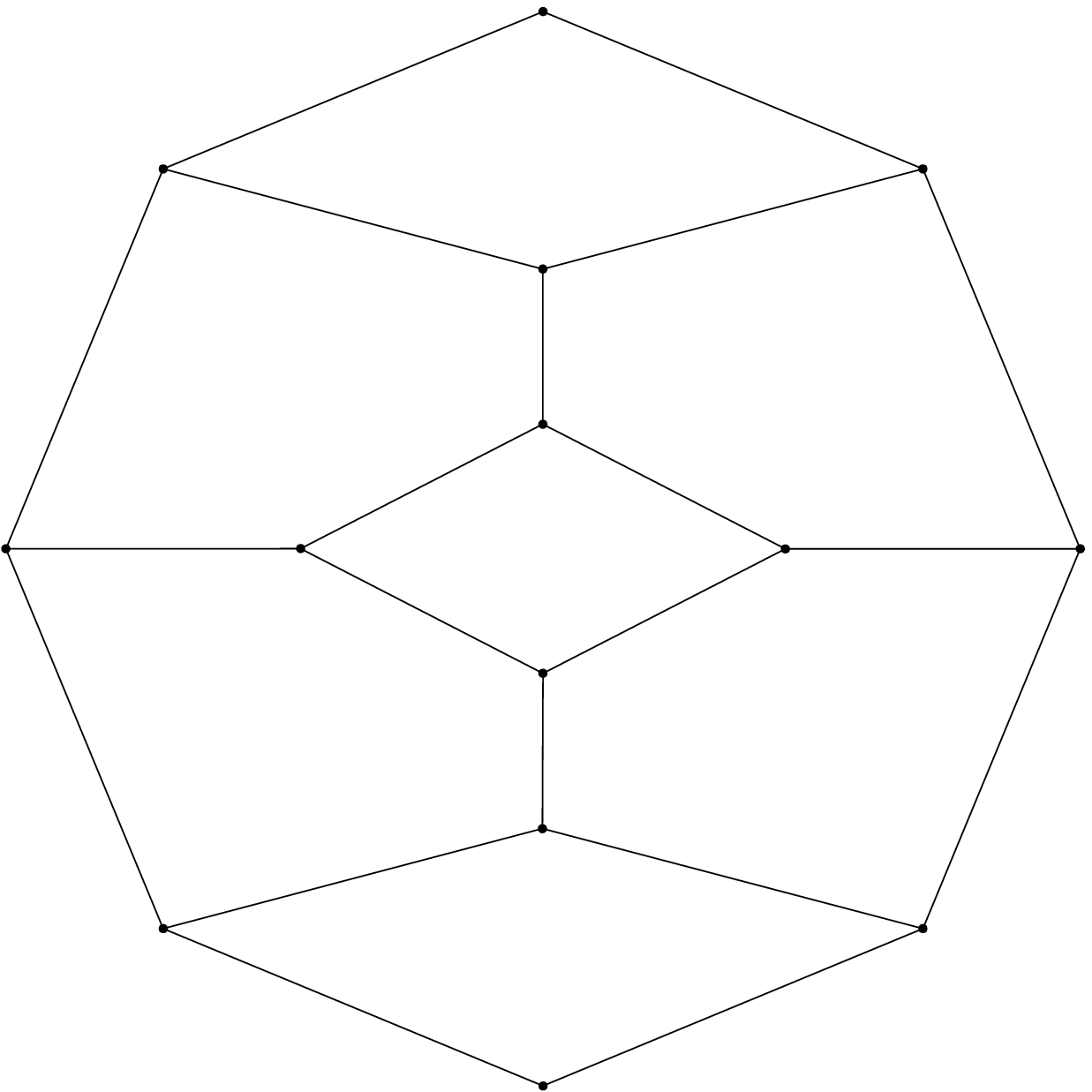}\par
$C_{2\nu}$
\end{minipage}
\begin{minipage}{3cm}
\centering
\epsfig{height=20mm, file=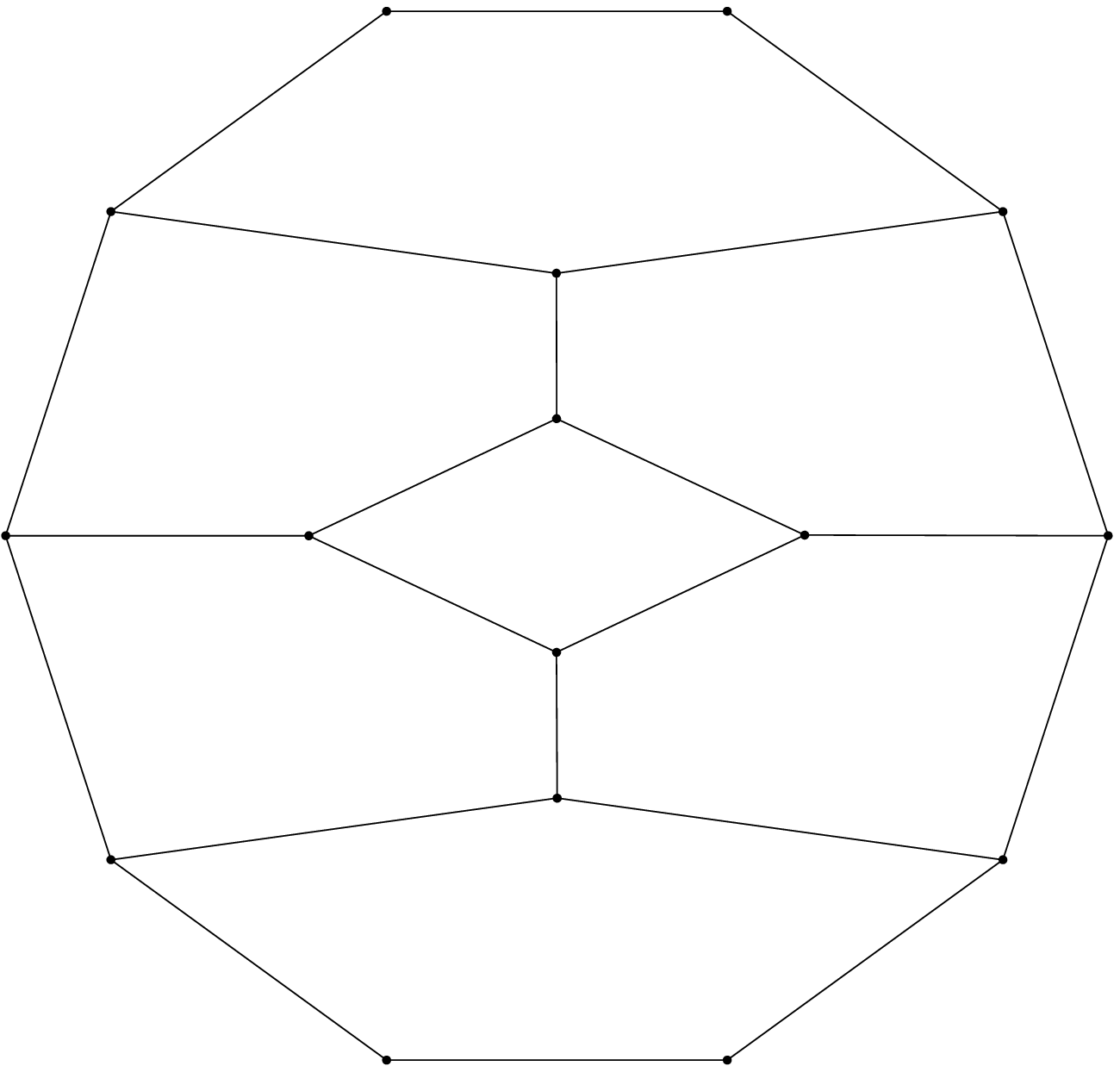}\par
$C_{2\nu}$
\end{minipage}
\begin{minipage}{3cm}
\centering
\epsfig{height=20mm, file=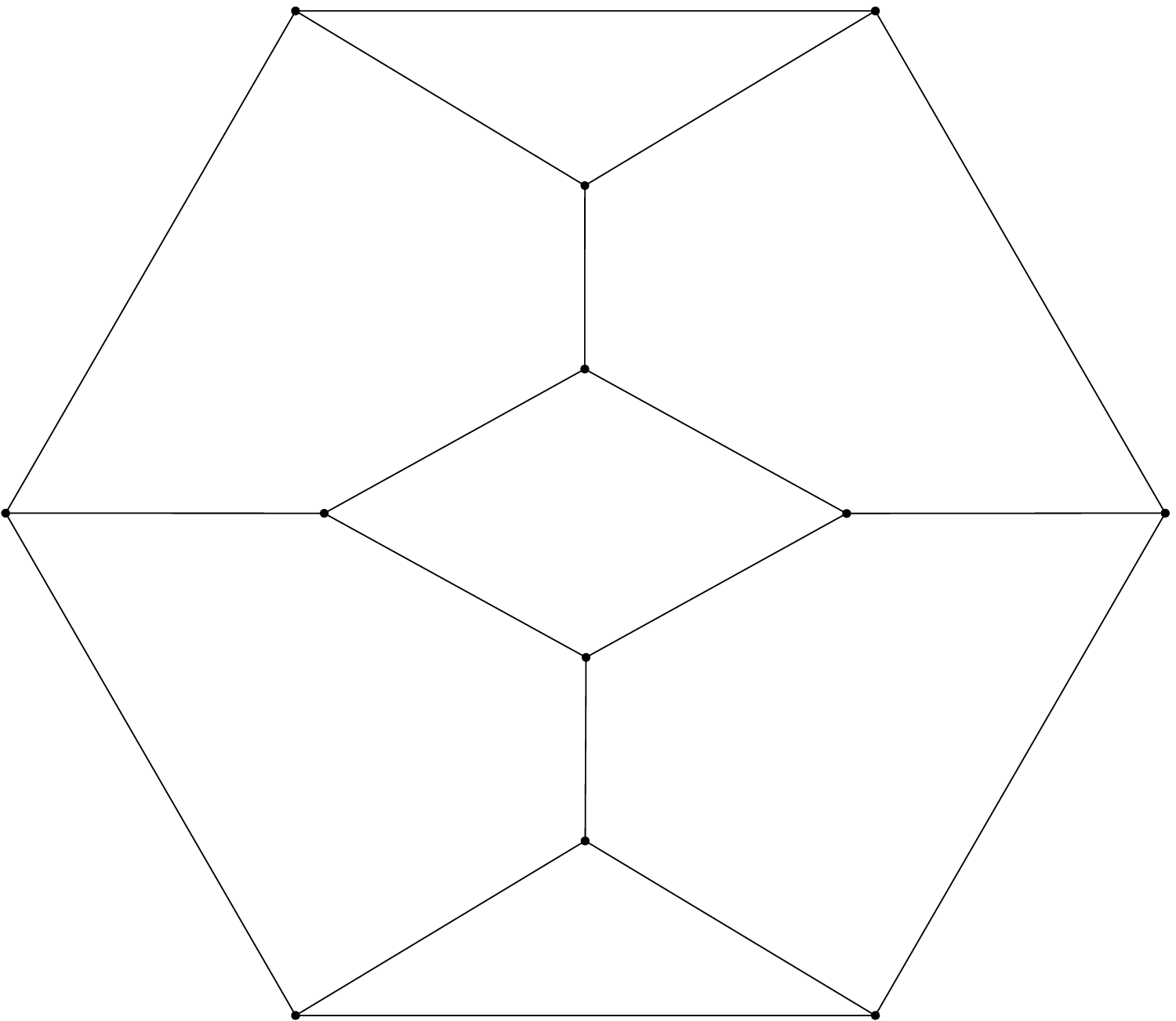}\par
$C_{2\nu}$, nonext.
\end{minipage}
\begin{minipage}{3cm}
\centering
\epsfig{height=20mm, file=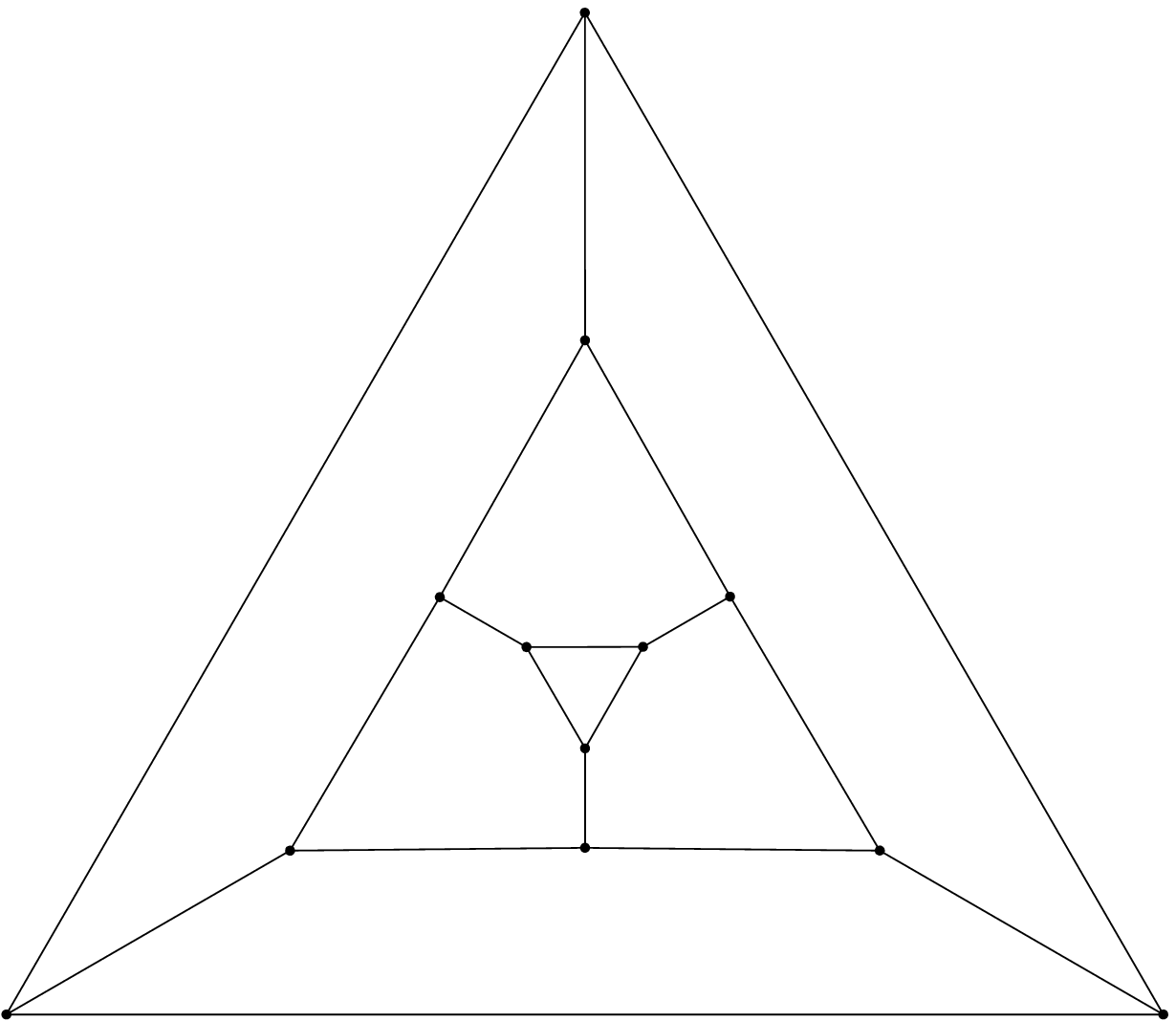}\par
$C_{3\nu}$, nonext.
%PAIR1
\end{minipage}

\end{center}
List of sporadic elementary $(\{3,4,5\},3)$-polycycles with $8$ faces:
\begin{center}
\begin{minipage}{3cm}
\centering
\epsfig{height=20mm, file=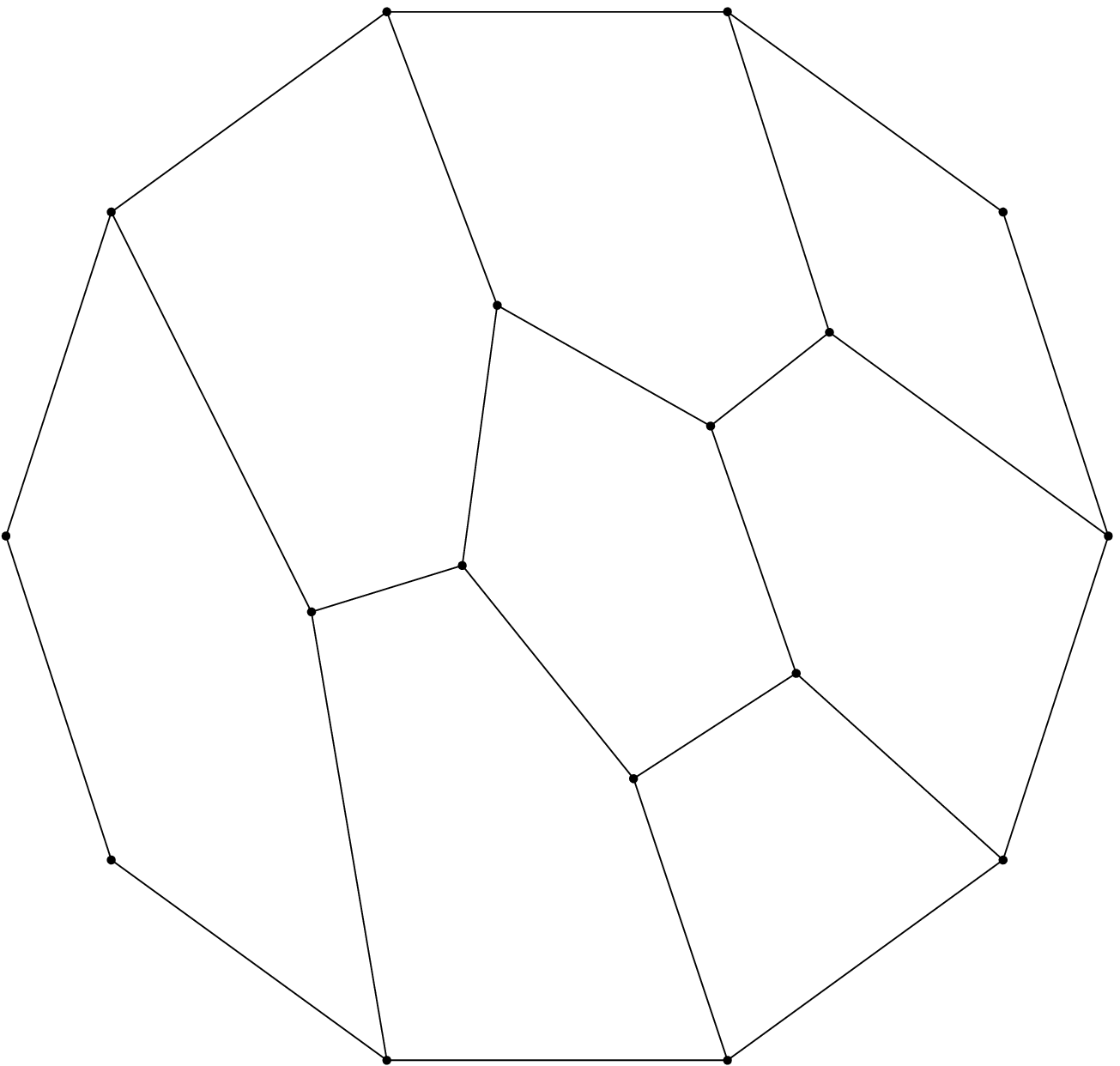}\par
$C_1$
\end{minipage}
\begin{minipage}{3cm}
\centering
\epsfig{height=20mm, file=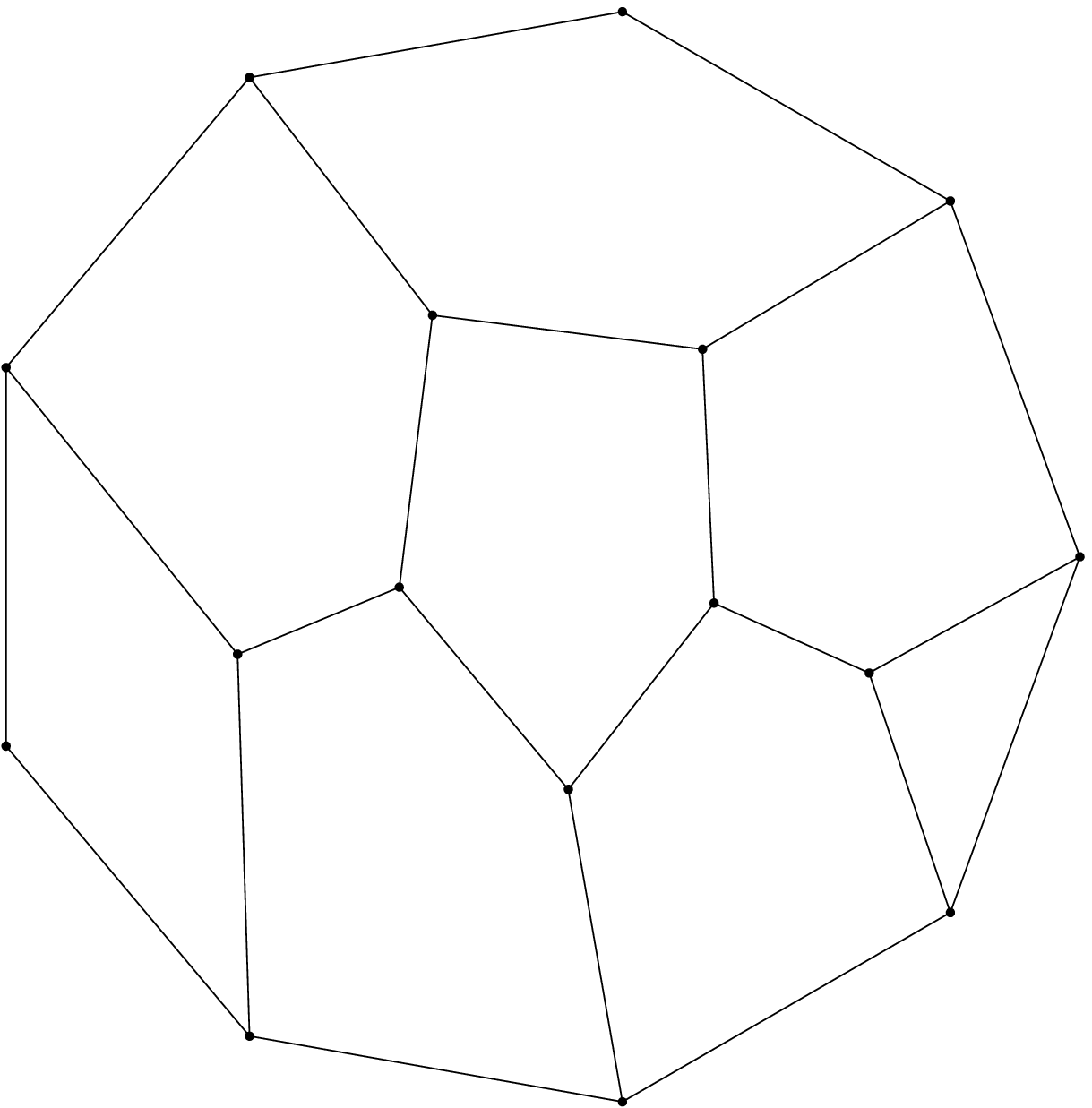}\par
$C_1$
\end{minipage}
\begin{minipage}{3cm}
\centering
\epsfig{height=20mm, file=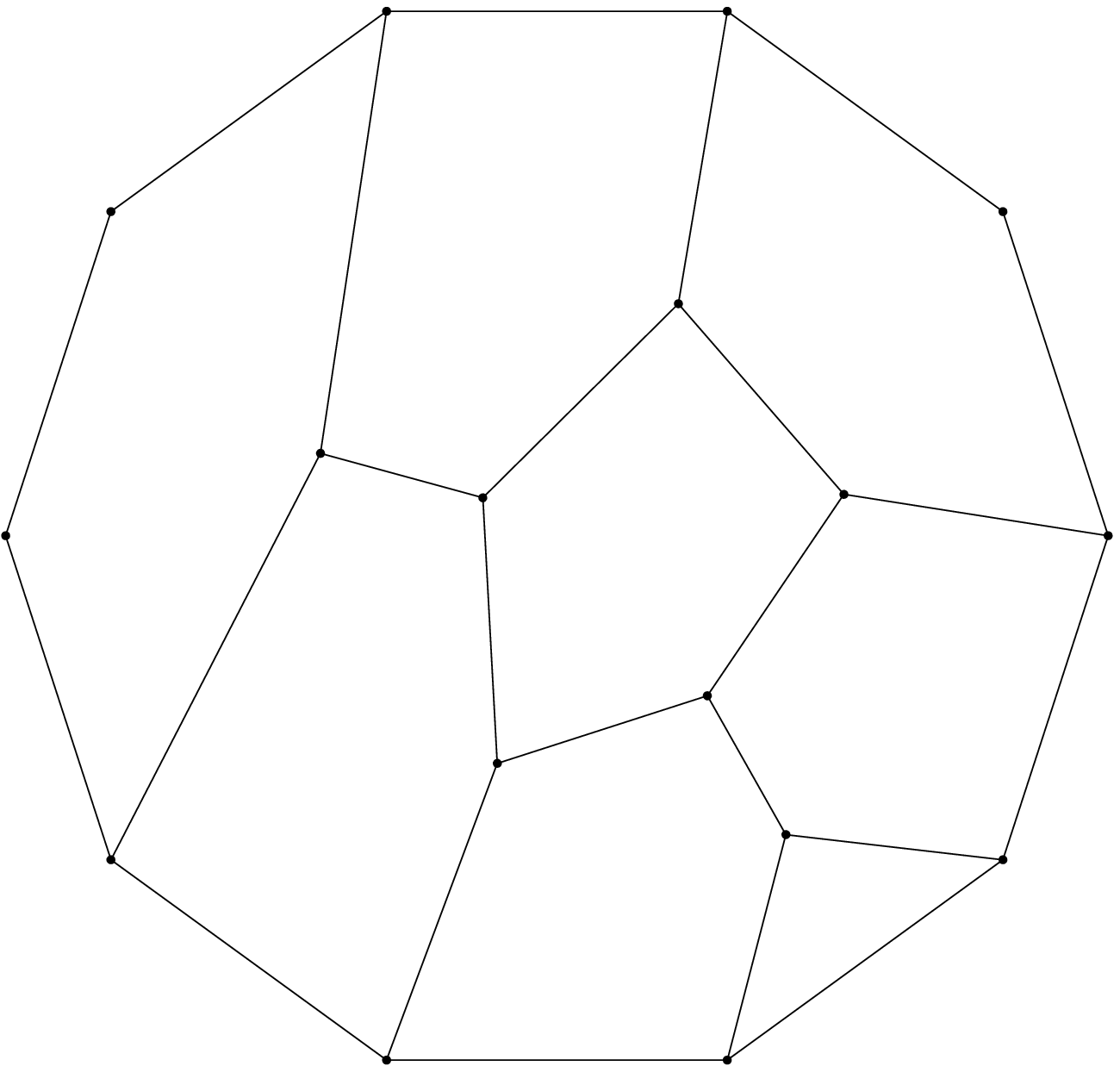}\par
$C_1$
\end{minipage}
\begin{minipage}{3cm}
\centering
\epsfig{height=20mm, file=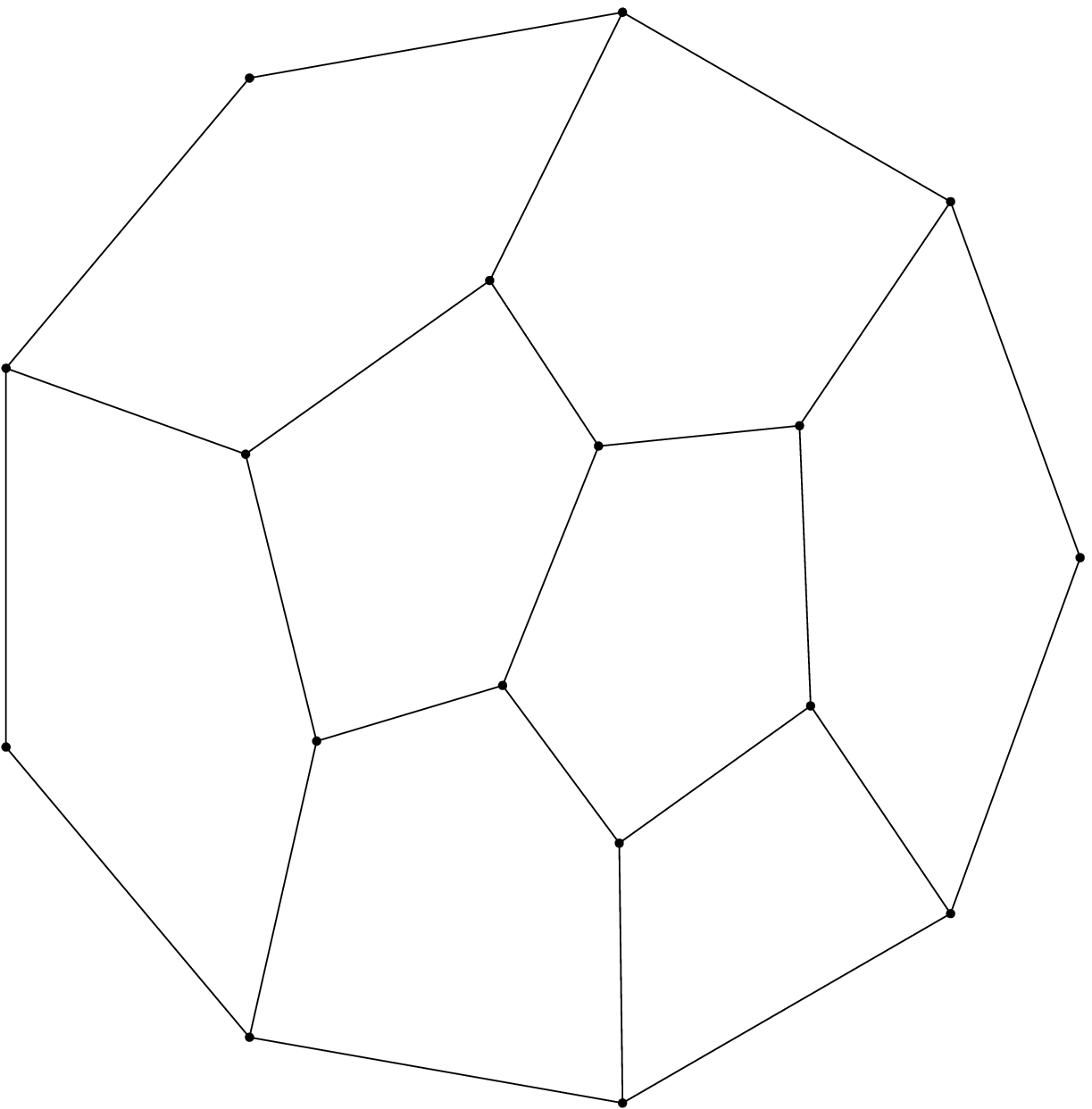}\par
$C_1$
\end{minipage}
\begin{minipage}{3cm}
\centering
\epsfig{height=20mm, file=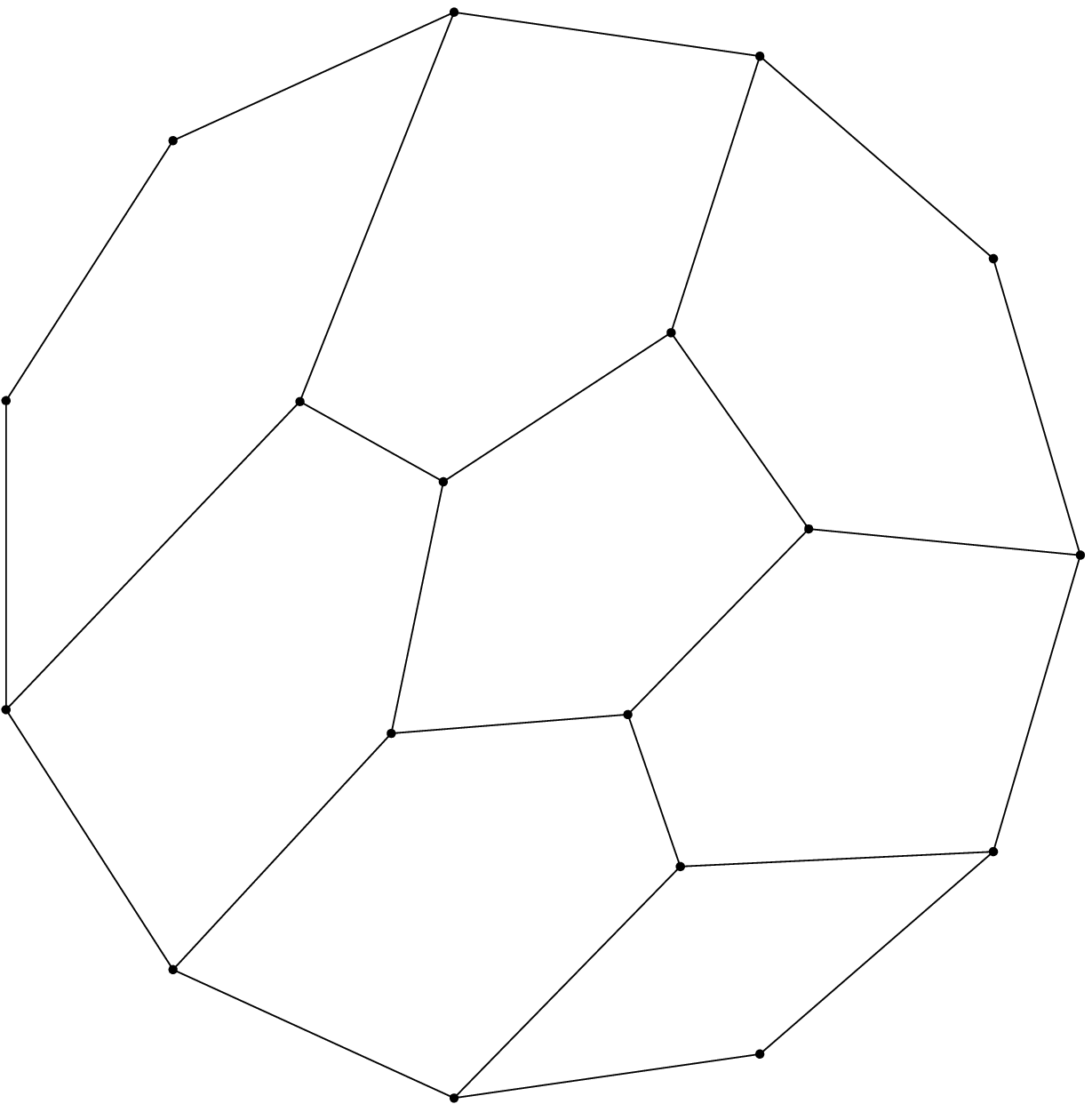}\par
$C_1$
\end{minipage}
\begin{minipage}{3cm}
\centering
\epsfig{height=20mm, file=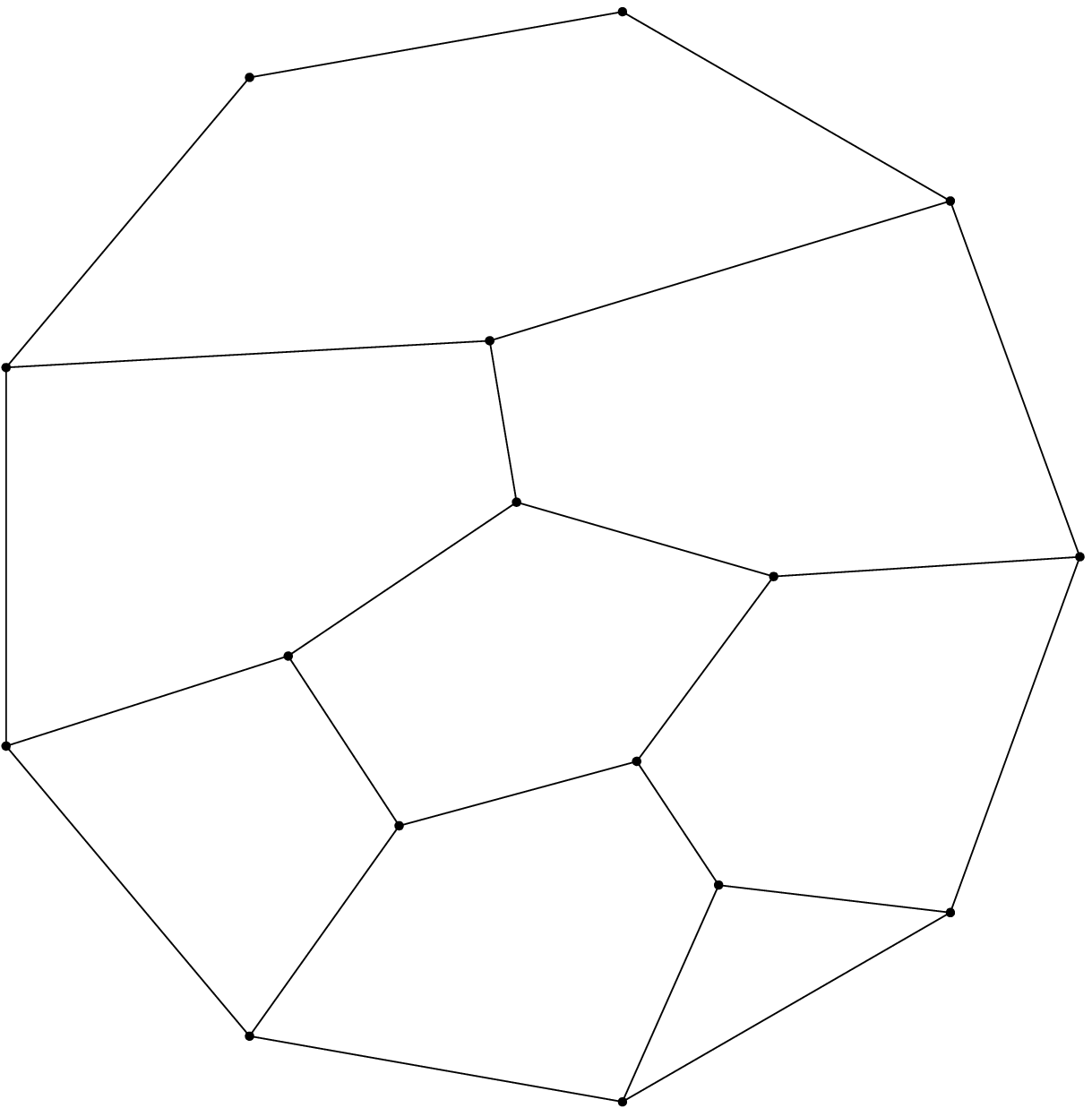}\par
$C_1$
\end{minipage}
\begin{minipage}{3cm}
\centering
\epsfig{height=20mm, file=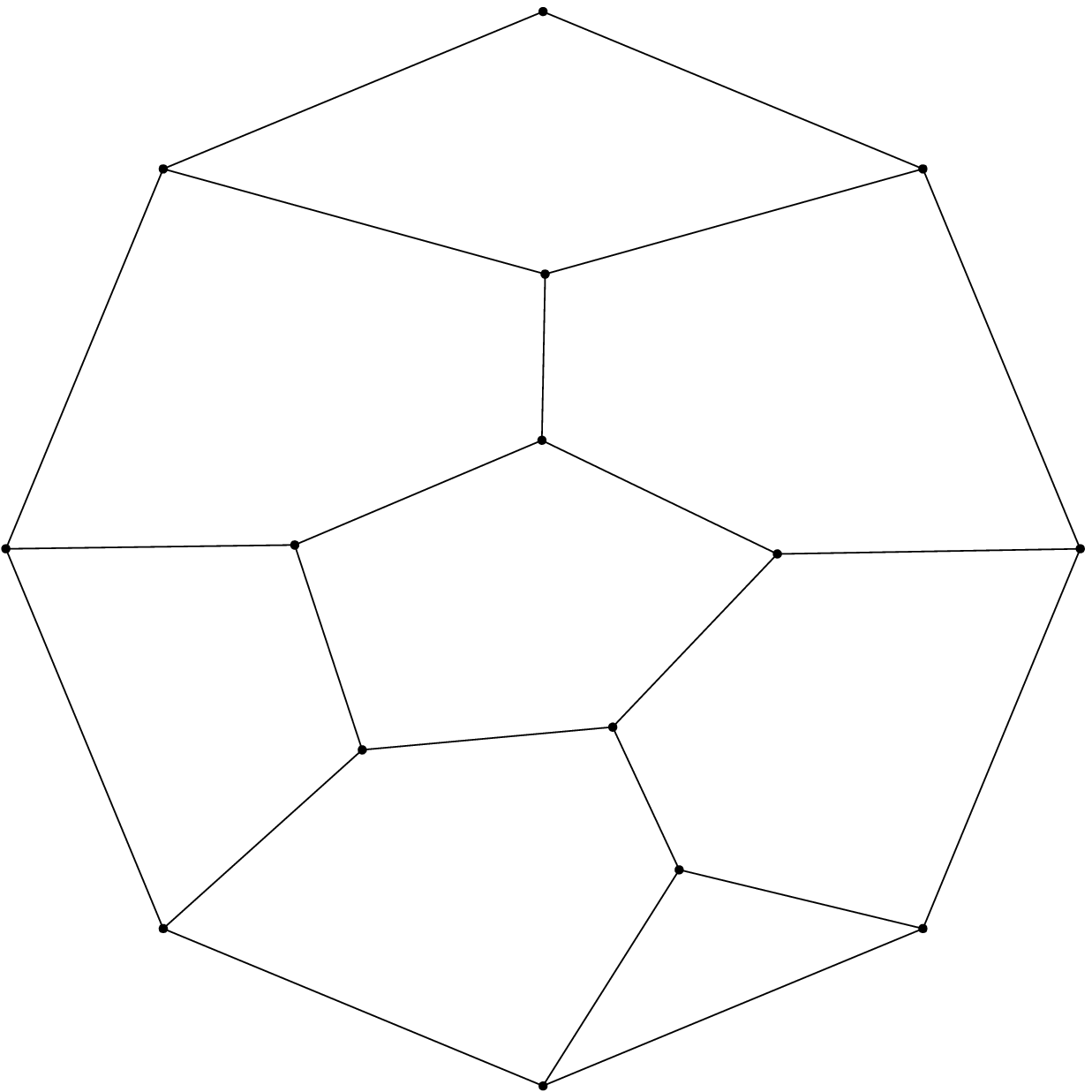}\par
$C_1$, nonext.
\end{minipage}
\begin{minipage}{3cm}
\centering
\epsfig{height=20mm, file=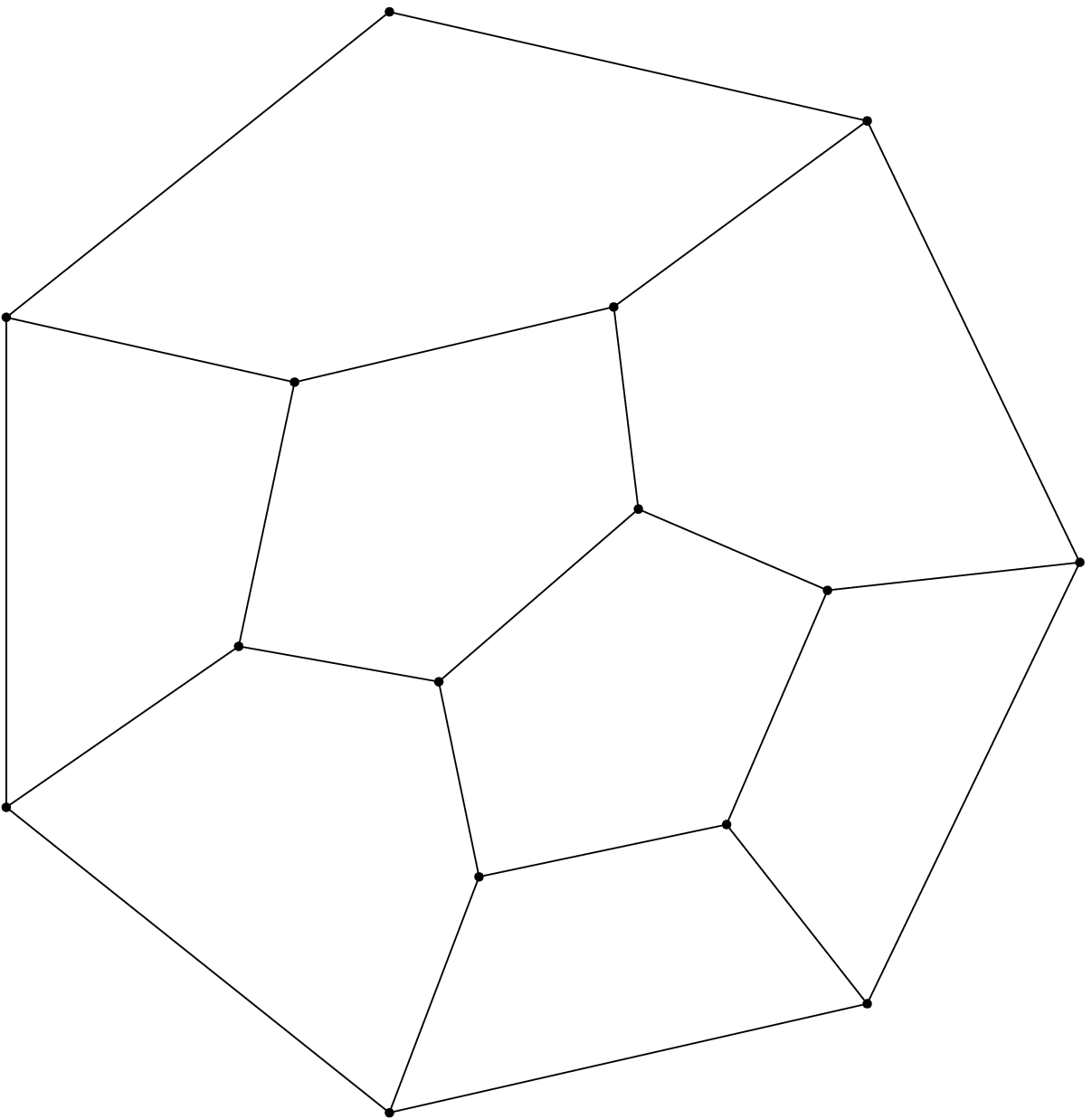}\par
$C_1$, nonext.
\end{minipage}
\begin{minipage}{3cm}
\centering
\epsfig{height=20mm, file=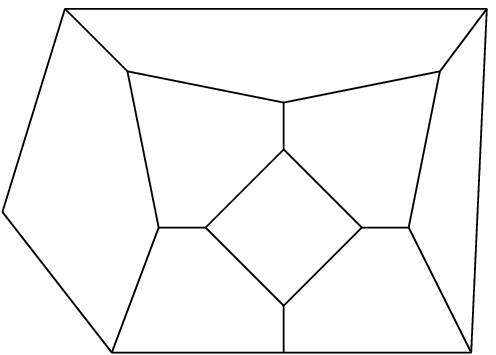}\par
$C_1$, nonext.
\end{minipage}
\begin{minipage}{3cm}
\centering
\epsfig{height=20mm, file=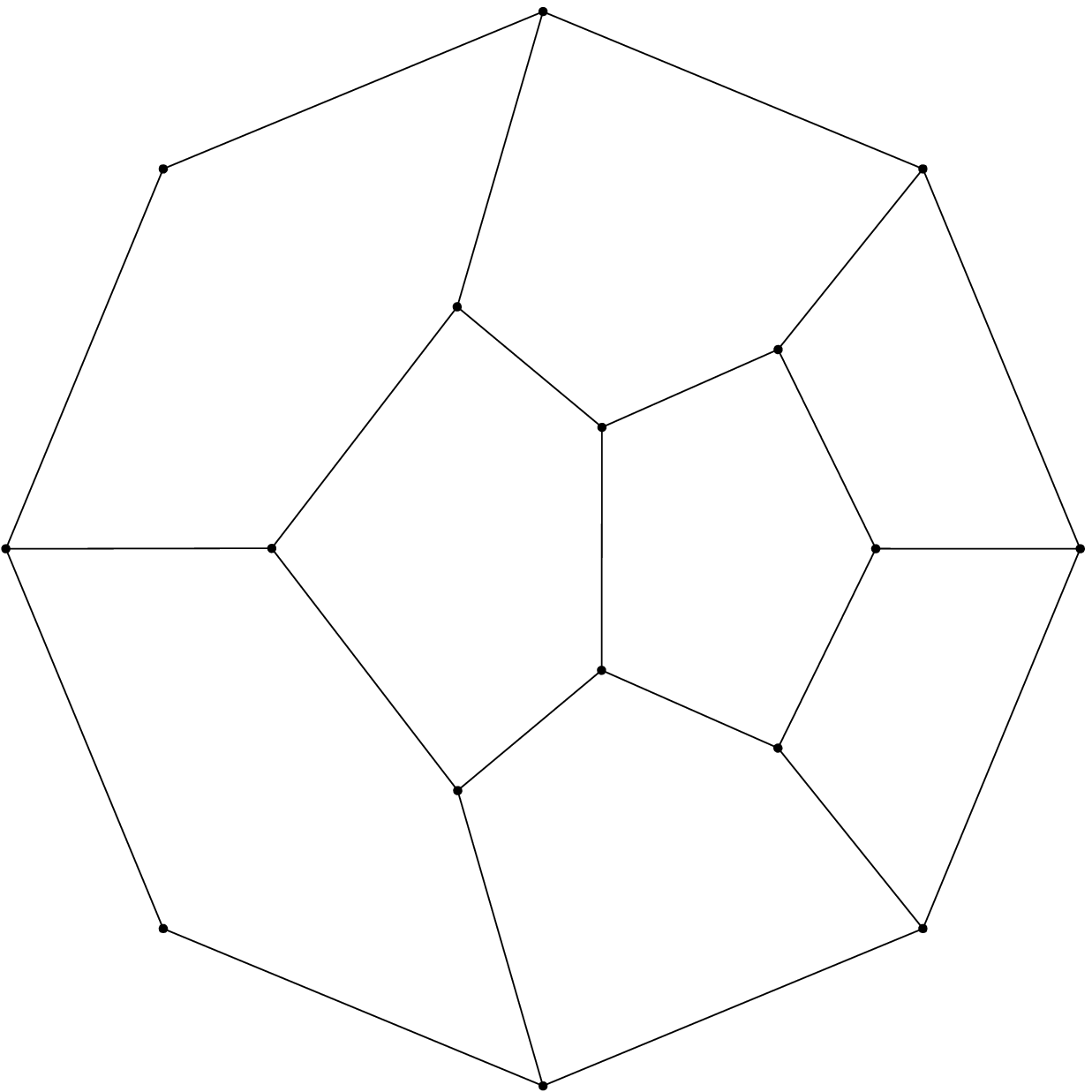}\par
$C_s$
\end{minipage}
\begin{minipage}{3cm}
\centering
\epsfig{height=20mm, file=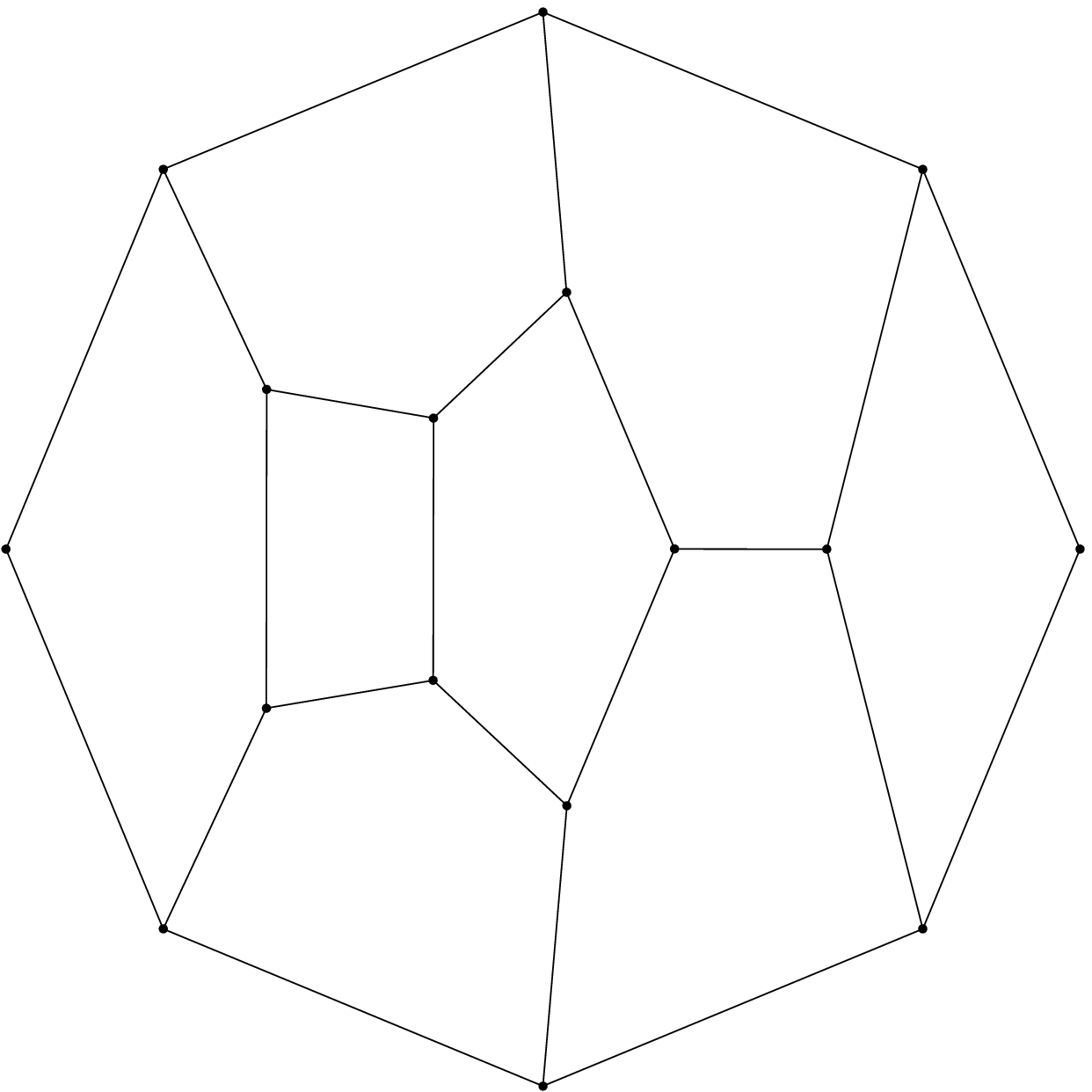}\par
$C_s$
\end{minipage}
\begin{minipage}{3cm}
\centering
\epsfig{height=20mm, file=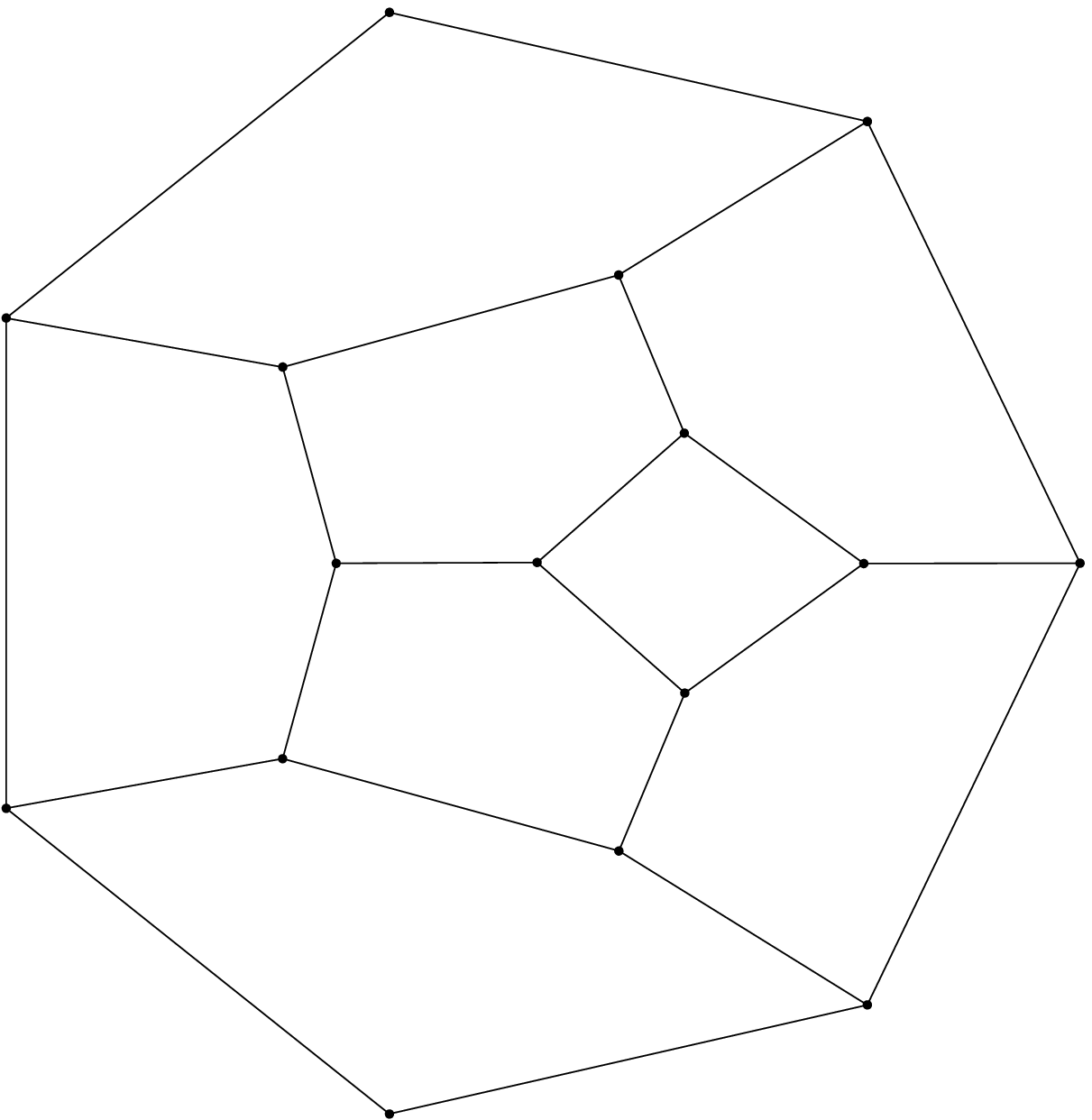}\par
$C_s$
\end{minipage}
\begin{minipage}{3cm}
\centering
\epsfig{height=20mm, file=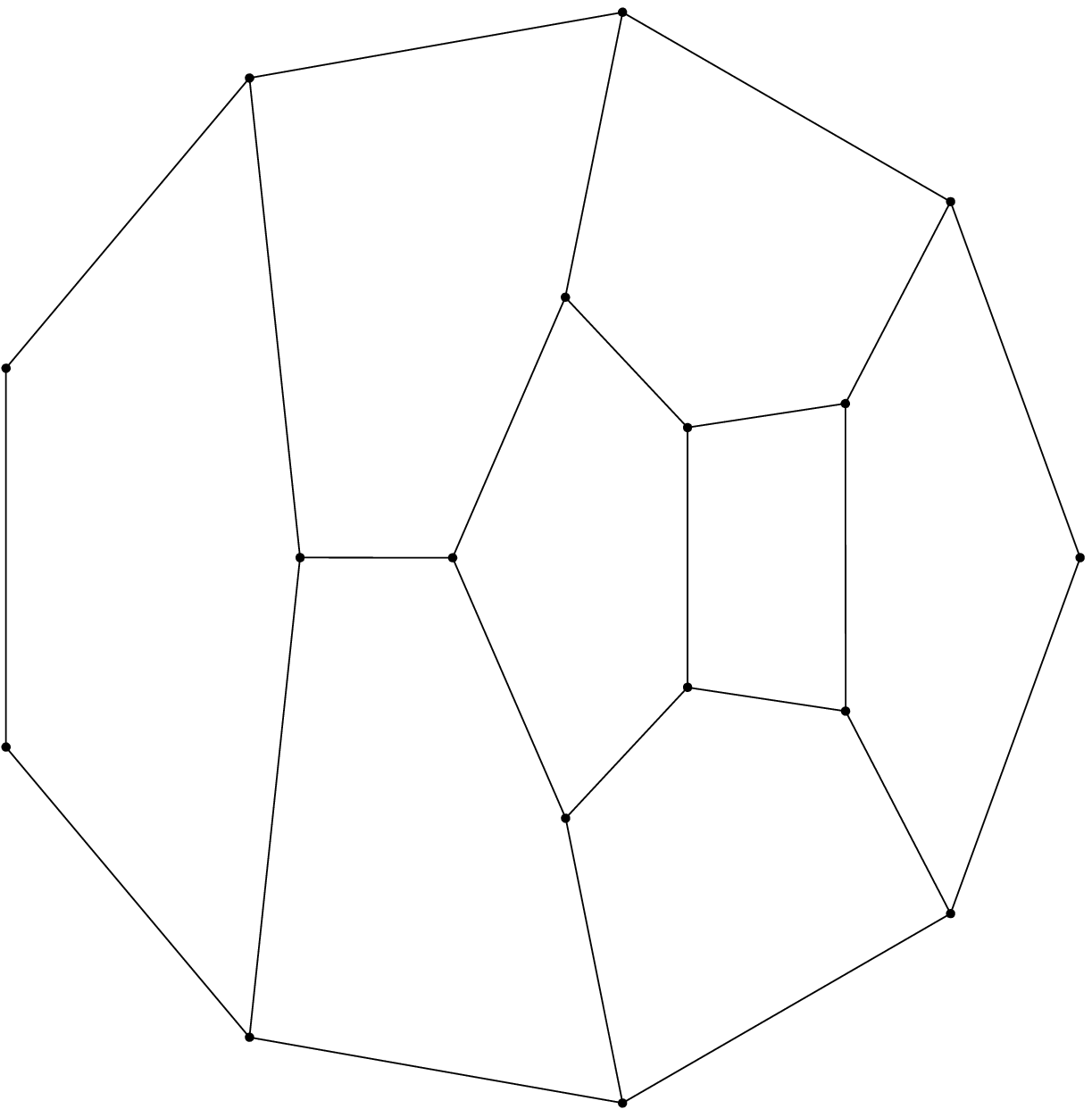}\par
$C_s$
\end{minipage}
\begin{minipage}{3cm}
\centering
\epsfig{height=20mm, file=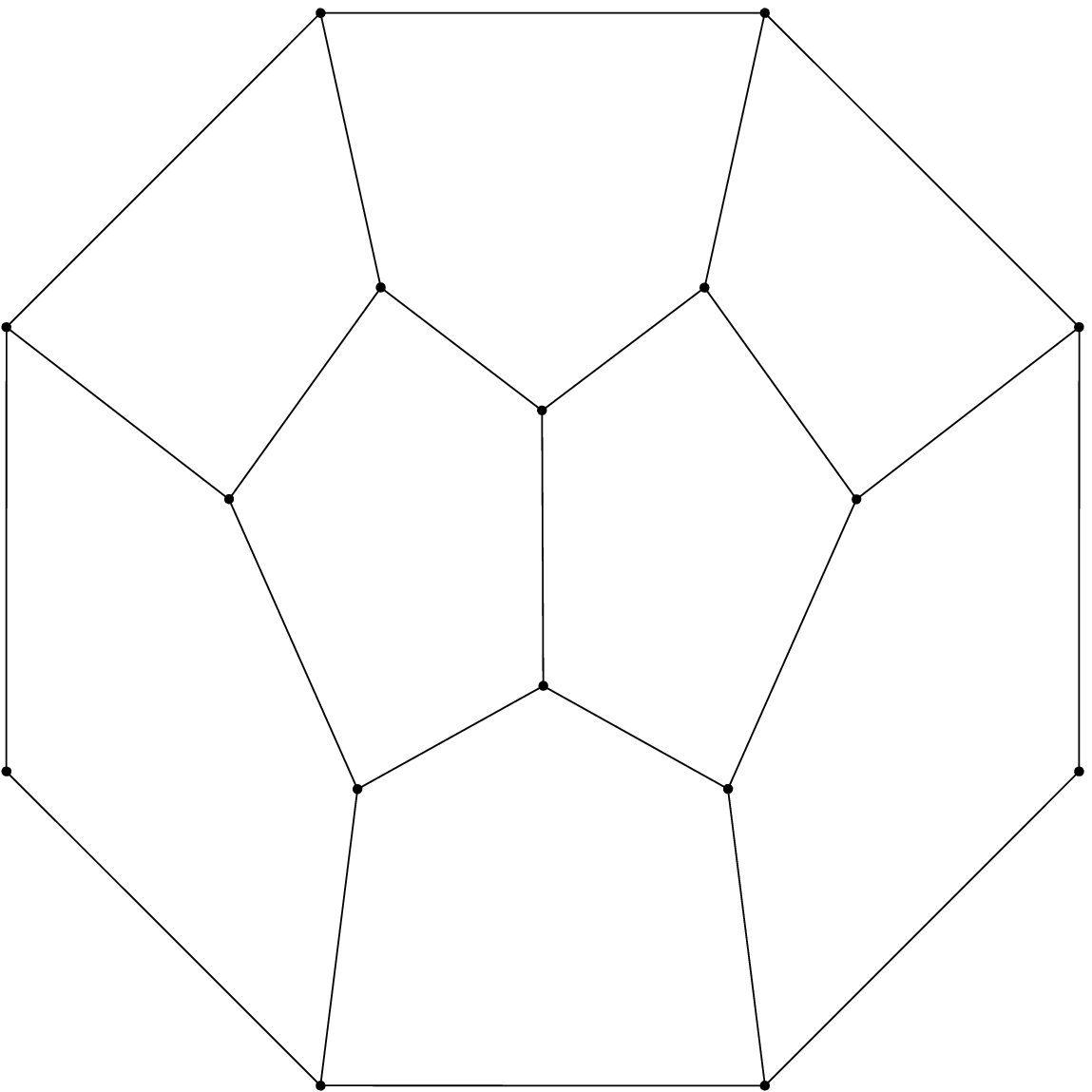}\par
$C_s$
\end{minipage}
\begin{minipage}{3cm}
\centering
\epsfig{height=20mm, file=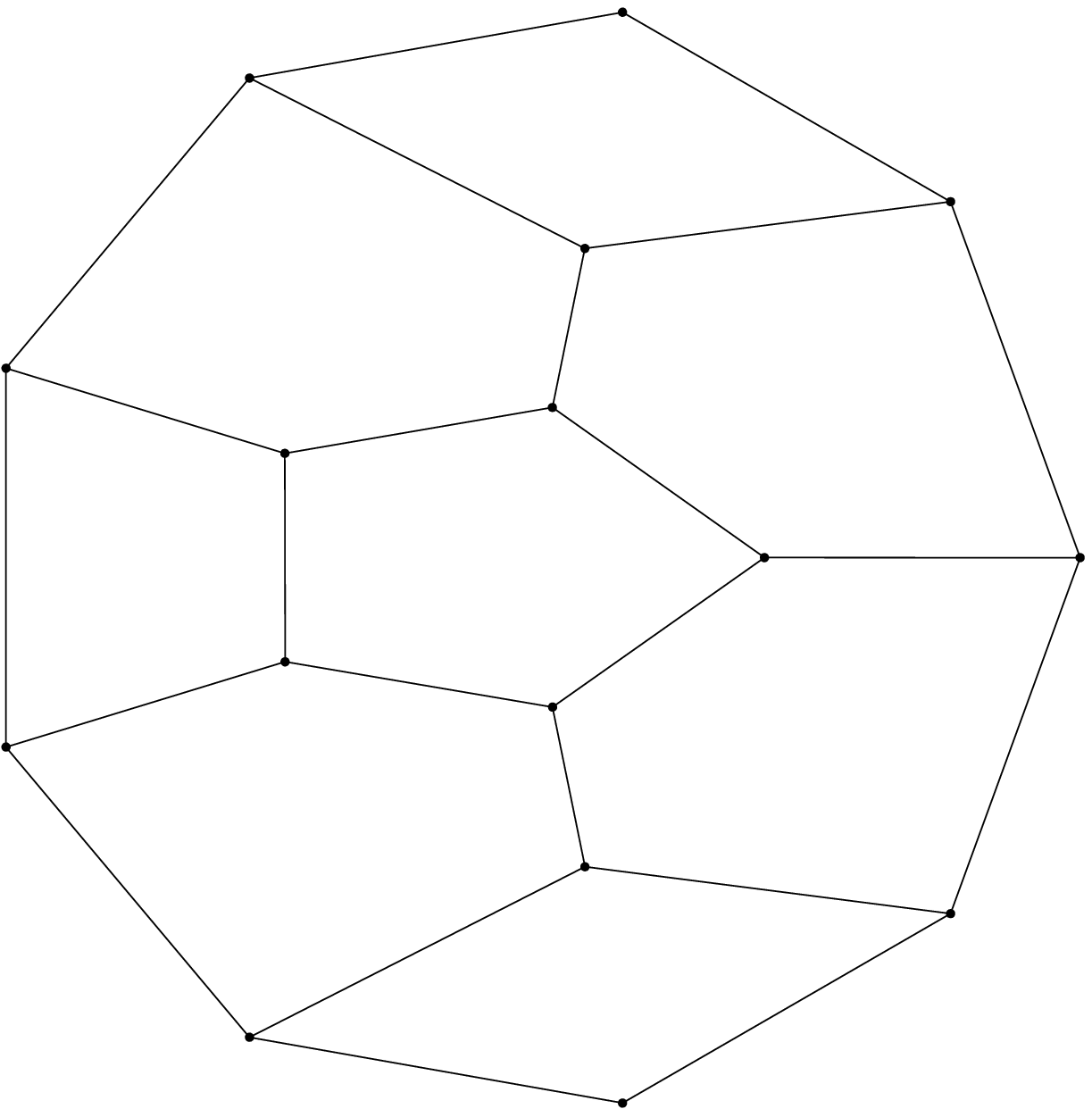}\par
$C_s$
\end{minipage}
\begin{minipage}{3cm}
\centering
\epsfig{height=20mm, file=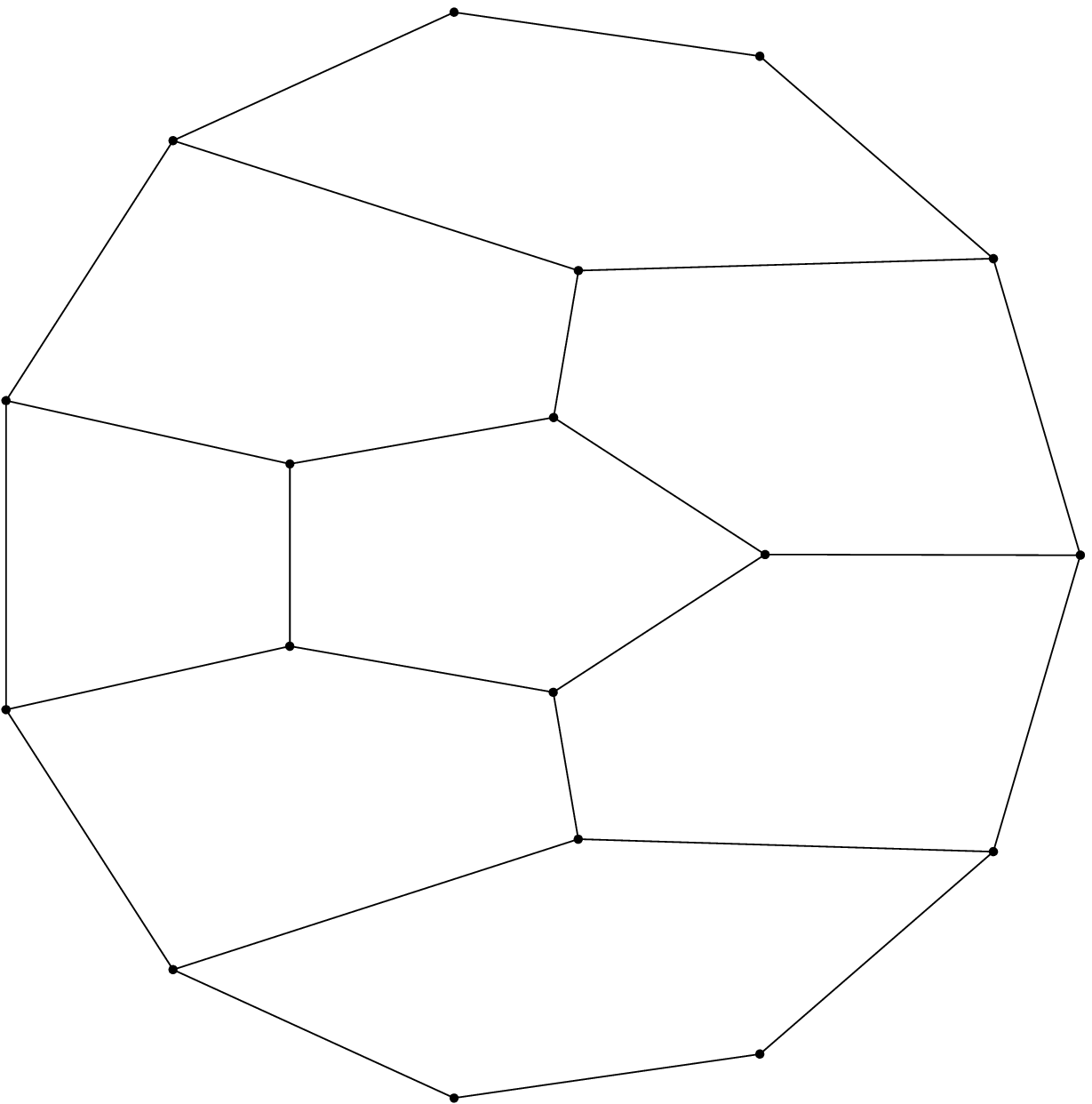}\par
$C_s$
\end{minipage}
\begin{minipage}{3cm}
\centering
\epsfig{height=20mm, file=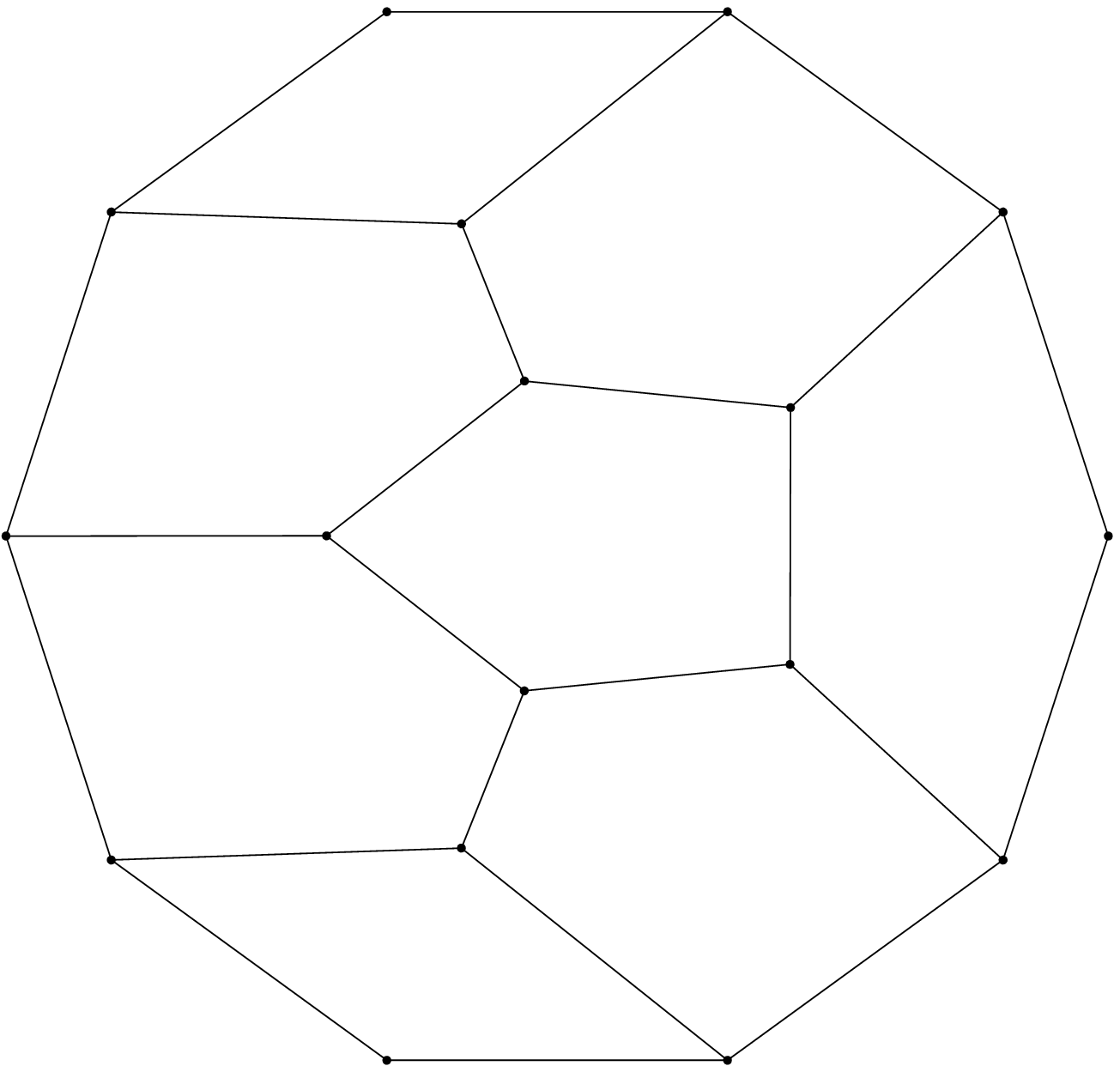}\par
$C_s$
\end{minipage}
\begin{minipage}{3cm}
\centering
\epsfig{height=20mm, file=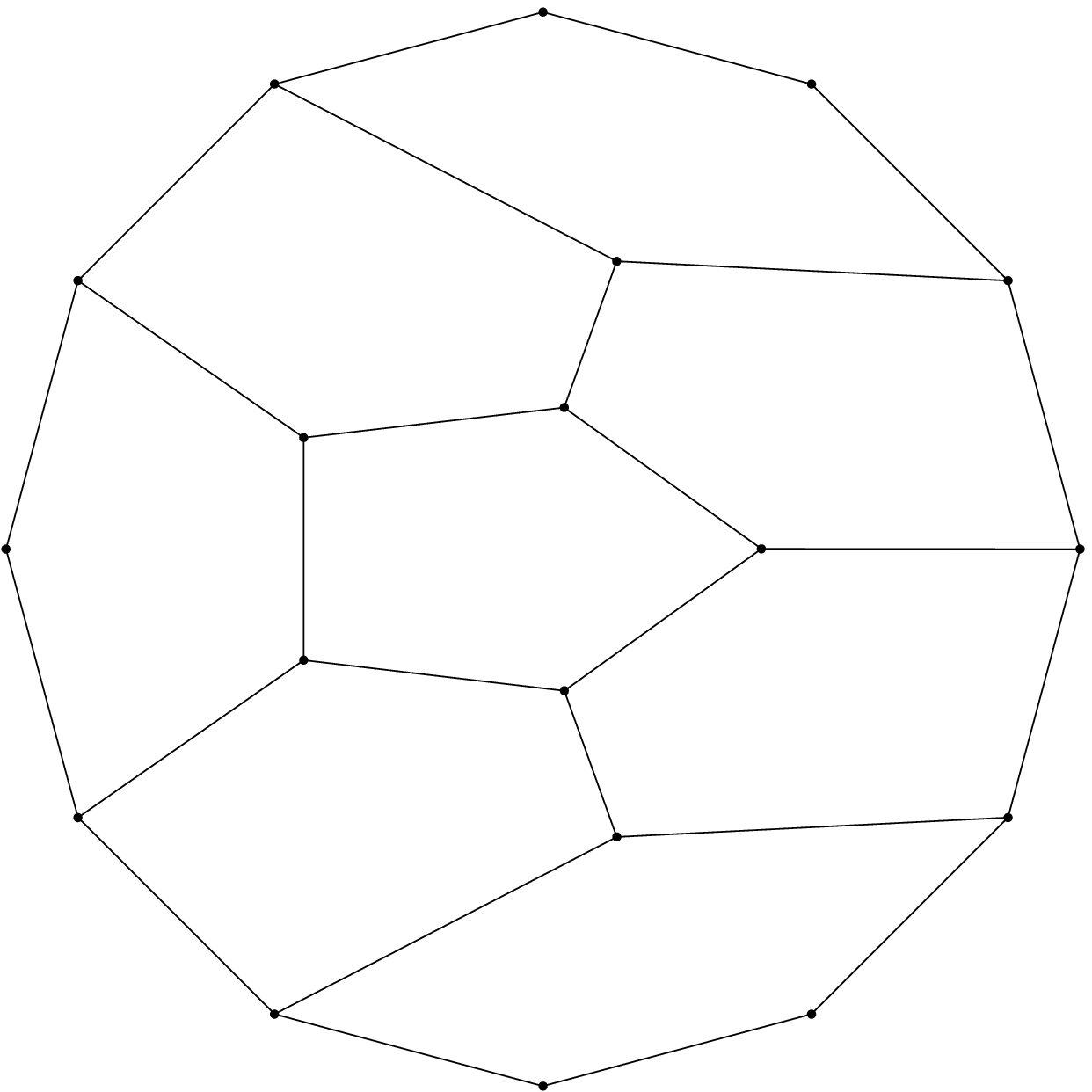}\par
$C_s$
\end{minipage}
\begin{minipage}{3cm}
\centering
\epsfig{height=20mm, file=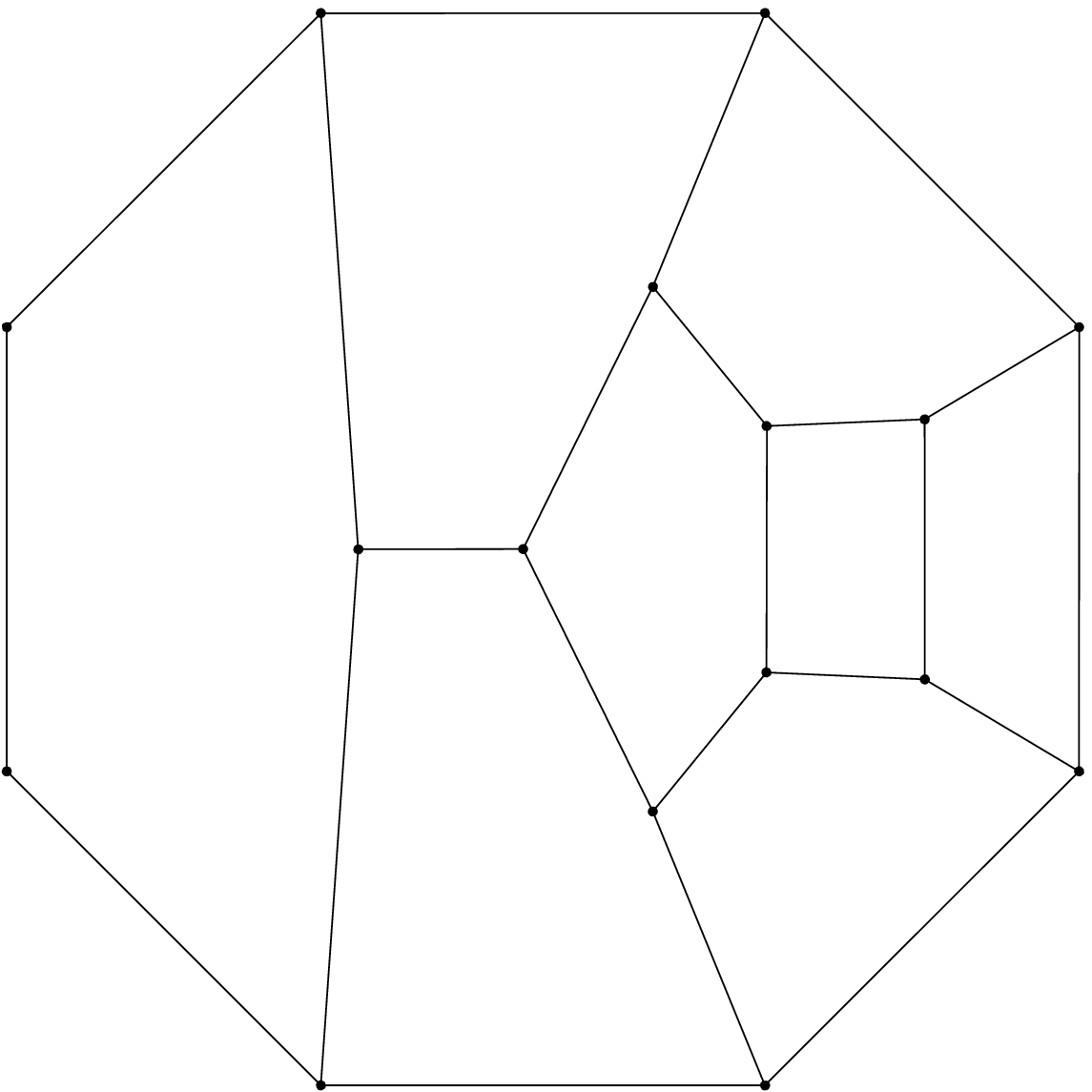}\par
$C_s$
\end{minipage}
\begin{minipage}{3cm}
\centering
\epsfig{height=20mm, file=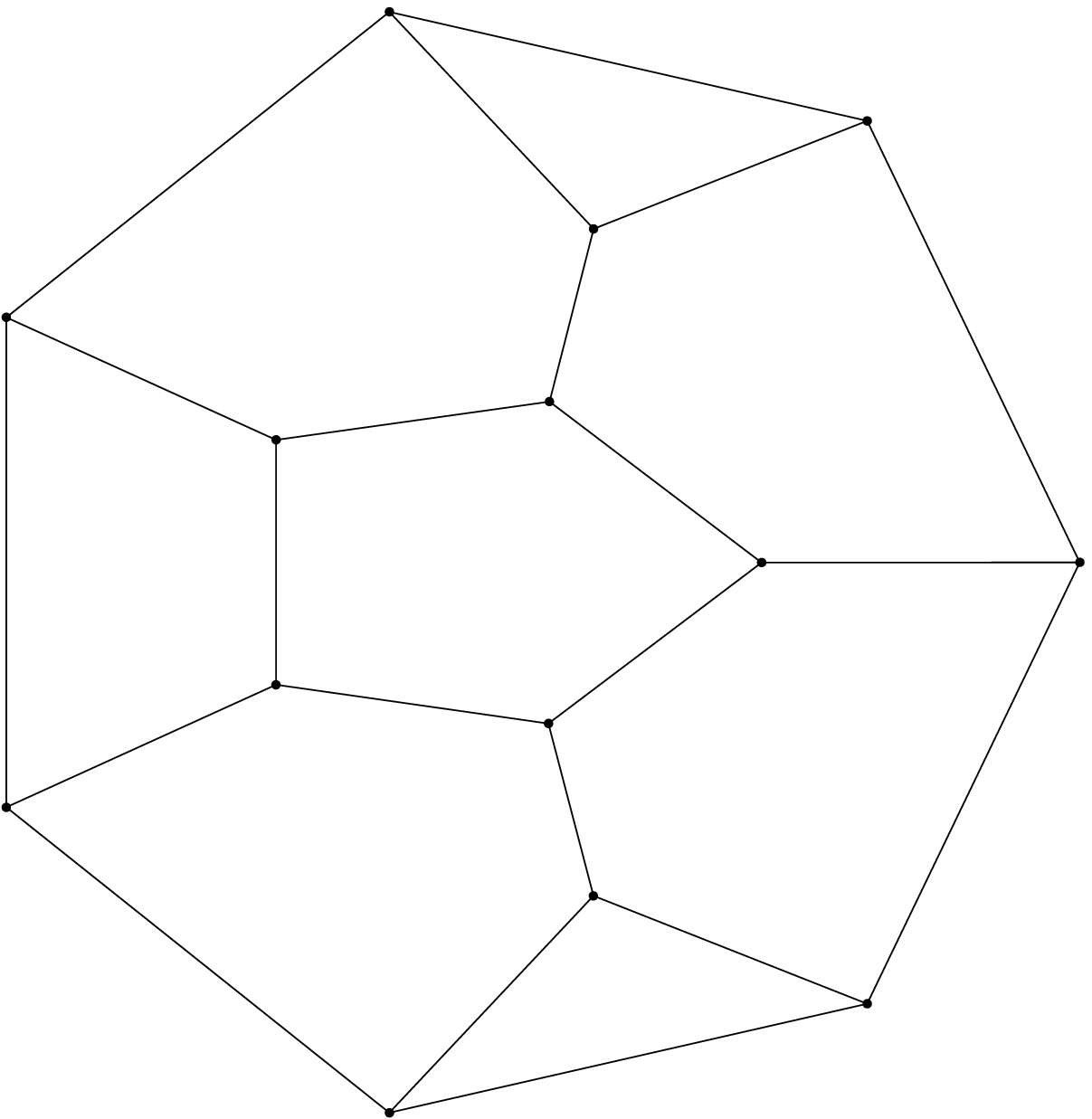}\par
$C_s$, nonext.
\end{minipage}
\begin{minipage}{3cm}
\centering
\epsfig{height=20mm, file=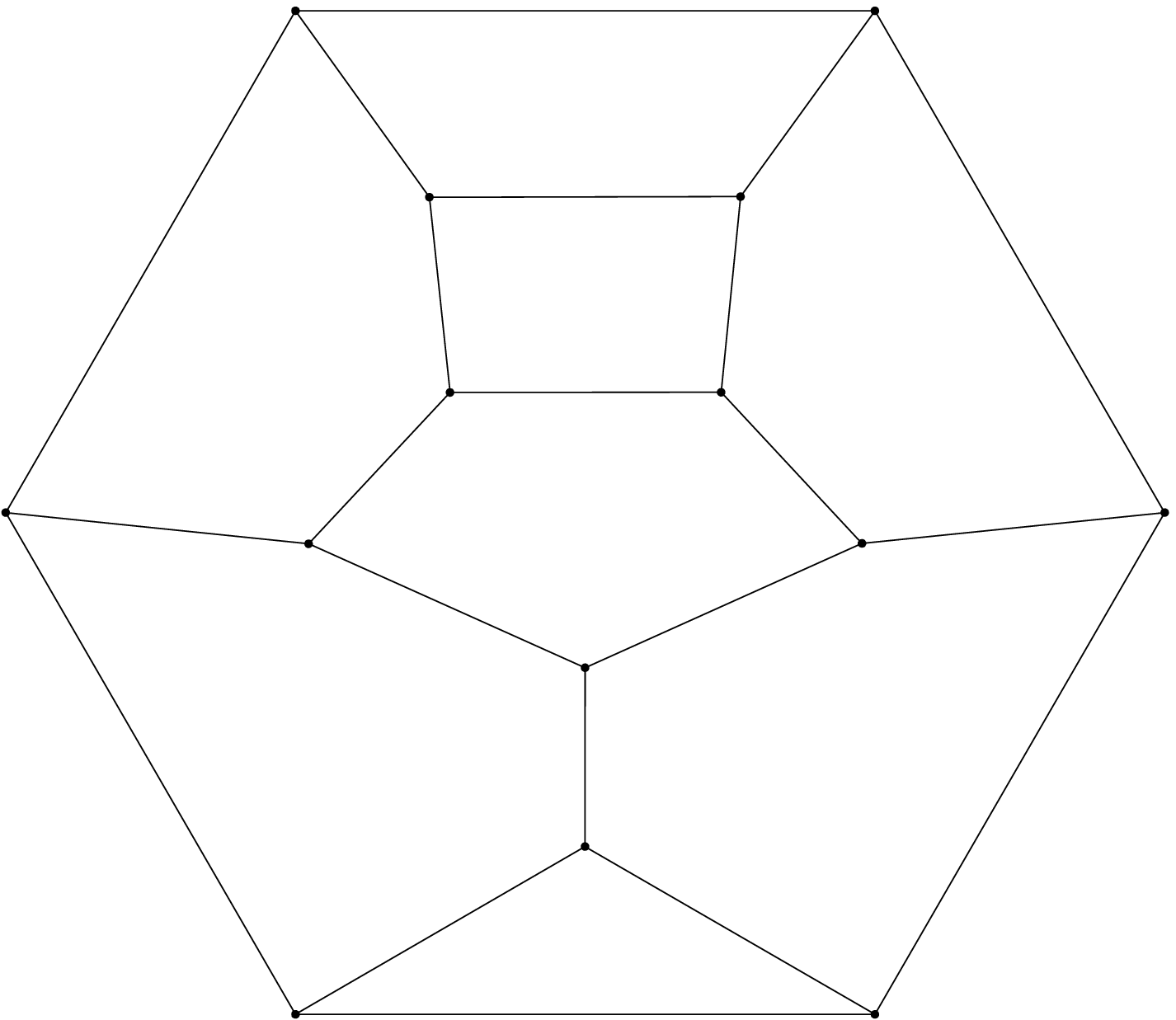}\par
$C_s$, nonext.
\end{minipage}
\begin{minipage}{3cm}
\centering
\epsfig{height=20mm, file=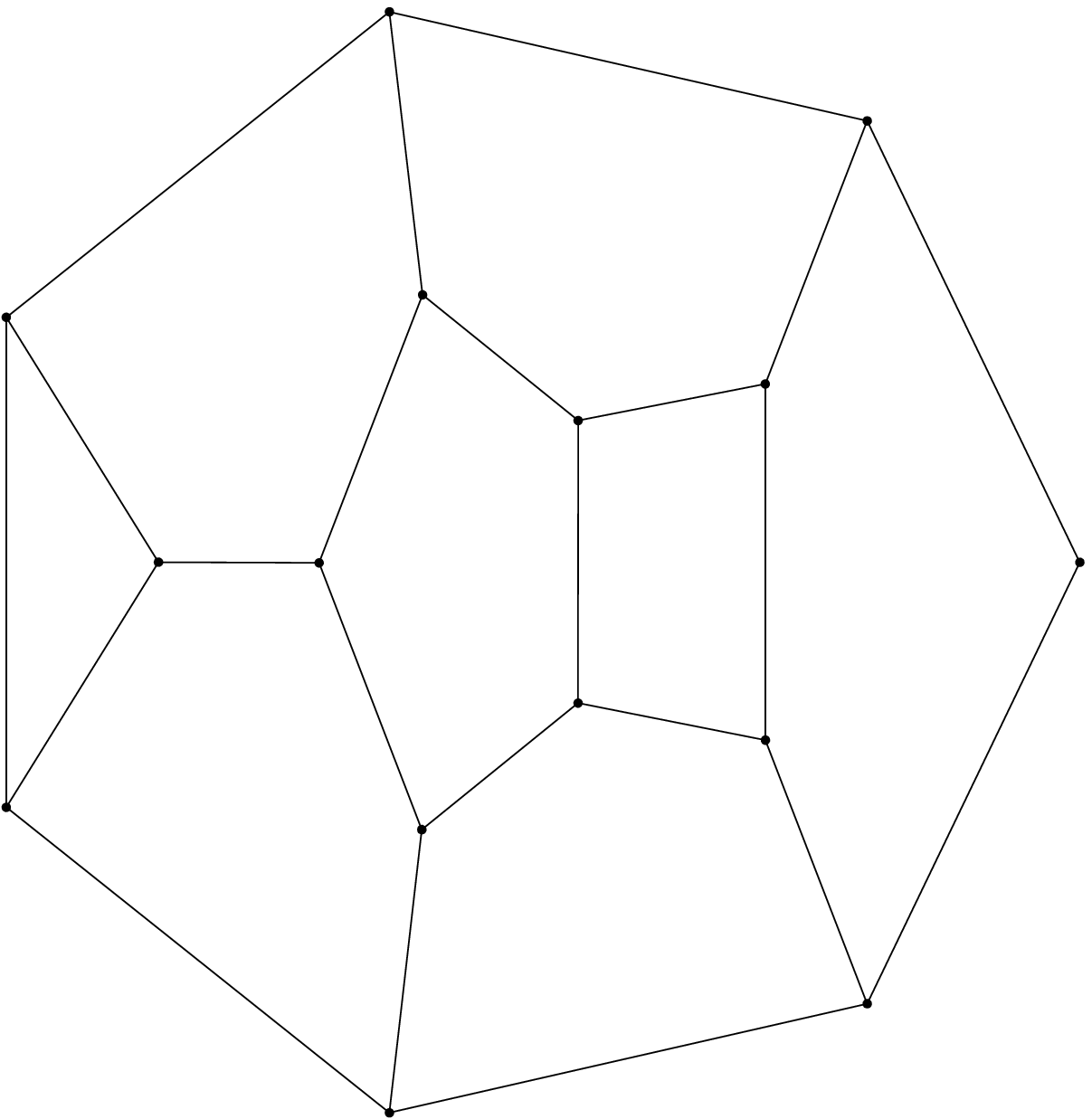}\par
$C_s$, nonext.
\end{minipage}
\begin{minipage}{3cm}
\centering
\epsfig{height=20mm, file=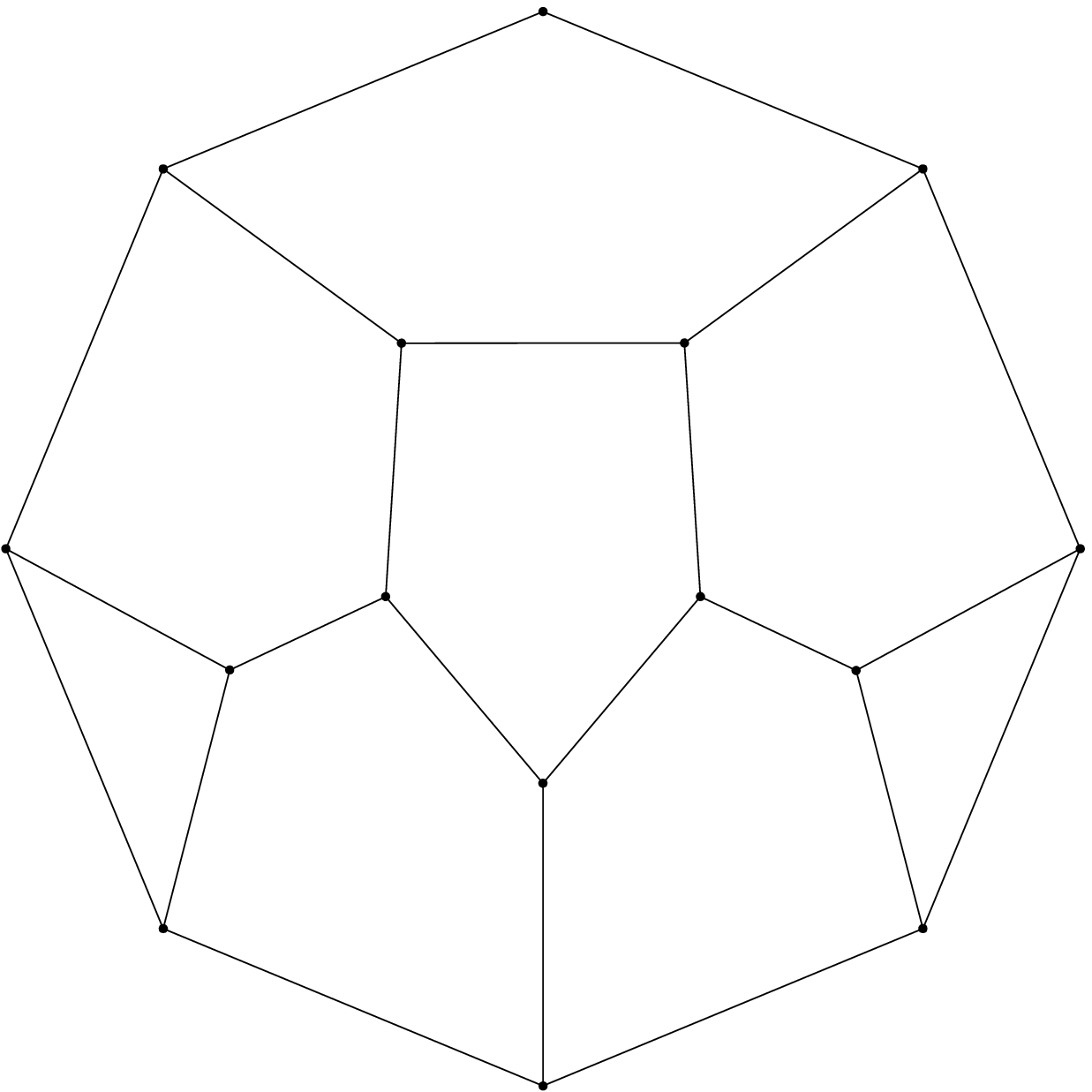}\par
$C_s$, nonext.
\end{minipage}
\begin{minipage}{3cm}
\centering
\epsfig{height=20mm, file=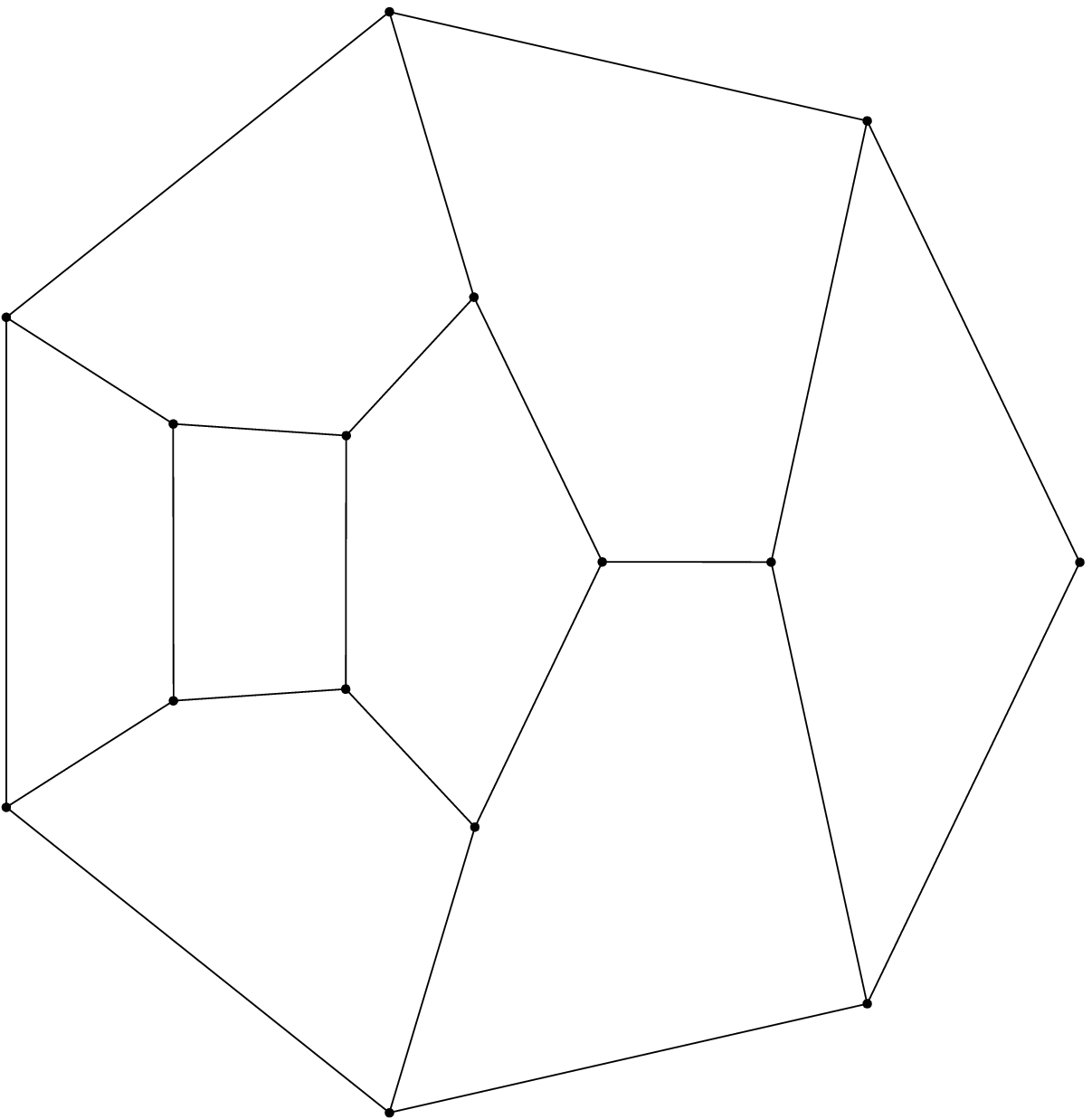}\par
$C_s$, nonext.
\end{minipage}
\begin{minipage}{3cm}
\centering
\epsfig{height=20mm, file=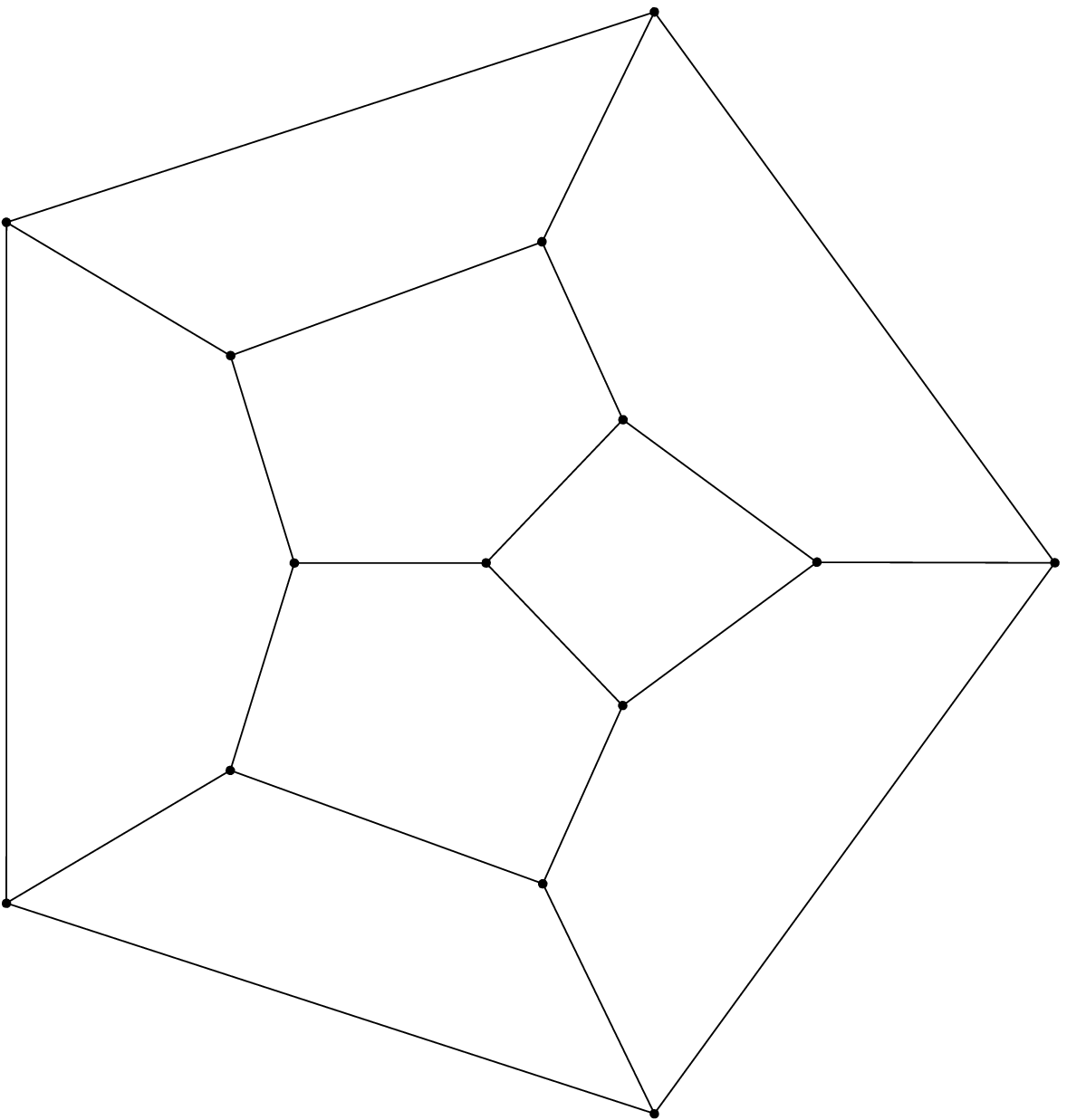}\par
$C_s$, nonext.
%PAIR1
\end{minipage}
\begin{minipage}{3cm}
\centering
\epsfig{height=20mm, file=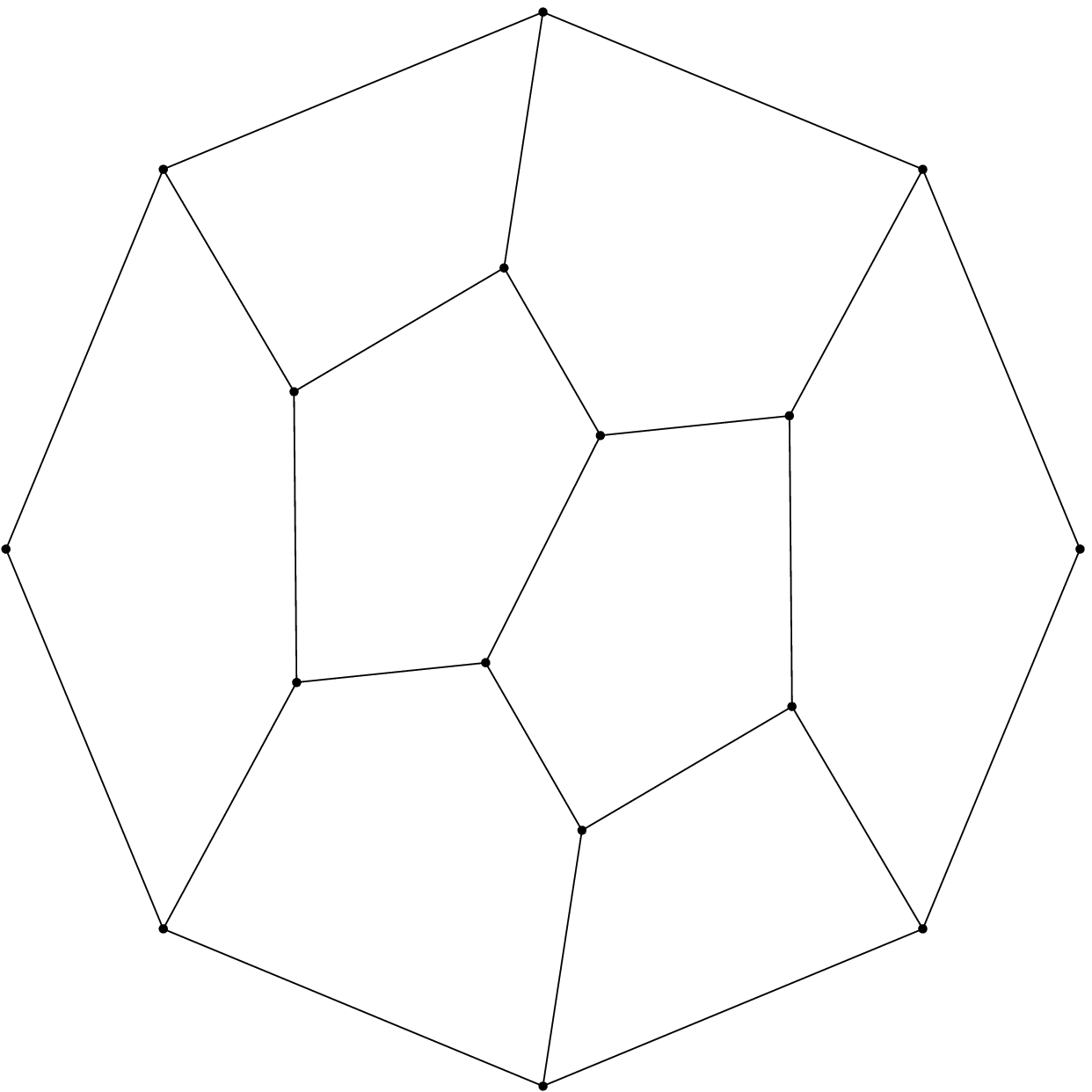}\par
$C_2$
\end{minipage}
\begin{minipage}{3cm}
\centering
\epsfig{height=20mm, file=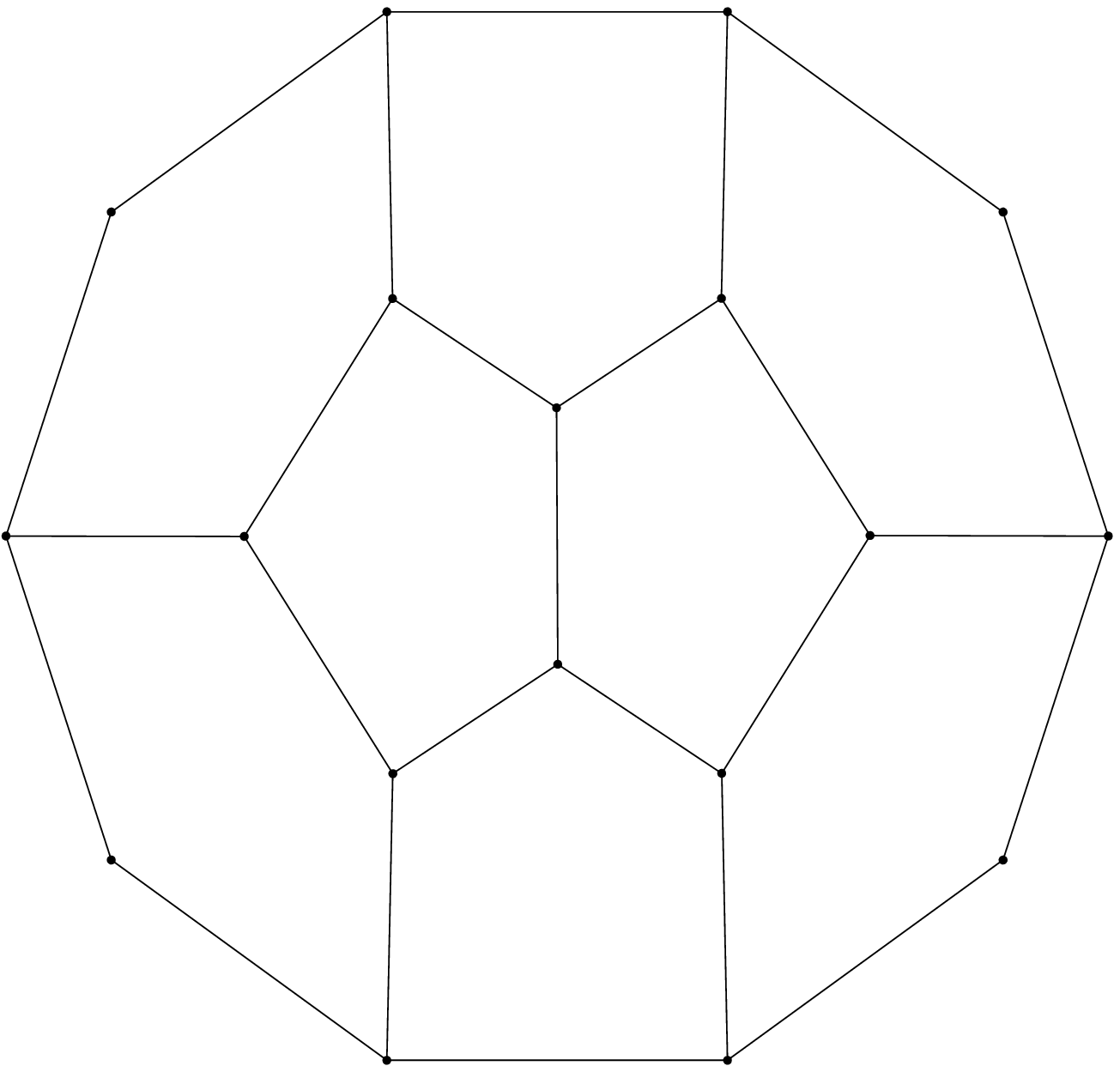}\par
$C_{2\nu}$
\end{minipage}
\begin{minipage}{3cm}
\centering
\epsfig{height=20mm, file=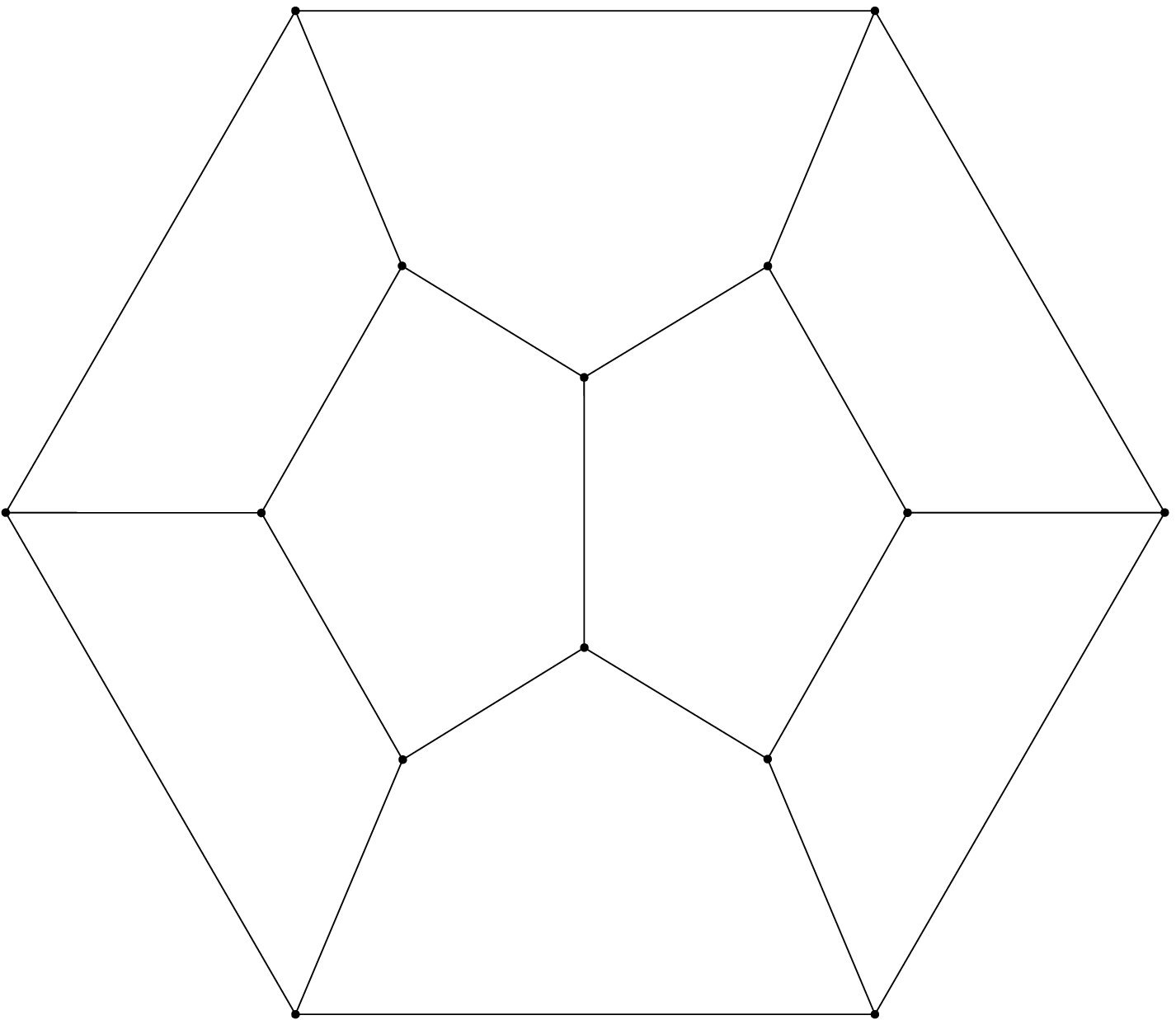}\par
$C_{2\nu}$, nonext.
\end{minipage}
\begin{minipage}{3cm}
\centering
\epsfig{height=20mm, file=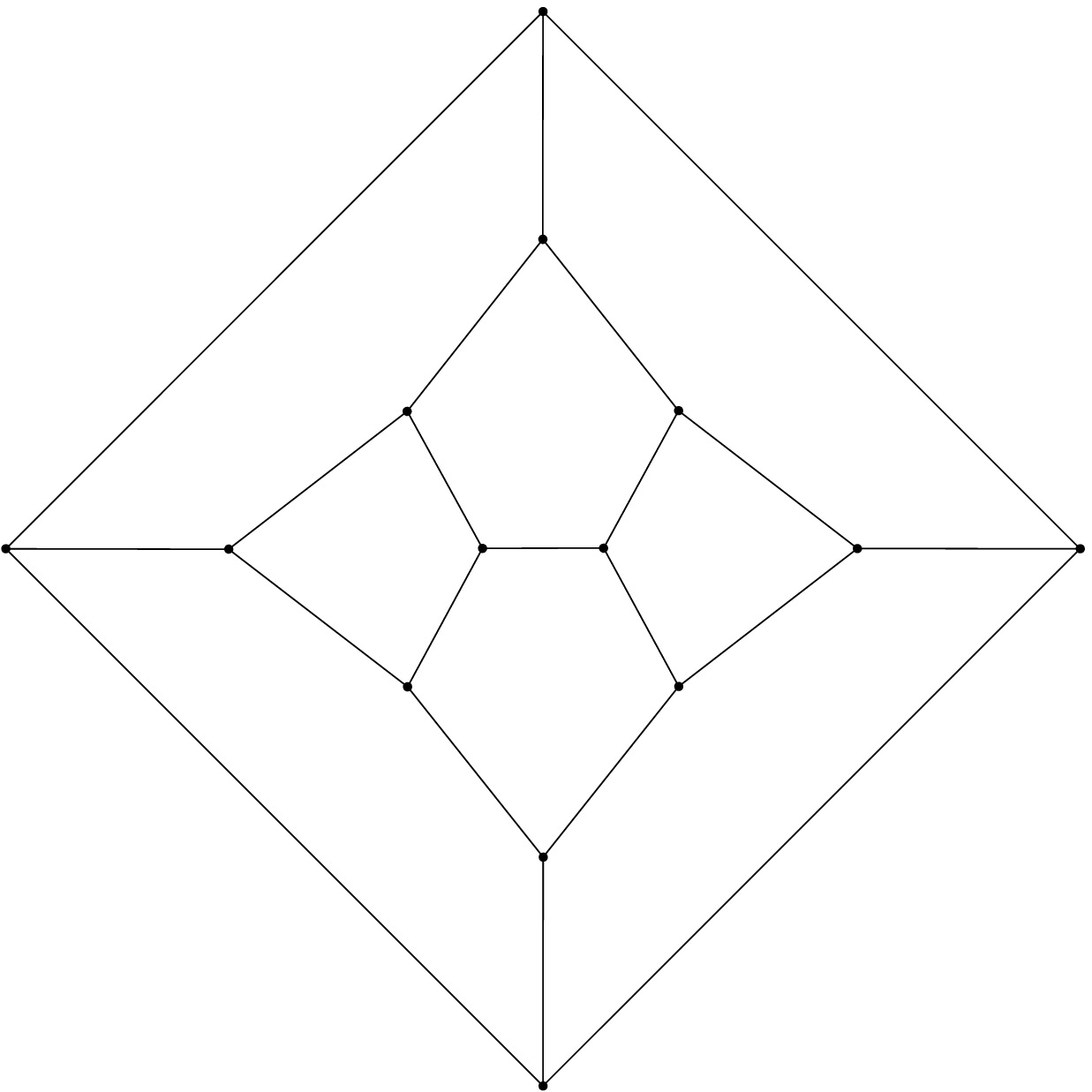}\par
$C_{2\nu}$, nonext.
%PAIR1
\end{minipage}

\end{center}
List of sporadic elementary $(\{3,4,5\},3)$-polycycles with $9$ faces:
\begin{center}
\begin{minipage}{3cm}
\centering
\epsfig{height=20mm, file=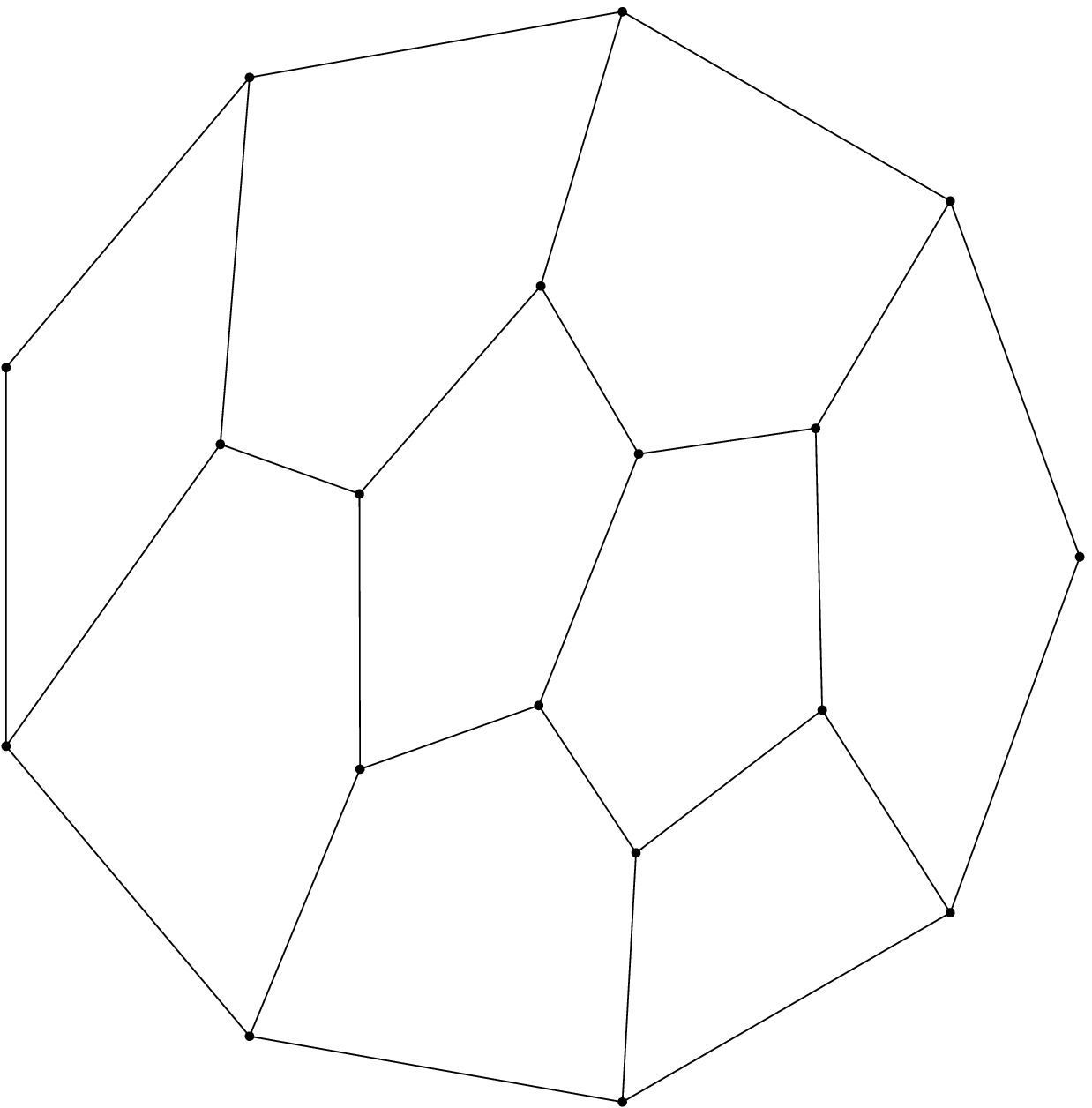}\par
$C_1$
\end{minipage}
\begin{minipage}{3cm}
\centering
\epsfig{height=20mm, file=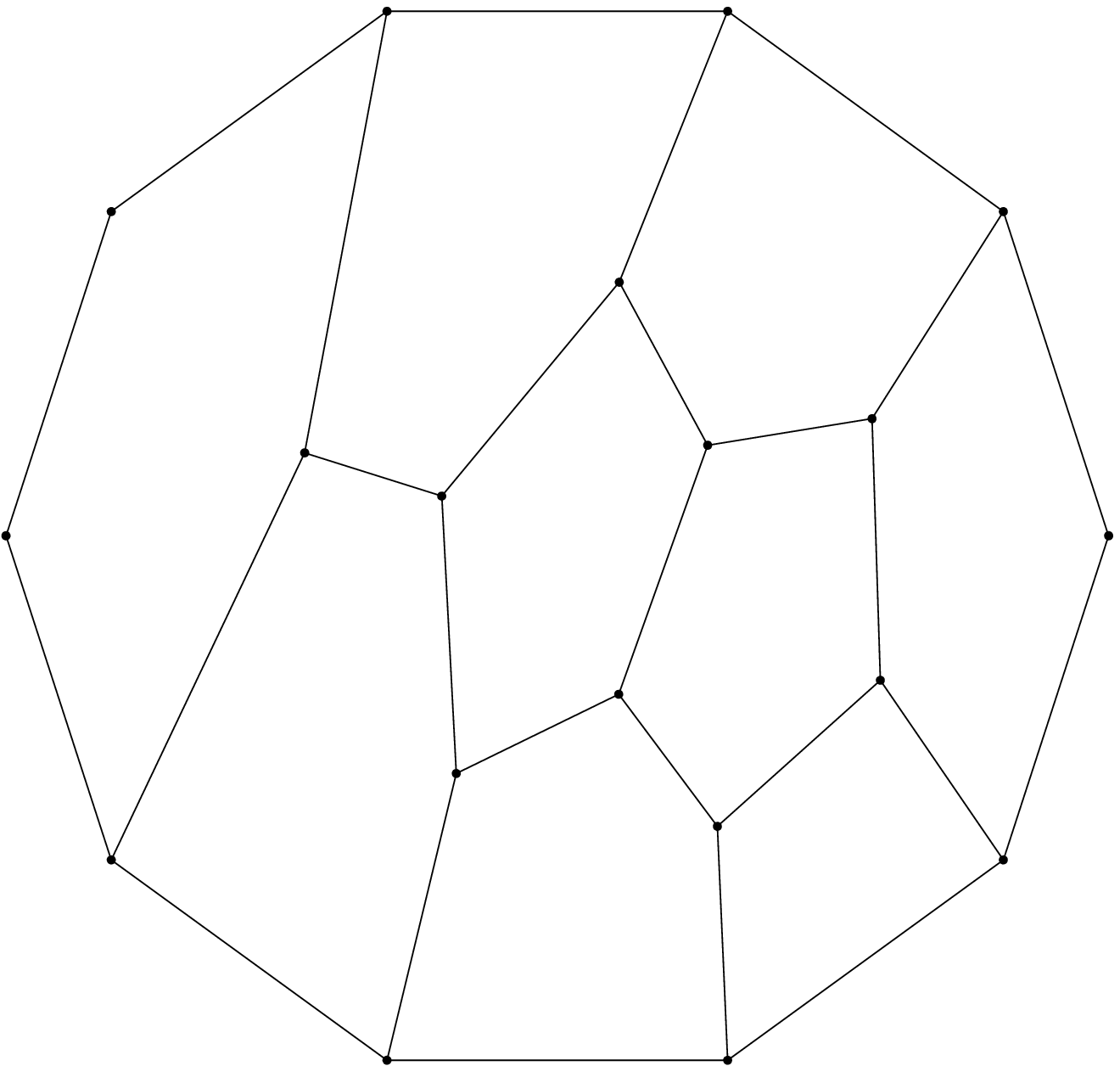}\par
$C_1$
\end{minipage}
\begin{minipage}{3cm}
\centering
\epsfig{height=20mm, file=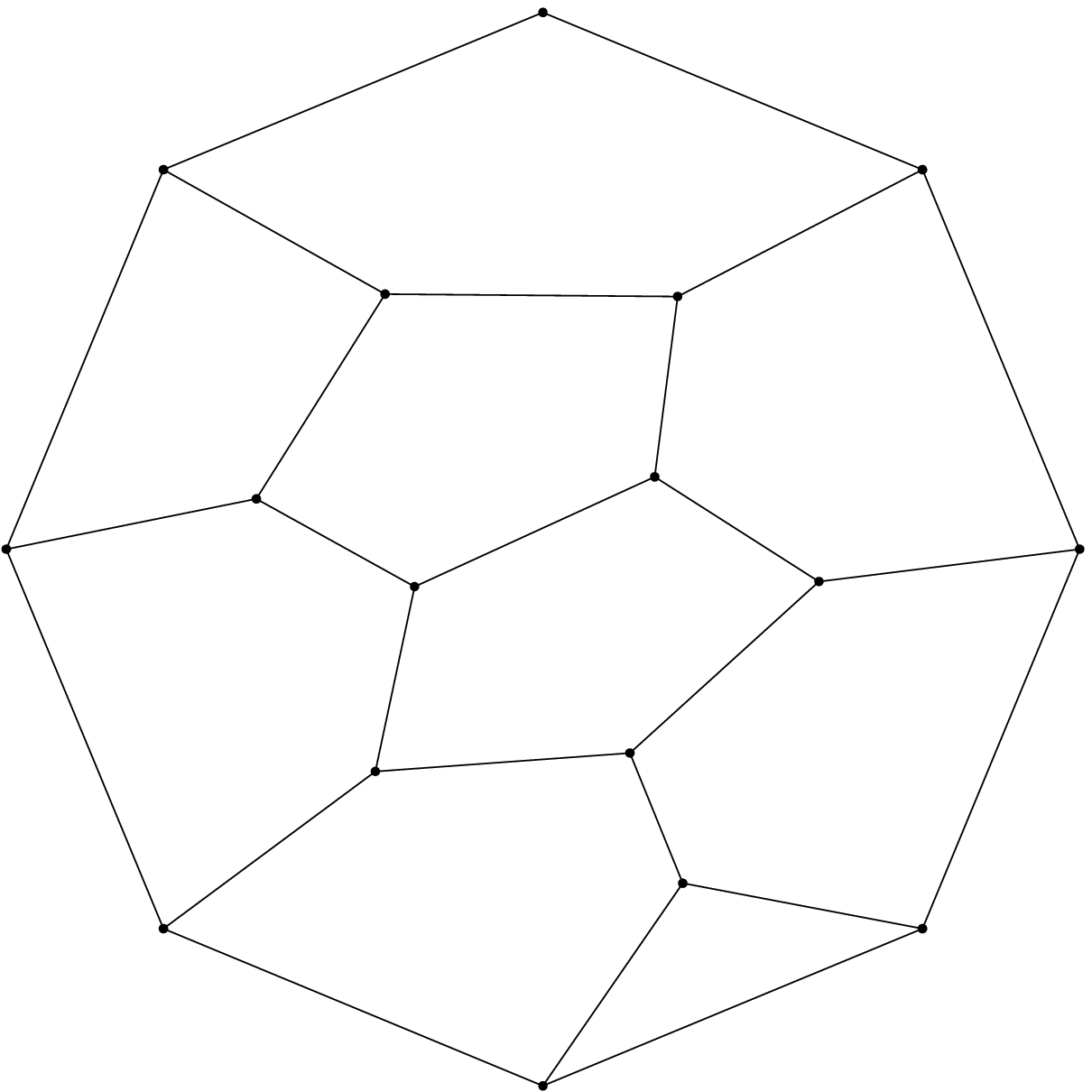}\par
$C_1$, nonext.
\end{minipage}
\begin{minipage}{3cm}
\centering
\epsfig{height=20mm, file=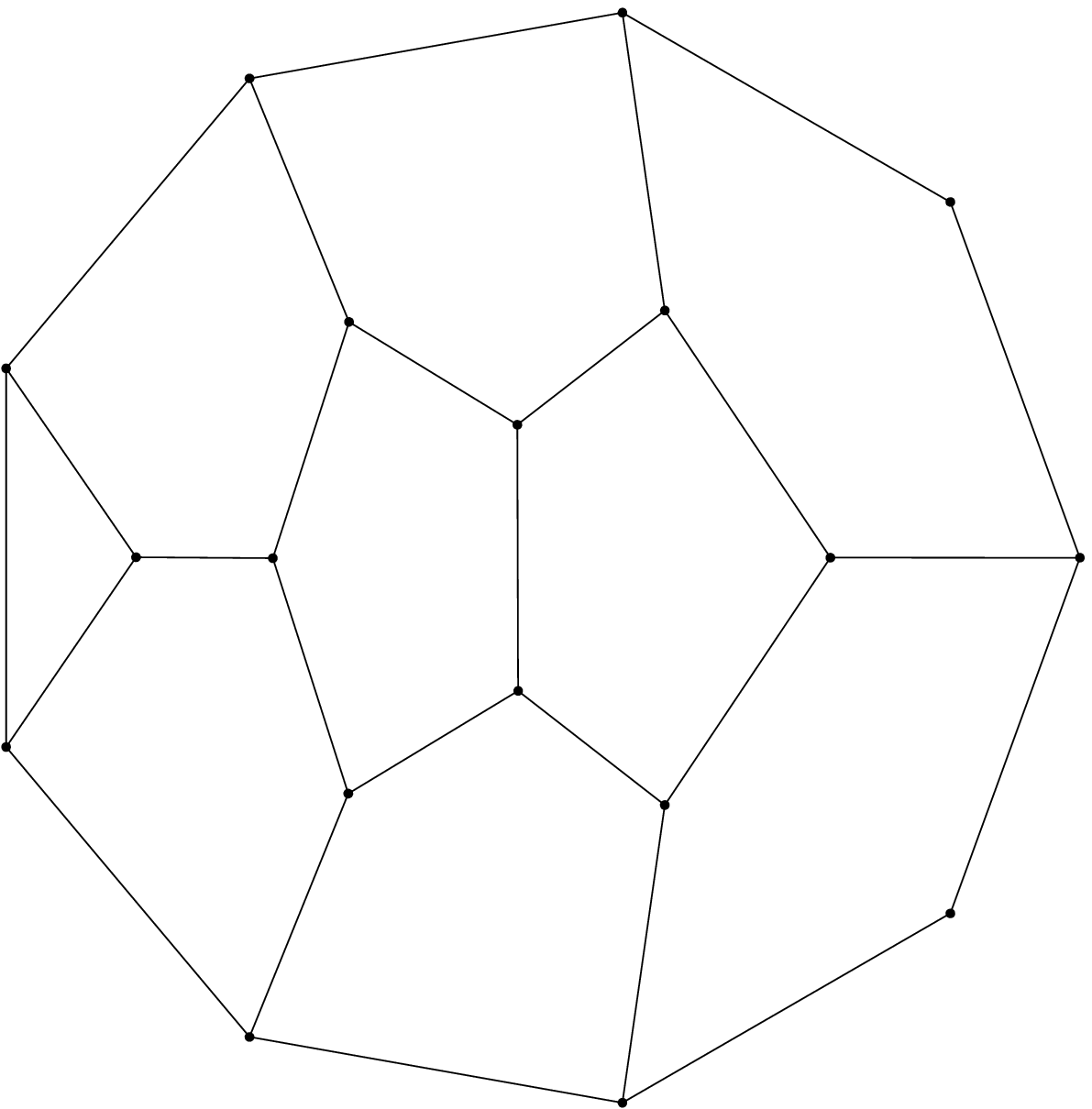}\par
$C_s$
\end{minipage}
\begin{minipage}{3cm}
\centering
\epsfig{height=20mm, file=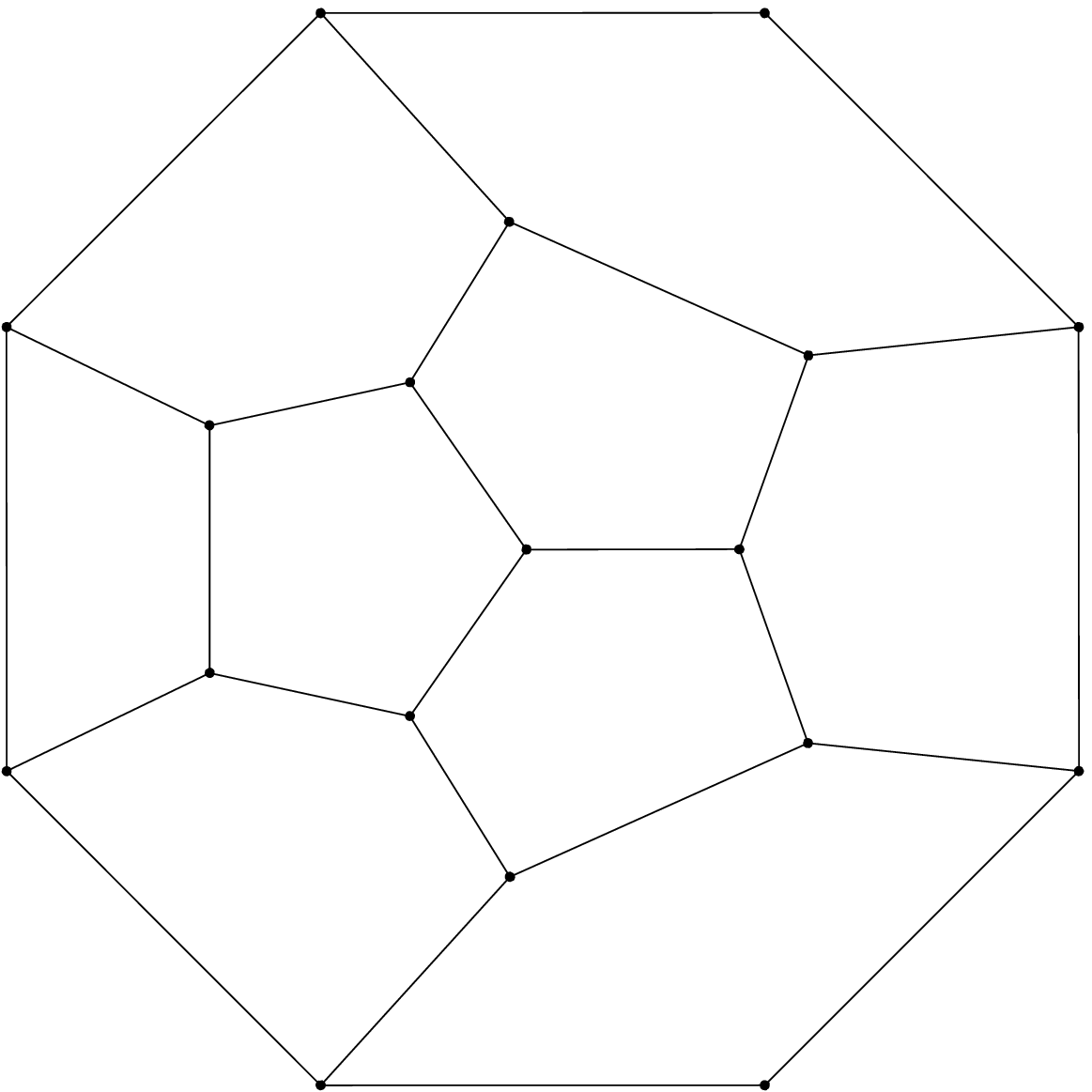}\par
$C_s$
\end{minipage}
\begin{minipage}{3cm}
\centering
\epsfig{height=20mm, file=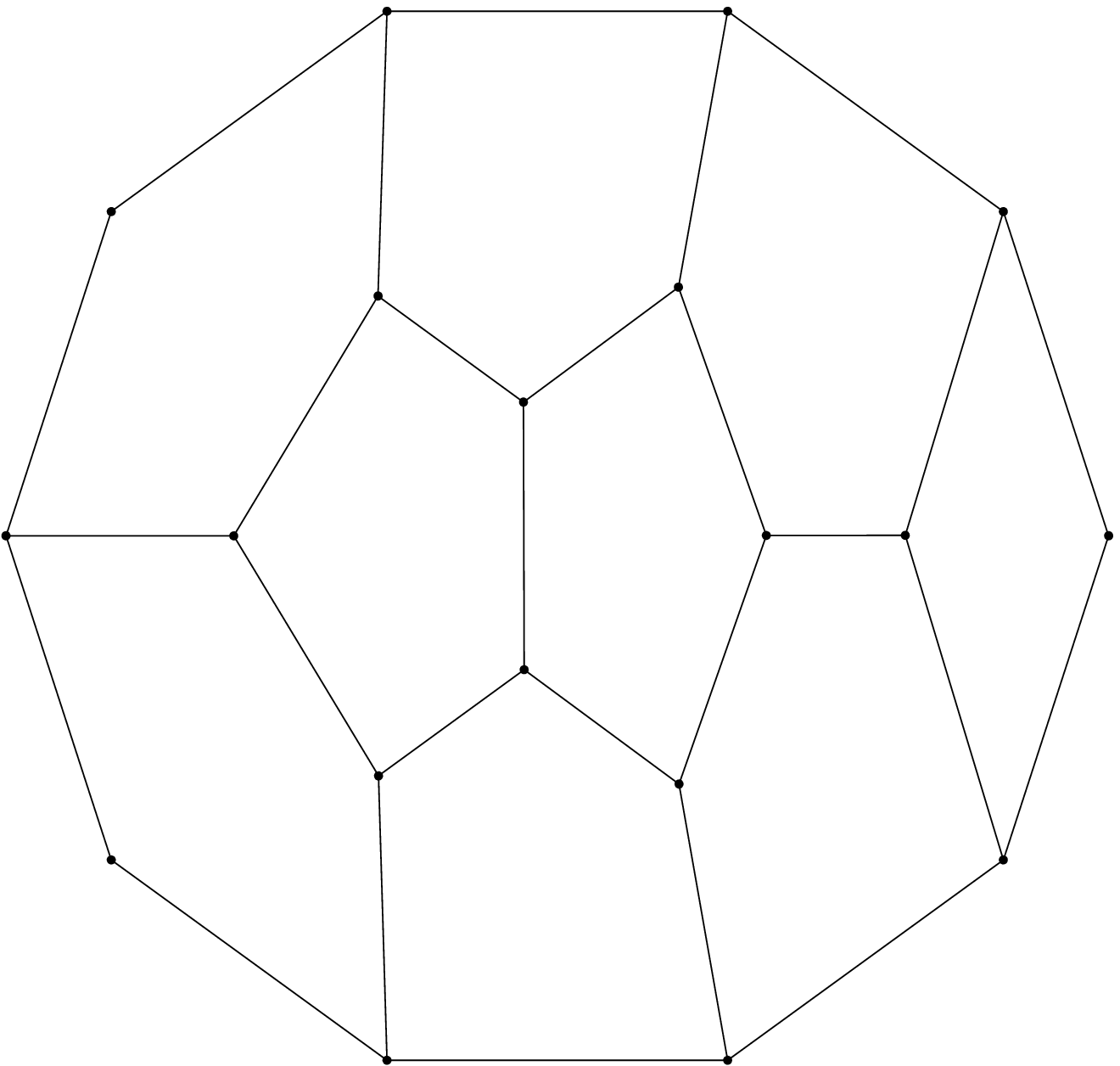}\par
$C_s$
\end{minipage}
\begin{minipage}{3cm}
\centering
\epsfig{height=20mm, file=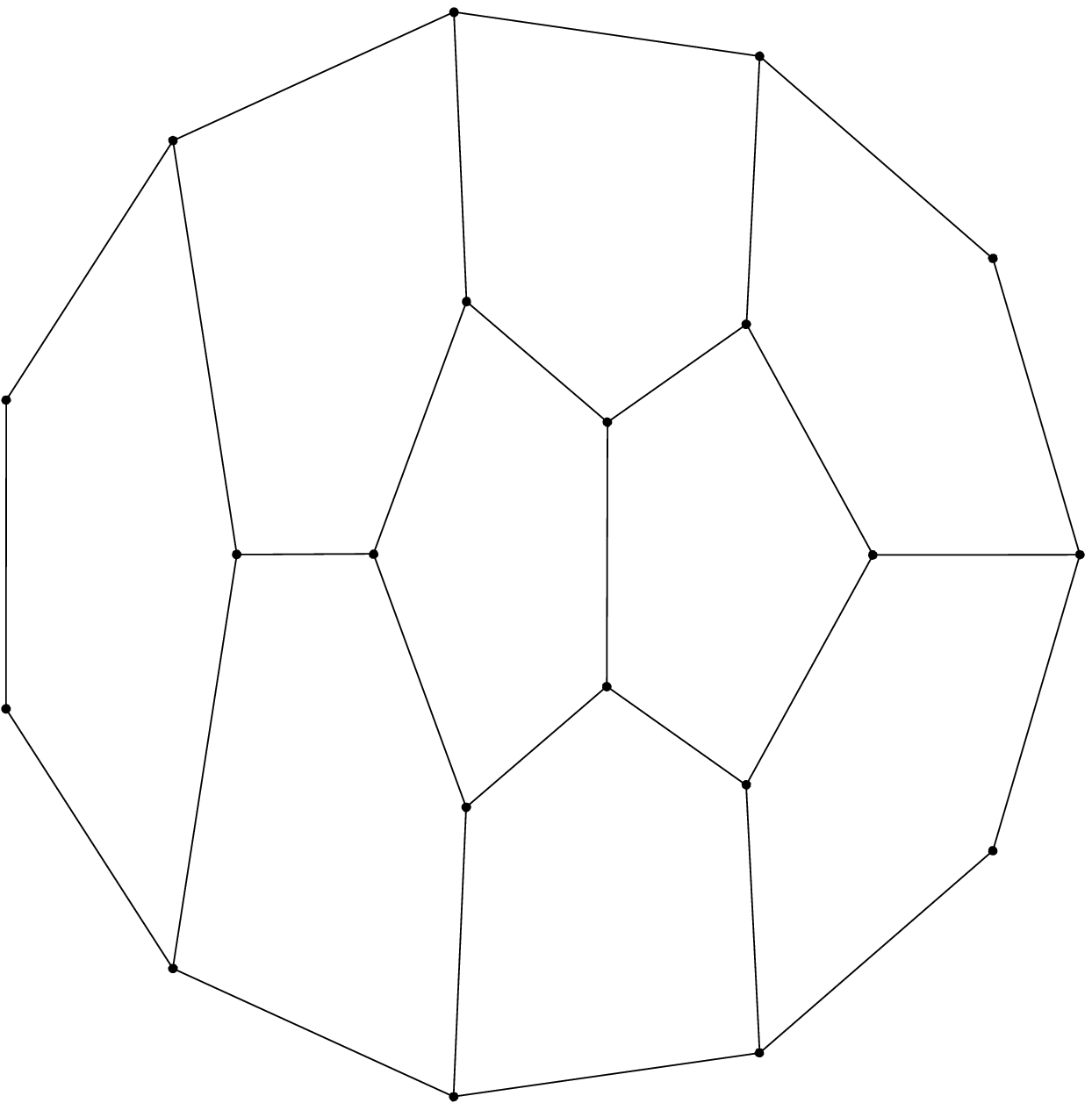}\par
$C_s$
\end{minipage}
\begin{minipage}{3cm}
\centering
\epsfig{height=20mm, file=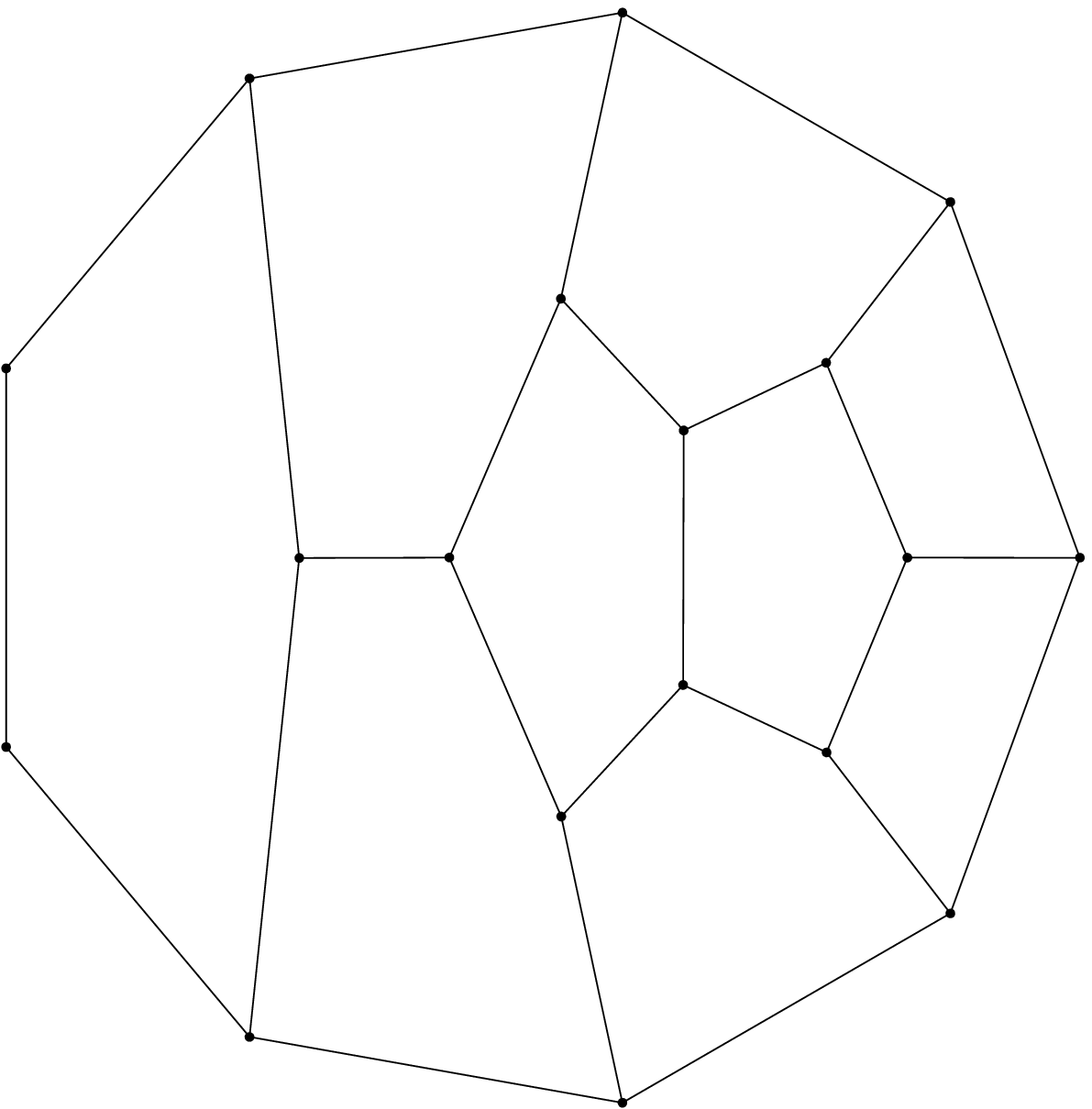}\par
$C_s$
\end{minipage}
\begin{minipage}{3cm}
\centering
\epsfig{height=20mm, file=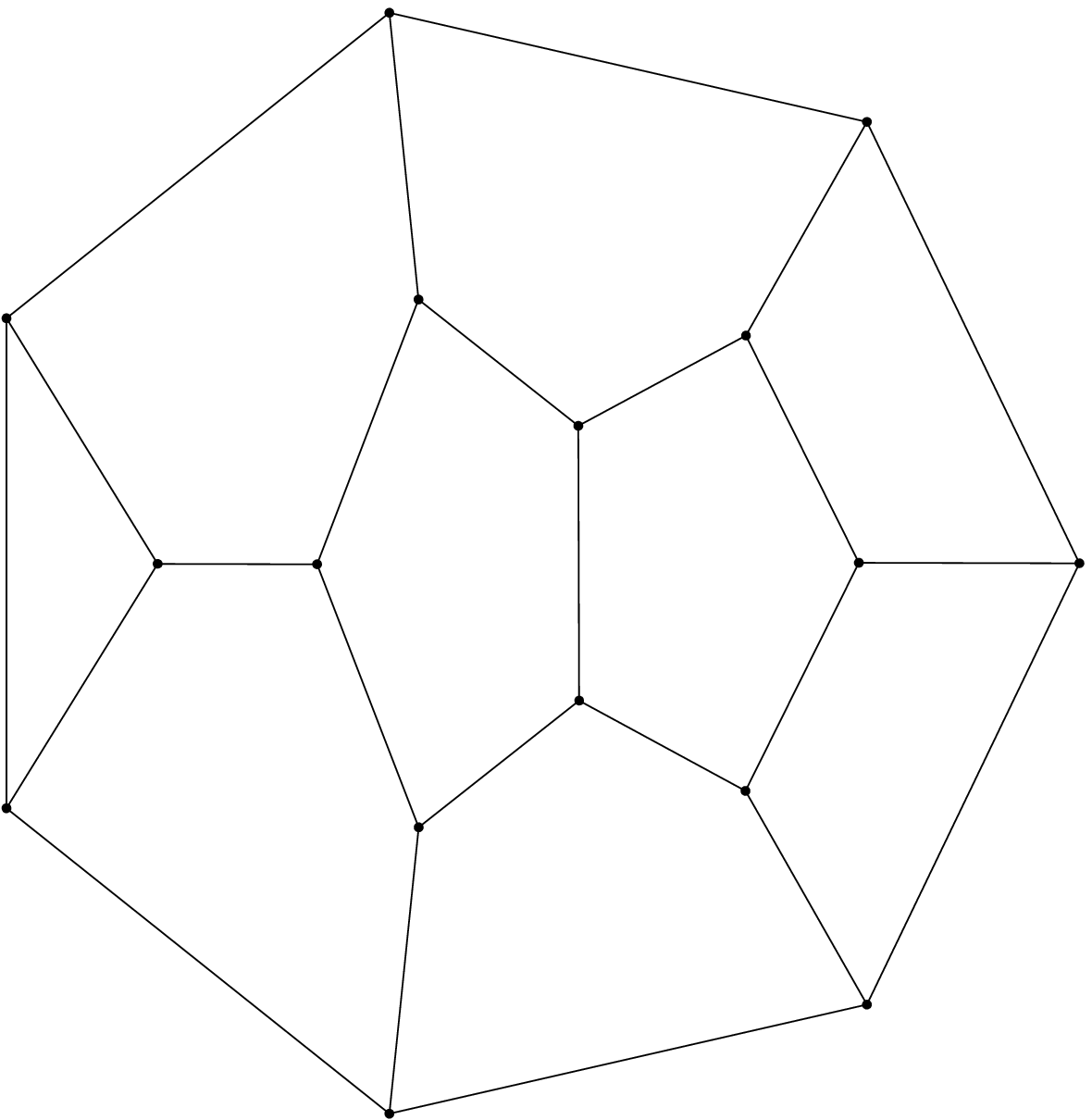}\par
$C_s$, nonext.
\end{minipage}
\begin{minipage}{3cm}
\centering
\epsfig{height=20mm, file=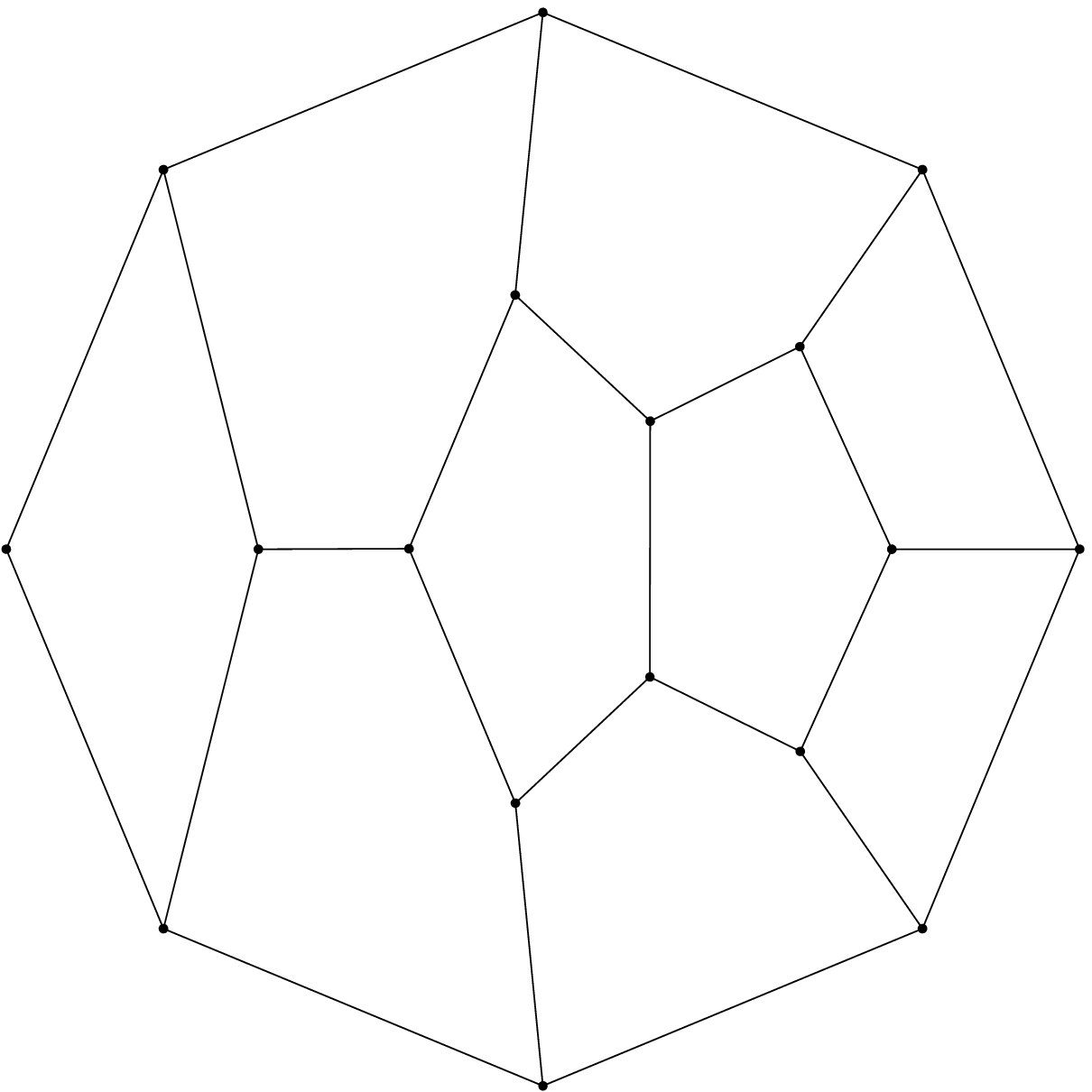}\par
$C_s$, nonext.
\end{minipage}
\begin{minipage}{3cm}
\centering
\epsfig{height=20mm, file=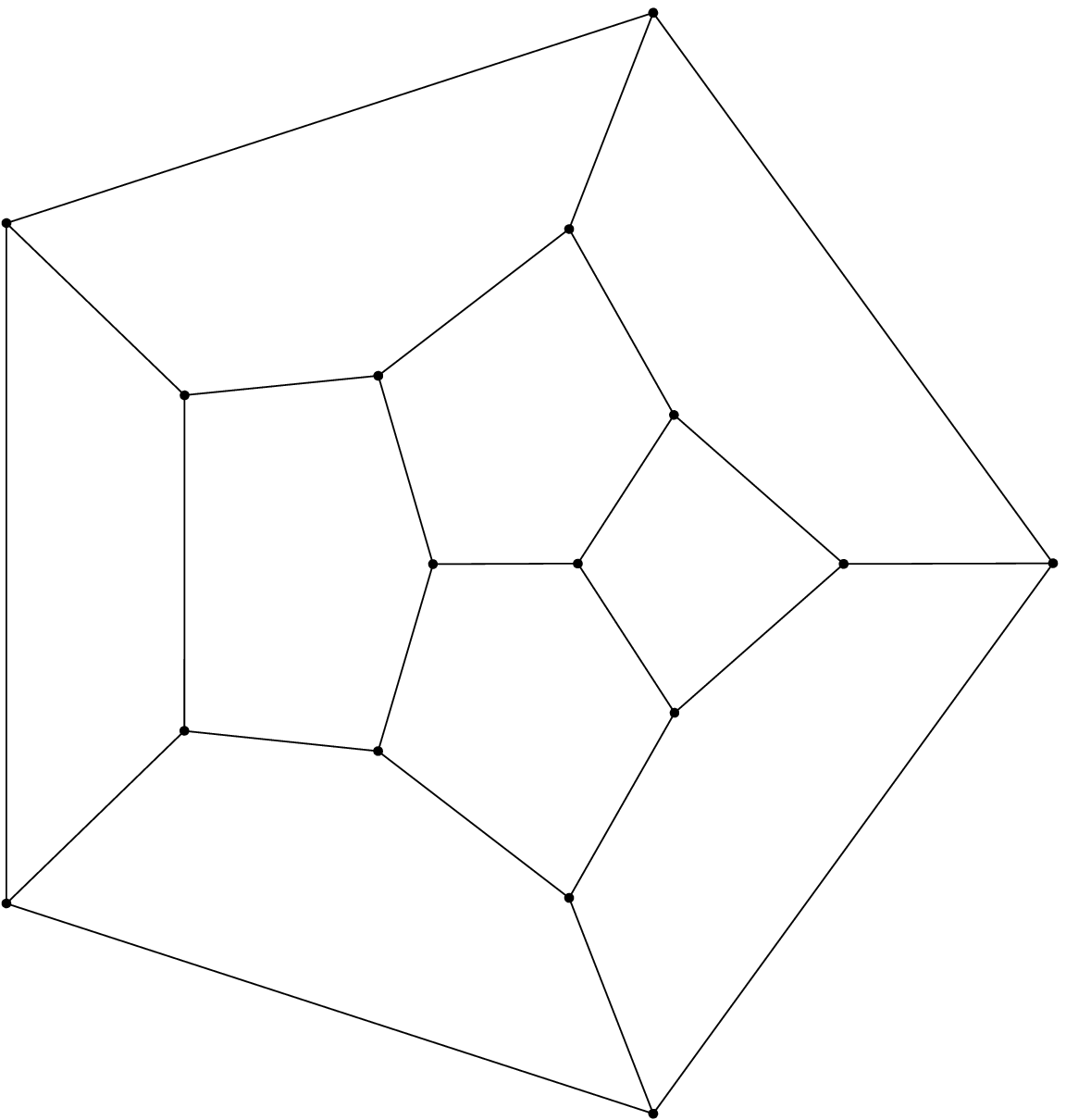}\par
$C_s$, nonext.
%PAIR1
\end{minipage}
\begin{minipage}{3cm}
\centering
\epsfig{height=20mm, file=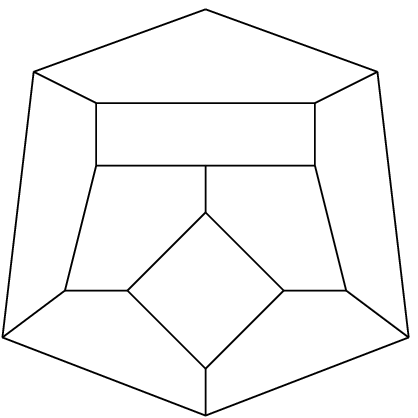}\par
$C_s$, nonext.
\end{minipage}
\begin{minipage}{3cm}
\centering
\epsfig{height=20mm, file=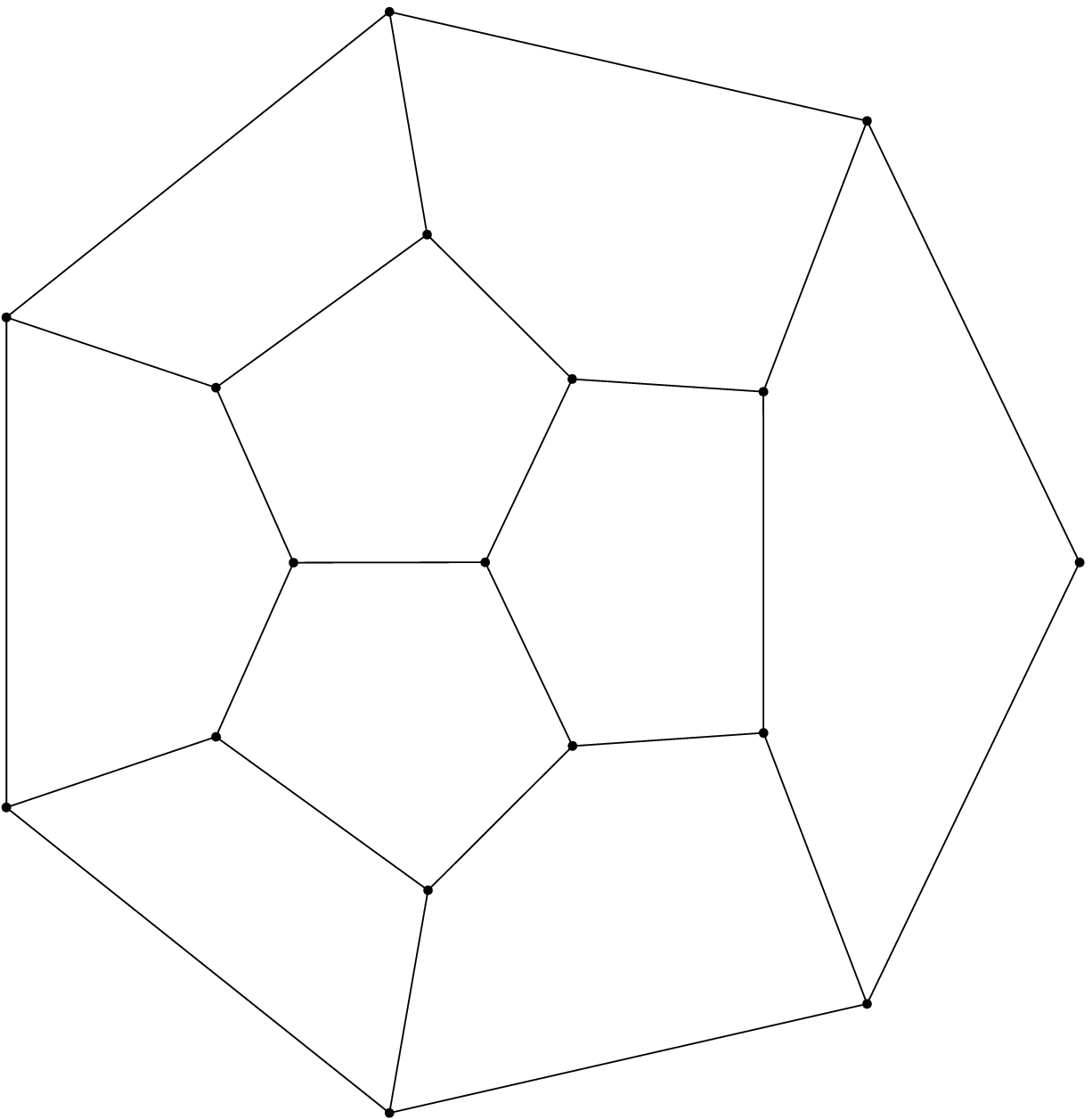}\par
$C_s$, nonext.
\end{minipage}
\begin{minipage}{3cm}
\centering
\epsfig{height=20mm, file=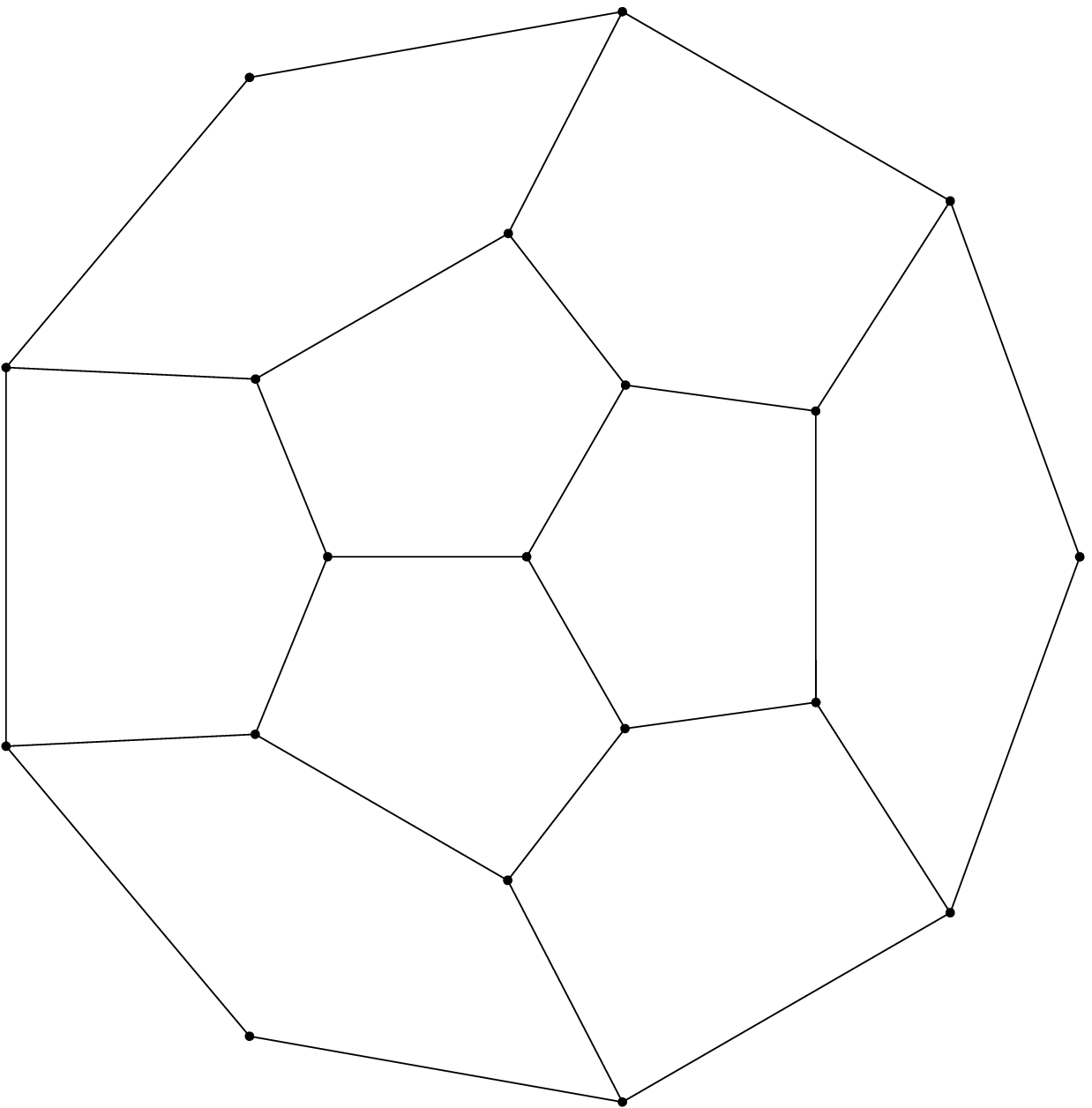}\par
$C_{3\nu}$
\end{minipage}
\begin{minipage}{3cm}
\centering
\epsfig{height=20mm, file=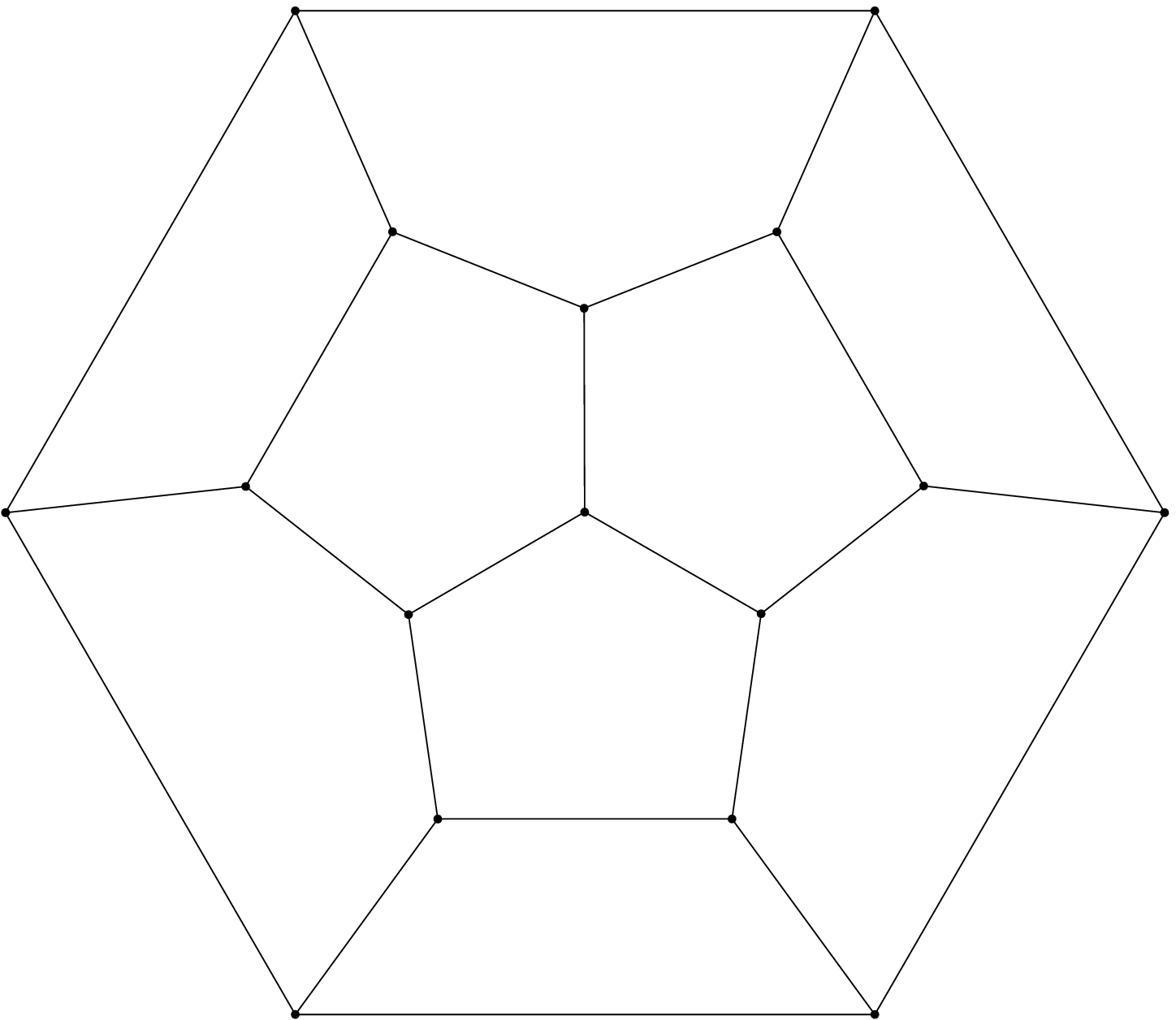}\par
$C_{3\nu}$, nonext.
\end{minipage}
\begin{minipage}{3cm}
\centering
\epsfig{height=20mm, file=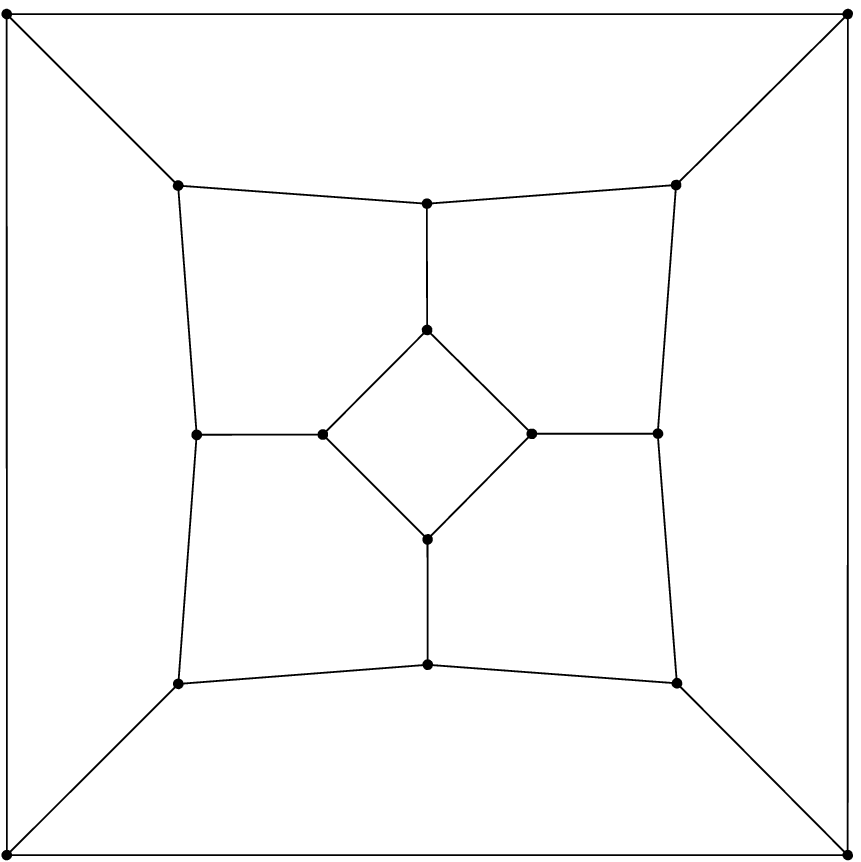}\par
$C_{4\nu}$, nonext.
%PAIR1
\end{minipage}

\end{center}
List of sporadic elementary $(\{3,4,5\},3)$-polycycles with $10$ faces:
\begin{center}
\begin{minipage}{3cm}
\centering
\epsfig{height=20mm, file=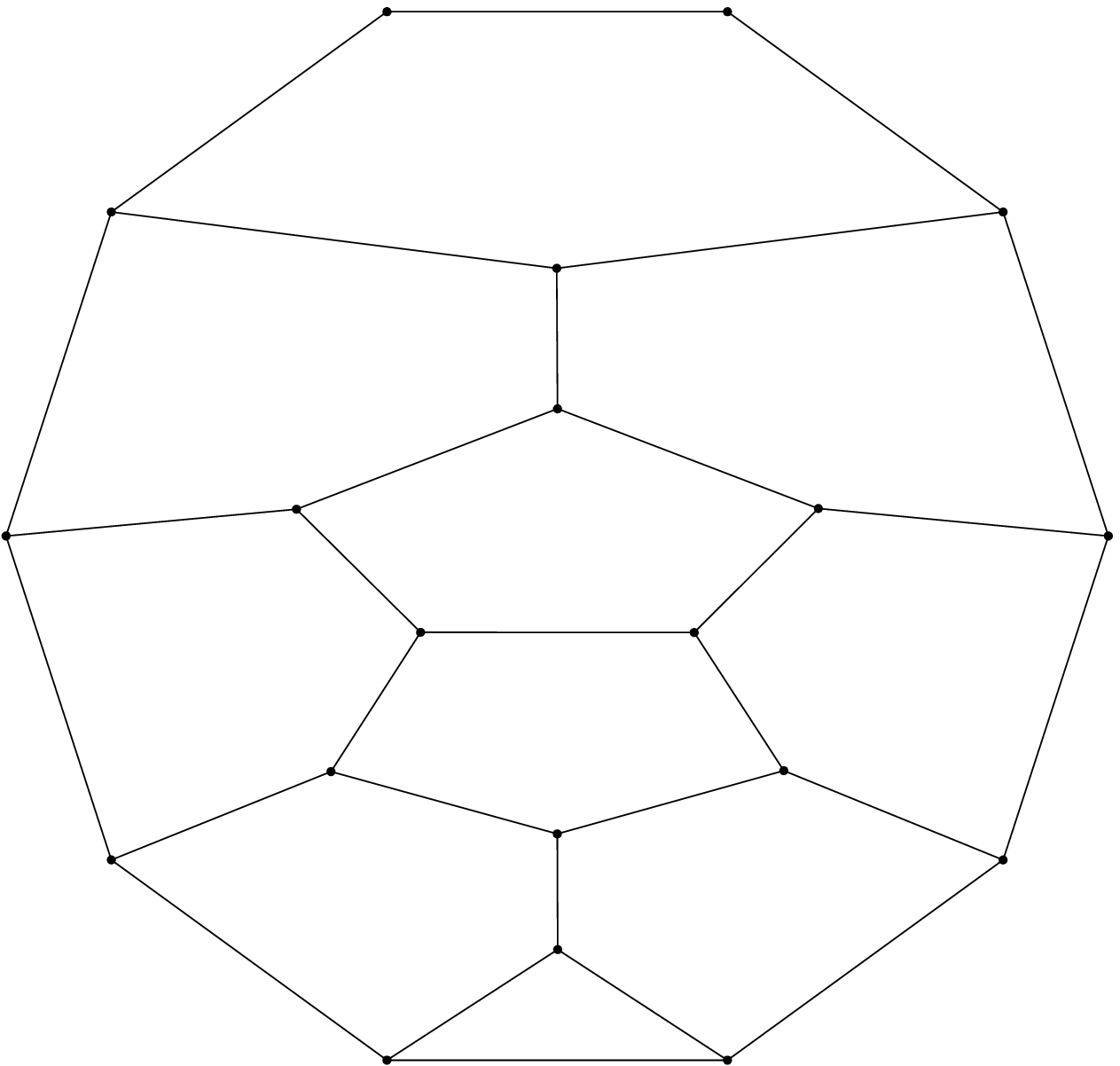}\par
$C_s$
\end{minipage}
\begin{minipage}{3cm}
\centering
\epsfig{height=20mm, file=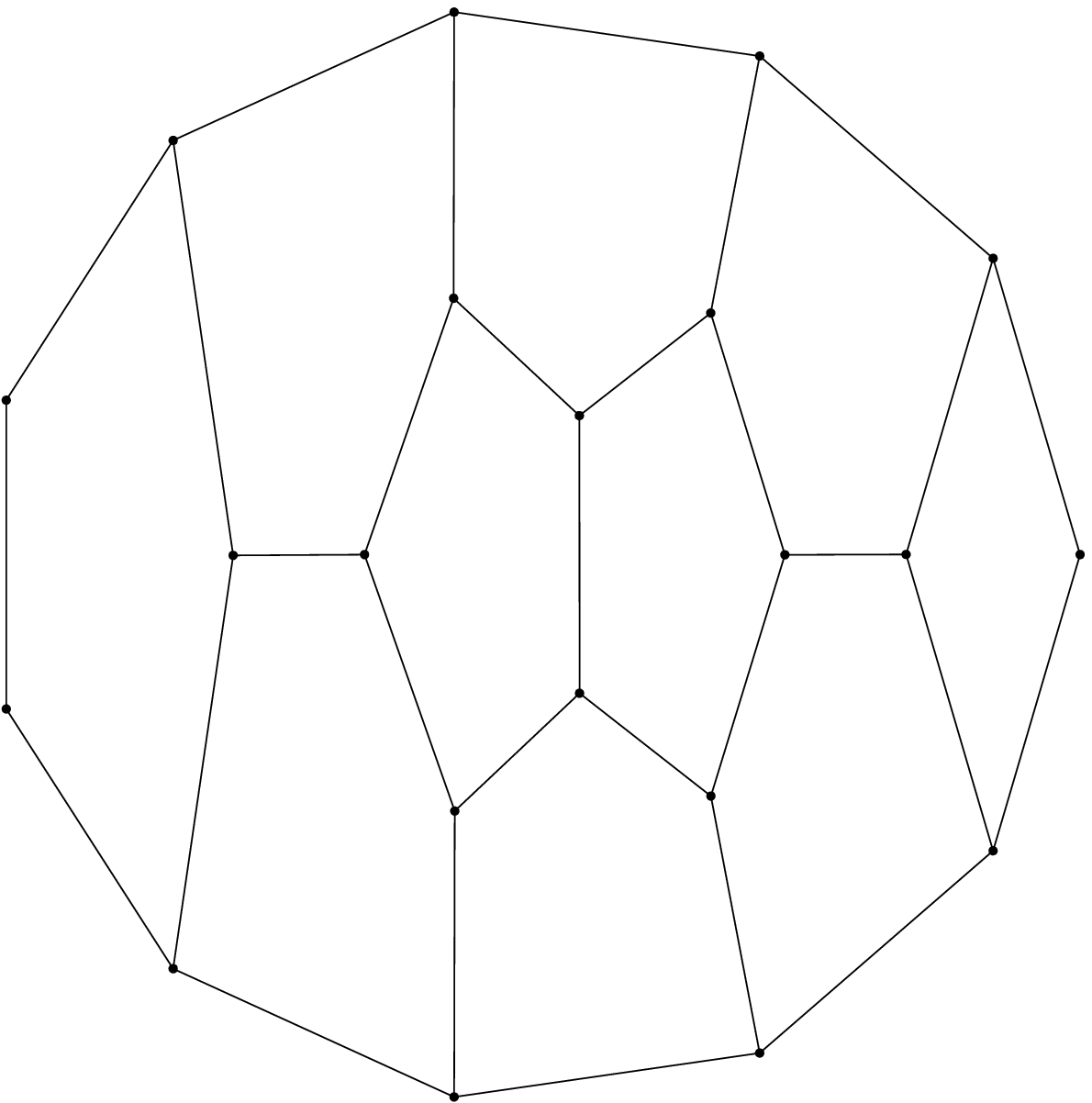}\par
$C_s$
\end{minipage}
\begin{minipage}{3cm}
\centering
\epsfig{height=20mm, file=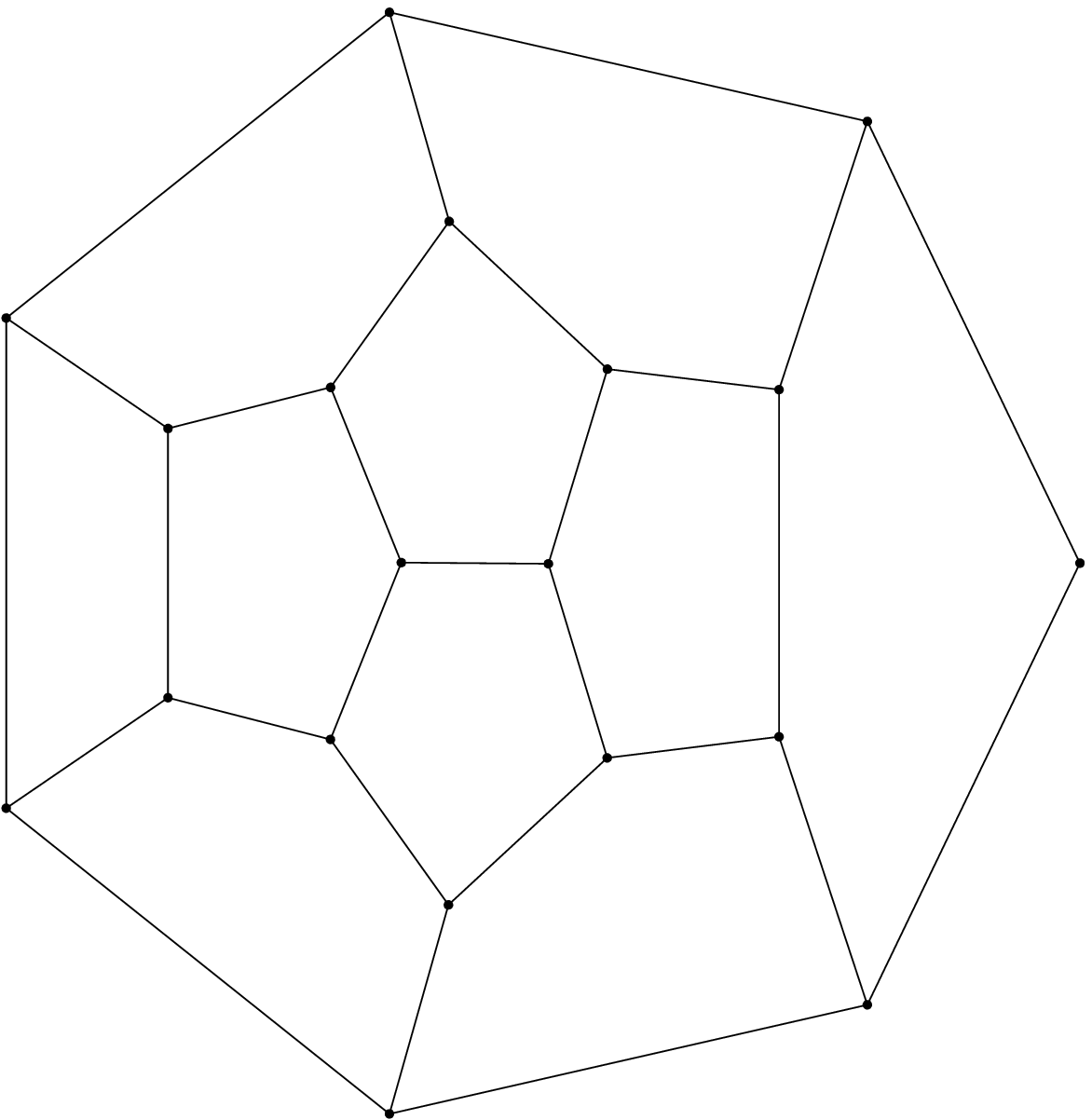}\par
$C_s$, nonext.
\end{minipage}
\begin{minipage}{3cm}
\centering
\epsfig{height=20mm, file=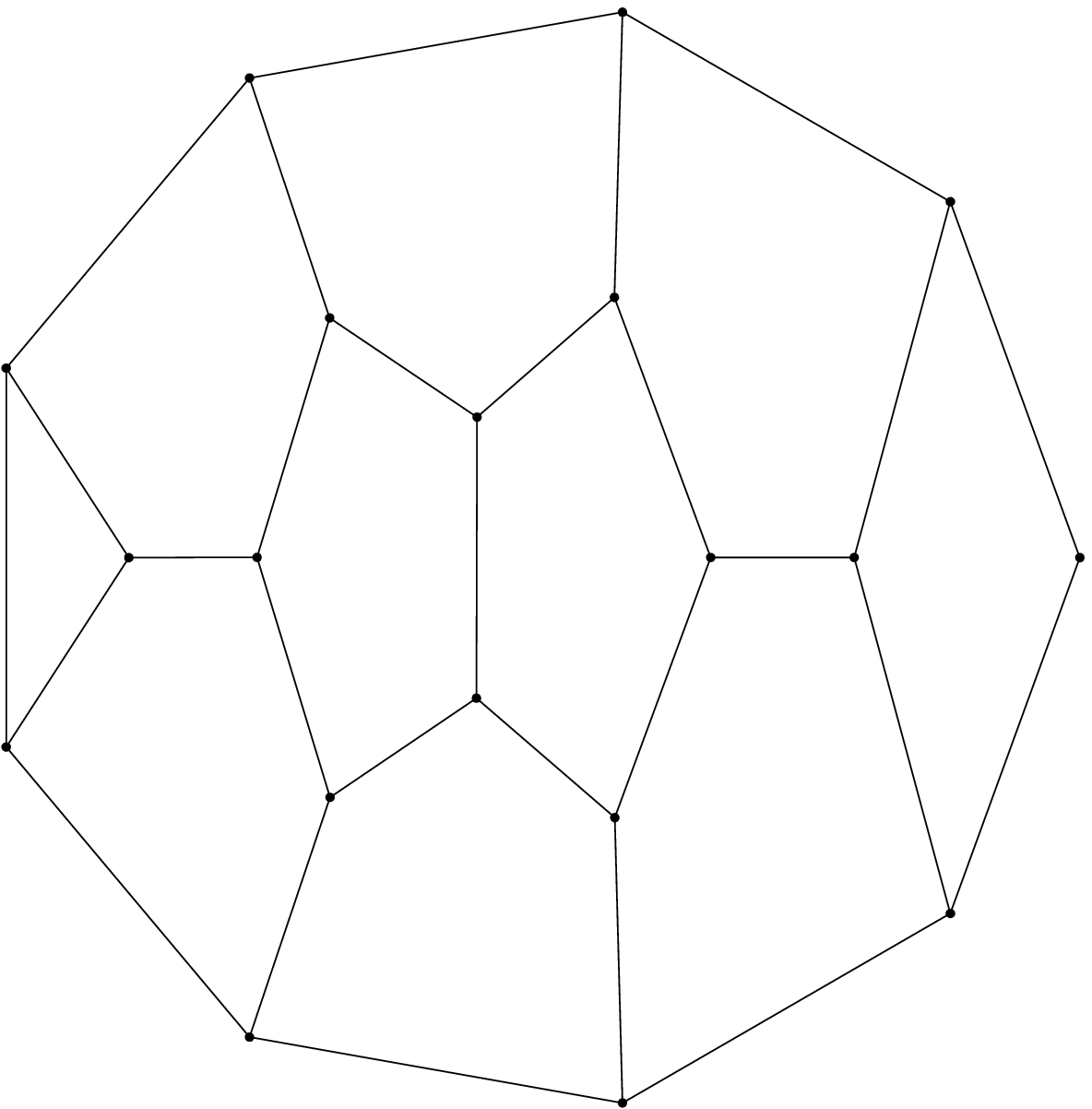}\par
$C_s$, nonext.
\end{minipage}
\begin{minipage}{3cm}
\centering
\epsfig{height=20mm, file=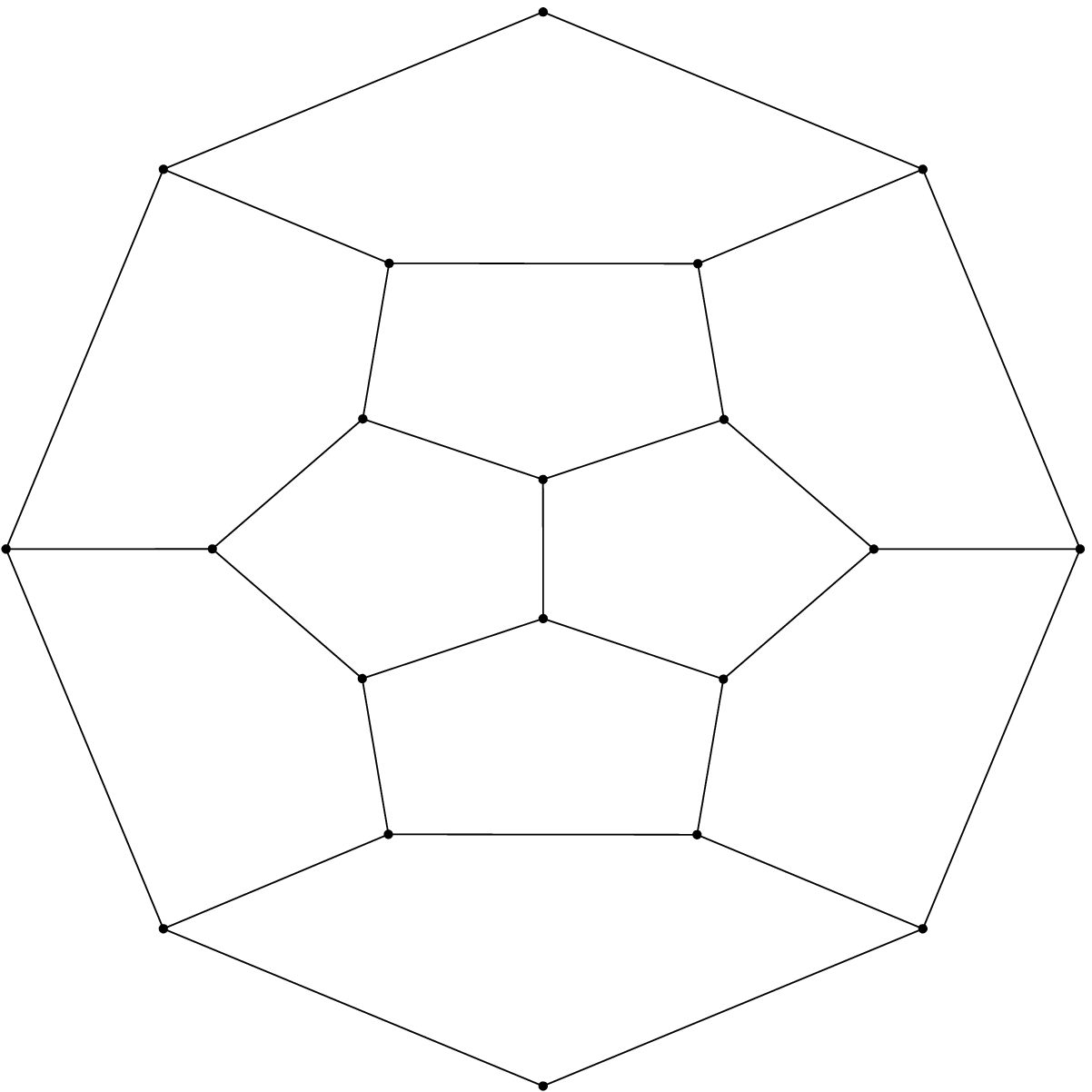}\par
$C_{2\nu}$
\end{minipage}
\begin{minipage}{3cm}
\centering
\epsfig{height=20mm, file=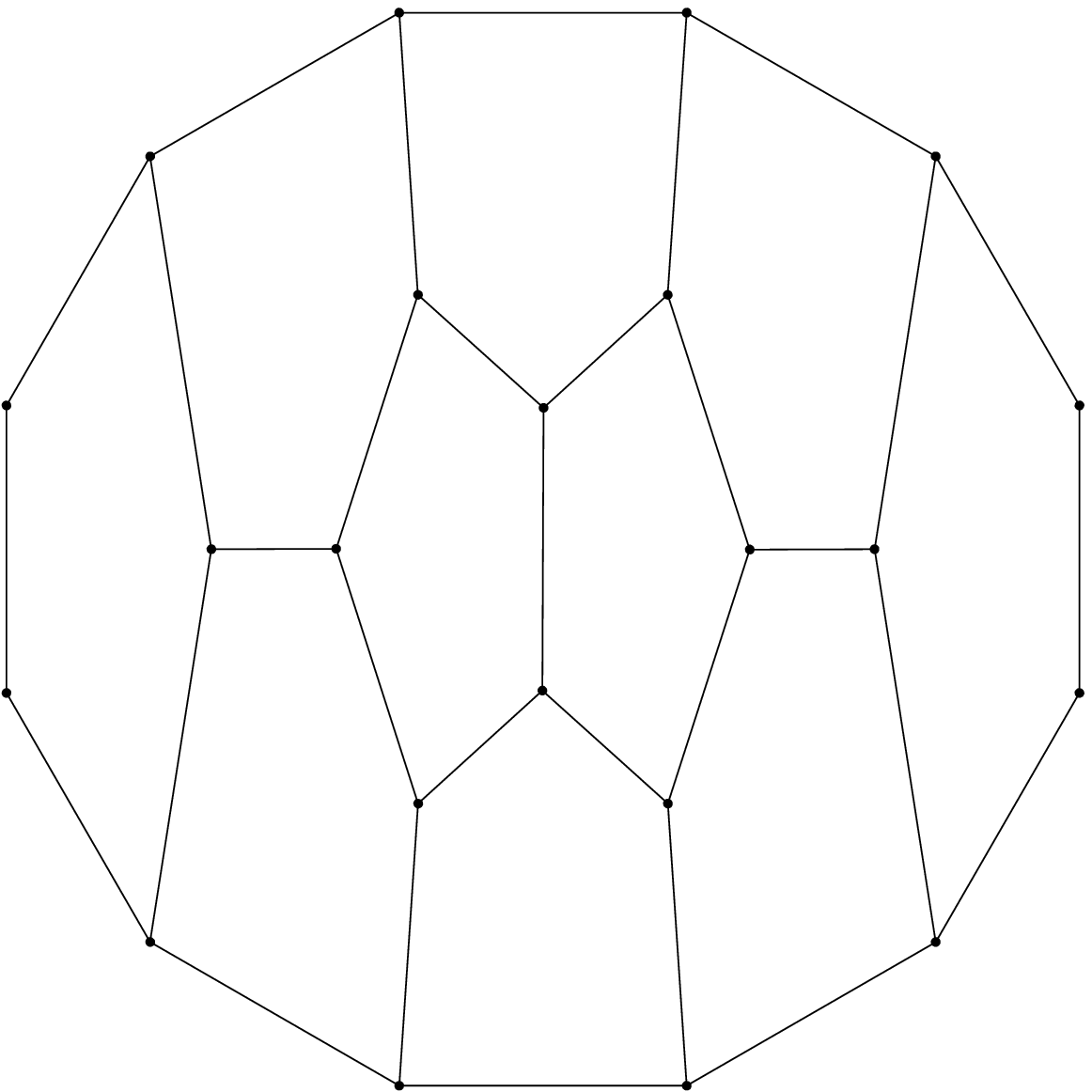}\par
$C_{2\nu}$
\end{minipage}
\begin{minipage}{3cm}
\centering
\epsfig{height=20mm, file=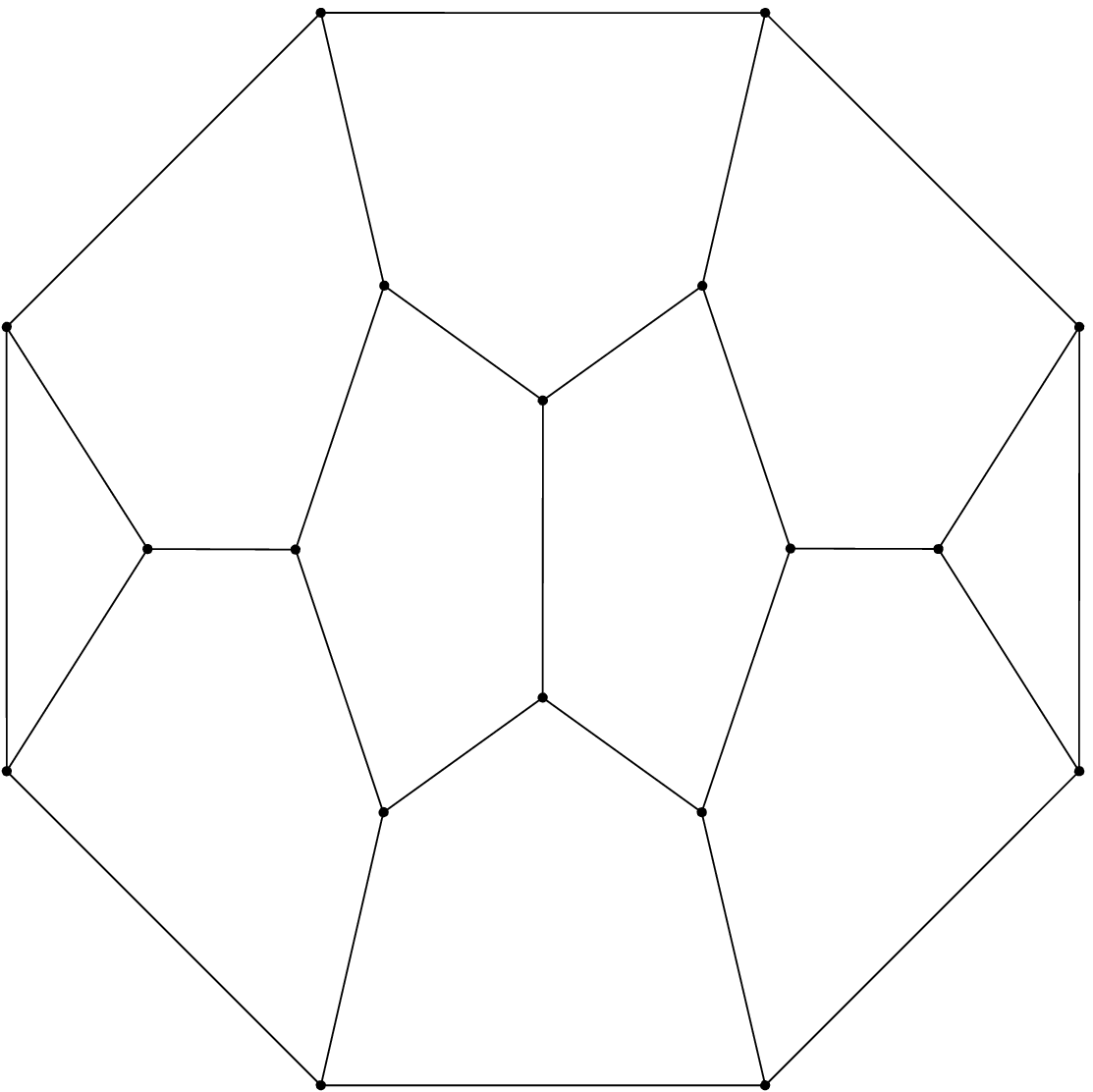}\par
$C_{2\nu}$, nonext.
\end{minipage}
\begin{minipage}{3cm}
\centering
\epsfig{height=20mm, file=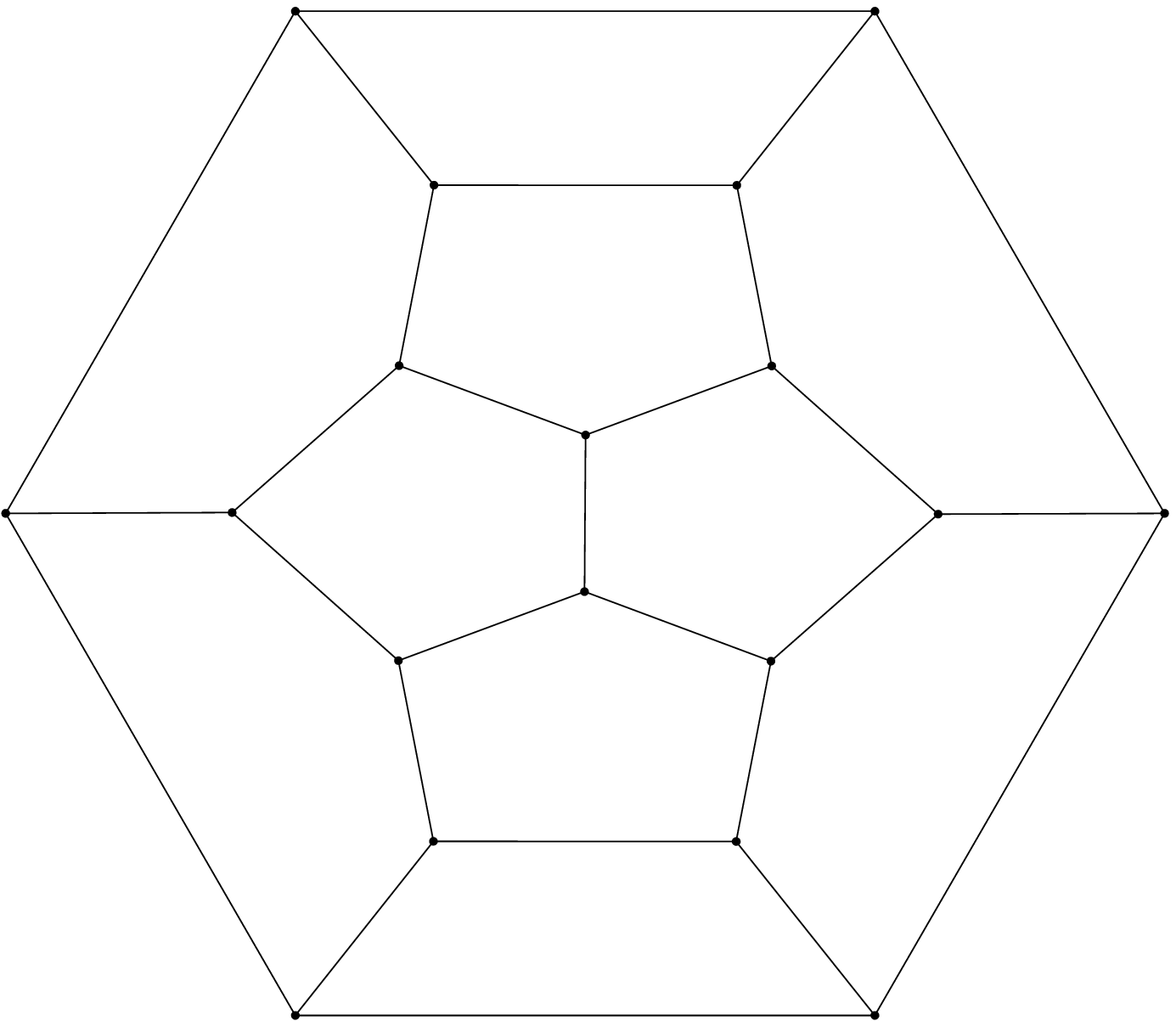}\par
$C_{2\nu}$, nonext.
\end{minipage}
\begin{minipage}{3cm}
\centering
\epsfig{height=20mm, file=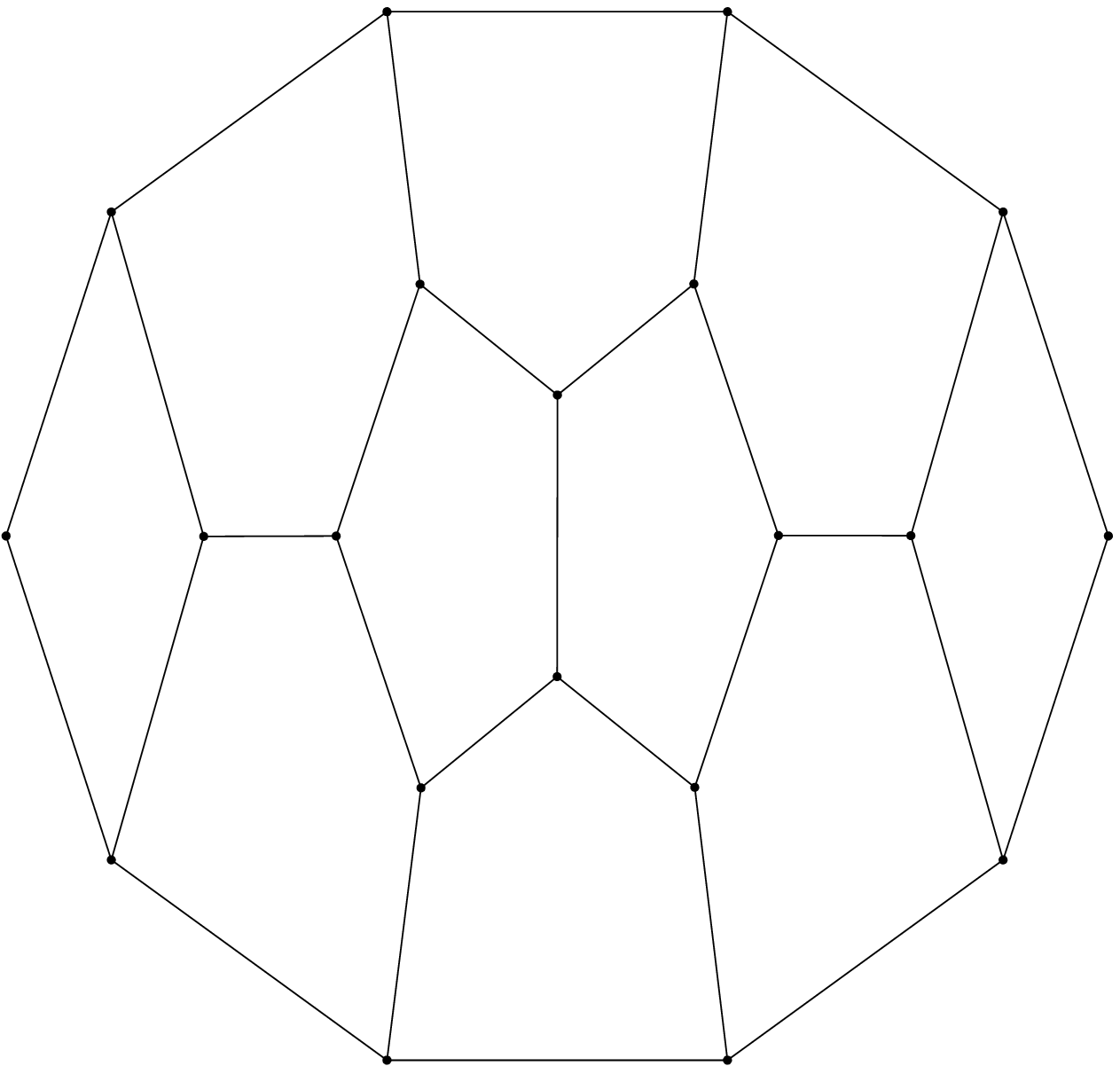}\par
$C_{2\nu}$, nonext.
\end{minipage}

\end{center}
List of sporadic elementary $(\{3,4,5\},3)$-polycycles with at least $11$ faces:
\begin{center}
\begin{minipage}{5cm}
\centering
\epsfig{height=20mm, file=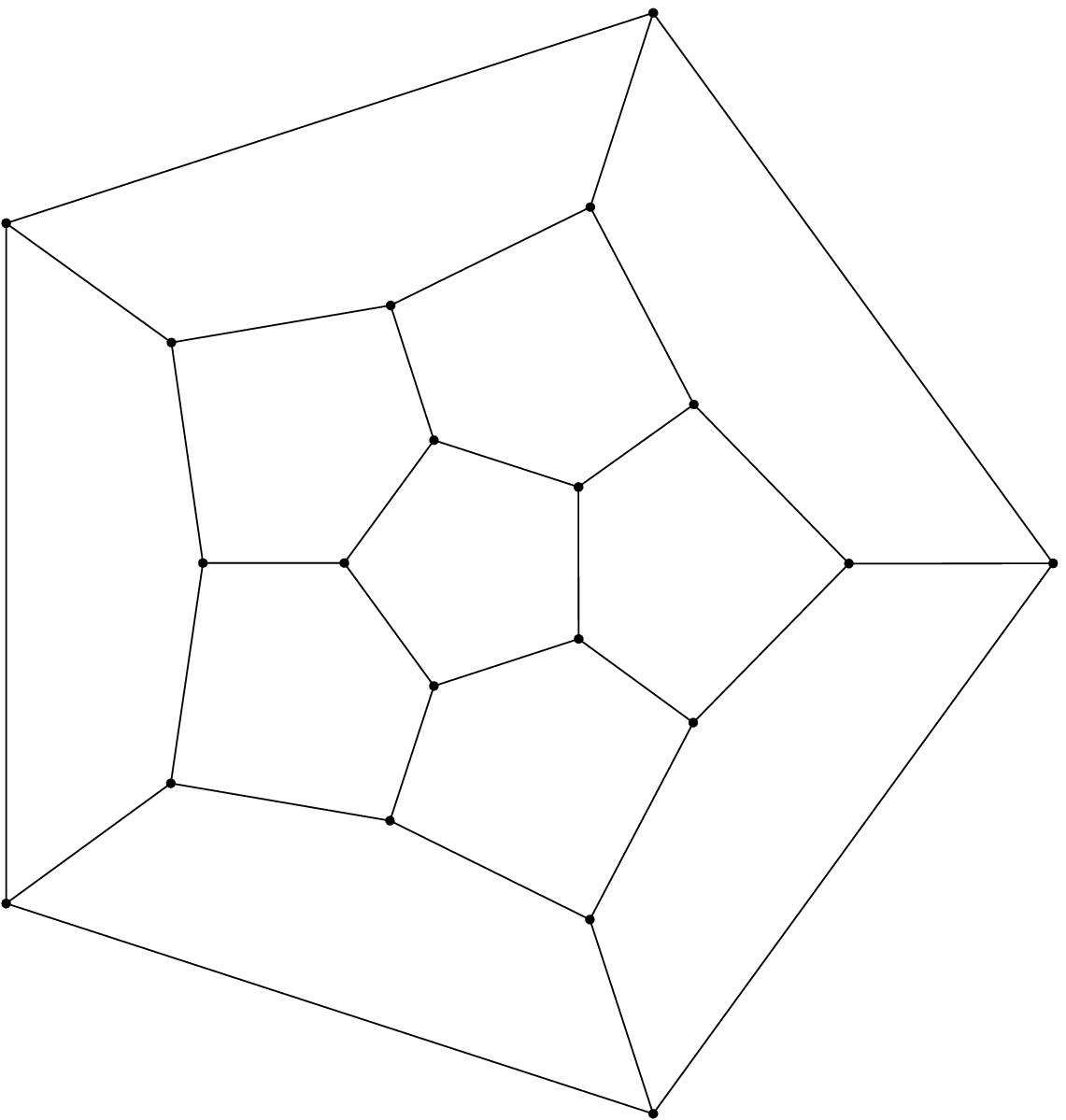}\par
$C_{5\nu}$, nonext. $(I_h)$
\end{minipage}

\end{center}

\section{Classification of elementary $(\{2,3\},4)$-polycycles}

\begin{theorem}\label{Theorem4valent}
Any elementary $(\{2,3\}, 4)$-polycycle is one of the following eight:

\begin{center}
\begin{minipage}{15cm}
\begin{minipage}{3.5cm}
\centering
\resizebox{2.0cm}{!}{\includegraphics{ElementaryDrawing/3gon.eps}}\par
$C_{3\nu}$ $(D_{3h})$
\end{minipage}
\begin{minipage}{3.5cm}
\centering
\epsfig{height=20mm, file=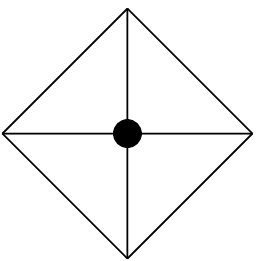}\par
$C_{4\nu}$
\end{minipage}
\begin{minipage}{3.5cm}
\centering
\epsfig{height=20mm, file=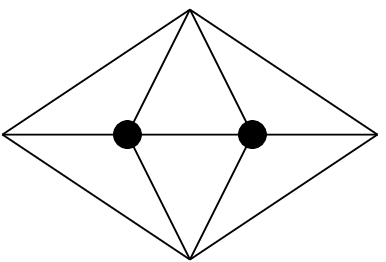}\par
$C_{2\nu}$
\end{minipage}
\begin{minipage}{3.5cm}
\centering
\resizebox{20mm}{!}{\rotatebox{0}{\includegraphics{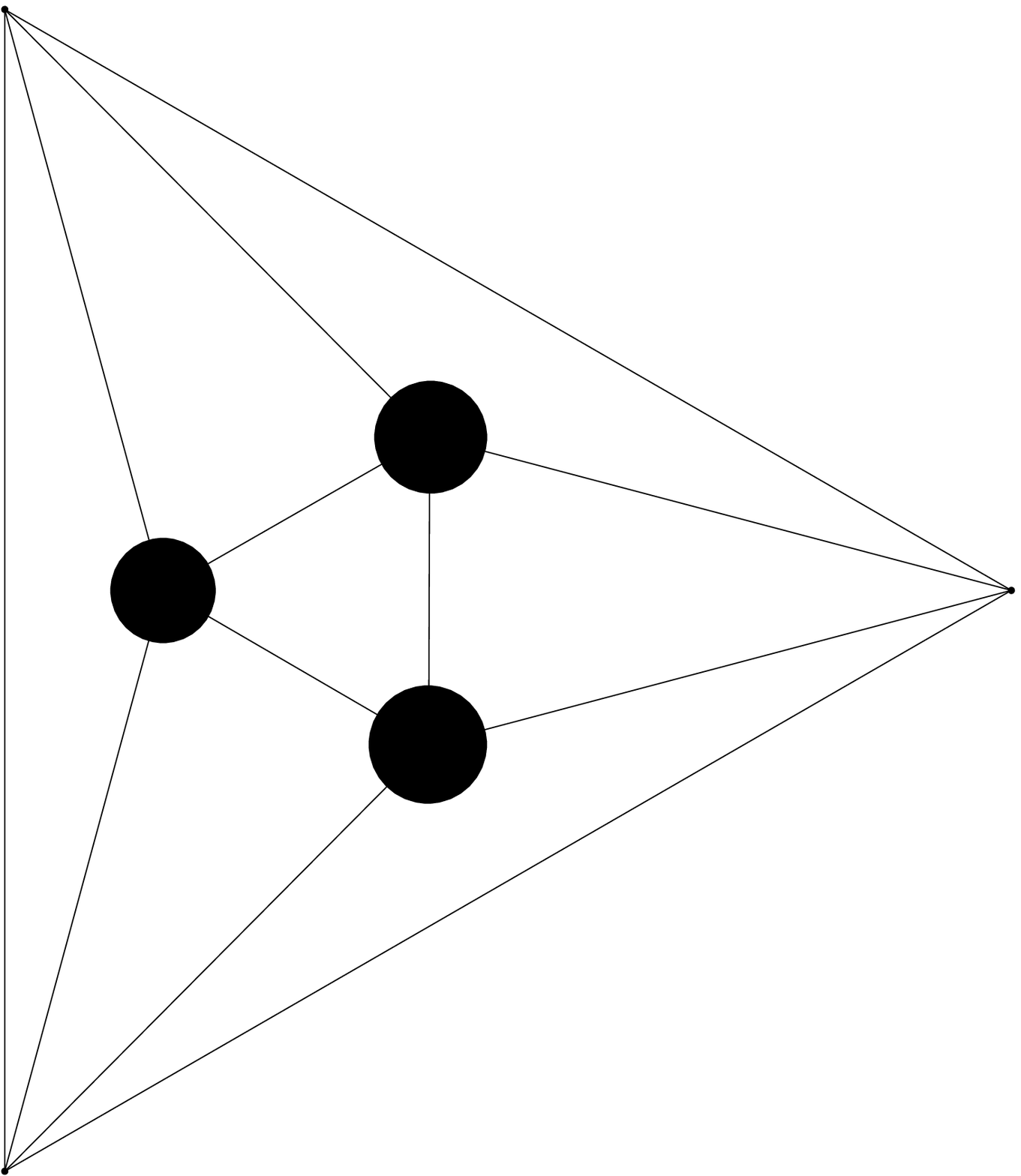}}}\par
$C_{3\nu}$, nonext. $(O_h)$
\end{minipage}
\begin{minipage}{3.5cm}
\centering
\epsfig{height=20mm, file=ElementaryDrawing/2gon.eps}\par
$C_{2\nu}$ $(D_{2h})$
\end{minipage}
\begin{minipage}{3.5cm}
\centering
\epsfig{height=20mm, file=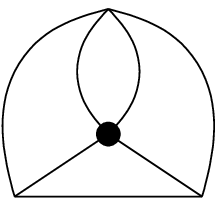}\par
$C_s$
\end{minipage}
\begin{minipage}{3.5cm}
\centering
\epsfig{height=20mm, file=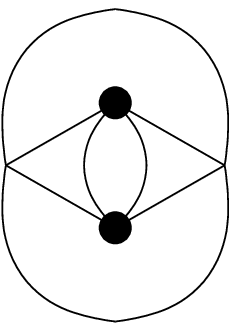}\par
$C_{2\nu}$,~nonext.~$(D_{2d})$
\end{minipage}
\begin{minipage}{3.5cm}
\centering
\epsfig{height=20mm, file=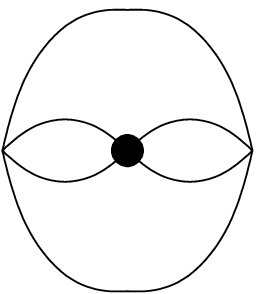}\par
$C_{2\nu}$,~nonext.~$(D_{3h})$
\end{minipage}

\end{minipage}

\end{center}

\end{theorem}
\proof The list of elementary $(\{3\},4)$-polycycles is determined in \cite{DDS2} and consists of the first four graphs of this theorem.
Let $P$ be a $(\{2,3\},4)$-polycycle, containing a $2$-gon.
If $|F_1|=1$, then it is
the $2$-gon. Clearly, the case, in which two $2$-gons share one edge,
is impossible. Assume that $P$ contains two $2$-gons, which share a
vertex. Then we should add triangle on both sides and so, obtain 
the second above polycycle.
If there is a $2$-gon, which does not share a vertex with a $2$-gon, then
$P$ contains the following pattern:
\begin{center}
\epsfig{height=20mm, file=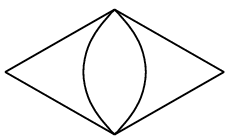}\par
\end{center}
So, clearly, $P$ is one of the last two possibilities above. \qed

Note that seventh and fourth polycycles in Theorem \ref{Theorem4valent}
are, respectively, $2$- and $3$-antiprisms; here the exterior face is
the unique hole. The $m$-antiprism for any $m\geq 2$ can also be
seen as $(\{2,3\},4)$-polycycle with $F_2$ consisting of the exterior
and interior $m$-gons; this polycycle is not elementary.

%Together with the $4$-valent plane graphs with $2$- or $3$-gonal
%faces, this form the complete list of elementary $(\{2,3\},4)$-polycycles.

\section{Classification of elementary $(\{2,3\},5)$-polycycles}

Let us consider an elementary $(\{2,3\},5)$-polycycle $P$.
Assume that $P$ is not an $i$-gon and has a $2$-gonal face $f$.
If $f$ is adjacent to a hole, then the polycycle is not
elementary.
So, holes are adjacent only to $3$-gons. If one remove such a
$3$-gon $t$, then the third vertex $v$ of $t$, which is necessarily
interior in $P$, becomes non-interior in $P-t$. The polycycle $P-t$ is 
not
necessarily elementary. Let us denote by $e_1,\dots, e_5$ the edges
incident to $v$ and assume that $e_1$, $e_2$ are edges of $t$.
The boundary is adjacent only to $3$-gons.
The potential bridges in $P-t$ are $e_3$, $e_4$ and $e_5$. Let
us check all five cases:
\begin{itemize}
\item If no edge $e_k$ is a bridge, then $P-t$ is
elementary.
\item If only $e_4$ is a bridge, then two cases can happen:
\begin{itemize}
\item This bridge goes from a hole to the same hole. This means
that $P$ is formed by the merging of two elementary $(\{2,3\},5)$-polycycles.
\item This bridge goes from a hole to another hole. This means
that $P$ has at least two holes and that $P-t$ is formed by
the merging of two open edges of an elementary
$(\{2,3\},5)$-polycycle, which has one hole less.
\end{itemize}
\item If $e_3$ or $e_5$ is a bridge, then $P-t$ is formed by the agglomeration
of an elementary $(\{2,3\},5)$-polycycle and a $i$-gon
with $i=2$ or $3$.
\item If $e_4$ is a bridge and $e_3$ or $e_5$ is a bridge, then $P-t$ is
formed by the agglomeration of an elementary $(\{2,3\},5)$-polycycle
and two $i$-gons with $i=2$ or $3$.
\item If all $e_k$ are bridges, then $P$ has only one interior vertex.
\end{itemize}

Given a hole of a $(R,q)$-polycycle, its {\em boundary sequence} is
the sequence of degrees of all consecutive vertices of the boundary of
this hole.

\begin{figure}
\begin{center}
\begin{center}
Infinite series $\alpha\alpha$ of elementary $(\{2,3\},5)$-polycycles:
\begin{center}
\begin{minipage}{2cm}
\centering
\epsfig{height=15mm, file=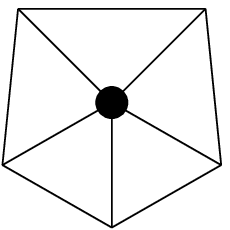}\par
$C_{5\nu}$
\end{minipage}
\begin{minipage}{2cm}
\centering
\epsfig{height=17mm, file=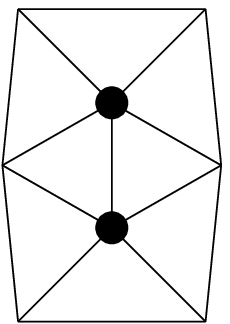}\par
$C_{2\nu}$
\end{minipage}
\begin{minipage}{3cm}
\centering
\epsfig{width=28mm, file=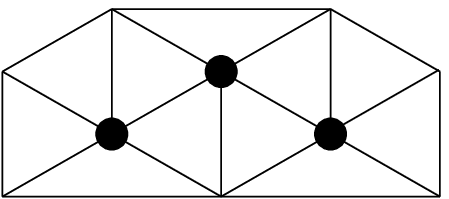}\par
$C_s$
\end{minipage}
\begin{minipage}{3cm}
\centering
\epsfig{width=28mm, file=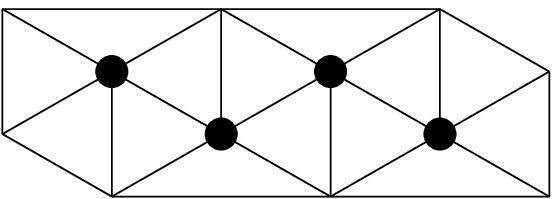}\par
$C_2$
\end{minipage}
\begin{minipage}{3cm}
\centering
\epsfig{width=28mm, file=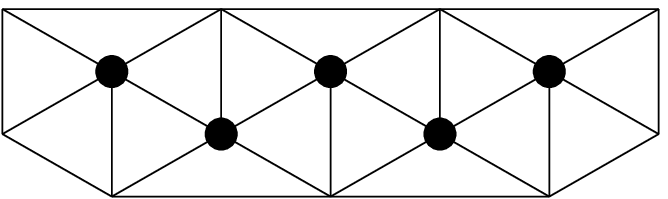}\par
$C_s$
\end{minipage}
\begin{minipage}{3cm}
\centering
\epsfig{width=28mm, file=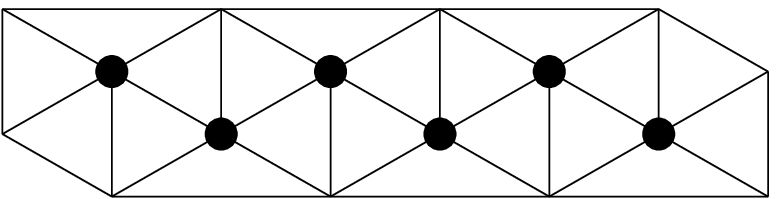}\par
$C_2$
\end{minipage}
%\begin{minipage}{3cm}
%\centering
%\epsfig{width=28mm, file=ElementaryDrawing/5_val7ver_2.eps}\par
%$C_s$
%\end{minipage}

\end{center}
\end{center}
\begin{center}
Infinite series $\alpha\beta$ of elementary $(\{2,3\},5)$-polycycles:
\begin{center}
\begin{minipage}{2cm}
\centering
\epsfig{height=16mm, file=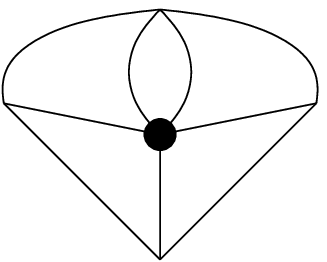}\par
$C_s$
\end{minipage}
\begin{minipage}{2cm}
\centering
\epsfig{height=16mm, file=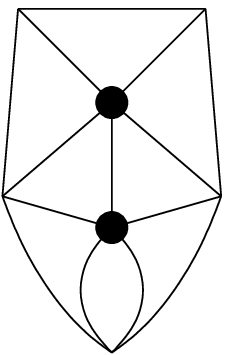}\par
$C_s$
\end{minipage}
\begin{minipage}{3cm}
\centering
\epsfig{width=28mm, file=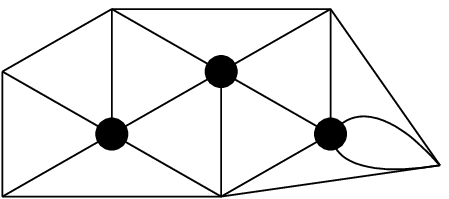}\par
$C_1$
\end{minipage}
\begin{minipage}{3cm}
\centering
\epsfig{width=28mm, file=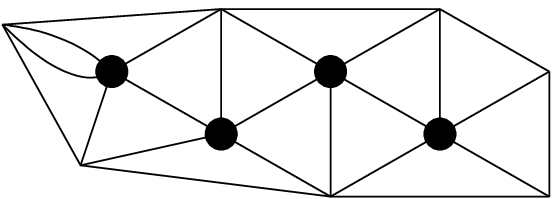}\par
$C_1$
\end{minipage}
\begin{minipage}{3cm}
\centering
\epsfig{width=28mm, file=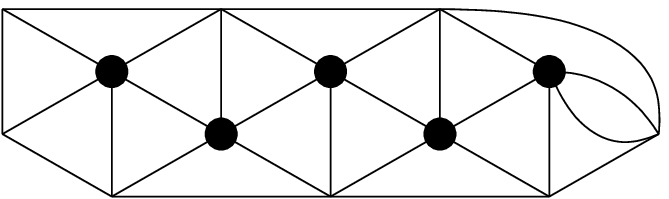}\par
$C_1$
\end{minipage}
\begin{minipage}{3cm}
\centering
\epsfig{width=28mm, file=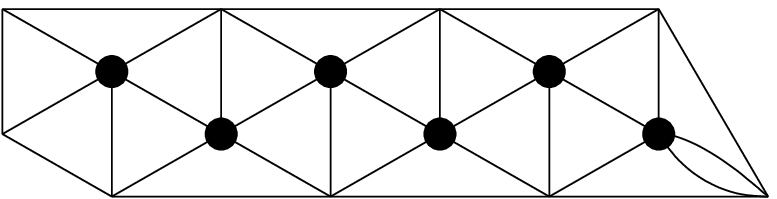}\par
$C_1$
\end{minipage}
%\begin{minipage}{3cm}
%\centering
%\epsfig{width=28mm, file=ElementaryDrawing/5_val7ver_3.eps}\par
%$C_1$
%\end{minipage}

\end{center}
\end{center}
\begin{center}
Infinite series $\alpha\gamma$ of elementary $(\{2,3\},5)$-polycycles:
\begin{center}
\begin{minipage}{3cm}
\centering
\epsfig{height=17mm, file=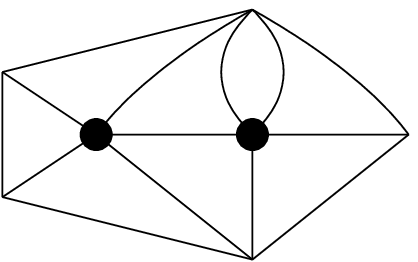}\par
$C_1$
\end{minipage}
\begin{minipage}{3cm}
\centering
\epsfig{width=28mm, file=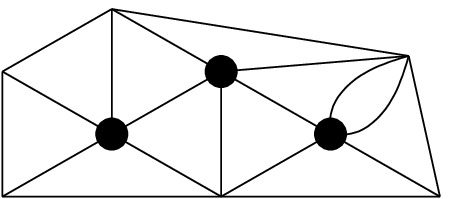}\par
$C_1$
\end{minipage}
\begin{minipage}{3cm}
\centering
\epsfig{width=28mm, file=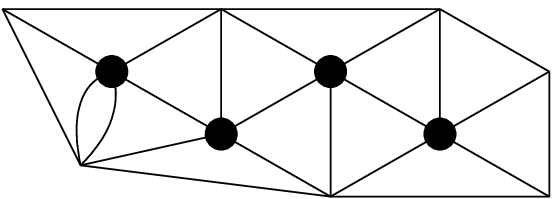}\par
$C_1$
\end{minipage}
\begin{minipage}{3cm}
\centering
\epsfig{width=28mm, file=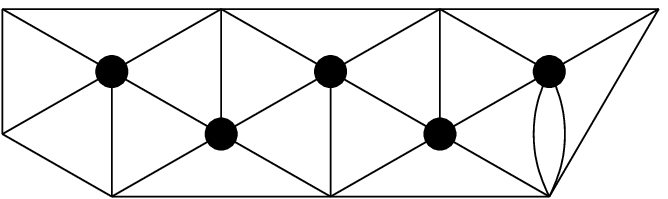}\par
$C_1$
\end{minipage}
\begin{minipage}{3cm}
\centering
\epsfig{width=28mm, file=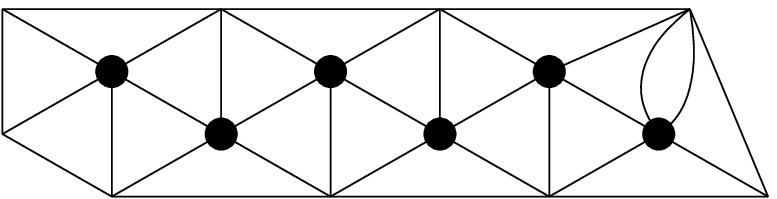}\par
$C_1$
\end{minipage}
%\begin{minipage}{3cm}
%\centering
%\epsfig{width=28mm, file=ElementaryDrawing/5_val7ver_4.eps}\par
%$C_1$
%\end{minipage}
\end{center}
\end{center}
\begin{center}
Infinite series $\beta\beta$ of elementary $(\{2,3\},5)$-polycycles:
\begin{center}
\begin{minipage}{3cm}
\centering
\epsfig{height=20mm, file=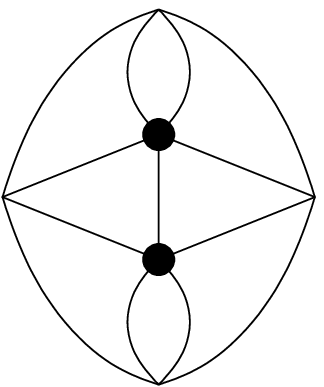}\par
$C_{2\nu}$
\end{minipage}
\begin{minipage}{3cm}
\centering
\epsfig{width=28mm, file=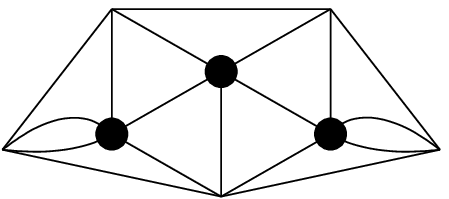}\par
$C_s$
\end{minipage}
\begin{minipage}{3cm}
\centering
\epsfig{width=28mm, file=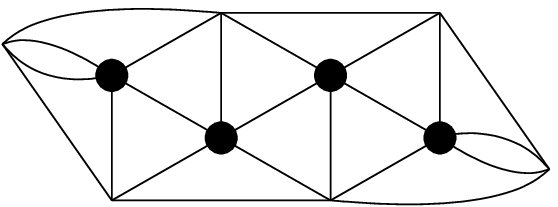}\par
$C_2$
\end{minipage}
\begin{minipage}{3cm}
\centering
\epsfig{width=28mm, file=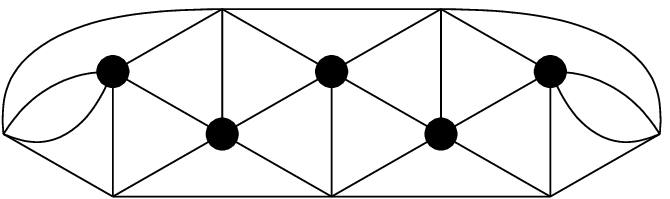}\par
$C_s$
\end{minipage}
\begin{minipage}{3cm}
\centering
\epsfig{width=28mm, file=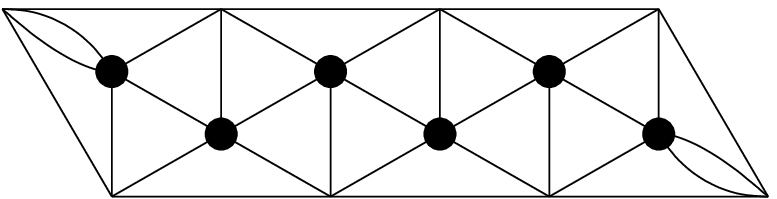}\par
$C_2$
\end{minipage}
%\begin{minipage}{3cm}
%\centering
%\epsfig{width=28mm, file=ElementaryDrawing/5_val7ver_5.eps}\par
%$C_s$, nonext.
%\end{minipage}

\end{center}
\end{center}
\begin{center}
Infinite series $\beta\gamma$ of elementary $(\{2,3\},5)$-polycycles:
\begin{center}
\begin{minipage}{3cm}
\centering
\epsfig{height=16mm, file=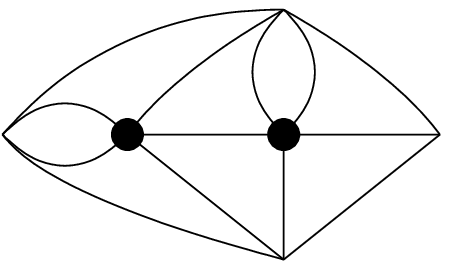}\par
$C_1$
\end{minipage}
\begin{minipage}{3cm}
\centering
\epsfig{width=28mm, file=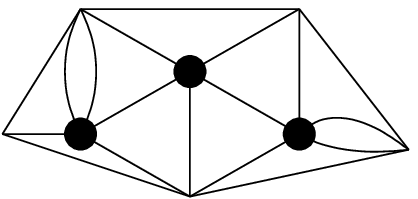}\par
$C_1$
\end{minipage}
\begin{minipage}{3cm}
\centering
\epsfig{width=28mm, file=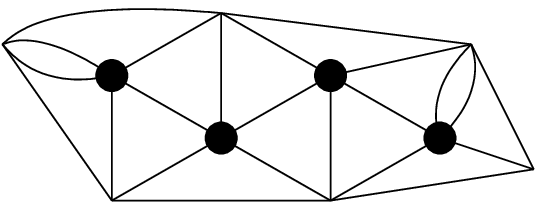}\par
$C_1$
\end{minipage}
\begin{minipage}{3cm}
\centering
\epsfig{width=28mm, file=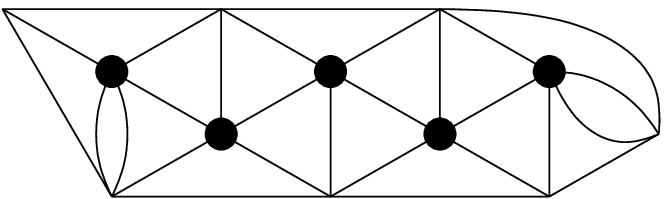}\par
$C_1$
\end{minipage}
\begin{minipage}{3cm}
\centering
\epsfig{width=28mm, file=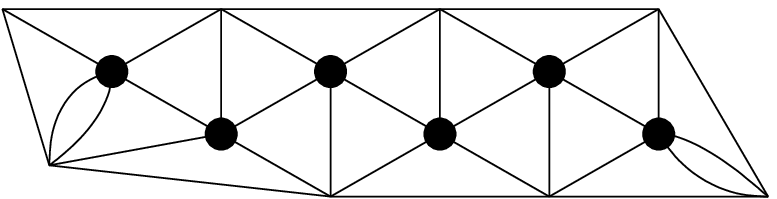}\par
$C_1$
\end{minipage}
%\begin{minipage}{3cm}
%\centering
%\epsfig{width=28mm, file=ElementaryDrawing/5_val7ver_6.eps}\par
%$C_1$, nonext.
%\end{minipage}

\end{center}
\end{center}
\begin{center}
Infinite series $\gamma\gamma$ of elementary $(\{2,3\},5)$-polycycles:
\begin{center}
\begin{minipage}{3cm}
\centering
\epsfig{height=16mm, file=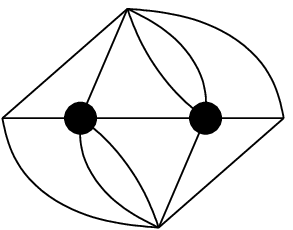}\par
$C_{2}$
\end{minipage}
\begin{minipage}{3cm}
\centering
\epsfig{width=24mm, file=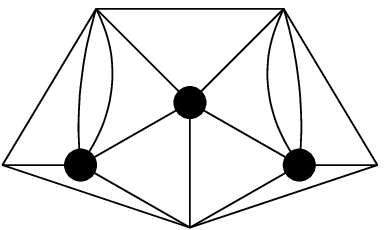}\par
$C_s$
\end{minipage}
\begin{minipage}{3cm}
\centering
\epsfig{width=28mm, file=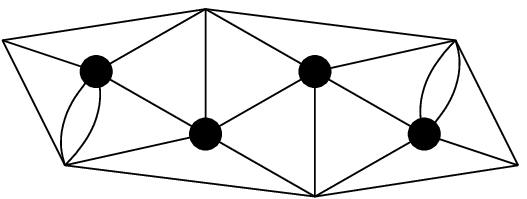}\par
$C_2$, nonext.
\end{minipage}
\begin{minipage}{3cm}
\centering
\epsfig{width=28mm, file=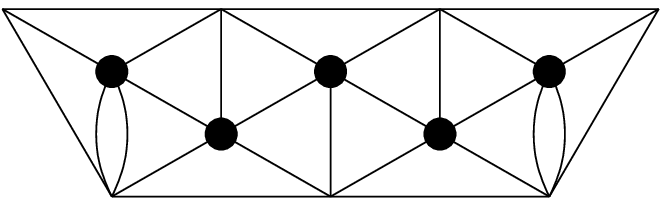}\par
$C_s$, nonext.
\end{minipage}
\begin{minipage}{3cm}
\centering
\epsfig{width=28mm, file=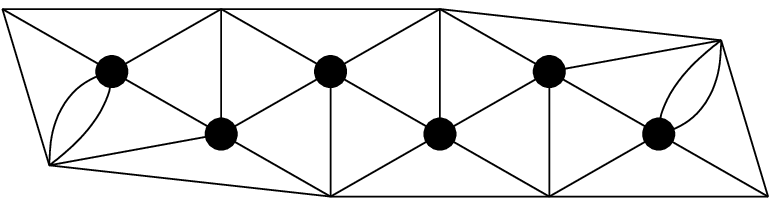}\par
$C_2$, nonext.
\end{minipage}
%\begin{minipage}{3cm}
%\centering
%\epsfig{width=28mm, file=ElementaryDrawing/5_val7ver_7.eps}\par
%$C_s$, nonext.
%\end{minipage}

\end{center}
\end{center}

\end{center}
\caption{The first $5$ members of the six infinite series of $(\{2,3\},5)$-polycycles}
\label{SixInfiniteSeries}
\end{figure}

\begin{figure}
\begin{center}
\begin{minipage}{3cm}
\centering
\resizebox{2.0cm}{!}{\includegraphics{ElementaryDrawing/3gon.eps}}\par
$C_{3\nu}$ $(D_{3h})$
\end{minipage}
\begin{minipage}{3cm}
\centering
\epsfig{height=20mm, file=ElementaryDrawing/2gon.eps}\par
$C_{2\nu}$ $(D_{2h})$
\end{minipage}
\begin{minipage}{3cm}
\centering
\epsfig{height=17mm, file=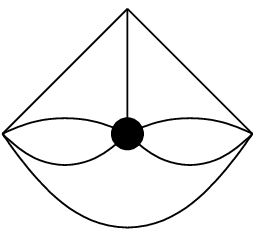}\par
$C_s$
% 1
\end{minipage}
\begin{minipage}{3cm}
\centering
\epsfig{height=20mm, file=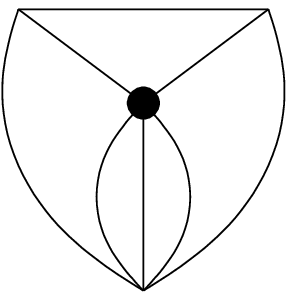}\par
$C_s$ ($C_{2v}$)
\end{minipage}
\begin{minipage}{3cm}
\centering
\epsfig{height=20mm, file=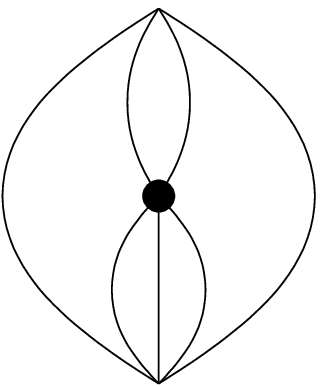}\par
$C_s$,~nonext.~$(C_{2\nu})$
% 3
\end{minipage}
\begin{minipage}{3cm}
\centering
\epsfig{height=20mm, file=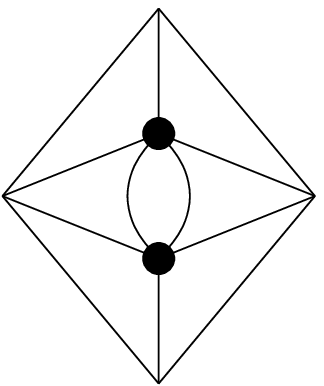}\par
$C_{2\nu}$
% 5
\end{minipage}
\begin{minipage}{3cm}
\centering
\epsfig{height=20mm, file=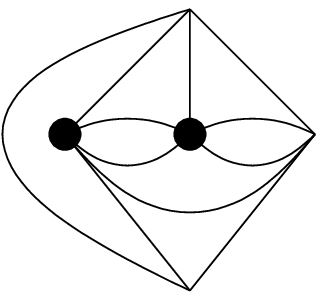}\par
$C_1$ $(C_s)$
% 2
\end{minipage}
\begin{minipage}{3cm}
\centering
\epsfig{height=20mm, file=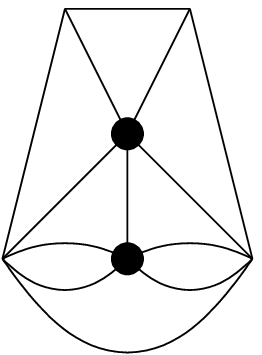}\par
$C_s$
% 4
\end{minipage}
\begin{minipage}{3cm}
\centering
\epsfig{height=20mm, file=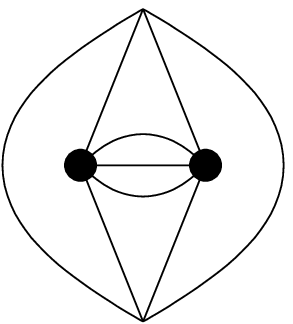}\par
$C_{2v}$, nonext.
\end{minipage}
\begin{minipage}{3cm}
\centering
\epsfig{height=20mm, file=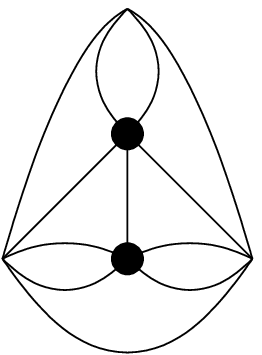}\par
$C_s$, nonext.
% 6
\end{minipage}
\begin{minipage}{3cm}
\centering
\epsfig{height=20mm, file=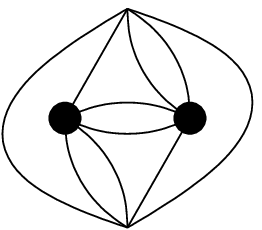}\par
$C_{2}$,~nonext.~($D_{2d}$)
\end{minipage}
\begin{minipage}{3cm}
\centering
\epsfig{width=25mm, file=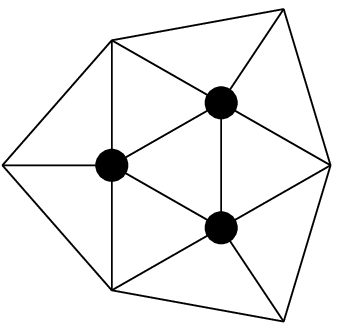}\par
$C_{3\nu}$
% 7
\end{minipage}
\begin{minipage}{3cm}
\centering
\epsfig{width=25mm, file=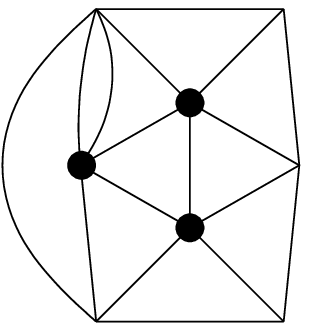}\par
$C_1$
% 8
\end{minipage}
\begin{minipage}{3cm}
\centering
\epsfig{width=28mm, file=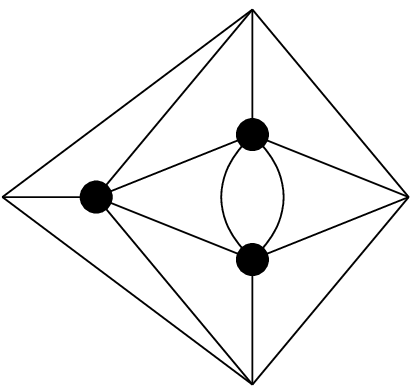}\par
$C_s$
% 9
\end{minipage}
\begin{minipage}{3cm}
\centering
\epsfig{width=25mm, file=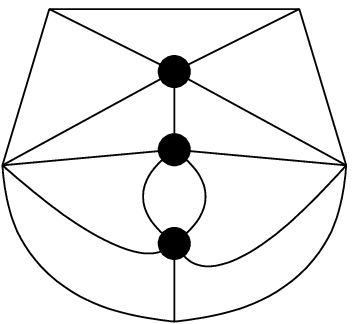}\par
$C_s$
% 10
\end{minipage}
\begin{minipage}{3cm}
\centering
\epsfig{width=28mm, file=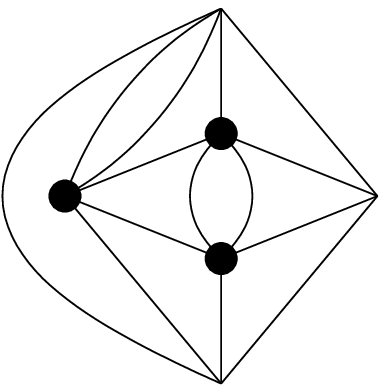}\par
$C_1$ $(C_2)$
% 11
\end{minipage}
\begin{minipage}{3cm}
\centering
\epsfig{width=25mm, file=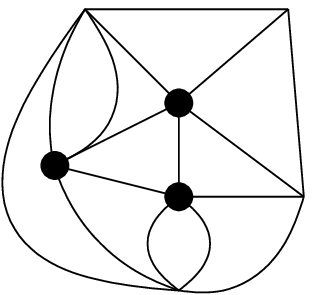}\par
$C_1$
% 12
\end{minipage}
\begin{minipage}{3cm}
\centering
\epsfig{width=28mm, file=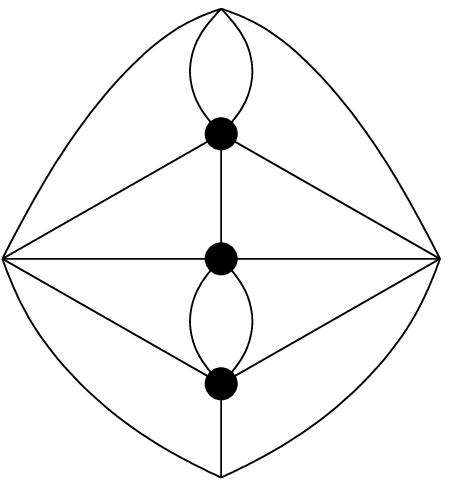}\par
$C_s$
% 15
\end{minipage}
\begin{minipage}{3cm}
\centering
\epsfig{width=28mm, file=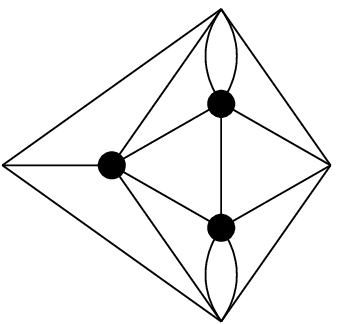}\par
$C_s$
% 14
\end{minipage}
%\begin{minipage}{3cm}
%\centering
%\epsfig{width=28mm, file=ElementaryDrawing/5_val3ver_15.eps}\par
%\end{minipage}
\begin{minipage}{3cm}
\centering
\epsfig{width=28mm, file=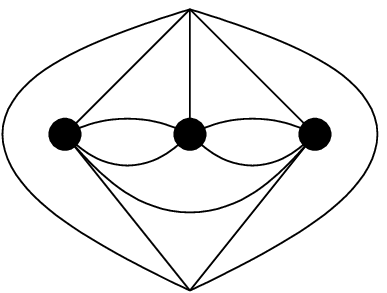}\par
$C_s$, nonext.
% 13
\end{minipage}
\begin{minipage}{3cm}
\centering
\epsfig{width=25mm, file=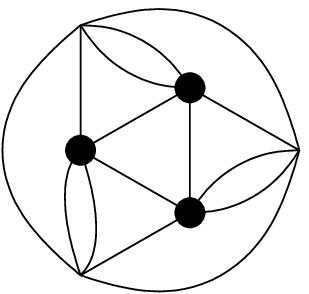}\par
$C_3$, nonext.
% 16
\end{minipage}
\begin{minipage}{3cm}
\centering
\epsfig{width=28mm, file=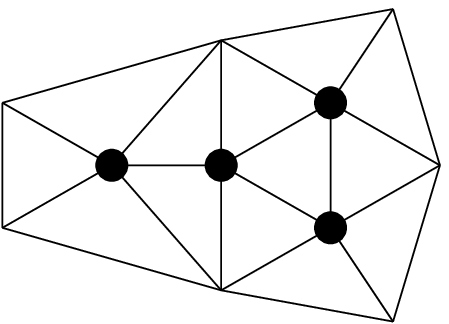}\par
$C_s$
% 18
\end{minipage}
\begin{minipage}{3cm}
\centering
\epsfig{width=28mm, file=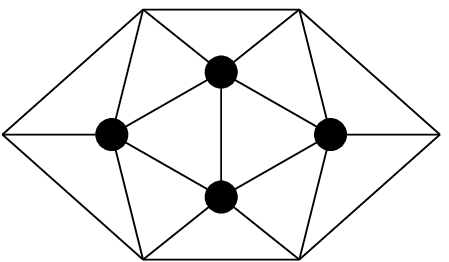}\par
$C_{2\nu}$
% 17
\end{minipage}
\begin{minipage}{3cm}
\centering
\epsfig{width=28mm, file=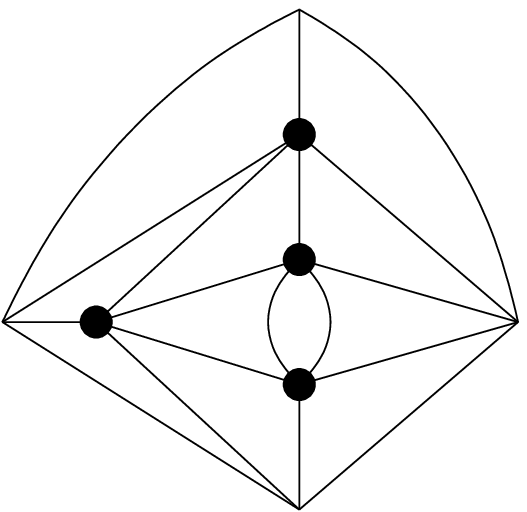}\par
$C_1$
% 19
\end{minipage}
\begin{minipage}{3cm}
\centering
\epsfig{width=28mm, file=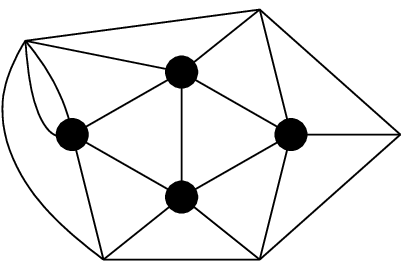}\par
$C_1$
% 20
\end{minipage}
\begin{minipage}{3cm}
\centering
\epsfig{width=28mm, file=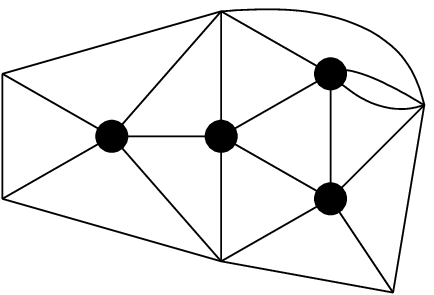}\par
$C_1$
% 21
\end{minipage}
\begin{minipage}{3cm}
\centering
\epsfig{width=28mm, file=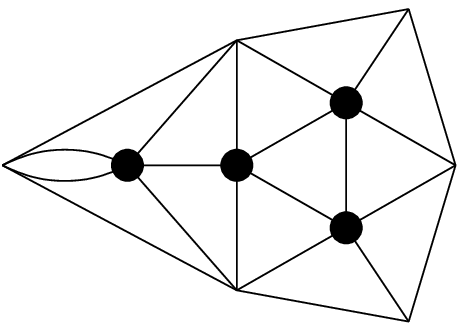}\par
$C_s$
% 22
\end{minipage}
\begin{minipage}{3cm}
\centering
\epsfig{width=28mm, file=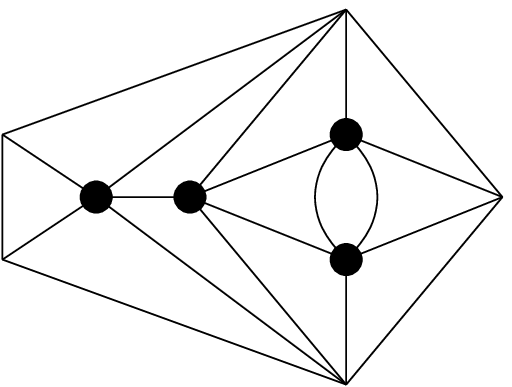}\par
$C_s$
% 23
\end{minipage}
\begin{minipage}{3cm}
\centering
\epsfig{width=28mm, file=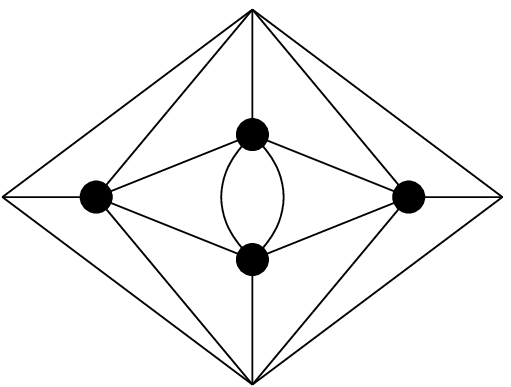}\par
$C_{2\nu}$
% 24
\end{minipage}
\begin{minipage}{3cm}
\centering
\epsfig{width=28mm, file=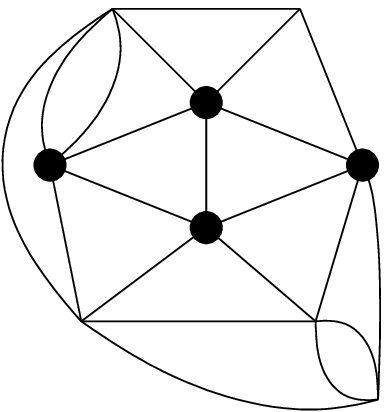}\par
$C_1$
% 25
\end{minipage}

\end{center}

\caption{Sporadic $5$-valent $(\{2,3\},5)$-polycycles (first part)}
\label{Sporadic5_valentFirst}
\end{figure}

\begin{figure}
\begin{center}
\begin{minipage}{3cm}
\centering
\epsfig{width=28mm, file=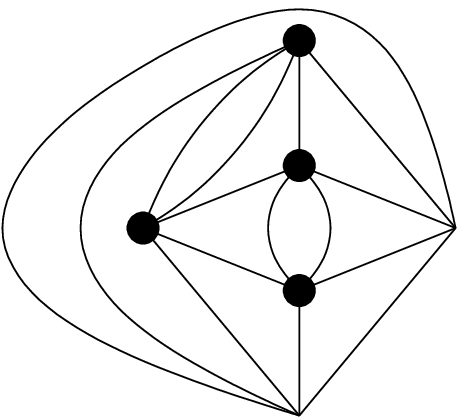}\par
$C_1$, nonext.
% 26
\end{minipage}
\begin{minipage}{3cm}
\centering
\epsfig{width=25mm, file=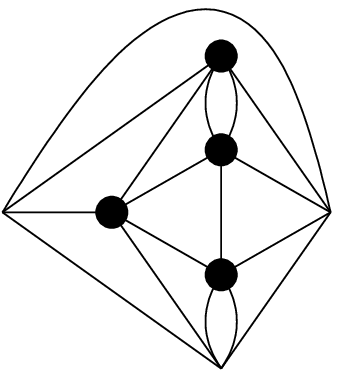}\par
$C_1$,~nonext.~$(C_{2\nu})$
% 27
\end{minipage}
\begin{minipage}{3cm}
\centering
\epsfig{width=28mm, file=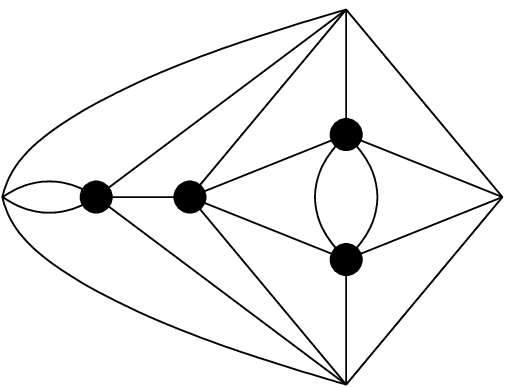}\par
$C_s$
% 28
\end{minipage}
\begin{minipage}{3cm}
\centering
\epsfig{width=28mm, file=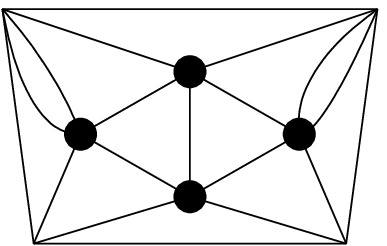}\par
$C_s$
% 29
\end{minipage}
\begin{minipage}{3cm}
\centering
\epsfig{width=25mm, file=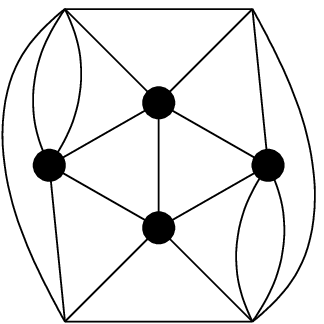}\par
$C_2$
% 30
\end{minipage}
\begin{minipage}{3cm}
\centering
\epsfig{width=28mm, file=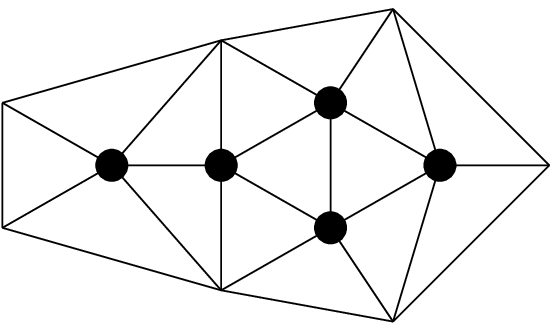}\par
$C_s$
\end{minipage}
\begin{minipage}{3cm}
\centering
\epsfig{width=28mm, file=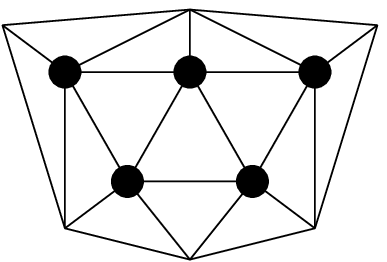}\par
$C_s$
\end{minipage}
\begin{minipage}{3cm}
\centering
\epsfig{width=28mm, file=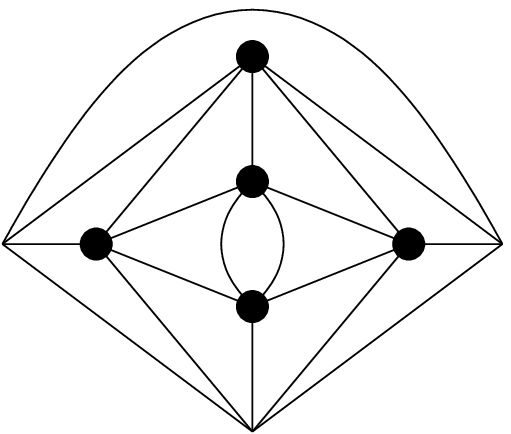}\par
$C_1$ $(C_{2\nu})$
\end{minipage}
\begin{minipage}{3cm}
\centering
\epsfig{width=28mm, file=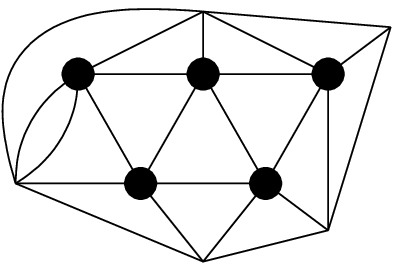}\par
$C_1$
\end{minipage}
\begin{minipage}{3cm}
\centering
\epsfig{width=28mm, file=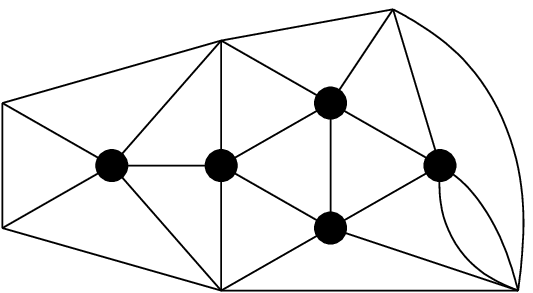}\par
$C_1$
\end{minipage}
\begin{minipage}{3cm}
\centering
\epsfig{width=28mm, file=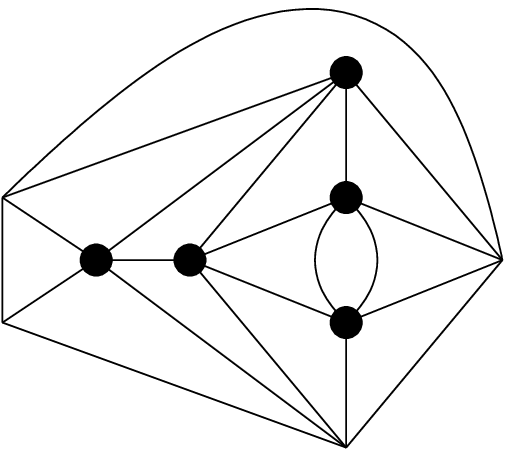}\par
$C_1$
\end{minipage}
\begin{minipage}{3cm}
\centering
\epsfig{width=28mm, file=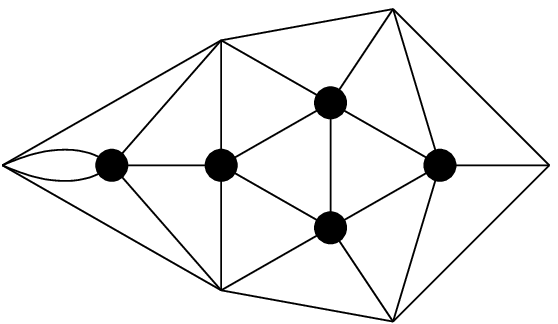}\par
$C_s$
\end{minipage}
\begin{minipage}{3cm}
\centering
\epsfig{width=28mm, file=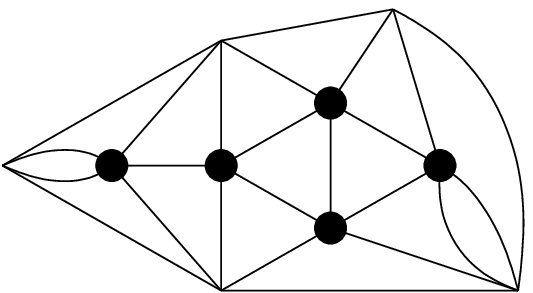}\par
$C_1$
\end{minipage}
\begin{minipage}{3cm}
\centering
\epsfig{width=28mm, file=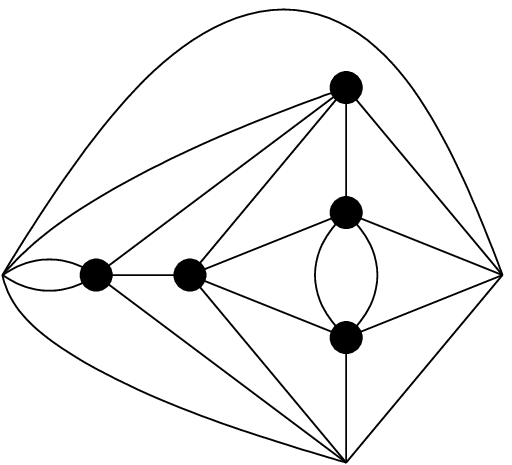}\par
$C_1$, nonext.
\end{minipage}
\begin{minipage}{3cm}
\centering
\epsfig{width=28mm, file=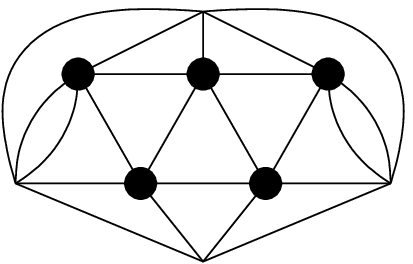}\par
$C_s$, nonext.
\end{minipage}
\begin{minipage}{3cm}
\centering
\epsfig{width=28mm, file=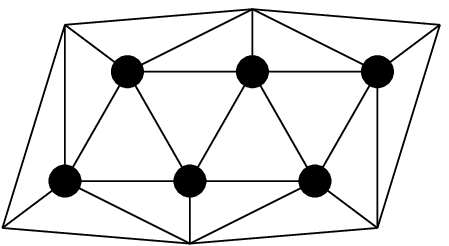}\par
$C_2$
\end{minipage}
\begin{minipage}{3cm}
\centering
\epsfig{width=28mm, file=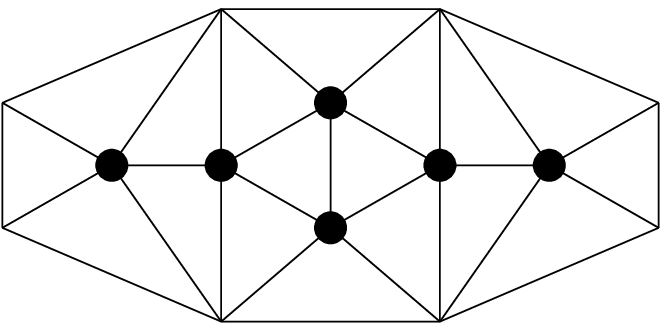}\par
$C_{2\nu}$
\end{minipage}
\begin{minipage}{3cm}
\centering
\epsfig{width=28mm, file=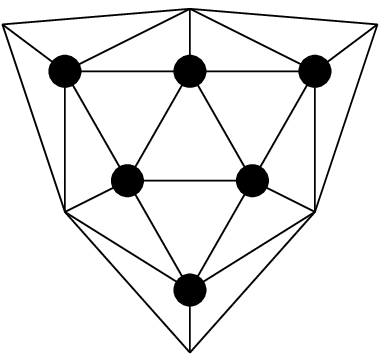}\par
$C_{3\nu}$
\end{minipage}
\begin{minipage}{3cm}
\centering
\epsfig{width=28mm, file=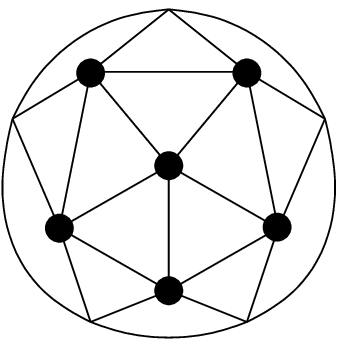}\par
$C_{5\nu}$
\end{minipage}
\begin{minipage}{3cm}
\centering
\epsfig{width=28mm, file=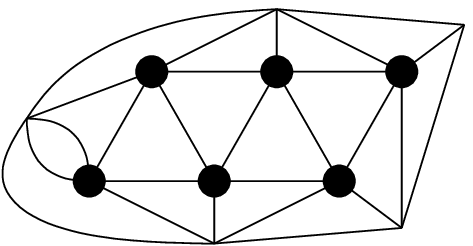}\par
$C_1$
\end{minipage}
\begin{minipage}{3cm}
\centering
\epsfig{width=28mm, file=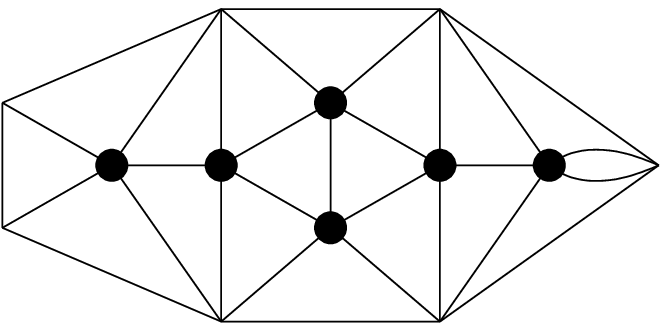}\par
$C_s$
\end{minipage}
\begin{minipage}{3cm}
\centering
\epsfig{width=28mm, file=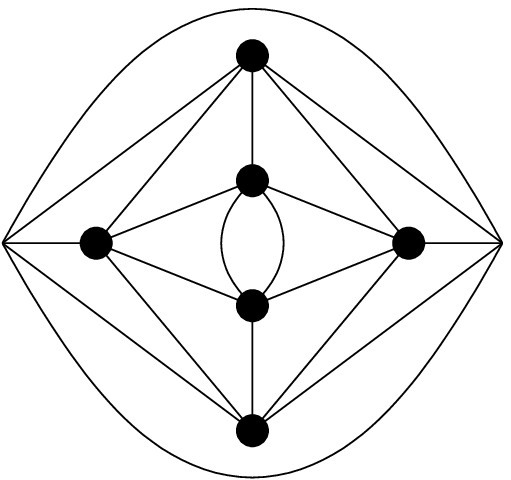}\par
$C_{2\nu}$, nonext.
\end{minipage}
\begin{minipage}{3cm}
\centering
\epsfig{width=28mm, file=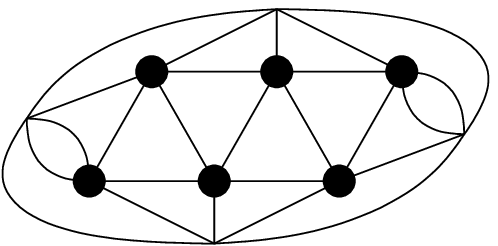}\par
$C_2$, nonext.
\end{minipage}
\begin{minipage}{3cm}
\centering
\epsfig{width=28mm, file=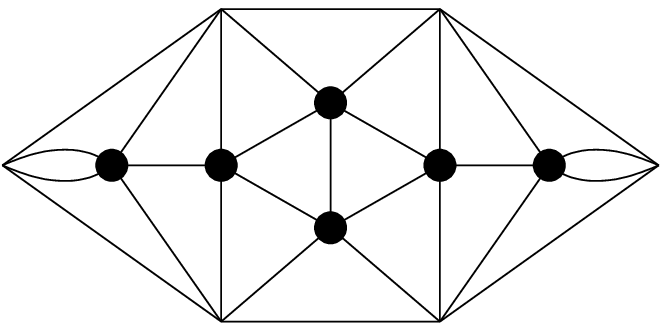}\par
$C_{2\nu}$, nonext.
\end{minipage}
\begin{minipage}{3cm}
\centering
\epsfig{width=23mm, file=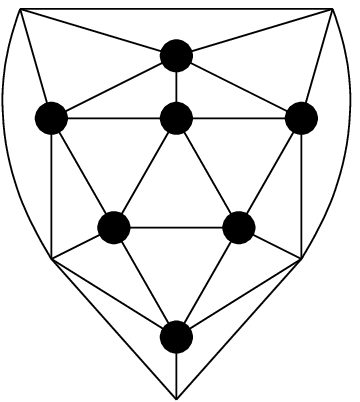}\par
$C_s$
\end{minipage}
\begin{minipage}{3cm}
\centering
\epsfig{width=28mm, file=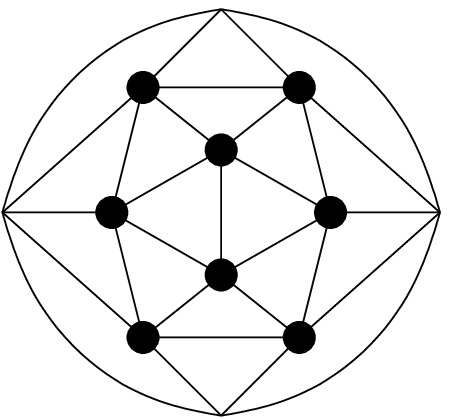}\par
$C_{2\nu}$
\end{minipage}
\begin{minipage}{3cm}
\centering
\epsfig{width=28mm, file=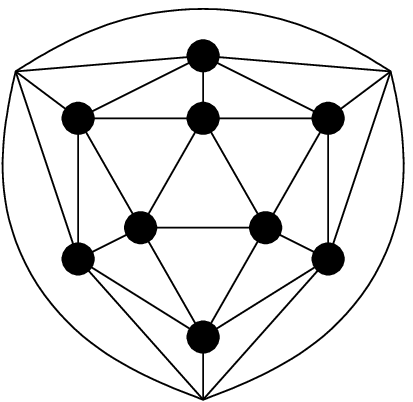}\par
$C_{3\nu}$,~nonext.~$(I_h)$
\end{minipage}

\end{center}

\caption{Sporadic $5$-valent $(\{2,3\},5)$-polycycles (second part)}
\label{Sporadic5_valentSecond}
\end{figure}

\begin{theorem}
The list of elementary $(\{2,3\},5)$-polycycles consists of:

\begin{enumerate}
\item[(i)] $57$ sporadic $(\{2,3\},5)$-polycycles given on
Figures \ref{Sporadic5_valentFirst} and \ref{Sporadic5_valentSecond},

\item[(ii)] three following infinite $(\{2,3\},5)$-polycycles:
\begin{center}
\begin{minipage}{8cm}
\centering
\epsfig{height=10mm, file=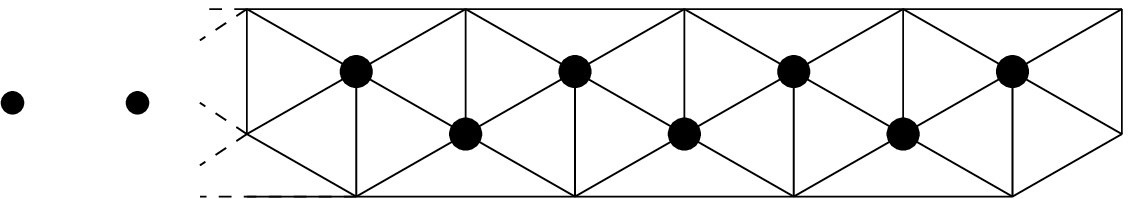}\par
$\alpha$: $C_1$
\end{minipage}
\begin{minipage}{8cm}
\centering
\epsfig{height=10mm, file=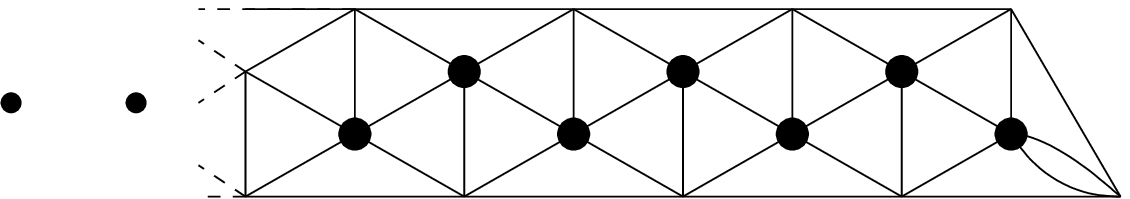}\par
$\beta$: $C_1$
\end{minipage}
\begin{minipage}{8cm}
\centering
\epsfig{height=10mm, file=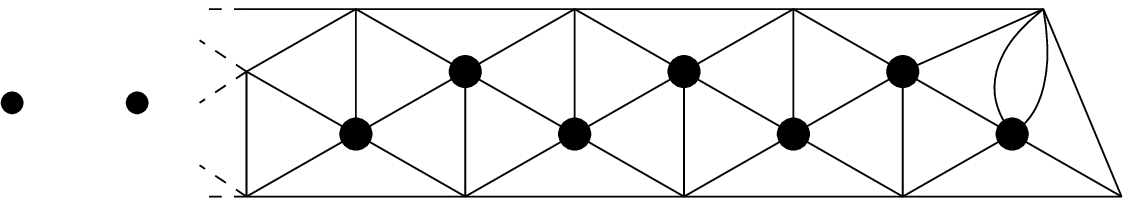}\par
$\gamma$: $C_1$, nonext.
\end{minipage}
\end{center}

\item[(iii)] six  infinite series of $(\{2,3\},5)$-polycycles with one hole
(they are obtained by concatenating endings of a pair of 
polycycles, given in (ii); see Figure \ref{SixInfiniteSeries} for the
first $5$ graphs),

\item[(iv)] the following $5$-valent doubly infinite
$(\{2,3\},5)$-polycycle, called {\em snub $\infty$-antiprism}:
\begin{center}
\begin{minipage}{8cm}
\centering
\epsfig{height=10mm, file=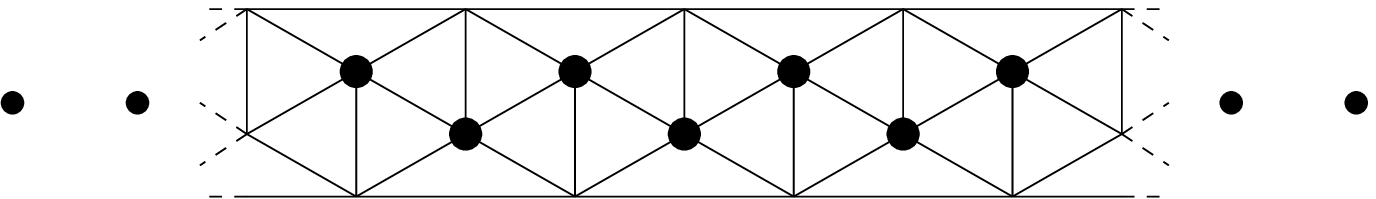}\par
pma2, nonext.
\end{minipage}
\end{center}

%NEED TO DO 33x y33

\item[(v)] the infinite series of 
{\em snub $m$-antiprisms}\footnote{{\em Snub $m$-antiprism} is a 
$5$-valent plane graph defined on page 119 of \cite{DGS}. For $m=4$ it 
is one of well-known {\em regular-faced polyhedra} (Johnson polyhedra).
For $m=3$ it has the 
skeleton of Icosahedron, but it is distinct from the 
corresponding sporadic one, having only one hole (exterior $3$-gonal face)}, 
$m\geq 2$ (two $m$-gonal holes, non-extensible for 
$m \ge 4$, with symmetry $D_{md}$), represented below for 
$m=2,3,4,5$:

\begin{center}
\begin{minipage}{3cm}
\centering
\epsfig{height=20mm, file=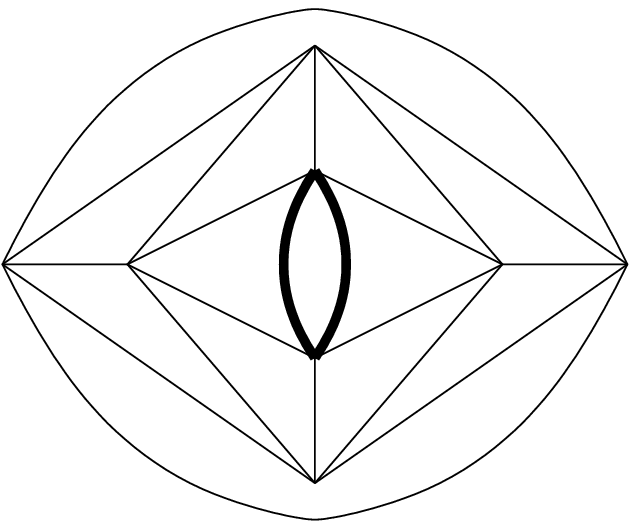}\par
$D_{2d}$
\end{minipage}
\begin{minipage}{3cm}
\centering
\epsfig{height=20mm, file=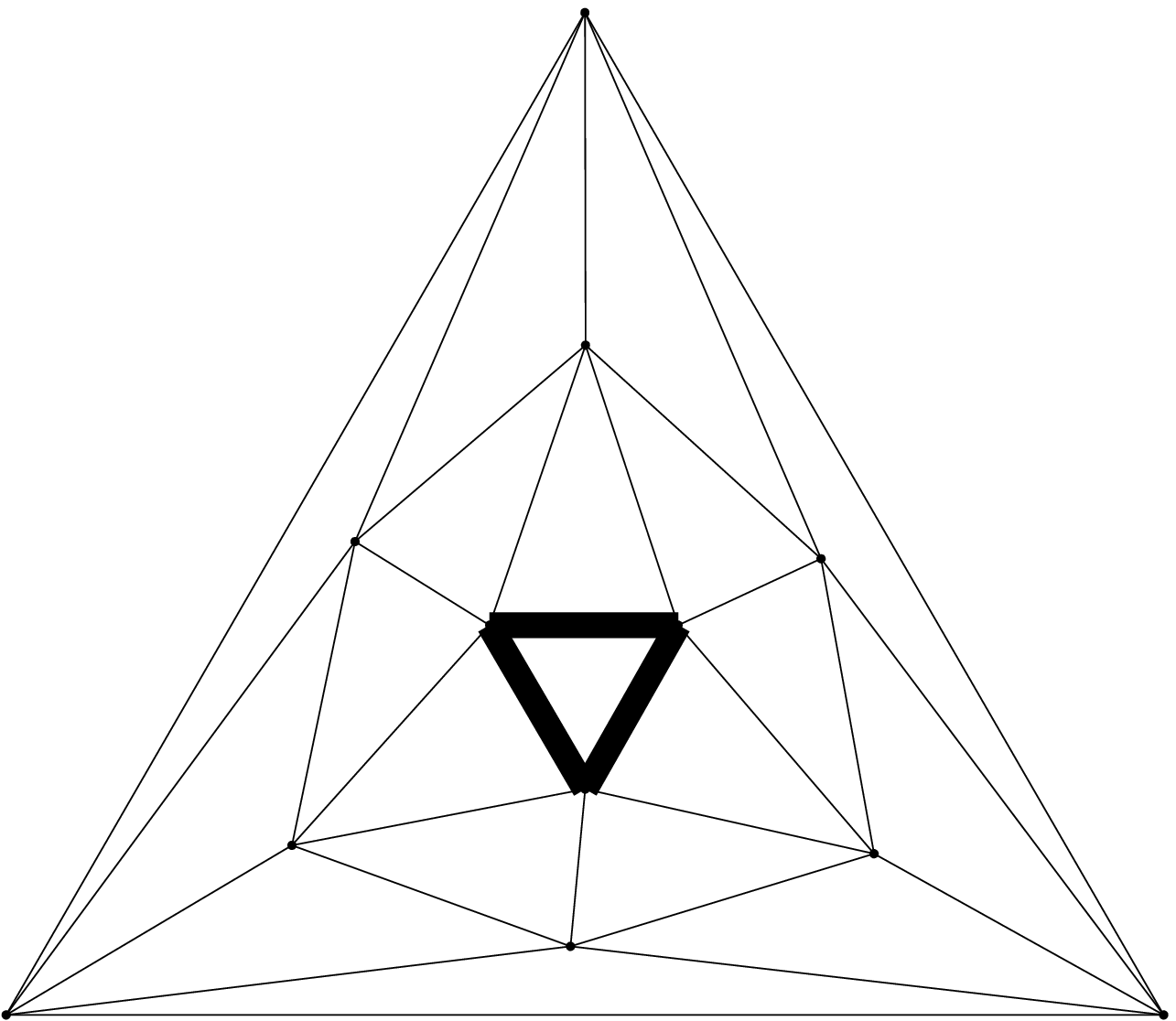}\par
$D_{3d}$ $(I_h)$
\end{minipage}
\begin{minipage}{3cm}
\centering
\epsfig{height=20mm, file=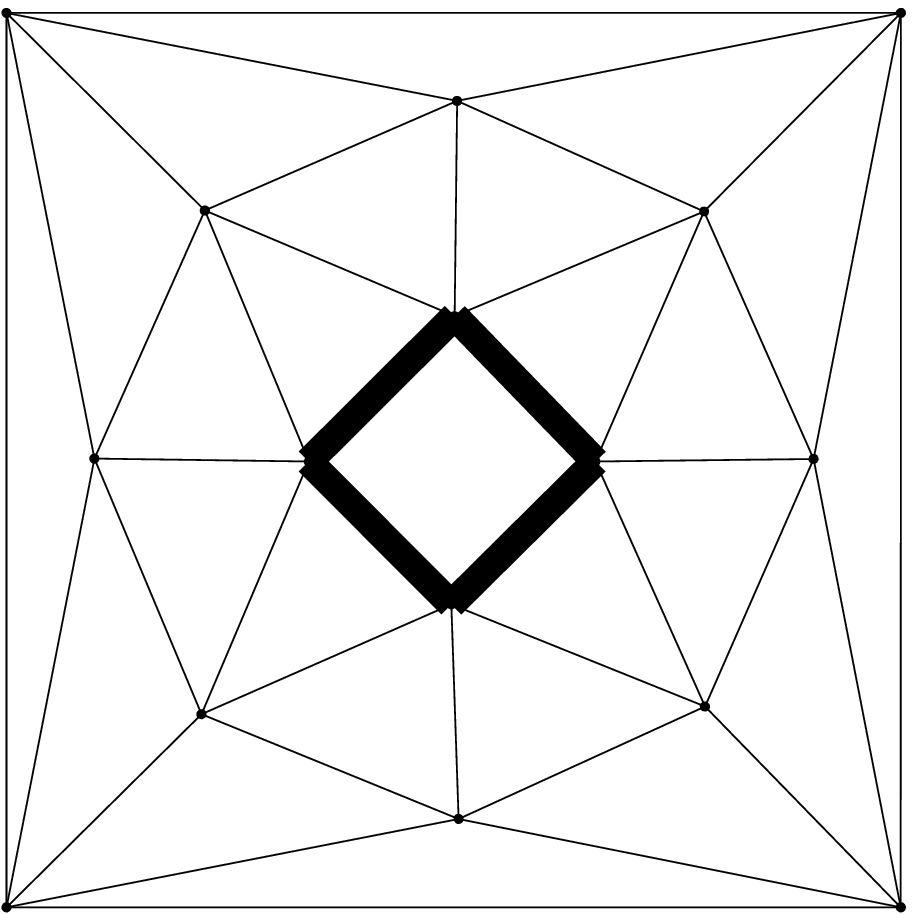}\par
$D_{4d}$
\end{minipage}
\begin{minipage}{3cm}
\centering
\epsfig{height=20mm, file=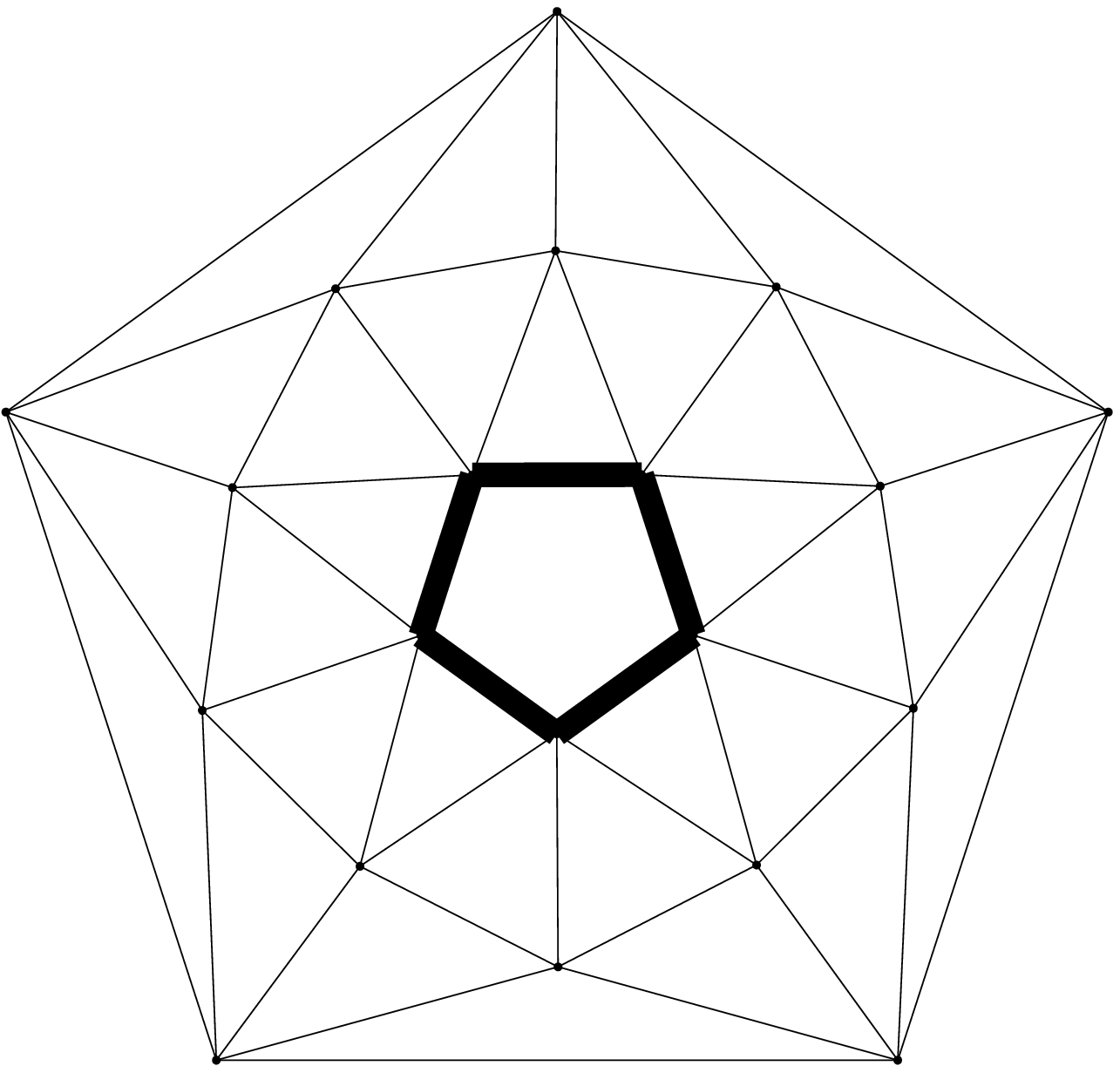}\par
$D_{5d}$
\end{minipage}
%\begin{minipage}{3cm}
%\centering
%\epsfig{height=20mm, file=ElementaryDrawing/5valSer6sec.eps}\par
%$D_{6d}$
%\end{minipage}

\end{center}

\end{enumerate}
\end{theorem}
\proof Let us take an elementary $(\{2,3\},5)$-polycycle, which
is finite. Then, by removing a triangle, which is adjacent to a boundary, 
one is led to the situation described above. Hence, the algorithm
for enumerating finite elementary $(\{2,3\},5)$ polycycles is 
the following.
\begin{enumerate}
\item Begin with isolated $i$-gons with $i=2$ or $3$.
\item For every vertex $v$ of an elementary polycycle with $n$ interior
vertices, consider all possibilities of adding $2$- and $3$-gons
incident to $v$, such that the obtained polycycle is elementary
and $v$ has become an interior vertex.
\item Reduce by isomorphism.
\end{enumerate}
The above algorithm first finds some sporadic $(\{2,3\},5)$-polycycles
and the first elements of the infinite series and then find only 
the elements of the infinite series.
In order to prove that this is the complete list of all finite
$(\{2,3\},5)$-polycycles with only one hole, one needs to consider the case,
in which only $e_4$ is a bridge going from a hole to the same hole.
So, we need 
to consider all possibilities, where the addition of two elementary
$(\{2,3\},5)$-polycycles and one $3$-gon make a larger 
elementary $(\{2,3\},5)$-polycycle. 
Given a sequence $a_1,\dots,a_n$, we say
that a sequence $b_1,\dots, b_p$ with $p<n$ is a {\em pattern} of that sequence
if, for some $n_0$, one has $a_{n_0+j-1}=b_j$ or $a_{n_0+1-j}=b_j$ with 
the
addition being modulo $n$.
The $(\{2,3\},5)$-polycycles,
used in that construction, should have the pattern $3,3,x$ with $x\leq 4$ 
in their boundary sequence. Only the polycycles, which belong to the
six infinite series, satisfy this and it is easy to see, that the result
of the operation is still one of the six infinite series.
So, the list of finite elementary $(\{2,3\},5)$-polycycles with one hole
is the announced one.

If a polycycle has more than one hole, then
it is obtained by the addition of $2$- and $3$-gonal faces to
a vertex, incident to the boundary
(or from another $(\{2,3\},5)$-polycycle
with a smaller number of holes) by the merging of two open edges
in a bridge and addition of a $3$-gon.
This suggest an iterative procedure, where we begin with a
$(\{2,3\},5)$-polycycle with one hole and consider all possibilities
of extension to a polycycle with two holes and then continue for more
complicated polycycles.

The polycycles, which are suitable for such a construction, should contain
the pattern $3,3,x$ and $y,3,3$ in their boundary sequence with $x\leq 4$
and $y\leq 3$. Amongst all finite $(\{2,3\},5)$-polycycles with one hole,
only members of the infinite series $\alpha\alpha$ of $(\{2,3\},5)$-polycycles
with one hole, are suitable ones and one obtains the infinite series of
the above theorem.
Furthermore, the members of this infinite series admit no extension.
This proves that
there are no other $(\{2,3\},5)$-polycycles with at least two holes.

Consider now an elementary infinite $(\{2,3\},5)$-polycycle $P$. 
Eliminate all $2$-gonal faces of $P$ and obtain another
$(\{3\},5)$-polycycle $P'$, which is not necessarily elementary.
We do a decomposition of $P'$ along its elementary components, which
are enumerated in \cite{DDS2}. If snub $\infty$-antiprism is one of
the components, then we are finished and $P=P'$
is the snub $\infty$-antiprism.
If $\alpha$ is one of the components, then one has two edges, along 
which to extend the polycycle; they are depicted below:
\begin{center}
\begin{minipage}{8cm}
\centering
\epsfig{height=15mm, file=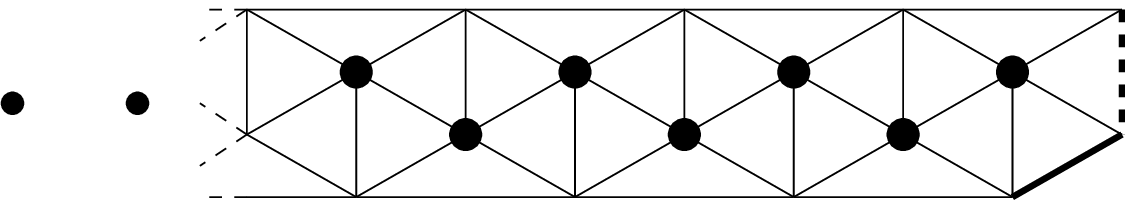}\par
\end{minipage}
\end{center}
Clearly, if we extend the polycycle along only one of those edges, then
the result is not an elementary polycycle. The consideration of all 
possibilities yields
$\beta$ and $\gamma$. Suppose now that $P'$ has no infinite components.
Then $P$ has at least one infinite path $f_0, \dots, f_i, \dots$, such 
that
$f_i$ is adjacent to $f_{i+1}$, but $f_{i-1}$ is not adjacent to 
$f_{i+1}$. The considerations, analogous to the $3$-valent case,
yield the result. \qed

%List of sporadic elementary $(\{2,3\},5)$-polycycles:

\newpage
\section{Possible extensions of the setting}

The classification, developped in this paper, is an extension of 
previous
work of Deza and Shtogrin, by allowing various possible 
number of sides of interior faces and several possible holes.
A natural question is, if one can further enlarge the class of 
polycycles.

We required $2$-connectivity and
that any two holes do not share a vertex.
If one removes those two hypothesis,
then many other graphs do appear.

Consider infinite series of $(\{2\},6)$-polycycles, {\em m-bracelets},
$m\geq 2$ (i.e. $m$-circle, but each edge is tripled). The central edge
is a bridge for those polycycles, for both $2$-gons of the triple of
edges. But if one removes those two digons, then the resulting 
plane graph 
has two holes sharing a face, i.e. violates the second of the crucial points
(i)--(iii) of the definition of $(R,q)$-polycycle. This shows that
our hypothesis were necessary.

For even $m$, each even edge (for some order $1,\dots, m$ of them)
can be duplicated $t$ times (for fixed $t$, $1\leq t \leq 5$), and each
odd edge duplicated $6-t$ times; so, all vertex-valencies will be still
$6$. On the other hand, two holes ($m$-gons inside and outside of the
$m$-bracelet), have common vertices; so, it is again not our 
polycycle.


\begin{thebibliography}{99}

\bibitem[DDS05a]{penthex}
M. Deza, M. Dutour and M.I. Shtogrin, {\em Filling of a given boundary by $p$-gons and related problems}, to appear in Proceedings of Int. Conference ``General Theory of Information Transfer and Combinatorics'' (Bielefeld, 2004), ed. by L. Baumer, in Electronic Notes in Discrete Mathematics (2005).

\bibitem[DDS05b]{DDS2}
M. Deza, M. Dutour and M.I. Shtogrin, {\em Elliptic polycycles with holes}, Uspechi Mat. Nauk. {\bf 60-2} (2005) 157--158 (in Russian). English translation in Russian Math. Surveys {\bf 60-2}.

\bibitem[DGS04]{DGS}
M. Deza, V.P. Grishukhin and M.I. Shtogrin, {\em Scale-Isometric Polytopal Graphs in Hypercubes and Cubic Lattices}, World Scientific and Imperial College Press, 2004.

%\bibitem[DSS05]{DSS1}
%M. Deza, S. Shpectorov and M.I. Shtogrin, {\em Non-extendible $(3,5)$-polycycles}, in preparation.

\bibitem[DS98]{DSp1}
M. Deza and M.I. Shtogrin, {\it Polycycles}, Voronoi Conference on Analytic Number Theory and Space Tilings (Kyiv, September 7--14, 1998), Abstracts, 19--23.

\bibitem[DS99]{DSp2}
M. Deza and M.I. Shtogrin, {\it Primitive polycycles and helicenes}, Uspechi Mat. Nauk. {\bf 54-6} (1999) 159--160 (in Russian). English translation in Russian Math. Surveys {\bf 54-6}, 1238--1239.

\bibitem[DS00a]{DSp3}
M. Deza and M.I. Shtogrin, {\it Infinite primitive polycycles}, Uspechi Mat. Nauk. {\bf 55-1} (2000) 179--180 (in Russian). English translation in Russian Math. Surveys {\bf 55-1}, 169--170.

\bibitem[DS00b]{DS5}
M. Deza and M.I.Shtogrin, {\em Embedding of chemical graphs into hypercubes}, 
Math. Zametki {\bf 68-3} (2000) 339--352 (in Russian). English translation 
in Mathematical Notes {\bf 68-3} 295--305.

\bibitem[DS00c]{7}
M. Deza and M.I. Shtogrin, {\it Polycycles: Symmetry and Embeddability}, Uspechi Mat. Nauk. {\bf 55-6} (2000) 129--130 (in Russian). English translation in Russian Math. Surveys {\bf 55-6}, 1146--1147.


%\bibitem[DSSt00b]{DSp4}
%M. Deza and M.I. Shtogrin, {\it Polycycles: symmetry and embedding}, Russian Math. Surveys {\bf 56-6} (2000) 159--160.


\bibitem[DS01]{DS8}
M. Deza and M.I. Shtogrin, {\em Clusters of Cycles}, Journal of Geometry
and Physics {\bf 40-3,4} (2001) 302--319.

\bibitem[DS02a]{DS11}
M.Deza and M.I.Shtogrin, {\em Criterion of embedding of $(r,q)$-polycycles},
Uspechi Mat. Nauk. {\bf 57-3} (2002) 149--150 (in Russian).
English translation in Russian Math. Surveys {\bf 57-3}, 589-591.

\bibitem[DS02b]{DS10}
M.Deza and M.I.Shtogrin, {\em Extremal and non-extendible polycycles}, 
Proceedings of Steklov Mathematical Institute, {\bf 239} (2002) 117--135. 
(Translated from Trudy of Steklov Math. Institut {\bf 239} (2002) 127--145).


%\bibitem{16}
%M.Deza and M.I. Shtogrin, {\em Embeddability Criterion for $(r,q)$-Polycycles}, Usp. Mat. Nauk, {\bf 57-3} (2002) 149--150.


\bibitem[DS04]{DS12}
M. Deza and M.I. Shtogrin, {\em Archimedean polycycles}, Uspechi Mat. Nauk. {\bf 59-3} (2004) 165--166 (in Russian). English translation in Russian Math. Surveys {\bf 59-3}, 564--566.

\bibitem[DS05a]{DS13a}
M. Deza and M.I. Shtogrin, {\em Metrics of constant curvature on polycycles}, to appear in Math. Zametki.

\bibitem[DS05b]{DS13b}
M. Deza and M.I. Shtogrin, {\em Types of polycycles}, in preparation.

\bibitem[DD05]{FaceRegular}
M. Dutour and M. Deza, {\em Face-regular $3$-valent two-faced spheres and tori}, in preparation.

\bibitem[Sh99]{S1}
M.I. Shtogrin, {\it Primitive polycycles: criterion}, Uspechi Mat. Nauk. {\bf 54-6} (1999) 177--178 (in Russian). English translation in Russian Math. Surveys {\bf 54-6}, 1261--1262.

\bibitem[Sh00]{S2}
M.I. Shtogrin, {\it Non-primitive polycycles and helicenes}, Uspechi Mat. Nauk. {\bf 55-2} (2000) 155--156 (in Russian). English translation in Russian Math. Surveys {\bf 55-2}, 358--360.



\end{thebibliography}
\end{document}